\documentclass[12pt]{amsart}

\title[Bordered HF modules for satellite operations using planar graphs]{Bordered Heegaard Floer modules for satellite operations using planar graphs}

\author{Shikhin Sethi}
\thanks{This work is licensed under Creative Commons Attribution-NonCommercial-ShareAlike 4.0 International.}

\usepackage{xcolor}
\definecolor{darkblue}{rgb}{0,0,0.4}
\usepackage[bookmarks=true, bookmarksopen=true,
bookmarksdepth=3,bookmarksopenlevel=2,
colorlinks=true,
linkcolor=darkblue,
citecolor=darkblue,
filecolor=darkblue,
menucolor=darkblue,
urlcolor=darkblue]{hyperref}

\hypersetup{pdftitle={Bordered Heegaard Floer Modules for Satellite Operations using Planar Graphs}}
\hypersetup{pdfauthor={shikhin}}

\usepackage[T1]{fontenc}
\usepackage[pdftex,lmargin=1in, rmargin=1in, tmargin=1in, bmargin=1in]{geometry}

\usepackage{setspace,caption}
\usepackage{amsmath, enumitem}
\usepackage{amssymb}
\usepackage{tikz}
\usepackage{subfig}
\usepackage{booktabs}
\usepackage{mathtools}
\usepackage{dsfont}
\usepackage{amsthm}

\usetikzlibrary{shapes, arrows, positioning, decorations.pathreplacing, patterns}
\usepackage{tikzit}
\usepackage{pgffor}

\tikzstyle{number-labels}=[fill=none, draw=none, shape=circle, font={\tiny}]
\tikzstyle{dot}=[fill=black, draw=black, shape=circle, inner sep=1 pt]
\tikzstyle{bigdot}=[fill=black, draw=black, shape=circle, inner sep=2 pt]
\tikzstyle{text-labels}=[font={\small}, fill=none]
\tikzstyle{surrounded}=[fill=white, draw=black, shape=circle, tikzit shape=circle, inner sep=1pt, font={\tiny}]
\tikzstyle{smalldot}=[fill=black, draw=black, shape=circle, inner sep=0.5 pt]

\tikzstyle{blue}=[-, draw=blue, thick]
\tikzstyle{red}=[-, draw=red, thick]
\tikzstyle{double-sided}=[<->, thick]
\tikzstyle{thickedge}=[-, draw=black, thick]
\tikzstyle{dashededged}=[-, dashed]
\tikzstyle{densedashedge}=[-, densely dashed]
\tikzstyle{denselydottedred}=[-, densely dotted, draw=red]
\tikzstyle{braceedge}=[-, decorate, decoration=brace]

\def\face#1{\hbox{$\bigcirc$\llap{\lower 2pt\hbox{\"{}\rlap{\lower 3pt
					\hbox{\kern -2.5pt$#1{}$}}}\kern 2.5pt}}}

\def\tongue{\hbox{$\bigcirc$\llap{\lower 2pt\hbox{\"{}\rlap{\lower 1pt
					\hbox{\kern -4.5pt${}^{{}_\sigma}$}}}\kern 2.5pt}}}

\DeclareMathAlphabet\mathbfcal{OMS}{cmsy}{b}{n}

\newcommand\numberthis{\addtocounter{equation}{1}\tag{\theequation}}

\newcommand{\defeq}{\vcentcolon=}
\newcommand{\eqdef}{=\vcentcolon}

\newcommand{\kk}{\mathds k}

\newcommand{\RR}{\mathbb R}
\newcommand{\CC}{\mathbb C}
\newcommand{\ZZ}{\mathbb Z}

\newcommand{\FF}{\mathbb F}

\newcommand{\bD}{\mathbb{D}}
\newcommand{\wt}{\widetilde}

\newcommand{\Filt}{\mathcal{F}}

\newcommand{\tensor}{\otimes}
\newcommand{\bdy}{\partial}

\newcommand{\spinc}{\mathfrak s}
\renewcommand{\Re}{\operatorname{Re}}

\DeclareMathOperator{\spin}{spin}
\newcommand{\SpinC}{\spin^c}

\DeclareMathOperator{\ind}{ind}
\DeclareMathOperator{\gr}{gr}
\newcommand{\gre}{\widetilde{\gr}}

\theoremstyle{plain}

\numberwithin{equation}{section}
\newtheorem{theorem}[equation]{Theorem}

\newtheorem{proposition}[equation]{Proposition}
\newtheorem{lemma}[equation]{Lemma}

\newtheorem{definition}[equation]{Definition}

\theoremstyle{definition}

\theoremstyle{remark}

\newtheorem{remark}[equation]{Remark}

\newcommand{\lab}[1]{$\scriptstyle #1$}

\newcommand{\HF}{\mathit{HF}}
\newcommand{\HFa}{\widehat {\HF}}

\newcommand{\HFm}{{\HF}^-}

\newcommand{\CFa}{\widehat {\mathit{CF}}}
\newcommand{\CFm}{\mathit{CF}^-}
\newcommand{\CFmm}{\mathbf{CF}^-}

\newcommand{\x}{\mathbf x}
\newcommand{\y}{\mathbf y}
\newcommand{\z}{\mathbf z}

\newcommand{\Ainf}{A_\infty}
\newcommand{\Algm}{\Alg_-}
\newcommand{\Algc}{\mathbfcal{A}_-}
\newcommand{\AlgmAs}{\Alg_-^{0,\text{as}}}
\newcommand{\Alg}{\mathcal{A}}

\newcommand{\alphas}{{\boldsymbol{\alpha}}}
\newcommand{\betas}{{\boldsymbol{\beta}}}

\newcommand{\cM}{\mathcal{M}}

\newcommand{\CFD}{\mathit{CFD}}
\newcommand{\CFDD}{\mathit{CFDD}}
\newcommand{\CFA}{\mathit{CFA}}

\newcommand{\CFDa}{\widehat{\CFD}}

\newcommand{\CFDm}{\CFD^-}
\newcommand{\CFAc}{\mathbf{CFA}^-}

\newcommand{\CFK}{\mathit{CFK}}
\newcommand{\CFKa}{\widehat{\CFK}}
\newcommand{\CFKm}{\CFK^-}

\newcommand{\CFDDm}{\CFDD^-}

\newcommand{\CFAa}{\widehat{\CFA}}

\newcommand{\fModule}{\mathfrak{M}}

\newcommand\Id{\mathbb{I}}

\newcommand\Gen{\mathfrak{S}}
\newcommand{\Heegaard}{\mathcal{H}}
\newcommand{\HD}{\Heegaard}
\newcommand{\CD}{\mathcal{C}}

\newcommand{\op}{\mathrm{op}}

\newcommand{\CFAm}{\CFA^-}

\newcommand{\AsDiag}{{\boldsymbol{\Gamma}}} 
\newcommand{\TrDiag}{\boldsymbol{\gamma}}
\newcommand{\TrMPrim}{\mathbf{p}}

\newcommand{\TrPMDiag}{\mathbf{p}}
\newcommand{\Yvar}{\mathcal{Y}}

\usepackage[
backend=biber,
style=alphabetic,
doi=false,isbn=false,url=false,
eprint=true
]{biblatex}

\begin{document}

\begin{abstract}
	Lipshitz, Ozsv\'ath, and Thurston extend the theory of bordered Heegaard Floer homology to compute $\CFmm$. Like with the \emph{hat} theory, their \emph{minus} invariants provide a recipe to compute knot invariants associated to satellite knots. We combinatorially construct the weighted $\Ainf$-modules associated to the $(p, 1)$-cable. The operations on these modules count certain classes of inductively constructed decorated planar graphs. This description of the weighted $\Ainf$-modules provides a combinatorial proof of the $\Ainf$ structure relations for the modules. We further prove a uniqueness property for the modules we construct: any weighted extensions of the unweighted $U = 0$ modules have isomorphic associated type D modules.	
\end{abstract}

\maketitle
\tableofcontents
\thispagestyle{empty}

\section{Introduction\label{ch:intro}}

At the start of the century, Ozsv\'ath and Szab\'o introduced the package of Heegaard Floer homology to study closed three-manifolds \cite{OSz04:HolomorphicDisks}. To a closed, oriented three-manifold $Y$, Heegaard Floer homology associates a graded chain complex $\CFm(Y)$ over the polynomial ring $\FF[U]$. After choosing a pointed Heegaard diagram representing $Y$, the chain complex is constructed as a variant of Lagrangian Floer homology on certain Lagrangians in a symmetric product arising from the Heegaard splitting. The chain homotopy type of $\CFm(Y)$ is an invariant of the manifold. There are different \emph{flavors} of the invariant; in particular, specializing with $U = 0$ provides a simpler and easier to compute invariant, denoted $\CFa(Y)$. The homologies $\HFm(Y)$ and $\HFa(Y)$ of the respective chain complexes are also invariants of the manifold.
	
The construction of Heegaard Floer homology as a version of Lagrangian Floer homology leaves it hard to compute; as defined, the differential, in particular, relies on counts of solutions to certain non-linear partial differential equations. In 2006, Sarkar and Wang provided an easier to compute combinatorial description of the \emph{hat} flavor $\HFa(Y)$ using especially \emph{nice} Heegaard diagrams \cite{Sarkar2010}. The \emph{hat} invariants suffice for significant applications to the study of three-manifolds. However, the Heegaard Floer package also captures invariants for smooth, closed, oriented four-manifolds with the \emph{mixed invariants} \cite{OSz06}. The mixed invariants rely intrinsically on the $U$-module structure of $\CFm$, and so it remains desirable to have tractable ways to compute $\CFm$.
		
Rasmussen \cite{Ras03} and, independently, Ozsv\'ath and Szab\'o \cite{OSz04:KnotInvariants} refined the Heegaard Floer package to an invariant of knots (or links) in three-manifolds. A knot $K$ in a three-manifold $Y$ can be described with a doubly-pointed Heegaard diagram, which gives rise to a chain complex $\CFKm(K)$ over the two-variable polynomial ring $\FF[U,V]$. The chain homotopy type of $\CFKm(K)$ is a knot invariant. As with the three-manifold invariant, the $V = 0$ specialization $\CFKa(K)$ provides an easier to compute flavor of the knot Floer complex. In particular, Sarkar's and Wang's combinatorial description extends to computing $\CFKa(K)$. Early computations of the knot Floer complex of knots relied upon analysis of appropriate Heegaard diagrams. For example, Goda, Matsuda, and Morifuji proved that $(1,1)$-knots admit a genus one Heegaard diagram from which the knot Floer homology groups can be combinatorially computed \cite{Goda2005}. Hedden computed some homology groups for certain cables of knots by finding appropriate Heegaard diagrams for the cables starting with a diagram of the pattern \cite{Hedden2005}.

As with the case for invariants of three-manifolds, specializing the chain complex $\CFKm(K)$ with $V = 0$ loses some information. The concordance invariant $\Upsilon$ defined by Ozsv\'ath, Stipsicz, and Szab\'o requires the full information contained in $\CFKm(K)$ \cite{O17}. The knot Floer complex of a knot can be used recover the Heegaard Floer homology of Dehn surgeries along the knot \cite{OSz08}; to compute $\CFm$ of the Dehn surgeries requires the full power of $\CFKm(K)$.
		
Lipshitz \cite{Lipshitz06:BorderedHF}, and subsequently Lipshitz, Ozsv\'ath, and Thurston \cite{LOT1}, extended Heegaard Floer homology to a theory associating invariants to manifolds with boundary. That is, cut your closed, oriented three-manifold $Y$ along an oriented surface $F$. The theory of \emph{bordered} Heegaard Floer homology associates invariants to each piece thus obtained, and provides a way to glue the invariants along the surface $F$ to recover $\CFa(Y)$. This provides a scheme to compute both Heegaard Floer invariants and knot Floer invariants by decomposing manifolds into simpler pieces. In particular, the knot Floer complexes of satellite knots can be better understood by such a decomposition.
	
Levine used bordered Floer homology to study generalized Whitehead doubles of knots \cite{Levine2012}; Petkova used it to generalize Hedden's result and computed the knot Floer homology groups of $(p, pn+1)$-cables of knots with \emph{thin} Floer homology \cite{Petkova2013}. Petkova's result relies on the computation of the bordered Floer invariant associated to the $(p, 1)$-cable in the solid torus. Hom provided a formula for a concordance invariant $\tau$ derived from the Heegaard Floer package for cables of knots using bordered Floer homology \cite{Hom2013}.
	
Recently, Lipshitz, Ozsv\'ath, and Thurston have been working on extending bordered Floer homology from the \emph{hat} to the \emph{minus} variant for manifolds split along a torus boundary \cite{LOT:torus-alg, LOT:torus-mod}. To the torus, they associate a \emph{weighted $\Ainf$-algebra} $\Algm$. $\Algm$ keeps track of the asymptotics of pseudoholomorphic curves near infinities that arise from the puncture---corresponding to the torus boundary---in the now \emph{bordered} Heegaard diagram. The construction of $\Algm$ is combinatorial: the $\Ainf$-operations on $\Algm$ count certain classes of planar graphs. To a three-manifold $Y$ with boundary identified with the torus, they associate a \emph{weighted $\Ainf$-module} $\CFAm(Y)$ over $\Algm$, with the operations given by counts of certain pseudoholomorphic curves. In upcoming work \cite{LOT:torus-pairing}, they establish a pairing between the invariants for two manifolds with torus boundary. This pairing yields the (completed) \emph{minus} invariants for the closed three-manifold obtained by gluing the two pieces along the common boundary.
	
In analogue with the \emph{hat} flavors, the new bordered Floer invariants lend themselves nicely to the computation of the unspecialized knot Floer complex $\CFKm(K)$ of satellite knots. In this paper, we study the weighted $\Ainf$-modules $\CFAm(\CD_p)$ associated to the $(p, 1)$-cable. This extends Petkova's work on computing $\CFAa(\CD_p)$. 

In Section~\ref{ch:background}, we recall preliminaries from the theory of bordered Floer homology, with an emphasis on the new technology developed for the \emph{minus} flavor. In Section~\ref{ch:cfa}, we provide a combinatorial description of the weighted $\Ainf$-structure
on the $\Ainf$-modules. The $\Ainf$-operations on the module count certain planar graphs, akin
to the algebra operations on $\Alg_-$. This explicit description of the
$\Ainf$-structure lets us verify the $\Ainf$-relations on the module. In Section~\ref{sec:11-patterns}, we remark on how a similar analysis would extend to other $(1,1)$-pattern knots. In Section~\ref{ch:algprop}, we study algebraic properties of the modules we construct. In particular, we prove a boundedness property needed to make sense of tensor products of the weighted $\Ainf$-modules. We further prove a certain uniqueness property, characterizing the weighted $\Ainf$-modules in terms of the specialized unweighted modules. In Section~\ref{ch:cfd}, we compute some sample tensor products using our modules.

Other methods have been developed to compute the knot Floer complex, and that assist in computation for satellite knots. Hanselman, Rasmussen, and Watson provide an interpretation of the bordered Floer invariants for the \emph{hat} flavor as immersed curves in the torus \cite{Hanselman2023}. In 2023, Hanselman constructed a refinement of the immersed curve invariant with additional decorations \cite{Ha23}. This theory sees the \emph{minus} flavor of Heegaard Floer homology; Chen and Hanselman build on top the immersed curve invariant to compute the $UV = 0$ specialization of the knot Floer complex by constructing analogues of the bordered Floer invariants for immersed Heegaard diagrams \cite{CH23}. For $(1,1)$-patterns, this computation is combinatorial. We rely on computations performed with the immersed curve invariants to check our modules.
	
	In a different vein, Zemke constructs a bordered Floer theory for the \emph{minus} flavor using input obtained from the link surgery formula \cite{Zemke2021}. Recent work of Chen, Zemke, and Zhou builds on top of this to provide a formula for the knot Floer complex of satellite knots when the link Floer complex of the $2$-component link associated to the pattern knot is an L-space link \cite{CZZ24}. There are $(1,1)$-patterns such that the associated $2$-component link is not an L-space link---say, the mirror of the Whitehead link---and their formulae do not extend to such pattern knots.

\section*{Acknowledgments}
This paper is an adapted rendition of the author's thesis at Princeton University. The author is grateful to their advisor, Peter Ozsv\'ath, for suggesting the project that lead to this paper, and for all their guidance in developing it. The author acknowledges Otte Hein\"avaara for mathematical conversations around this paper, and for the graph theoretic arguments in Section~\ref{sec:char-mod-tiling}.

\section{Bordered Heegaard Floer homology}
\label{ch:background}

In this section, we recall the preliminaries of the bordered Heegaard Floer theory. We assume familiarity with Heegaard Floer homology and knot Floer homology. The \emph{hat} flavor of the theory was developed in \cite{Lipshitz06:BorderedHF} by Lipshitz and subsequently in \cite{LOT1} by Lipshitz, Ozsv\'ath, and Thurston. The \emph{minus} flavor of the theory was developed in \cite{LOT:torus-alg,LOT:torus-mod,LOT:torus-pairing} by Lipshitz, Ozsv\'ath, and Thurston. An interested reader is referred to these texts for details; in particular, we only consider bordered manifolds with torus boundary, while the \emph{hat} flavor is developed for more general manifolds.

\subsection{Algebraic underpinnings}
\label{sec:alg-underpinnings}

The algebraic objects defined in this section are studied in depth in \cite{LOT:abstract}. Let $\kk$ be a unital ring of characteristic two.

\begin{definition}[$\Ainf$-algebra]
	An $\Ainf$-algebra $\Alg$ over $\kk$ is a $\kk$-module $A$
	with $n$-input algebra operations 
	$\mu_n: A^{\otimes_\kk n} \to A$ for all $n \geq 1$, such that the $\mu_n$ are $\kk$-module maps and they satisfy the $\Ainf$ structure relations
		\[ \sum_{i+j+l = n}
			\mu_{i + 1 + l}\left(a_1, \cdots, a_{i}, \mu_{j}(a_{i+1}, \dots, a_{i+j}),
			a_{i+j+1}, \dots, a_{i+j+l}\right) = 0\]
		for all $n \geq 1$ and $a_1, \dots, a_n \in A$.
		
	The $\Ainf$-algebra $\Alg$ is \emph{strictly unital} if there is an element $1 \in A$ such that $\mu_2(a, 1) = \mu_2(1, a) = a$ for all $a$; and if some $a_i = 1$, then $\mu_n(a_1, \dots, a_n) = 0$ for any $n \neq 2$.
		\end{definition}

An operation $\mu_n$ can be depicted graphically by an oriented tree, each edge oriented downwards, with exactly one vertex with $n$ incoming edges and one outgoing edge. This tree we call an $n$-corolla. Natural composition of such trees yields a planar tree representing a composition of operations in $\Alg$.
To obtain the $\Ainf$ relation corresponding to $n$, start with the $n$-corolla and consider all planar trees $T$ with two vertices such that contracting any edge in $T$ results in the $n$-corolla. The sum of the compositions given by all such trees $T$ is what equals zero in the $\Ainf$ relations.
	
\begin{center}
	\begin{tikzpicture}
	\begin{pgfonlayer}{nodelayer}
		\node [style=none] (0) at (2.5, 1) {};
		\node [style=smalldot] (1) at (2.5, -0.4) {};
		\node [style=none] (2) at (2.5, -1) {};
		\node [style=text-labels] (3) at (11.25, 0) {};
		\node [style=text-labels] (4) at (11.25, 0) {$=$};
		\node [style=smalldot] (5) at (2.5, 0.4) {};
		\node [style=none] (6) at (-8.5, 0) {$n = 1$:};
		\node [style=none] (7) at (-8.5, -3.5) {$n = 2$:};
		\node [style=none] (8) at (-2.15, -2.25) {};
		\node [style=smalldot] (9) at (-1.5, -3) {};
		\node [style=none] (10) at (-0.85, -2.25) {};
		\node [style=smalldot] (11) at (-1.5, -3.5) {};
		\node [style=none] (12) at (-1.5, -4.15) {};
		\node [style=none] (13) at (0.5, -3.5) {};
		\node [style=none] (14) at (0.5, -3.5) {$+$};
		\node [style=none] (15) at (1.75, -2.25) {};
		\node [style=smalldot] (16) at (2.15, -2.9) {};
		\node [style=none] (17) at (3.25, -2.25) {};
		\node [style=smalldot] (18) at (2.5, -3.475) {};
		\node [style=none] (19) at (2.5, -4.125) {};
		\node [style=none] (20) at (5.75, -2.25) {};
		\node [style=none] (22) at (7.25, -2.25) {};
		\node [style=smalldot] (23) at (6.5, -3.475) {};
		\node [style=none] (24) at (6.5, -4.15) {};
		\node [style=none] (25) at (4.5, -3.5) {};
		\node [style=none] (26) at (4.5, -3.5) {$+$};
		\node [style=text-labels] (27) at (11.25, -3.5) {};
		\node [style=text-labels] (28) at (11.25, -3.5) {$=$};
		\node [style=smalldot] (29) at (6.875, -2.9) {};
		\node [style=none] (30) at (-8.5, -7) {$n = 3$:};
		\node [style=none] (31) at (-4.5, -6) {};
		\node [style=smalldot] (32) at (-4.15, -6.55) {};
		\node [style=smalldot] (33) at (-3.75, -7.25) {};
		\node [style=none] (34) at (-3.75, -6) {};
		\node [style=none] (35) at (-3, -6) {};
		\node [style=none] (36) at (-3.75, -8) {};
		\node [style=none] (37) at (-2, -6) {};
		\node [style=smalldot] (39) at (-1.25, -7.25) {};
		\node [style=none] (40) at (-1.25, -6) {};
		\node [style=none] (41) at (-0.5, -6) {};
		\node [style=none] (42) at (-1.25, -8) {};
		\node [style=none] (43) at (0.5, -6) {};
		\node [style=smalldot] (45) at (1.25, -7.25) {};
		\node [style=none] (46) at (1.25, -6) {};
		\node [style=none] (47) at (2, -6) {};
		\node [style=none] (48) at (1.25, -8) {};
		\node [style=none] (55) at (3, -6) {};
		\node [style=smalldot] (57) at (3.75, -6.75) {};
		\node [style=none] (58) at (3.75, -6) {};
		\node [style=none] (59) at (4.5, -6) {};
		\node [style=none] (60) at (3.75, -8) {};
		\node [style=none] (61) at (5.5, -6) {};
		\node [style=smalldot] (63) at (5.95, -6.75) {};
		\node [style=none] (64) at (6.25, -6) {};
		\node [style=none] (65) at (7, -6) {};
		\node [style=none] (66) at (6.25, -8) {};
		\node [style=none] (67) at (12.25, 0) {$0$};
		\node [style=none] (68) at (12.25, -3.5) {$0$};
		\node [style=none] (69) at (8, -6) {};
		\node [style=smalldot] (71) at (9.05, -6.75) {};
		\node [style=none] (72) at (8.75, -6) {};
		\node [style=none] (73) at (9.5, -6) {};
		\node [style=none] (74) at (8.75, -8) {};
		\node [style=none] (75) at (-2.5, -7) {$+$};
		\node [style=none] (76) at (0, -7) {$+$};
		\node [style=none] (77) at (2.5, -7) {$+$};
		\node [style=none] (78) at (5, -7) {$+$};
		\node [style=none] (79) at (7.5, -7) {$+$};
		\node [style=text-labels] (80) at (11.25, -7) {};
		\node [style=text-labels] (81) at (11.25, -7) {$=$};
		\node [style=none] (82) at (12.25, -7) {$0$};
		\node [style=smalldot] (83) at (-1.25, -6.5) {};
		\node [style=smalldot] (84) at (1.675, -6.55) {};
		\node [style=smalldot] (85) at (3.75, -7.3) {};
		\node [style=smalldot] (86) at (6.25, -7.35) {};
		\node [style=smalldot] (88) at (8.75, -7.35) {};
		\node [style=none] (89) at (-12, 1) {};
		\node [style=smalldot] (90) at (-12, 0) {};
		\node [style=none] (91) at (-12, -1) {};
		\node [style=none] (92) at (-12.5, -2.5) {};
		\node [style=smalldot] (93) at (-12, -3.5) {};
		\node [style=none] (94) at (-12, -4.5) {};
		\node [style=none] (95) at (-11.5, -2.5) {};
		\node [style=none] (96) at (-12.75, -6) {};
		\node [style=smalldot] (97) at (-12, -7) {};
		\node [style=none] (98) at (-12, -8) {};
		\node [style=none] (99) at (-11.25, -6) {};
		\node [style=none] (100) at (-12, -6) {};
	\end{pgfonlayer}
	\begin{pgfonlayer}{edgelayer}
		\draw (0.center) to (1);
		\draw (1) to (2.center);
		\draw (8.center) to (9);
		\draw (10.center) to (9);
		\draw (9) to (11);
		\draw (11) to (12.center);
		\draw (15.center) to (16);
		\draw (16) to (18);
		\draw (18) to (19.center);
		\draw (17.center) to (18);
		\draw (23) to (24.center);
		\draw (29) to (22.center);
		\draw (29) to (23);
		\draw (23) to (20.center);
		\draw (31.center) to (32);
		\draw (32) to (33);
		\draw (33) to (34.center);
		\draw (33) to (35.center);
		\draw (33) to (36.center);
		\draw (39) to (41.center);
		\draw (39) to (42.center);
		\draw (45) to (46.center);
		\draw (45) to (48.center);
		\draw (57) to (58.center);
		\draw (57) to (59.center);
		\draw (63) to (64.center);
		\draw (71) to (72.center);
		\draw (37.center) to (39);
		\draw (40.center) to (83);
		\draw (83) to (39);
		\draw (47.center) to (84);
		\draw (84) to (45);
		\draw (45) to (43.center);
		\draw (55.center) to (57);
		\draw (57) to (85);
		\draw (85) to (60.center);
		\draw (65.center) to (86);
		\draw (63) to (86);
		\draw (86) to (66.center);
		\draw (71) to (88);
		\draw (88) to (74.center);
		\draw (69.center) to (88);
		\draw (61.center) to (63);
		\draw (73.center) to (71);
		\draw (89.center) to (90);
		\draw (90) to (91.center);
		\draw (92.center) to (93);
		\draw (93) to (94.center);
		\draw (95.center) to (93);
		\draw (96.center) to (97);
		\draw (97) to (98.center);
		\draw (99.center) to (97);
		\draw (100.center) to (97);
	\end{pgfonlayer}
\end{tikzpicture}\\
\end{center}

The $n = 1$ relation says that $\mu_1$  is a differential; the $n = 2$ relation is the Leibniz rule; the $n = 3$ relation says that associativity holds up to a homotopy.

\begin{definition}[Weighted (curved) $\Ainf$-algebra]
			A weighted (curved) $\Ainf$-algebra over $\kk$ is a $\kk$-module $A$ with $n$-input algebra operations 
			$\mu_n^w: A^{\otimes_\kk n} \to A$ for all $n,w \in \ZZ_{\geq0}$, such that the $\Ainf$-relations
			\[\sum_{\substack{i+j+l = n \\ w_1 + w_2 = w}}
			\mu_{i + 1 + l}^{w_1}\left(a_1, \cdots, a_{i}, \mu_{j}^{w_2}(a_{i+1}, \dots, a_{i+j}),
			a_{i+j+1}, \dots, a_{i+j+l}\right) = 0\]
			are satisfied for all $n, w \geq 0$ and $a_1, \dots, a_n \in A$. Here, $w$ is called the \emph{weight} of the operation $\mu^w_n$.
			
				The weighted $\Ainf$-algebra $\Alg$ is \emph{strictly unital} if there is an element $1 \in A$ such that $\mu_2^0(a, 1) = \mu_2^0(1, a) = a$ for all $a$; and if some $a_i = 1$, then  $\mu_n^w(a_1, \dots, a_n) = 0$ for all $(n,w) \neq (2,0)$.
			
			The $0$-input operations $\mu_0^w: \kk \to A$ provide the \emph{curvature} $\mu_0^w \defeq \mu_0^w(1)$ of $\Alg$. We require that $\mu_0^0 = 0$.
\end{definition}

A weighted operation $\mu_n^w$ can be represented graphically with an $n$-corolla where the vertex is marked to be of weight $w$.
To obtain the $\Ainf$ relation corresponding to $n$, start with the weight $w$ $n$-corolla, and consider all planar trees $T$ with two vertices such that contracting any edge in $T$ results in the weight $w$ $n$-corolla. Here, contracting an edge between two weighted vertices adds their weight. The sum of the compositions given by all such trees $T$ is what equals zero in the $\Ainf$ relations.
	
\begin{center}
	\begin{tikzpicture}
	\begin{pgfonlayer}{nodelayer}
		\node [style=surrounded] (1) at (4.5, 0) {$0$};
		\node [style=none] (2) at (4.5, -1) {};
		\node [style=text-labels] (3) at (18, 0) {};
		\node [style=text-labels] (4) at (18, 0) {$=$};
		\node [style=surrounded] (5) at (4.5, 1) {$1$};
		\node [style=none] (6) at (-9.5, 0) {$(n,w) = (0,1)$:};
		\node [style=none] (7) at (-9.5, -3.5) {$(n,w) = (1,1)$:};
		\node [style=text-labels] (27) at (18, -3.5) {};
		\node [style=text-labels] (28) at (18, -3.5) {$=$};
		\node [style=none] (30) at (-9.5, -7) {$(n,w) = (2,1)$:};
		\node [style=none] (31) at (-6.25, -6) {};
		\node [style=surrounded] (32) at (-5.9, -6.55) {$1$};
		\node [style=surrounded] (33) at (-5.5, -7.25) {$0$};
		\node [style=none] (35) at (-4.75, -6) {};
		\node [style=none] (36) at (-5.5, -8) {};
		\node [style=none] (37) at (9.5, -6) {};
		\node [style=surrounded] (39) at (10.25, -7.25) {$0$};
		\node [style=none] (41) at (11, -6) {};
		\node [style=none] (42) at (10.25, -8) {};
		\node [style=none] (43) at (-1.25, -6) {};
		\node [style=surrounded] (45) at (-0.5, -7.25) {$0$};
		\node [style=none] (47) at (0.25, -6) {};
		\node [style=none] (48) at (-0.5, -8) {};
		\node [style=none] (55) at (3.75, -6) {};
		\node [style=surrounded] (57) at (4.5, -6.7) {$0$};
		\node [style=none] (59) at (5.25, -6) {};
		\node [style=none] (60) at (4.5, -8) {};
		\node [style=none] (61) at (12.25, -6) {};
		\node [style=none] (65) at (13.75, -6) {};
		\node [style=none] (66) at (13, -8) {};
		\node [style=none] (67) at (19, 0) {$0$};
		\node [style=none] (68) at (19, -3.5) {$0$};
		\node [style=none] (69) at (14.75, -6) {};
		\node [style=none] (73) at (16.25, -6) {};
		\node [style=none] (74) at (15.5, -8) {};
		\node [style=none] (75) at (-4.25, -7) {$+$};
		\node [style=none] (76) at (8.5, -7) {$+$};
		\node [style=none] (77) at (3.25, -7) {$+$};
		\node [style=none] (78) at (11.75, -7) {$+$};
		\node [style=none] (79) at (14.25, -7) {$+$};
		\node [style=text-labels] (80) at (18, -7) {};
		\node [style=text-labels] (81) at (18, -7) {$=$};
		\node [style=none] (82) at (19, -7) {$0$};
		\node [style=surrounded] (84) at (-0.075, -6.55) {$1$};
		\node [style=surrounded] (85) at (4.5, -7.375) {$1$};
		\node [style=surrounded] (86) at (13, -7.25) {$0$};
		\node [style=surrounded] (88) at (15.5, -7.25) {$0$};
		\node [style=surrounded] (100) at (1.25, -3.75) {$0$};
		\node [style=none] (101) at (1.25, -4.5) {};
		\node [style=surrounded] (102) at (1.25, -3) {$1$};
		\node [style=none] (103) at (1.25, -2.25) {};
		\node [style=surrounded] (104) at (3.5, -3.75) {$0$};
		\node [style=none] (105) at (3.5, -4.5) {};
		\node [style=surrounded] (106) at (3.5, -3) {$1$};
		\node [style=none] (107) at (3.5, -2.25) {};
		\node [style=none] (108) at (2.375, -3.5) {$+$};
		\node [style=surrounded] (109) at (5.75, -3.5) {$0$};
		\node [style=none] (110) at (5.75, -4.5) {};
		\node [style=none] (111) at (5.75, -2.25) {};
		\node [style=surrounded] (112) at (6.25, -3) {$1$};
		\node [style=surrounded] (113) at (8.75, -3.5) {$0$};
		\node [style=none] (114) at (8.75, -4.5) {};
		\node [style=none] (115) at (8.75, -2.25) {};
		\node [style=surrounded] (116) at (8.25, -3) {$1$};
		\node [style=none] (117) at (4.625, -3.5) {$+$};
		\node [style=none] (118) at (7.25, -3.5) {$+$};
		\node [style=surrounded] (119) at (9.5, -6.75) {$1$};
		\node [style=surrounded] (120) at (16.25, -6.75) {$1$};
		\node [style=surrounded] (121) at (13, -6.5) {$1$};
		\node [style=none] (122) at (-3.75, -6) {};
		\node [style=surrounded] (123) at (-3.4, -6.55) {$0$};
		\node [style=surrounded] (124) at (-3, -7.25) {$1$};
		\node [style=none] (125) at (-2.25, -6) {};
		\node [style=none] (126) at (-3, -8) {};
		\node [style=none] (127) at (-1.75, -7) {$+$};
		\node [style=none] (128) at (1.25, -6) {};
		\node [style=surrounded] (129) at (2, -7.25) {$1$};
		\node [style=none] (130) at (2.75, -6) {};
		\node [style=none] (131) at (2, -8) {};
		\node [style=surrounded] (132) at (2.425, -6.55) {$0$};
		\node [style=none] (133) at (0.75, -7) {$+$};
		\node [style=none] (134) at (6.25, -6) {};
		\node [style=surrounded] (135) at (7, -6.7) {$1$};
		\node [style=none] (136) at (7.75, -6) {};
		\node [style=none] (137) at (7, -8) {};
		\node [style=surrounded] (138) at (7, -7.375) {$0$};
		\node [style=none] (139) at (5.75, -7) {$+$};
		\node [style=surrounded] (141) at (-13, 0.5) {$1$};
		\node [style=none] (142) at (-13, -0.5) {};
		\node [style=none] (143) at (-13, -2.5) {};
		\node [style=surrounded] (144) at (-13, -3.5) {$1$};
		\node [style=none] (145) at (-13, -4.5) {};
		\node [style=none] (146) at (-13.5, -6) {};
		\node [style=surrounded] (147) at (-13, -7) {$1$};
		\node [style=none] (148) at (-13, -8) {};
		\node [style=none] (149) at (-12.5, -6) {};
	\end{pgfonlayer}
	\begin{pgfonlayer}{edgelayer}
		\draw (1) to (2.center);
		\draw (31.center) to (32);
		\draw (32) to (33);
		\draw (33) to (35.center);
		\draw (33) to (36.center);
		\draw (39) to (41.center);
		\draw (39) to (42.center);
		\draw (45) to (48.center);
		\draw (57) to (59.center);
		\draw (37.center) to (39);
		\draw (47.center) to (84);
		\draw (84) to (45);
		\draw (45) to (43.center);
		\draw (55.center) to (57);
		\draw (57) to (85);
		\draw (85) to (60.center);
		\draw (65.center) to (86);
		\draw (86) to (66.center);
		\draw (88) to (74.center);
		\draw (69.center) to (88);
		\draw (5) to (1);
		\draw (100) to (101.center);
		\draw (102) to (100);
		\draw (103.center) to (102);
		\draw (104) to (105.center);
		\draw (106) to (104);
		\draw (107.center) to (106);
		\draw (109) to (110.center);
		\draw (111.center) to (109);
		\draw (112) to (109);
		\draw (113) to (114.center);
		\draw (115.center) to (113);
		\draw (116) to (113);
		\draw (61.center) to (86);
		\draw (73.center) to (88);
		\draw (86) to (121);
		\draw (88) to (120);
		\draw (119) to (39);
		\draw (122.center) to (123);
		\draw (123) to (124);
		\draw (124) to (125.center);
		\draw (124) to (126.center);
		\draw (129) to (131.center);
		\draw (130.center) to (132);
		\draw (132) to (129);
		\draw (129) to (128.center);
		\draw (135) to (136.center);
		\draw (134.center) to (135);
		\draw (135) to (138);
		\draw (138) to (137.center);
		\draw (141) to (142.center);
		\draw (143.center) to (144);
		\draw (144) to (145.center);
		\draw (146.center) to (147);
		\draw (147) to (148.center);
		\draw (149.center) to (147);
	\end{pgfonlayer}
\end{tikzpicture}\\
\end{center}

We note that some of the terms in the $\Ainf$ relations include the curvature popping off the weight $w$ vertex.

\begin{definition}[Weighted $\Ainf$-module] Given a weighted $\Ainf$-algebra $\Alg$ over a ring $\kk$,
			a weighted (right) $\Ainf$-module $\fModule$ over $\Alg$ is a right $\kk$-module $M$; along with $n$-input module operations
			\[ m_n^w: M \otimes_\kk \overbrace{A \otimes_\kk \dots \otimes_\kk A}^n \to M \]
			for all $n,w \in \ZZ_{\geq0}$ such that the maps $m^w_{1+n}$ are $\kk$-module maps, and they satisfy the structure relations
			\begin{align*}
			\sum_{\substack{i+j+l = n \\ w_1 + w_2 = w}}
			&m_{2 + i + l}^{w_1}\left(\x, a_1, \cdots, a_{i}, \mu_{j}^{w_2}(a_{i+1}, \dots, a_{i+j}),
			a_{i+j+1}, \dots, a_{i+j+l}\right) \\
			&+ \sum_{\substack{i + j = n \\ w_1 + w_2 = w}}
			m^{w_1}_{1+j}\left(m^{w_2}_{1+i}\left(\x, a_1, \cdots, a_{i}\right), a_{i+1}, \dots, a_{i+j}\right) = 0
			\end{align*}
			for all $n, w \geq 0$ and $\x \in M, a_1, \dots, a_n \in \Alg$. Here, $w$ is called the \emph{weight} of the operation $m^w_{1+n}$.
			
			A weighted $\Ainf$-module $\fModule$ over a strictly unital weighted $\Ainf$-algebra $\Alg$ is \emph{strictly unital} if $m_2^0(\x, 1) = \x$ for all $\x \in M$; and if some $a_i = 1$, then $m_n^w(\x, a_1, \dots, a_n) = 0$ for all $(n,w) \neq (2,0)$.
\end{definition}

There is a graphical representation of the weighted module operations, entirely  like that of weighted algebra operations, except the left-most edge is marked to represent module elements. The $\Ainf$ structure relations can then be constructed similarly.

\begin{definition}[Weighted type D module]
			Given a weighted $\Ainf$-algebra $\Alg$ over a ring $\kk$,
			a weighted (left) type D module $\mathcal{P}$ over $\Alg$ consists of a $\kk$-module $P$ and a $\kk$-module map $\delta^1:P \to \Alg \tensor_{\kk} P$ satisfying the structure relation
			\[ \sum_{n, w \geq 0} (\mu^w_n \tensor \mathbb{I}) \circ \delta^n = 0. \]
			Here, $\delta^n: P \to \Alg^{\tensor_{\kk} n} \tensor_{\kk} P$ is the map obtained by iterating $\delta^1$ $n$ times.
\end{definition}

\subsection{The torus algebra}
\label{sec:torus-algebra}

This section recalls the construction of the torus algebra $\Algm$, introduced
by Lipshitz, Ozsv\'ath, and Thurston in \cite{LOT:torus-alg}. 
The construction of the torus algebra $\Algm$ begins with the algebra
$\widehat \Alg = \Alg(T^2)$, which is used to define $\HFa$
for 3-manifolds with torus boundary \cite{LOT1}. $\widehat \Alg$ has a natural
extension to an algebra $\AlgmAs$. Over $\FF_2$, $\widehat \Alg$ has basis given by
the idempotents $\iota_0, \iota_1$, and elements of the form 
$\rho_{i} \rho_{i+1}
\cdots \rho_{i+k} \eqdef \rho_{i,\dots,i+k}$ 
for $i \in \ZZ/4\ZZ$ and $k \in \ZZ_{\geq 0}$. These are called
the \emph{basic elements} of the algebra. The non-idempotent elements $\rho_{i, \dots, i_k}$ are called the \emph{Reeb elements}.

$\AlgmAs$
can be seen as the path algebra over $\FF_2$ of
\begin{center}
	\begin{tikzpicture}
		\node at (0,0) (iota0) {$\iota_0$};
		\node at (4,0) (iota1) {$\iota_1$};
		\draw[->, bend right=15] (iota0) to node[above]{$\rho_3$} (iota1);
		\draw[->, bend left=45] (iota0) to node[above]{$\rho_1$} (iota1);
		\draw[->, bend right=15] (iota1) to node[above]{$\rho_2$} (iota0);
		\draw[->, bend left=45] (iota1) to node[above]{$\rho_4$} (iota0);
	\end{tikzpicture}
\end{center} 
with the relations that any of decreasing products
$\rho_2 \rho_1, \rho_3 \rho_2, \rho_4 \rho_3$, and $\rho_1 \rho_4$ vanish.
We note that the quotient of $\AlgmAs$ by the ideal generated by $\rho_4$ gives us the algebra $\widehat A$.

The algebra $\widehat A$ keeps track of the Reeb chords associated to pseudoholomorphic disks that do not cross a basepoint $z$ near the Reeb chord $\rho_4$. An algebra that performs a similar bookkeeping for disks that are allowed to cross the basepoint need weighted algebra operations to account for boundary degenerations.

For the remainder of this paper, our ground ring $\kk$ will be the 
sub-algebra generated by $\iota_0$ and
$\iota_1$ adjoined with the variable $U$. The torus algebra $\Algm$ is given by a (curved) weighted $\Ainf$-algebra structure on
$\AlgmAs[U]$ over the ground ring $\kk$.

There is non-zero curvature for $\Algm$, all concentrated in weight $1$. It is given by the sum of all length four Reeb chords.
\[\mu_0^1 = \rho_{1234} + \rho_{2341} + \rho_{3412} + \rho_{4123}.\]
The differential $\mu_1^0$ and the higher weight $1$-input operations are all trivial.

As the operations $\mu_n^w$ on $\Algm$ are $U$-multilinear, it suffices to define them
on basic algebra elements.
For two inputs, the operation $\mu_2^0$ is given by multiplication. For any $(n, w) \neq (2, 0)$, if any of the $a_i$ are idempotents,
$\mu_n^w(a_1, \dots, a_n) = 0$. That is, $\Algm$ is strictly unital.

For Reeb elements $a_1, \dots, a_n$, the operations $\mu_n^w(a_1, \dots, a_n)$ are defined
by counting certain classes of planar graphs that have to them associated the chord sequence $a_1 \otimes \cdots \otimes a_n$ and are of weight $w$. We call these tiling patterns. In Section~\ref{ch:cfa}, we make use of similar tiling patterns to define our type A modules. We then distinguish the two as algebra or module tiling patterns. 

\subsubsection{Planar graphs and tiling patterns}

Our tiling patterns are connected planar graphs
$\Gamma$, and come with an embedding into the disk $\bD$. The leaves of $\Gamma$ are precisely where it intersects the boundary of the disk. The tiling patterns come equipped with a choice of a leaf as the root. As a convention, we use \emph{vertex} to refer to vertices of our planar graphs that are not leaves. In particular, the root is a not a vertex.

\begin{definition}
	A \emph{centered tiling pattern} is a connected,
	rooted, planar $\Gamma$ such that every vertex in $\Gamma$ is
	$4$-valent, and every component of $\bD \setminus \Gamma$ either meets
	$\bdy \bD$ or is bound by four edges. We call components of $\bD \setminus \Gamma$ that do not meet
	the boundary \emph{internal faces}, and the boundaries of internal faces \emph{short cycles}.
	
	A tiling pattern comes equipped with a valid labeling around the vertices.
	For any vertex $v$, let $Q_v$ be the set of faces that meet
	$v$. A \emph{valid labeling} on $\Gamma$ consists of a
	labeling $\Lambda_v:
	Q_v \to \{1,2,3,4\}$ for every vertex $v$ such that:
	\begin{itemize}
		\item Let $Q_v = \{r_1, r_2, r_3, r_4\}$ be the set of faces around $v$, labelled in clockwise order. Then, up to a cyclic reordering, $(\Lambda_v(r_1), \Lambda_v(r_2), \Lambda_v(r_3), \Lambda_v(r_4))$ is a cyclic reordering of $(1,2,3,4)$.
		\item If $e : v_1 \to v_2$ is an edge oriented such that a face $f$
		is to the right of $e$, then $\Lambda_{v_1}(f) + 1 \equiv \Lambda_{v_2}(f) \mod 4$.
	\end{itemize}
\end{definition}

Along with the centered tiling patterns, we considered an extended class of tiling patterns called the left-extended and right-extended patterns. These are obtained from centered tiling patterns by adding $2$-valent vertices that respect the labelings.

\begin{definition}
	Consider any centered tiling pattern. There is a distinguished
	edge $e$ adjacent to the root. Consider it with the orientation pointing
	towards the root. All \emph{left-extended} (respectively \emph{right-extended})
	\emph{tiling patterns} are obtained by adding $k$ new
	$2$-valent vertices along
	$e$, and labeling the new vertices
	to the left (respectively right) of $e$ such that the property
	\begin{itemize}
	\item If $e : v_1 \to v_2$ is an edge oriented such that a face $f$
		is to the right of $e$, then $\Lambda_{v_1}(f) + 1 \equiv \Lambda_{v_2}(f) \mod 4$.
	\end{itemize}
	is maintained.
\end{definition}

\begin{figure}
	\centering
	\begin{tikzpicture}[scale=1.25]
	\begin{pgfonlayer}{nodelayer}
		\node [style=none] (0) at (-4, 0) {};
		\node [style=none] (1) at (4, 0) {};
		\node [style=bigdot] (2) at (-3.175, 1.875) {};
		\node [style=none] (3) at (2.625, 2.25) {};
		\node [style=none] (4) at (-1.725, 2.65) {};
		\node [style=none] (5) at (1.275, 2.8) {};
		\node [style=none] (6) at (2.125, -2.5) {};
		\node [style=none] (7) at (-1.95, -2.575) {};
		\node [style=none] (8) at (-3.35, -1.7) {};
		\node [style=none] (9) at (3.625, -1.325) {};
		\node [style=dot] (10) at (-1.625, 1.6) {};
		\node [style=dot] (11) at (1.35, 1.775) {};
		\node [style=dot] (12) at (-1.75, -1.4) {};
		\node [style=dot] (13) at (1.825, -1.2) {};
		\node [style=number-labels] (14) at (-2, 2.125) {$1$};
		\node [style=number-labels] (15) at (-1.25, 2.175) {$2$};
		\node [style=number-labels] (16) at (-1.25, 1.25) {$3$};
		\node [style=number-labels] (17) at (-2, 1.25) {$4$};
		\node [style=number-labels] (18) at (1, 2.25) {$1$};
		\node [style=number-labels] (19) at (1.75, 2.25) {$2$};
		\node [style=number-labels] (20) at (1.75, 1.5) {$3$};
		\node [style=number-labels] (21) at (1, 1.425) {$4$};
		\node [style=number-labels] (22) at (1.425, -0.8) {$1$};
		\node [style=number-labels] (23) at (2.25, -0.8) {$2$};
		\node [style=number-labels] (24) at (2.25, -1.5) {$3$};
		\node [style=number-labels] (25) at (1.5, -1.5) {$4$};
		\node [style=number-labels] (26) at (-2.15, -1) {$1$};
		\node [style=number-labels] (27) at (-1.275, -1.025) {$2$};
		\node [style=number-labels] (28) at (-1.325, -1.875) {$3$};
		\node [style=number-labels] (29) at (-2.15, -1.875) {$4$};
		\node [style=none] (30) at (13.5, 0) {};
		\node [style=bigdot] (31) at (7.5, 0) {};
		\node [style=none] (32) at (10.5, 1.775) {};
		\node [style=none] (33) at (10.5, -1.775) {};
		\node [style=dot] (34) at (10.5, 0) {};
		\node [style=dot] (35) at (9.25, 0) {};
		\node [style=number-labels] (36) at (9.25, 0.5) {$2$};
		\node [style=number-labels] (37) at (10, 0.5) {$1$};
		\node [style=number-labels] (38) at (11, 0.5) {$2$};
		\node [style=number-labels] (39) at (11, -0.5) {$3$};
		\node [style=number-labels] (40) at (10, -0.5) {$4$};
		\node [style=dot] (41) at (8.5, 0) {};
		\node [style=number-labels] (42) at (8.5, 0.5) {$3$};
	\end{pgfonlayer}
	\begin{pgfonlayer}{edgelayer}
		\draw [bend left=90, looseness=1.25] (0.center) to (1.center);
		\draw [bend right=90, looseness=1.25] (0.center) to (1.center);
		\draw [in=-150, out=-15, looseness=0.75] (2) to (3.center);
		\draw [bend left, looseness=0.25] (8.center) to (9.center);
		\draw [bend left=15, looseness=0.50] (4.center) to (7.center);
		\draw [bend right=15, looseness=0.25] (5.center) to (6.center);
		\draw [bend left=90] (31) to (30.center);
		\draw [bend right=90] (31) to (30.center);
		\draw [in=180, out=0] (31) to (30.center);
		\draw (32.center) to (33.center);
	\end{pgfonlayer}
\end{tikzpicture}
	\caption[Examples of tiling patterns.]
	{\textbf{Examples of tiling patterns.} The root on the boundary is marked with the black dot. On the left is a centered tiling pattern. On the right, a left-extended one.}
	\label{fig:alg-tiling-pattern}
\end{figure}
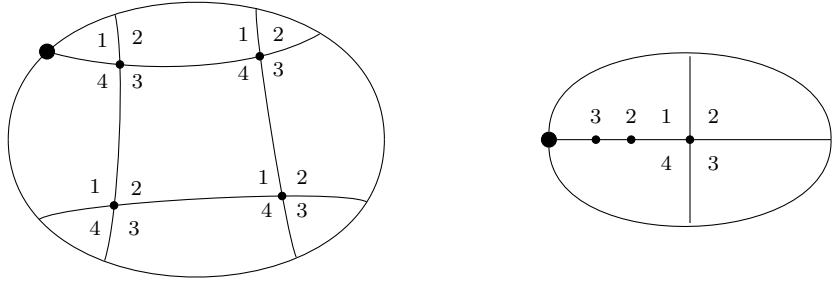

To each tiling pattern $\Gamma$, we
associate a weight $w \in \ZZ_{\geq0}$, a chord sequence $a_1 \tensor \cdots \otimes a_n$ where each $a_i \in \Algm$,
and an output $o \in \Algm$. The operations on our weighted $\Ainf$-algebra count such tiling patterns. That is, the operation $\mu_n^w(a_1, \dots, a_n)$ is the sum of the output $o$ of all tiling patterns $\Gamma$ that have weight $w$ and chord sequence $a_1 \otimes \cdots \otimes a_n$.

\begin{definition}
	Let $\Gamma$ be a tiling pattern. The \emph{weight} $w$ of $\Gamma$ is defined as the number of short cycles in $\Gamma$.
	Equivalently,
	the first homology $H_1(\Gamma; \ZZ)$ of the underlying planar graph has dimension $w$.
\end{definition}

As an example, the tiling patterns on the left in Figure~\ref{fig:alg-tiling-pattern} has weight 1. The tiling pattern on the right has weight 0.

\begin{definition} \label{def:alg-chord-seq} Let $\Gamma$ be a tiling pattern.
	The \emph{chord sequence} of $\Gamma$ is given by reading the labels visible
	along the faces while traversing $\bdy \bD$ (with its chosen orientation), starting
	at the root.
	That is, let $f_1, \dots, f_k$ be the faces in order along the boundary.
	For any $1 \leq i \leq k$, let $v^i_1, \dots, v^i_{j_i}$ be the vertices in order
	along the internal boundary of the face $f_i$. Then the chord sequence
	of $\Gamma$ is defined as
	\[\left(\prod_{l = 1}^{j_1} \rho_{\Lambda_{v^1_l}(f_1)}\right)\tensor \dots \tensor 
	\left(\prod_{l = 1}^{j_i} \rho_{\Lambda_{v^i_l}(f_i)}\right)\tensor \dots \tensor
	\left(\prod_{l = 1}^{j_k} \rho_{\Lambda_{v^k_l}(f_k)}\right).\]
\end{definition}

As an example, the tiling patterns on the left in Figure~\ref{fig:alg-tiling-pattern} has associated chord sequence $\rho_{41} \tensor \rho_4 \tensor \rho_{34} \tensor \rho_3 \tensor \rho_{23} \tensor \rho_2 \tensor \rho_{12} \tensor \rho_1$. The tiling pattern on the right has associated chord sequence $\rho_4 \tensor \rho_3 \tensor \rho_2 \tensor \rho_{123}$.

\begin{definition}
	Let $\Gamma$ be a tiling pattern.
	Let $d$ be the number of 4-valent vertices in $\Gamma$. Order
	the faces and vertices as in Definition~\ref{def:alg-chord-seq}.
	
	If $\Gamma$ is a centered tiling pattern, the \emph{output} is $U^d$ times
	the left idempotent of the first label; namely $\rho_{\Lambda_{v^1_l}(f_1)}$.
	Equivalently, it is
	$U^d$ times the right idempotent of the last label; namely
	$\rho_{\Lambda_{v^k_{j_k}}(f_k)}$.
	
	If $\Gamma$ is an extended tiling pattern, the \emph{output} is $U^d$ times
	the product of the labels along the 2-valent vertices.
\end{definition}

\sloppy So, the operation corresponding to the tiling pattern on the left in Figure~\ref{fig:alg-tiling-pattern} is $\mu^1_8(\rho_{41}, \rho_4,  \rho_{34}, \rho_3, \rho_{23}, \rho_2, \rho_{12}, \rho_1) = U^4 \iota_1$. The operation corresponding to the tiling pattern on the right is $\mu^0_4(\rho_4, \rho_3, \rho_2, \rho_{123}) = U \rho_{23}$.

\subsection{The type A module} 
Let $Y$ be an oriented $3$-manifold with torus boundary, along with an identification of its boundary with the torus. Such a $3$-manifold can be represented by a bordered Heegaard diagram.

\begin{definition}[Bordered Heegaard diagram]
	A \emph{bordered Heegaard diagram} for a bordered $3$-manifold with torus boundary is given by $\HD = (\Sigma, \alphas, \betas, z)$, where
	\begin{itemize}
	\item $\Sigma$ is a compact, oriented surface of genus $g$ with one boundary component;
	\item $\alphas = \{\alpha_1^c, \alpha_2^c, \alpha_1, \dots, \alpha_{g-1}\}$ are pairwise disjoint curves in $\Sigma$, where the $\alpha_1^c$ and $\alpha_2^c$ are arcs in $\Sigma$ with boundary on $\bdy \Sigma$, and the $\alpha_1, \dots, \alpha_{g-1}$ lie in the interior of $\Sigma$;
	\item $\betas = \{\beta_1, \dots, \beta_g\}$ are pairwise disjoint circles in $\Sigma$; and
	\item $z$ is a basepoint in $\bdy \Sigma \setminus (\alphas \cap \Sigma)$; such that
	\end{itemize}
	the $\alphas$ and $\betas$ intersect transversely, and $\Sigma \setminus \alphas$ and $\Sigma \setminus \betas$ are both connected.
\end{definition}
Construction~4.6 in \cite{LOT1} spells out how such a bordered Heegaard diagram gives rise to a bordered $3$-manifold. Lemma~4.9 of \cite{LOT1} demonstrates that every bordered $3$-manifold with torus boundary arises from a bordered Heegaard diagram.

The two $\alpha$-arcs divide the boundary $\bdy \Sigma$ into four segments, with the basepoint $z$ in one segment. With respect to the boundary orientation, the segments are labeled with the Reeb chords $\rho_4, \rho_3, \rho_2,$ and $\rho_1$, with $z$ in the region labeled $\rho_4$. The Reeb elements of the torus algebra track asymptotics of pseudoholomorphic curves corresponding to the Reeb chords on $\bdy \Sigma$. As a convention, with respect to the boundary orientation, the first $\alpha$-arc after the basepoint $z$ is $\alpha_1^c$.

Let $\Gen(\HD)$ be the set of $g$-tuples of points $\x = (x_1, \dots, x_g)$ on $\Sigma$ such that exactly one $x_i$ lies on each $\beta$-circle; exactly one $x_i$ lies on each $\alpha$-circle; and one $x_i$ lies on one of the two $\alpha$-arcs. Recall that $\Gen(\HD)$ generates the type A module $\CFAa(\HD)$, and will generators the module $\CFAm(\HD)$.

The operations on the type A modules count pseudoholomorphic curves in $\Sigma \times [0, 1] \times \RR$ for appropriately chosen families of almost complex structures $J$, and with certain boundary conditions. Each curve from a generator $\x$ to $\y \in \Gen(\HD)$ corresponds to a homology class in $\pi_2(\x, \y)$. A \emph{domain} $B \in \pi_2(\x, \y)$ is a linear combination of the connected components of $\Sigma \setminus (\alphas \cup \betas)$, such that its boundary lying in $\betas$ is a linear chain from $\x$ to $\y$. That is, if $\bdy^\beta B$ denotes the part of $\bdy B$ lying in $\betas$, then $\bdy (\bdy^\beta B) = \x - \y$. A domain is \emph{positive} if it is a positive linear combination, and only such domains admit pseudoholomorphic representatives.

We recall a few more definitions before sketching a construction of the type A module $\CFAm(\HD)$. A \emph{periodic domain} is a domain in $\pi_2(\x, \x)$ such that its multiplicity at the basepoint $z$ is zero. In notation, $n_z(B) = 0$. A \emph{provincial domain} is a domain such that its boundary does not lie in $\bdy \Sigma$. In notation, $\bdy^\bdy B = 0$. A bordered Heegaard diagram is called \emph{admissible} if every periodic domain has both positive and negative coefficients, and \emph{provincially admissible} if every provincial periodic domain has both positive and negative coefficients.

Consider a provincially admissible bordered Heegaard diagram $\HD = (\Sigma, \alphas, \betas, z)$ for a bordered 3-manifold with torus boundary. 

Let $\vec{a}$ be a sequence of Reeb elements and $w$ a non-negative integer.
For generators $\x$ and $\y$,
and a homology class $B \in \pi_2(\x, \y)$ between them, let
$\cM^B(\x, \y; \vec{a}; w)$ denote the moduli space of embedded pseudoholomorphic
maps
$$u: (S, \bdy S) \to (\Sigma \times [0,1] \times \RR, (\alphas \times \{1\} \times \RR) \cup (\betas \times \{0\} \times \RR))$$
that are asymptotic to $\x \times [0,1]$ at $-\infty$, to
$\y \times [0, 1]$ at $+\infty$, and to the chord
sequence $\vec{a}$ and $w$ simple Reeb orbits at \emph{east infinity}. A simple Reeb orbit is a map $S^1 \to \bdy \Sigma \times [0,1] \times \RR$ which maps to $[0, 1] \times \RR$ by a constant map, and winds around $\bdy \RR$ once. The precise technical constraints are considered in \cite[Section~3]{LOT:torus-mod}. Let
$\ind(B,\vec{a},w)$ denote one more than the expected dimension of this moduli space, such that when $\ind(B,\vec{a},w) = 1$ then $\cM^B(\x, \y; \vec{a}; w)$ has finitely many points.

The construction of $\CFAm(\HD)$ proceeds by first constructing a
weighted $\Ainf$-module $\CFAm_{nu}(\HD)$. The weighted $\Ainf$-module $\CFAm_{nu}(\HD)$ as a left $\kk$-module is generated by $\Gen(\HD)$. The idempotents $\iota_0$ and $\iota_1$ correspond to $\alpha^a_1$ and $\alpha^a_2$ such that $m^0_2(\x, \iota_i) = x$ if $\x \cap \alpha^a_{i+1} \neq \emptyset$, and $m^0_2(\x, \iota_i) = 0$ otherwise. Further, if any $a_i$ is an idempotent, $m^w_{1+n}(\x, a_1, \dots, a_n) = 0$ whenever $(n, w) \neq (2,0)$; that is, $\CFAm_{nu}(\HD)$ is strictly unital.
For a sequence of Reeb elements $a_1, \dots, a_n$, the weighted operations are defined by counting pseudoholomorphic curves.
$$m^w_{1+n}(\x, a_1, \dots, a_n) =
\sum_{\y \in \mathfrak{S}(\HD)} \sum_{\substack{B \in \pi_2(\x,\y)\\ \ind(B,\vec{a}, w) = 1}}
\#\left({\cM}^B(\x, \y; a_1, \dots, a_n; w)\right) U^{n_z(B)}  \y.$$

 The operations on $\CFAm_{nu}(\HD)$ are not $U$-equivariant. That is, $\CFAm_{nu}(\HD)$ is defined over the ground ring $\FF_2$. The complications arise in choosing appropriate families of almost complex structures \cite[Remark~6.29]{LOT:torus-mod}. 
To construct $\CFAm(\HD)$, Lipshitz, Ozsv\'ath, and Thurston extend $\CFAm_{nu}(\HD)$ to a weighted $\Ainf$-module $\CFAm_g(\HD)$ defined over a
larger algebra
which is quasi-isomorphic to $\Algm$.
Restricting scalars to $\Algm$ produces $\CFAm(\HD)$. The weighted
$\Ainf$-module $\CFAm(\HD)$ has $U$-equivariant module operations and is
homotopy equivalent to $\CFAm_{nu}(\HD)$.

In this paper, we focus on the case when genus $g = 1$.  In this case, there is a choice of complex structure that avoids the complications of \cite[Remark~6.29]{LOT:torus-mod} (c.f. \cite[Section~4.2]{LOT:torus-mod}).  That is, we avoid needing to pass to $\CFAm_g(\HD)$.

\subsection{More algebraic underpinnings}
\label{sec:more-algebra}

The type D module in the \emph{minus} bordered Floer theory is defined by  tensoring the type A module with a \emph{dualizing} bimodule. To do so, we must define the tensor products of weighted $\Ainf$-algebras and modules. We make use of these notions in Section~\ref{ch:algprop}, and our definitions are merely precise enough to carry us through. A deeper dive is in \cite{LOT:abstract}, while \cite[Sections~1.3~and~1.5]{LOT:torus-mod} present a synopsis relevant to the constructions in bordered Floer theory. 

In this section, we work over $\FF_2[\Yvar_1, \Yvar_2]$. In our constructions, $\Yvar_1$ will act by 1, and $\Yvar_2$ will act by $U$. Let $\Alg$ and $\mathcal{B}$ be weighted $\Ainf$ algebra over ground rings $\kk_1$ and $\kk_2$ over $\FF_2[\Yvar_1, \Yvar_2]$. 

Given a weighted type A module $\fModule$ and weighted type D structure $D$ over $\Alg$, we can define their \emph{box tensor product} $\fModule \boxtimes_{\Alg} D$. This is a chain complex, with the underlying module $M \tensor_{\FF_2[\Yvar_1, \Yvar_2]} D$ and differential
\[ \partial = \sum_{n, w \geq 0} (m^w_{1+n} \tensor \mathbb{I}_{\mathcal{P}}) \circ (\mathbb{I}_{\fModule} \tensor \delta^n). \]

In Section~\ref{sec:alg-underpinnings}, we made use of weighted trees to represent operations on our algebras and modules. The root of the tree is its \emph{output}, and a certain subset of the leaves are marked as the \emph{inputs}. Every other vertex is \emph{internal}, and is marked with a non-negative integral weight. The \emph{dimension} of a weight $w$ $n$-input tree $T$ is defined to be $\dim(T) = n + 2w - v - 1$, where $v$ is the number of internal vertices. Such a tree induces a map $\mu(T): A^{\tensor n} \to A$, as we saw in Section~\ref{sec:alg-underpinnings}. The trees where the left-most leaf is an input are called the \emph{module trees}, and induce maps on the weighted $\Ainf$-modules. A weighted tree is said to be stable if any internal vertex with valence less than three has positive weight.

A weighted algebra diagonal $\AsDiag$ consists, for all $n, w \geq 0$ (and $n + 2w \geq 2$), of linear combinations $\TrDiag^{n,w}$ of pairs of stably-weighted $n$-input trees with weight at most $w$ and dimension $n + 2w - 2$, and marked with appropriate powers of $\Yvar_1$ and $\Yvar_2$. These linear combinations satisfy additional constraints of \emph{compatibility} and \emph{non-degeneracy} \cite[Definition~6.6]{LOT:abstract}. It is these pairs of trees that tell us what the higher operations on a tensor product of weighted $\Ainf$-algebras should be. That is, the tensor product $\Alg \tensor_{\AsDiag} \mathcal{B}$ has higher operations
\[ \mu^w_n = \sum_{n_{S,T}(S \tensor T) \in \TrDiag^{n,w}} n_{S,T} \mu_{\Alg}(T) \tensor \mu_{\mathcal B}(T), \]
	with the underlying module $A \tensor_{\FF_2[\Yvar_1, \Yvar_2]} B$.
	In \cite{LOT:abstract} they show that a weighted algebra diagonal exists, and this construction gives $\Alg \tensor_{\AsDiag} \mathcal{B}$ the structure of a weighted $\Ainf$-algebra. We will fix as our weighted algebra diagonal the diagonal they construct, and whose first few terms are listed in \cite[Appendix~A]{LOT:abstract}.
	
	\begin{definition}[Weighted type DD bimodules]
	A weighted type DD bimodule $P$ over $(\Alg, \mathcal B)$ is a weighted type D module over $\Alg \tensor_{\AsDiag} \mathcal B$.
	\end{definition}
	
	There is a similar definition of a weighted module diagonal that allows us to define the tensor product of weighted $\Ainf$-modules.
	
	To define the tensor product of a weighted $\Ainf$-module with a type DD bimodule, we require pairs of stably-weighted \emph{module} trees and trees. A weighted module diagonal primitive consists, for all $n \geq 1$ and $w \geq 0$, of linear combinations $\TrMPrim^{n,w}$ of pairs of stably-weighted $n$-input module trees and $n-1$ input trees, with weight at most $w$ and dimension $n + 2w - 2$. These too are marked with appropriate powers of $\Yvar_1$ and $\Yvar_2$, and satisfy certain compatibility and non-degeneracy constraints (that depend on the weighted module diagonal $\AsDiag$) \cite[Definition~6.29]{LOT:abstract}. \cite{LOT:abstract} construct such a weighted module diagonal primitive, and we list a few of the weight $0$ terms we rely upon in Section~\ref{ch:algprop}.
	
\begin{center}
  \includegraphics[scale=.75]{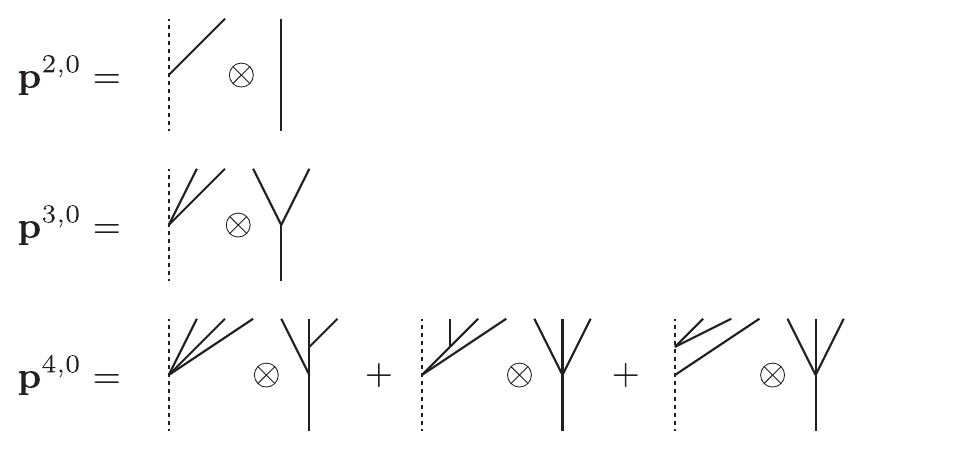}
\end{center}

Given a weighted type A module $\fModule$ over $\Alg$ and a type DD bimodule $P$ over $(\Alg, \mathcal B)$, we define the \emph{box tensor product} $\fModule \boxtimes_{\Alg} P$ to be the weighted type D structure with underlying module $M \tensor_{\kk} P$ and differential
\begin{equation}
	\label{eqn:m-dd-tensor}
	\delta^1(x \tensor y) = m^0_1(x) \tensor y \tensor 1 + \sum_{\substack{n_{S,T}(S \tensor T) \in \TrMPrim^{n,w} \\ n,w \geq 0}} n_{S,T}(m^M(S) \tensor \Id_P \tensor \mu^{\mathcal{B}^{\op}} (T^\op)) \circ (\Id_{M} \tensor \delta_P).
\end{equation}

To ensure the structure relations of the type D and type DD modules make sense---and to further ensure the box tensor product is well-defined---we require boundedness properties on our type A and type D modules, and our type DD bimodules. We defer the discussion of the boundedness properties for the type D and type DD bimodules to \cite[Section~7.2]{LOT:torus-mod}, and discuss the case of the type A modules. 

An $\Ainf$-algebra $\Alg$ is said to admit a filtration $A = \Filt^0 A \supset \Filt^1 A \supset \dots$ if the operations $\mu^w_n$ map
$\Filt^{m_1} A \tensor \cdots \tensor \Filt^{m_n} A$ to
$\Filt^{m_1 + \cdots + m_n + kw} A$ for a positive number $k$. $\Alg$ is complete with respect to $\Filt$ if it is complete as vector spaces with respect to the filtration. The definition extends to that of a filtration of a weighted $\Ainf$-module over a filtered, weighted $\Ainf$-algebra.

$\Alg$ is \emph{bonsai} if there exists an $N$ such that for any stably weighted tree $T$ with dimension $\dim(T) > N$, the operation $\mu(T) = 0$. $\Alg$ is \emph{filtered bonsai} with respect to a filtration $\Filt$ if it is complete with respect to $\Filt$ and the quotients $A/\Filt^m A$ are bonsai. A filtered, weighted $\Ainf$-module $\fModule$ over a filtered bonsai $\Alg$ is filtered bonsai if $\fModule$ is complete with respect to $\Filt$ and the quotients $M/\Filt^m M$ are bonsai.

For our type D and type DD modules, it will suffice that the type A modules we construct are filtered bonsai with respect to a chosen filtration. Proposition~7.7 of \cite{LOT:torus-mod} provides a summary of the boundedness constraints on the modules at play. These notions are further studied in Section~9 of \cite{LOT:abstract}.

\subsection{The type D module}

The dualizing bimodule $\CFDDm(\Id)$ is the type DD bimodule over $(\Algm^{U=1}, \Algm)$ with generators $\iota_0 \otimes \iota_0$ and $\iota_1 \otimes \iota_1$, and differential $\delta^1$ given by
\begin{equation}
	\begin{aligned}
	\delta^1(\iota_0 \tensor \iota_0) &= (\rho_1 \tensor \rho_3 + \rho_3 \tensor \rho_1 + \rho_{123} \tensor \rho_{123} + \rho_{341} \tensor \rho_{341}) \tensor (\iota_1 \tensor \iota_1) \\
	\delta^1(\iota_1 \tensor \iota_1) &= (\rho_2 \tensor \rho_2 + \rho_4 \tensor \rho_4 + \rho_{234} \tensor \rho_{412} + \rho_{412} \tensor \rho_{234}) \tensor (\iota_0 \tensor \iota_0) \\
	\end{aligned}
\end{equation}

For a provincially admissible bordered Heegaard diagram $\HD$ for
a bordered 3-manifold with torus boundary, Lipshitz, Ozsv\'ath, and Thurston define the weighted type D module $\CFDm(\HD)$
by tensoring $\CFAm_{U = 1}(\HD)$ with $\CFDDm(\Id)$. For computational purposes,
it is however more convenient to find $\CFDm(\HD)$ as a suitably defined
extension of the \emph{hat} invariant $\CFDa(\HD)$ using the following result.

\begin{proposition}[{\cite[Proposition~7.12]{LOT:torus-mod}}]
	\label{prop:type-d-structure}
	Consider a bordered $3$-manifold $Y$ with torus boundary. Let $P$ be a type D structure over $\Alg = \Alg(T^2)$ that is homotopy
	equivalent to $\CFDa(Y)$. We assume it to be reduced; that is, that the image of $\delta^1_P$ has no terms
	of the form $1 \otimes x$.
	
	Then, $P$ can be extended to a model for $\CFDm(Y)$. Namely, there exist operations $D_\rho: P \to P$ for Reeb chords $\rho$ in
	$\Algm$ but not in $\Alg$, such that the type D structure $Q$ over $\Algm$
	defined over the same vector space as $P$ with the differential
	\[\delta^1_Q = \delta^1_P + \sum_{\substack{\text{chords }\rho \in \Algm \\ \rho \not \in \Alg}} \rho \otimes D_\rho.\]
	is homotopy equivalent to $\CFDm(Y)$.
\end{proposition}
\subsection{Gradings}
\label{sec:gradings}

We first briefly recall what it means to grade a weighted $\Ainf$-algebra $\Alg$ by a non-commutative group, a weighted $\Ainf$-module $\fModule$ by a set with a right-action of this group, and a weighted type D module $D$ by a set with a left-action of this group.
Let $H$ be our non-commutative group, and consider central elements $\lambda_d, \lambda_w \in H$. The element $\lambda_w$ is the \emph{weight grading}; in the case when it is the identity, we will abbreviate $\lambda_d$ as $\lambda$. The grading on $\Alg$ satisfies the relation
\[ \gr(\mu^k_n(a_1, \dots, a_n)) = \lambda_d^{n-2} \lambda^k_w \gr(a_1) \cdots \gr(a_n). \]
The grading on $\fModule$ satisfies the relation
\begin{equation} 
\label{eqn:module-grading}
\gr(m^k_{1+n})(\x, a_1, \dots, a_n) = \gr(x) \cdot \lambda^{n-1}_d \lambda^k_w \gr(a_1) \cdots \gr(a_n).
\end{equation}
The grading on $D$ satisfies the relation
\[ \gr(\delta^1(x)) = \lambda^{-1}_d \gr(x). \]

We recall two groups by which $\Algm$ is graded: the \emph{big grading group} $G'$ and the \emph{intermediate grading group} $G$.

Consider the group $\frac{1}{2} \ZZ \times \ZZ^4$ with multiplication given by
\begin{align*}
(m; a, &b, c, d) \cdot (m'; a', b', c', d')  \\
	&= (m + m' + \frac{1}{2} \begin{vsmallmatrix}a & b\\a' & b'\end{vsmallmatrix}
	+ \frac{1}{2} \begin{vsmallmatrix}b & c\\b' & c'\end{vsmallmatrix}
	+ \frac{1}{2} \begin{vsmallmatrix}c & d\\c' & d'\end{vsmallmatrix}
	+ \frac{1}{2} \begin{vsmallmatrix}d & a\\d' & a'\end{vsmallmatrix};
	a + a', b + b', c + c', d + d')
\end{align*}
The big grading group $G'$ is the subgroup generated by the elements
\begin{align*}
\gr'(\rho_1) &= (-1/2; 1, 0, 0, 0) \qquad \qquad \gr'(\rho_2) = (-1/2; 0, 1, 0, 0) \\
\gr'(\rho_3) &= (-1/2; 0, 0, 1, 0) \qquad \qquad \gr'(\rho_4) = (-1/2; 0, 0, 0, 1) \\
\lambda &= (1; 0, 0, 0, 0)
\end{align*}
The first entry is called the \emph{Maslov component} of the grading; the remaining the \emph{$\SpinC$ component}. The special central elements are $\lambda_d = \lambda$ and $\lambda_w = (1; 1, 1, 1, 1)$. These, along with $\gr'(U) = (-1; 1, 1, 1, 1)$ determine a grading on $\Algm$.

The intermediate grading group $G$ is given by
\[ G = \left\{ (m; a, b) \in \left(\frac{1}{2} \ZZ \right)^3 ~\middle|~ a + b \in \ZZ, m + \frac{(2a + 1)(a + b + 1) + 1}{2} \in \ZZ \right \} .\] 
The multiplication on $G$ is given by
\[ (m; a, b) \cdot (m'; a', b') = (m + m' + \begin{vsmallmatrix}a & b \\ a' & b'\end{vsmallmatrix}; a + a'; b + b'). \]
There is a homomorphism $\pi: G' \to G$, and letting $\gr(a) = \pi(\gr'(a))$ determines a grading on $\Algm$ by $G$. We leave the details of this homomorphism to \cite[Section~4.2]{LOT:torus-alg}, and instead record $\gr$ for some elements. The special central elements are $\lambda_d = \lambda = (1; 0, 0)$ and $\lambda_w = (0; 0, 0)$. On Reeb chords, it takes the values
\begin{align*}
\gr(\rho_1) &= (-1/2; 1/2, -1/2) \qquad \qquad \gr(\rho_2) = (-1/2; 1/2, 1/2) \\
\gr(\rho_3) &= (-1/2; -1/2, 1/2) \qquad \qquad \gr(\rho_4) = (-3/2; -1/2, -1/2) \\
\gr(U) &= (-2; 0, 0)
\end{align*}
Now we turn to grading $\CFAm(\HD)$ for a bordered Heegaard diagram $\HD$. Given generators $\x$ and $\y$ in $\Gen(\HD)$, we can associate to each domain $B \in \pi(\x, \y)$ a grading in $G'$. That is,
\begin{equation}
\label{eqn:grading-domain}
 \gr'(B) = (-e(B) - n_{\x}(B) - n_{\y}(B); \bdy^\bdy B),	
\end{equation}
where $e(B)$ is the Euler measure of B \cite[Section~5.2]{LOT1} and $n_{\x}(B)$ and $n_{\y}(B)$ are the average of the local multiplicities of $B$ in the four regions near the points in $\x$ and $\y$ \cite[Section~5.7]{LOT1}. The grading set for $\CFAm(\HD)$ is the quotient of $G'$ by the gradings of the periodic domains. That is, for a fixed $\SpinC$-structure $\spinc$ and a generator $\x$ with $\spinc(\x) = \spinc$, let 
\[ P'(\x) = \langle \gr'(B) \mid B \in \pi_2(\x, \x), n_z(B) = 0 \rangle. \]
Then the grading set is $S'(\HD, \spinc) = P'(\x) \backslash G'$. Given a generator $\y$ with $\spinc(\y) = \spinc$, choosing any $B_y \in \pi_2(\x, \y)$, we can define
\[ \gr'(\y) = P'(\x) \gr(B_y). \]
This definition is independent of the choice of domain $B_y$ \cite[Proposition~8.2]{LOT:torus-mod}.

In a similar fashion, we can grade $\CFAm(\HD)$ with the analogous right $G$-set associated to the intermediate grading group.
We remark that the gradings on $\CFAm(\HD)$ are a direct extension of the gradings on $\CFAa(\HD)$.

\subsection{Embedded index formula}

In this section, we record some useful results about the expected dimension of the moduli spaces. The first of these proves a combinatorial formula for the expected dimension, and is known as the embedded index formula.

\begin{proposition}[{\cite[Proposition 3.30]{LOT:torus-alg}}] \label{prop:embedded-index-formula}
	Let $B$ be a homology class in $\pi_2(\x,\y)$ and $\vec a$ a sequence of Reeb chords compatible with $B$. Let $w$ be the number of simple Reeb orbits. Then,
	\begin{equation} \label{eq:embeedded-index-formula}
		\ind(u) = \ind(B, \vec a, w) = e(B) + n_{\x}(B) + n_{\y}(B) + |\vec a| + \iota(\vec a) + w.
	\end{equation}
\end{proposition}

Here, $e(B)$ is the Euler measure. $n_{\x}(B)$ and the $n_{\y}(B)$ are the average multiplicities of $B$ at the points in $\x$ and $\y$. We shall also call these the point measures. We refer the reader to \cite[Formula~(3.28)]{LOT:torus-mod} for a definition of $\iota(\vec{a})$. Instead, we rely on the following result to compute it.

\begin{lemma}[{\cite[Lemma 3.29]{LOT:torus-mod}}] \label{lem:iota-maslov}
	Let $\vec a = (a_1, \dots, a_m)$ be a sequence of Reeb chords such that the right idempotent of $a_i$ is equal to the left idempotent of $a_{i+1}$ for all $i$. Then $\iota(\vec a)$ is the Maslov component of
	$$\gr'(a_1) \cdots  \gr'(a_m).$$
\end{lemma}

The embedded index formula leads to a proof of the additivity of the expected dimensions of the moduli spaces.

\begin{proposition}[{\cite[Proposition 3.36]{LOT:torus-alg}}]
	\label{prop:index-additivity}
	Let $B \in \pi_2(\x, \y)$ and $B' \in \pi_2(\y, \z)$ be homology classes, and $\vec a$ and $\vec a'$ be sequences of Reeb chords compatible with $B$ and $B'$ respectively. Then $B + B'$ is a homology class in $\pi_2(\x, \z)$, and the concatenation $(\vec a, \vec a')$ of the sequences $\vec a$ and $\vec a'$ is compatible with it. The index is additive, namely
	$$\ind(B + B', (\vec a, \vec a'), w + w') = \ind(B, \vec a, w) + \ind(B, \vec a', w').$$
\end{proposition}

\subsection{Pairing theorems}
\label{sec:pairing-theorems}

The dream of bordered Heegaard Floer theory is to recover the Heegaard Floer invariant of a manifold by pairing the type A and type D invariants of the pieces. At the time of writing this paper, pairing theorems for the \emph{minus} variant remain under construction in upcoming work of Lipshitz, Ozsv\'ath, and Thurston \cite{LOT:torus-pairing}. We provide a brief sketch of the nature of this pairing.

There is a filtration on $\Algm$ given by the length of its elements, where each Reeb chord $\rho_i$ has length $1$, and $U$ has length $4$. We let $\Algc$ denote the completion of $\Algc$ with respect to this filtration, and remark that it is filtered bonsai. Let $\CFAc(\HD)$ denote the completion of $\CFAc(\HD)$ with respect to the $U$-power filtration, such that $\CFAc(\HD) = \CFAm(\HD) \tensor_{\FF_2[U]} \FF_2[[U]]$.

\begin{theorem}[{\cite[Theorem~1.36]{LOT:torus-mod}}, \cite{LOT:torus-pairing}]
	Given bordered $3$-manifolds $Y_1$ and $Y_2$ with $-\bdy Y_1 = T^2 = \bdy Y_2$, there is a homotopy equivalence
	\[ \CFmm(Y_2 \cup_{T^2} Y_1) \simeq \CFAc(Y_2) \boxtimes \CFDm(Y_1). \]
\end{theorem}

For the tensor products to be well-defined, we require boundedness properties on the type A and type D modules. In particular, we require that $\CFAc(\HD)$ is filtered bonsai. This is guaranteed when the bordered Heegaard diagram $\HD$ is admissible \cite[Lemma~7.5, Corollary~7.8]{LOT:torus-mod}. There is a graded version of the pairing theorem as well \cite[Theorem~8.11]{LOT:torus-mod}, \cite{LOT:abstract}. In practice, the graded computations work exactly as in the case for the \emph{hat} flavor of bordered Heegaard Floer.

\section{Description of \texorpdfstring{$\CFAm$}{CFA minus} with planar graphs}
\label{ch:cfa}

In this section, we construct the type A modules $\CFAm(\HD)$ for certain Heegaard diagrams $\HD$ corresponding to pattern knots for satellite operations. The operations on the type A modules will count certain classes of decorated planar graphs. This is analogous to the construction of the torus algebra, as described in Section~\ref{sec:torus-algebra}. A similar calculus of planar graphs has found use in other aspects of bordered knot Floer theory, to explicitly construct weighted $\Ainf$ algebra extensions of familiar objects: Ozsv\'ath and Szab\'o construct weighted deformations of their \emph{bordered knot algebras} \cite{OSzPong22, OSzPong23}; and Khan constructs weighted $\Ainf$ algebras associated to star-shaped Heegaard diagrams, and proves they are Koszul dual \cite{Khan2024}.

\subsection{Motivation}

\begin{figure}
	\centering
	\begin{tikzpicture}[scale=1.75]
		\draw[dashed] (.25,0) arc[radius=.25, start angle = 0, end angle = 90];
		\draw[dashed] (0,1.75) arc[radius=.25, start angle = 270, end angle = 360];
		\draw[dashed] (1.75,0) arc[radius=.25, start angle = 180, end angle = 90];
		\draw[dashed] (1.75,2) arc[radius=.25, start angle = 180, end angle = 270];
		\draw[color = red,thick] (.25,0) to (1.75,0);
		\draw[color = red,thick] (0,0.25) to (0,1.75);
		\draw[color = red,thick] (.25,2) to (1.75,2);
		\draw[color = red,thick] (2,0.25) to (2,1.75);
		\node at (.25,1.7) (r4) {$\rho_4$};
		\node at (1.7,1.7) (r1) {$\rho_1$};
		\node at (1.7,.3) (r2) {$\rho_2$};
		\node at (.3,.3) (r3) {$\rho_3$};
		\node at (.4,1.9) (z) {$z$};
		\draw[color = blue,thick] (1,0) to[in=270,out=90] (1,2);
		\node at (1,-.2) {$a$};
	\end{tikzpicture}
	\caption[A bordered Heegaard
	diagrams \texorpdfstring{$\HD_{st}$}{HDst} for the $0$-framed solid torus.]
	{\textbf{A bordered Heegaard
	diagrams $\HD_{st}$ for the $0$-framed solid torus.} The intersection of the $\alpha$-arc and the $\beta$-circle is $a$. There is a marked point $z$ near the boundary with Reeb chord $\rho_4$.}
	\label{fig:solid-torus}
\end{figure}
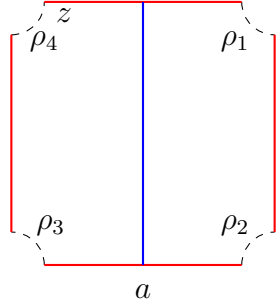

As motivation, we will construct the type A module for the $0$-framed solid torus. The bordered Heegaard diagram $\HD_{st}$ for the $0$-framed solid torus is pictured in Figure~\ref{fig:solid-torus}. There is one homologically non-trivial $\beta$-circle, which intersects an $\alpha$-arc at $a$. As a left $\FF_2[U]\langle \iota_0, \iota_1 \rangle$-module, $\CFAm(\HD_{st})$ is freely generated by $a$.

\begin{remark}
	Our construction of $\CFAm(\HD_{st})$ using certain planar graphs will readily generalize to our construction of other type A modules. We use $\x, \y, \z, \dots$ to denote generators of the type A module in constructions and arguments where specifics of the Heegaard diagram $\HD_{st}$ do not feature. To specialize to $\HD_{st}$, read all of these as $a$.
\end{remark}

As our diagrams $\HD$ are of genus $g = 1$, instead of looking at positive domains $B \in \pi_2(\x, \y)$, it is equivalent to look at immersions of the disk into the torus with appropriate boundary conditions. That is, let $\bD$ denote the standard unit disk in $\CC$. Consider the bordered Heegaard diagram to be drawn on the torus $T^2$. Let the $\alpha$-arcs intersect transversely at a point $p$, which represents the boundary component of $\Sigma$. We can label the quadrants around $p$ with $1, 2, 3, 4$ in a clockwise orientation; this records the Reeb chords. We further let the $\beta$-circle intersect the $\alpha$-arcs transversely too; in the case of $\HD_{st}$, in a single point $a$. The immersions we are interested are maps $u: \bD \to T^2$ such that:
\begin{itemize}
	\item If $|z| = 1$ and $\Re(z) < 0$, then $u(z)$ lies on the $\alpha$-arcs; if $|z| = 1$ and $\Re(z) > 0$, then $u(z)$ lies on the $\beta$-circle. That is, the left boundary of the disk is mapped to the $\alpha$-arcs; the right, to the $\beta$-circle.
	\item $u(-i) = \x$ and $u(i) = \y$.
	\item The map $u$ is an immersion everywhere, except for pre-images of the point $p$ that lie on the boundary $\bdy \bD$.
\end{itemize}

In slight overload of terminology, we shall say such an immersion $u$ is an immersed disk from $\x$ to $\y$.

The presence of such an immersion does not guarantee a contribution to an operation on the type A module. We further require a count of the number of pseudoholomorphic representatives of any such immersion. It is a fact of complex analysis that any such immersion $u$ that maps the \emph{corners} at $-i$ and $i$ to \emph{acute} corners admits a unique pseudoholomorphic representative.

\begin{proposition}
	\label{prop:acute-angled-disk}
	Let $u: \bD \to T^2$ be a disk from $\x$ to $\y$, and suppose further that it maps the corners at $-i$ and $i$ to \emph{acute} intersections of the $\alpha$-arcs and $\beta$-circles.
	
	Let $B$ be the positive domain in $\pi_2(\x, \y)$ corresponding to the image of $u$; that is, each elementary domain is included with multiplicity by which $u$ covers it. The labeling of the quadrants around $p$ determines a sequence of Reeb chords $\vec a$ associated to $u$, as encountered in order along $\bdy \bD$. The number of simple Reeb orbits $w$ is the count of pre-images of $p$ in the interior of $\bdy D$.
	
	Then, $\# \left(\cM^B(\x, \y; \vec a; w)\right) \equiv 1 \mod 2$.
\end{proposition}
\begin{proof}
	This is a classical result in Heegaard Floer theory. See \cite[Proposition 2.7.3]{Hans14} for one exposition of a proof.
\end{proof}

\subsubsection{Planar graphs and tiling patterns}
\label{sec:tiling-patterns}

Instead of working with immersions of the disk, we associate to each immersion a decorated planar graph. These decorated planar graphs we call \emph{module tiling patterns}, and they uniquely specify an immersion of the disk into the torus. The description of the operations on $\CFAm$ is in terms of these module tiling patterns. Not only does this make the description of the modules combinatorial, but also aids in proving the structure relations on $\CFAm$, as the operations on the algebra $\Algm$ are constructed with similar classes of planar graphs.

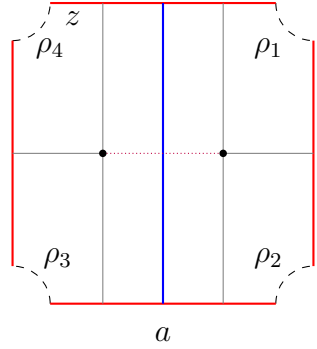
\begin{figure}
	\centering
	\begin{tikzpicture}[scale=2]
		\draw[dashed] (.25,0) arc[radius=.25, start angle = 0, end angle = 90];
		\draw[dashed] (0,1.75) arc[radius=.25, start angle = 270, end angle = 360];
		\draw[dashed] (1.75,0) arc[radius=.25, start angle = 180, end angle = 90];
		\draw[dashed] (1.75,2) arc[radius=.25, start angle = 180, end angle = 270];
		\draw[color = red,thick] (.25,0) to (1.75,0);
		\draw[color = red,thick] (0,0.25) to (0,1.75);
		\draw[color = red,thick] (.25,2) to (1.75,2);
		\draw[color = red,thick] (2,0.25) to (2,1.75);
		\node at (.25,1.7) (r4) {$\rho_4$};
		\node at (1.7,1.7) (r1) {$\rho_1$};
		\node at (1.7,.3) (r2) {$\rho_2$};
		\node at (.3,.3) (r3) {$\rho_3$};
		\node at (.4,1.9) (z) {$z$};
		\draw[color = blue,thick] (1,0) to[in=270,out=90] (1,2);
		\node at (1,-.2) {$a$};
		
		\draw[color = gray] (.6,1) to (0,1);
		\draw[color = gray] (.6,1) to (0.6,0);
		\draw[color = gray] (.6,1) to (0.6,2);
		\draw[color = purple,densely dotted] (.6,1) to (1.4,1);
		\draw[color = gray] (1.4,1) to (2,1);
		\draw[color = gray] (1.4,1) to (1.4,0);
		\draw[color = gray] (1.4,1) to (1.4,2);

		\node at (.6,1) [circle,fill,inner sep=1pt]{};
		\node at (1.4,1) [circle,fill,inner sep=1pt]{};
	\end{tikzpicture}
	\caption[The dual graph for the $0$-framed solid torus.]
	{\textbf{The dual graph for the $0$-framed solid torus.} The bordered Heegaard
	diagrams is given by the thick red arcs and the thick blue curve. The black
	nodes are the vertices in $\Gamma(\HD_{st})$.
	The edges in $\Gamma^\beta(\HD_{st})$ are dotted purple, and the remaining edges
	in $\Gamma(\HD_{st})$ are grey.}
	\label{fig:solid-torus-dual}
\end{figure}

Given a bordered Heegaard diagram $\HD$, we can associate a graph $\Gamma(\HD)$ that is \emph{dual} to the $\alpha$-arcs and the $\beta$-circle: To each component of $\Sigma \setminus (\alpha_1 \cup \alpha_2 \cup \beta)$, associate a vertex. There is an edge connecting any two components across an $\alpha$-arc or a $\beta$-circle if the two components share the arc or circle as a boundary. Figure~\ref{fig:solid-torus-dual} illustrates this for a bordered Heegaard diagram of the solid torus. We denote the edges that cross the $\beta$-circle by $\Gamma^\beta(\HD)$.

An immersion $u$ from $\x$ to $\y$ naturally gives rise to a (necessarily) planar graph in the disk, by taking the pre-image of $\Gamma(\HD)$. These planar graphs are connected and come with an embedding into the disk $\bD$. We will mark the boundary of the disk with red where it maps to the $\alpha$-arcs, and with blue where it maps to the $\beta$-circles. The pre-image of $\Gamma(\HD)$ meets the boundary of the disk in \emph{leaves}, and these do not correspond to vertices of $\Gamma(\HD)$. As a convention, we use \emph{vertex} to refer to vertices of our planar graphs that are not leaves.

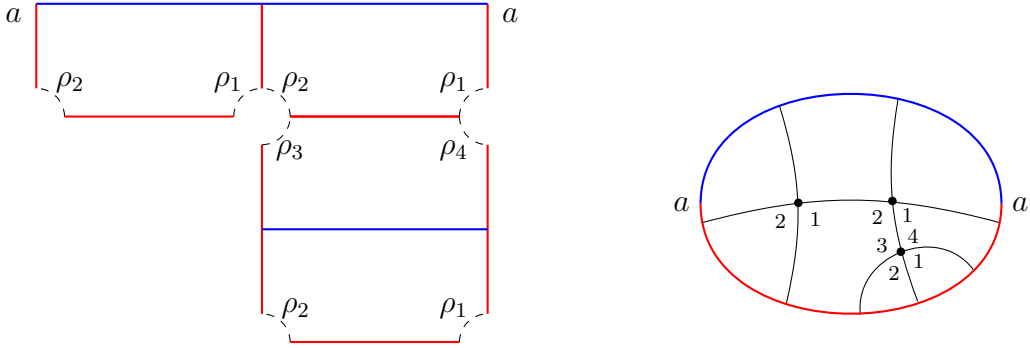
\begin{figure}
	\centering
		\begin{tikzpicture}[scale=1.5]
			\begin{scope}[shift={(0,0)}]
				\node at (.25,1.7) (r4) {$\rho_3$};
				\node at (1.7,1.7) (r1) {$\rho_4$};
				\node at (1.7,.3) (r2) {$\rho_1$};
				\node at (.3,.3) (r3) {$\rho_2$};
				\draw[dashed] (.25,0) arc[radius=.25, start angle = 0, end angle = 90];
				\draw[dashed] (0,1.75) arc[radius=.25, start angle = 270, end angle = 360];
				\draw[dashed] (1.75,0) arc[radius=.25, start angle = 180, end angle = 90];
				\draw[dashed] (1.75,2) arc[radius=.25, start angle = 180, end angle = 270];
				\draw[color = red,thick] (.25,0) to (1.75,0);
				\draw[color = red,thick] (0,0.25) to (0,1.75);
				\draw[color = red,thick] (.25,2) to (1.75,2);
				\draw[color = red,thick] (2,0.25) to (2,1.75);
		
				\draw[color = blue,thick] (0,1) to (2,1);
			\end{scope}
			\begin{scope}[shift={(0,2)}]
				\node at (2.2,0.9) (a2) {$a$};
				\node at (1.7,.3) (r2) {$\rho_1$};
				\node at (.3,.3) (r3) {$\rho_2$};
				\draw[dashed] (.25,0) arc[radius=.25, start angle = 0, end angle = 90];
				\draw[dashed] (1.75,0) arc[radius=.25, start angle = 180, end angle = 90];
				\draw[color = red,thick] (.25,0) to (1.75,0);
				\draw[color = red,thick] (0,0.25) to (0,1);
				\draw[color = red,thick] (2,0.25) to (2,1);
		
				\draw[color = blue,thick] (0,1) to (2,1);
			\end{scope}
			\begin{scope}[shift={(-2,2)}]
				\node at (-0.2,0.9) (a2) {$a$};
				\node at (1.7,.3) (r2) {$\rho_1$};
				\node at (.3,.3) (r3) {$\rho_2$};
				\draw[dashed] (.25,0) arc[radius=.25, start angle = 0, end angle = 90];
				\draw[dashed] (1.75,0) arc[radius=.25, start angle = 180, end angle = 90];
				\draw[color = red,thick] (.25,0) to (1.75,0);
				\draw[color = red,thick] (0,0.25) to (0,1);
				\draw[color = red,thick] (2,0.25) to (2,1);
		
				\draw[color = blue,thick] (0,1) to (2,1);
			\end{scope}
		\end{tikzpicture}	
		\qquad \qquad
		\begin{tikzpicture}
	\begin{pgfonlayer}{nodelayer}
		\node [style=none] (0) at (-4.5, 3.5) {};
		\node [style=none] (1) at (3.5, 3.5) {};
		\node [style=none] (2) at (-4.45, 3) {};
		\node [style=none] (3) at (3.45, 3) {};
		\node [style=none] (4) at (-2.4, 6.1) {};
		\node [style=none] (5) at (-2.225, 0.825) {};
		\node [style=none] (6) at (-5, 3.5) {$a$};
		\node [style=none] (7) at (4, 3.5) {$a$};
		\node [style=number-labels] (8) at (-2.375, 3.025) {$2$};
		\node [style=number-labels] (9) at (-1.45, 3.1) {$1$};
		\node [style=dot] (10) at (-1.9, 3.525) {};
		\node [style=none] (11) at (0.75, 6.275) {};
		\node [style=none] (12) at (1.275, 0.9) {};
		\node [style=none] (13) at (-0.25, 0.55) {};
		\node [style=none] (14) at (2.75, 1.75) {};
		\node [style=number-labels] (15) at (0.2, 3.15) {$2$};
		\node [style=number-labels] (16) at (0.325, 2.4) {$3$};
		\node [style=number-labels] (17) at (0.65, 1.65) {$2$};
		\node [style=number-labels] (18) at (1.3, 1.975) {$1$};
		\node [style=number-labels] (19) at (1.15, 2.625) {$4$};
		\node [style=number-labels] (20) at (1, 3.2) {$1$};
		\node [style=dot] (21) at (0.6, 3.575) {};
		\node [style=dot] (22) at (0.825, 2.225) {};
	\end{pgfonlayer}
	\begin{pgfonlayer}{edgelayer}
		\draw [bend left=15] (2.center) to (3.center);
		\draw [bend left=15] (4.center) to (5.center);
		\draw [bend right=15] (11.center) to (12.center);
		\draw [bend left=75, looseness=1.25] (13.center) to (14.center);
		\draw [style=blue, bend left=90, looseness=1.25] (0.center) to (1.center);
		\draw [style=red, bend right=90, looseness=1.25] (0.center) to (1.center);
	\end{pgfonlayer}
\end{tikzpicture}
		\caption [An immersed disk and its corresponding module tiling pattern.] {\textbf{An immersed disk and its corresponding module tiling pattern.}}
		\label{fig:solid-torus-illustration}
\end{figure}

The pre-image of $\Gamma(\HD)$ may contain pre-images of edges from $\Gamma^\beta(\HD)$. We contract these edges. If this results in parallel edges between vertices, we delete the parallel copies. We call this resulting graph, associated to our immersion u, $\Gamma(u)$. Figure~\ref{fig:solid-torus-illustration} illustrates the planar graph corresponding to an immersed disk for $\HD_{st}$. There are four vertices before we contract the edges in the pre-image of $\Gamma^\beta(\HD_{st})$, but three vertices post contraction.

\begin{definition} Let $\Gamma = \Gamma(u)$ for an immersion $u$ from $\x$ to $\y$.
	The faces of the planar graph $\Gamma$ are the components of $\bD \setminus \Gamma$.
	If a face intersects the blue boundary, we call it a \emph{blue face}. If it otherwise
	intersects the red
	boundary, we call it a \emph{red face}. 
	The faces that are neither blue faces nor red faces are the \emph{internal faces}.
	
	We call a vertex with an edge to
	the blue boundary a \emph{blue vertex}. We call every other vertex a \emph{red vertex}.
	
	The internal faces in
	our graphs
	are bound by four edges. We refer to these boundaries as \emph{short cycles}.
\end{definition}

The red vertices correspond to regions that have as boundary solely the $\alpha$-arcs. Equivalently, they correspond to gluing in a copy of the whole torus cut along the $\alpha$-arcs. We further observe that every vertex in a planar graphs $\Gamma(u)$ meets at most four faces that are not blue.

A module tiling pattern consists of such a planar graph $\Gamma(u)$ along with additional
data. Along the intersections of the red and the blue boundaries of the disk, we mark $\x$ and $\y$ for $u$ an immersion from $\x$ to $\y$. Along with each vertex $v$, we record a power of $U$ to keep track of the multiplicity of intersection with $z$; that is, $U_v = U^i$ for $i \in \{0, 1\}$. The red vertices represent an additional copy of the torus; therefore, $U_v = U$. When we illustrate the module tiling patterns, we omit recording $U_v$. 

The module tiling patterns come equipped with a labeling around the vertices to keep track of the Reeb chords.

\begin{definition}
	Let $\Gamma$ be a planar graph arising from an immersion. For any vertex $v$, let $Q_v$ be the set of non-blue faces that meet
	$v$. A \emph{valid labeling} on $\Gamma$ consists of a
	labeling $\Lambda_v:
	Q_v \to \{1,2,3,4\}$ for every vertex $v$ such that:
	\begin{itemize}
		\item Let $Q_v = \{r_1, \dots, r_k\}$ be the non-blue faces around $v$, labelled in clockwise order. We observed previously that $k$ is at most four. Then, up to a cyclic reordering, $\Lambda_v(r_1), \dots, \Lambda_v(r_k)$ is a contiguous subsequence of $(1,2,3,4)$.
		\item If $e : v_1 \to v_2$ is an edge oriented such that a non-blue face $f$
		is to the right of $e$, then $\Lambda_{v_1}(f) + 1 \equiv \Lambda_{v_2}(f) \mod 4$.
	\end{itemize}
\end{definition}

This is analogous to the labelings on the algebra tiling patterns. Analogously, a valid labeling is uniquely determined by its value on any single non-blue face. An immersion $u$ from $\x$ to $\y$ gives rise to a valid labeling on $\Gamma(u)$ using the labeling of the quadrants near $p$ in $T^2$. Figure~\ref{fig:solid-torus-illustration} illustrates one such valid labeling that arises from an immersion.

The above discussion constructs from an immersion $u$ a module tiling pattern $\Gamma(u)$. Conversely, given a module tiling pattern $\Gamma$ that coheres with the $\beta$-circle, we can construct an immersion $u$: the dual graph provides a recipe to take regions of $T^2 \setminus (\alpha_1 \cup \alpha_2 \cup \beta)$ and identify them along shared edges. We do not spell out the conditions near the $\beta$-circle; our module tiling operations will be constructed iteratively from simpler module tiling operations.

As a convention, we orient the boundary of the disk on which we draw our planar graphs counterclockwise.
We draw the module tiling patterns with the blue boundary on the top, and the
red boundary on the bottom. If $\Gamma$ is a module tiling operation from $\x$ to $\y$, we list $\x$ on the left
blue-red intersection, and $\y$ on the right intersection.

Immersions correspond to operations on the type A module, provided they admit pseudoholomorphic representatives. We spell out this correspondence with operations on the module for the module tiling patterns.

\begin{definition} Let $\Gamma$ be a module tiling pattern.
	The \emph{weight} $w$ of $\Gamma$ is defined as the number of short cycles in $\Gamma$.
	Equivalently,
	the first homology $H_1(\Gamma; \ZZ)$ of the underlying planar graph has dimension $w$.
\end{definition}

Recall that for an immersion $u$, the number of Reeb orbits equals the count of the pre-image of $p$ in the interior of $\bD$. This is equal to the count of the short cycles in $\Gamma(u)$.

\begin{definition} Let $\Gamma$ be a module tiling pattern.
	The \emph{chord sequence} of $\Gamma$ is given by reading the labels visible
	along the red faces while traversing $\bdy \bD$ (with its chosen orientation).

	That is, let $f_1, \dots, f_k$ be the red faces in order along the boundary,
	starting from any blue face.
	For any $1 \leq i \leq k$, let $v^i_{1}, \dots, v^i_{j_i}$ be the vertices in order
	along the internal boundary of the face $f_i$.
	Then, the chord sequence corresponding to $\Gamma$ is defined as
	\[ \left(\prod_{l = 1}^{j_1} \rho_{\Lambda_{v^1_{l}}(f_1)}\right)\tensor \dots \tensor 
	  \left(\prod_{l = 1}^{j_i} \rho_{\Lambda_{v^i_{l}}(f_i)}\right)\tensor \dots \tensor
	  \left(\prod_{l = 1}^{j_k} \rho_{\Lambda_{v^k_{l}}(f_k)}\right).\]
\end{definition}

\begin{lemma}[c.f. {\cite[Lemma 3.6]{LOT:torus-alg}}]
	Let $\Gamma$ be a module tiling pattern, and $a_1 \tensor \cdots \tensor a_n$ its chord sequence. Then, $a_1 \tensor \cdots \tensor a_n$ is non-zero; and for any $1 \leq i < n$, we have $a_i a_{i+1} = 0$.
\end{lemma}
\begin{proof}
	For any $1 \leq i < n$, factorize $a_i$ as $a_i = a \rho_j$. This means the last vertex of the red face $f_i$ carries the label $j$. This vertex must be the first for the red face $f_{i+1}$, and its corresponding label is forced to be $j - 1 \mod 4$ by virtue of the labeling conventions. So, $a_{i+1} = \rho_{j - 1} a'$. This ensures that $a_i a_{i + 1} = 0$; and that the idempotents of $a_i$ and $a_{i+1}$ agree and the chord sequence is non-zero.
\end{proof}

Let $\Gamma$ be a tiling pattern from $\x$ to $\y$ with weight $w$
and corresponding chord sequence
$a_1 \tensor \dots \tensor a_n$. The module operation associated to $\Gamma$ is
\[m^w_{1+n}(\x, a_1, \dots, a_n) = \left(\prod_{\text{vertex }v \in \Gamma} U_v\right) \y.\]
For example, the module operation associated to the module tiling pattern in Figure~\ref{fig:solid-torus-illustration} is $m^0_{1+5}(a, \rho_2, \rho_{123}, \rho_2, \rho_1, \rho_{41}) = Ua$.

\subsubsection{Moves to construct new module tiling patterns}
\label{sec:three-moves}

In \cite[Section 9.2]{LOT:torus-mod}, Lipshitz, Ozsv{\'a}th, and Thurston sketch three moves to obtain new immersions from a starting set of immersions of disks into the torus. We recall these moves now, and spell out their analogues in terms of the tiling patterns. All operations on our type A modules will admit a simple inductive description based on these moves.

Move (1): Consider two immersions of the disk corresponding to the operations $m^{w}_{1+n}(\x, a_1, \dots, a_n) = U^k \y$ and $m^{w'}_{1+m}(\y, a'_1, \dots, a'_m) = U^l \z$ such that $a_n a'_1 \neq 0$. The fact that $a_n a'_1 \neq 0$ means that we can glue the two disks along the $\alpha$-arc between $\y$ and $a_n$ of the first disk, and between $\y$ and $a'_1$ of the second disk. This yields a new immersion from $\x$ to $\z$. The weights add, and the corresponding operation is
		\[m^{w+w'}_{n+m}(\x, a_1, \dots, a_n a'_1, \dots, a'_m) = U^{k+l} \z.\]

\begin{figure}
	\centering
		\begin{tikzpicture}[scale=1.5]
			\begin{scope}[shift={(2,0)}]
				\node at (.25,1.7) (r4) {$\rho_3$};
				\node at (1.7,1.7) (r1) {$\rho_4$};
				\node at (1.7,.3) (r2) {$\rho_1$};
				\node at (.3,.3) (r3) {$\rho_2$};
				\draw[dashed] (.25,0) arc[radius=.25, start angle = 0, end angle = 90];
				\draw[dashed] (0,1.75) arc[radius=.25, start angle = 270, end angle = 360];
				\draw[dashed] (1.75,0) arc[radius=.25, start angle = 180, end angle = 90];
				\draw[dashed] (1.75,2) arc[radius=.25, start angle = 180, end angle = 270];
				\draw[color = red,thick] (.25,0) to (1.75,0);
				\draw[color = red,thick] (0,0.25) to (0,1.75);
				\draw[color = red,thick] (.25,2) to (1.75,2);
				\draw[color = red,thick] (2,0.25) to (2,1.75);
		
				\draw[color = blue,thick] (0,1) to (2,1);
			\end{scope}
			\begin{scope}[shift={(2,2)}]
				\node at (2.2,0.9) (a2) {$a$};
				\node at (-0.2,0.9) (a2) {$a$};
				\node at (1.7,.3) (r2) {$\rho_1$};
				\node at (.3,.3) (r3) {$\rho_2$};
				\draw[dashed] (.25,0) arc[radius=.25, start angle = 0, end angle = 90];
				\draw[dashed] (1.75,0) arc[radius=.25, start angle = 180, end angle = 90];
				\draw[color = red,thick] (.25,0) to (1.75,0);
				\draw[color = red,thick] (0,0.25) to (0,1);
				\draw[color = red,thick] (2,0.25) to (2,1);
		
				\draw[color = blue,thick] (0,1) to (2,1);
				
				\draw [line width = 4pt, draw=purple, opacity=0.4]
				(0,0.25) to (0,1);
			\end{scope}
			\begin{scope}[shift={(-2,2)}]
				\node at (2.2,0.9) (a2) {$a$};
				\node at (-0.2,0.9) (a2) {$a$};
				\node at (1.7,.3) (r2) {$\rho_1$};
				\node at (.3,.3) (r3) {$\rho_2$};
				\draw[dashed] (.25,0) arc[radius=.25, start angle = 0, end angle = 90];
				\draw[dashed] (1.75,0) arc[radius=.25, start angle = 180, end angle = 90];
				\draw[color = red,thick] (.25,0) to (1.75,0);
				\draw[color = red,thick] (0,0.25) to (0,1);
				\draw[color = red,thick] (2,0.25) to (2,1);
		
				\draw[color = blue,thick] (0,1) to (2,1);
				
				\draw [line width = 4pt, draw=purple, opacity=0.4]
				(2,0.25) to (2,1);
			\end{scope}
		\end{tikzpicture}
		\caption [An illustration of disks where move (1) can be performed.] {\textbf{An illustration of disks where move (1) can be performed.} The two disks can be glued along the shaded edge. This results in the disk from Figure~\ref{fig:solid-torus-illustration}.}
		\label{fig:move-1-illustration}
\end{figure}
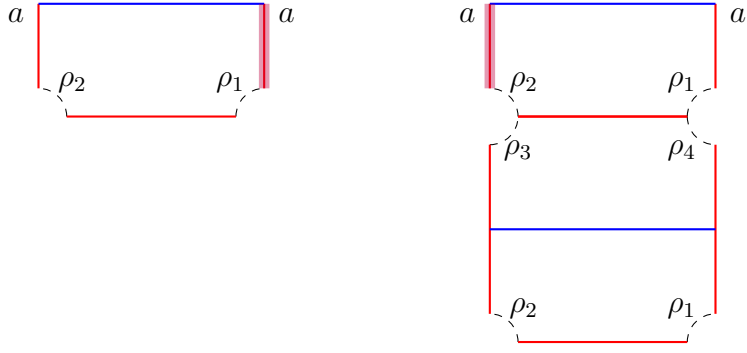

	\sloppy Figure~\ref{fig:move-1-illustration} illustrates an example of disks that can be glued as such. Here, $m^0_{3}(a, \rho_2, \rho_1) = a$ and $m^0_{5}(a, \rho_{23}, \rho_2, \rho_1, \rho_{41}) = Ua$ yield the new operation
	$m^0_{6}(a, \rho_2, \rho_{123}, \rho_2, \rho_1, \rho_{41}) = Ua$.

\begin{figure}
	\centering
	\begin{tikzpicture}[scale=1.25]
	\begin{pgfonlayer}{nodelayer}
		\node [style=none] (0) at (-6.5, 0) {};
		\node [style=none] (1) at (-8.5, 2) {};
		\node [style=none] (2) at (-8.5, -2) {};
		\node [style=none] (3) at (-4, 0) {};
		\node [style=none] (4) at (-2, 2) {};
		\node [style=none] (5) at (-2, -2) {};
		\node [style=dot] (6) at (-7.75, 0) {};
		\node [style=none] (7) at (-7.25, 0.75) {};
		\node [style=none] (8) at (-8.5, -0.25) {};
		\node [style=dot] (9) at (-2.75, 0) {};
		\node [style=none] (10) at (-3.25, 0.75) {};
		\node [style=none] (11) at (-2, -0.25) {};
		\node [style=none] (12) at (-8.25, 0.75) {};
		\node [style=none] (13) at (-2.25, 0.75) {};
		\node [style=none] (14) at (-6.575, -0.5) {};
		\node [style=none] (15) at (-3.925, -0.5) {};
		\node [style=none] (17) at (2.25, 2) {};
		\node [style=none] (18) at (2.25, -2) {};
		\node [style=none] (20) at (6.25, 2) {};
		\node [style=none] (21) at (6.25, -2) {};
		\node [style=dot] (22) at (3, 0) {};
		\node [style=none] (23) at (3.5, 0.75) {};
		\node [style=none] (24) at (2.25, -0.25) {};
		\node [style=dot] (25) at (5.5, 0) {};
		\node [style=none] (26) at (5, 0.75) {};
		\node [style=none] (27) at (6.25, -0.25) {};
		\node [style=none] (28) at (2.5, 0.75) {};
		\node [style=none] (29) at (6, 0.75) {};
		\node [style=none] (30) at (4.175, -0.4) {};
		\node [style=none] (31) at (4.325, -0.4) {};
		\node [style=none] (32) at (-0.5, 0) {};
		\node [style=none] (33) at (0.5, 0) {};
		\node [style=number-labels] (34) at (-7.75, -0.525) {$i$};
		\node [style=number-labels] (35) at (-2.75, -0.6) {$i+1$};
		\node [style=number-labels] (36) at (5.5, -0.6) {$i+1$};
		\node [style=number-labels] (37) at (3, -0.525) {$i$};
		\node [style=none] (40) at (-6, 0) {$\mathbf{x}$};
		\node [style=none] (41) at (-4.5, 0) {$\mathbf{x}$};
	\end{pgfonlayer}
	\begin{pgfonlayer}{edgelayer}
		\draw [style=blue, bend left=45] (1.center) to (0.center);
		\draw [style=blue, bend left=45] (3.center) to (4.center);
		\draw [bend right=15] (12.center) to (6);
		\draw [bend left=15] (7.center) to (6);
		\draw [bend left=15] (6) to (8.center);
		\draw [bend right=15] (9) to (13.center);
		\draw [bend right=15] (10.center) to (9);
		\draw [bend right=15] (9) to (11.center);
		\draw [bend right=15] (6) to (14.center);
		\draw [bend right=15] (15.center) to (9);
		\draw [style=red, bend left=45] (0.center) to (2.center);
		\draw [style=red, bend right=45] (3.center) to (5.center);
		\draw [bend right=15] (28.center) to (22);
		\draw [bend left=15] (23.center) to (22);
		\draw [bend left=15] (22) to (24.center);
		\draw [bend right=15] (25) to (29.center);
		\draw [bend right=15] (26.center) to (25);
		\draw [bend right=15] (25) to (27.center);
		\draw [bend right=15] (22) to (30.center);
		\draw [bend right=15] (31.center) to (25);
		\draw (30.center) to (31.center);
		\draw [style=blue, bend right=45, looseness=0.75] (17.center) to (20.center);
		\draw [style=red, in=135, out=45, looseness=0.75] (18.center) to (21.center);
		\draw [style=double-sided] (32.center) to (33.center);
	\end{pgfonlayer}
\end{tikzpicture}
	\caption[Move (1) on module tiling patterns.]
	{\textbf{Move (1) on module tiling patterns.} The labeling near the vertex is forced by the condition that $a_n a'_1 \neq 0$. The move identifies the two edges to the red boundary on the left.}
	\label{fig:move-1}
\end{figure}
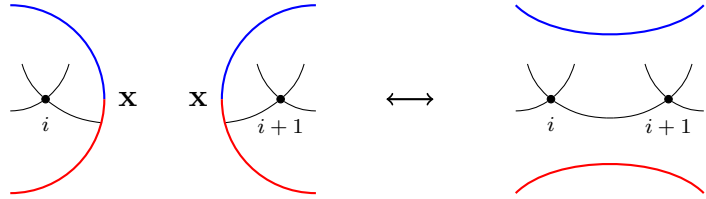

In terms of the module tiling patterns, the new module tiling pattern is obtained by gluing the right blue-red intersection of the first pattern with the left blue-red intersection of the second pattern, in such a way that the last edge to the red boundary of the first pattern is identified with the first edge to the red boundary of the first pattern. This is illustrated in Figure~\ref{fig:move-1}.

Move (2): Consider an immersion of the disk corresponding to the operation $m^w_{1+n}(\x, a_1, \dots, a_n) = U^k \y$. For any $1 \leq i < n$ with $a_i a_{i+1} = 0$, write $a_i$ as $a' \rho_j$. We can glue in a copy of the torus along an $\alpha$-arc to the edge between $a_i$ and $a_{i+1}$. This yields a new operation
\[m^w_{3+n}(\x, a_1, \dots, a_i \rho_{j+1}, \rho_j, \rho_{j-1},
    \rho_{j-2} a_{i+1}, \dots, a_n) = U^{k+1} \y.\]
    
    Move (2) on the immersed disk on the left in Figure~\ref{fig:move-1-illustration} with $i = 1$ yields the immersed disk on the right. As another example, $m^0_3(a, \rho_4, \rho_3) = Ua$ yields $m^0_5(a, \rho_{41}, \rho_4, \rho_3, \rho_{23}) = U^2 a$.
    
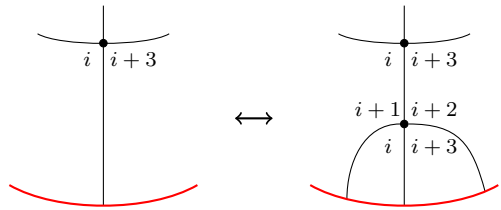
\begin{figure}
	\centering
	\begin{tikzpicture}
	\begin{pgfonlayer}{nodelayer}
		\node [style=none] (32) at (-0.5, 0) {};
		\node [style=none] (33) at (0.5, 0) {};
		\node [style=none] (34) at (-1.5, -1.775) {};
		\node [style=none] (35) at (-6.5, -1.775) {};
		\node [style=none] (36) at (-4, 3) {};
		\node [style=none] (37) at (-4, -2.325) {};
		\node [style=number-labels] (40) at (-4.425, 1.575) {$i$};
		\node [style=number-labels] (41) at (-3.2, 1.55) {$i+3$};
		\node [style=none] (42) at (6.5, -1.775) {};
		\node [style=none] (43) at (1.5, -1.775) {};
		\node [style=none] (44) at (4, 3) {};
		\node [style=none] (45) at (4, -2.325) {};
		\node [style=number-labels] (48) at (3.575, 1.575) {$i$};
		\node [style=number-labels] (49) at (4.8, 1.55) {$i+3$};
		\node [style=number-labels] (56) at (3.3, 0.225) {$i+1$};
		\node [style=number-labels] (57) at (3.575, -0.725) {$i$};
		\node [style=number-labels] (58) at (4.8, -0.775) {$i+3$};
		\node [style=number-labels] (60) at (4.8, 0.225) {$i+2$};
		\node [style=none] (62) at (6.15, -1.95) {};
		\node [style=none] (63) at (2.475, -2.15) {};
		\node [style=dot] (64) at (4, -0.15) {};
		\node [style=none] (65) at (-5.75, 2.25) {};
		\node [style=none] (66) at (-2.25, 2.25) {};
		\node [style=none] (67) at (2.25, 2.25) {};
		\node [style=none] (68) at (5.75, 2.25) {};
		\node [style=dot] (69) at (-4, 2) {};
		\node [style=dot] (70) at (4, 2) {};
	\end{pgfonlayer}
	\begin{pgfonlayer}{edgelayer}
		\draw [style=double-sided] (32.center) to (33.center);
		\draw (36.center) to (37.center);
		\draw (44.center) to (45.center);
		\draw [in=90, out=180] (64) to (63.center);
		\draw [in=105, out=0, looseness=1.25] (64) to (62.center);
		\draw [style=red, bend right, looseness=0.75] (35.center) to (34.center);
		\draw [style=red, bend right, looseness=0.75] (43.center) to (42.center);
		\draw [bend right, looseness=0.50] (65.center) to (66.center);
		\draw [bend right, looseness=0.50] (67.center) to (68.center);
	\end{pgfonlayer}
\end{tikzpicture}
	\caption[Move (2) on module tiling patterns.]
	{\textbf{Move (2) on module tiling patterns.} The move adds a new vertex with labelings as illustrated.}
	\label{fig:move-2}
\end{figure}

    In terms of the module tiling patterns, the move adds a new red vertex to the pattern. This is illustrated in Figure~\ref{fig:move-2}.

Move (3): Consider an immersion of the disk corresponding to the operation $m^w_{1+n}(\x, a_1, \dots, a_n) = U^k \y$ such that there is an
	$1 < i < n$ with $|a_i| = 4$. There is an $\alpha$-arc that precedes and follows the chord $a_i$, and we can glue the two together to obtain a new disk. This increases the number of simple Reeb orbits by one, and yields the operation
	\[m^{w+1}_{n-1}(\x, a_1, \dots, a_{i-2}, a_{i-1}a_{i+1}, a_{i+2}, \dots, a_n)
	= U^k \y.\]

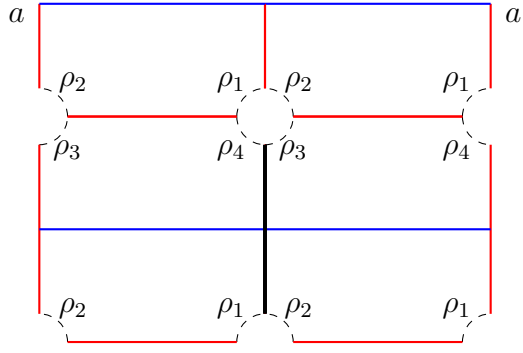
\begin{figure}
	\centering
		\begin{tikzpicture}[scale=1.5]
			\begin{scope}[shift={(-2,0)}]
				\node at (.25,1.7) (r4) {$\rho_3$};
				\node at (1.7,1.7) (r1) {$\rho_4$};
				\node at (1.7,.3) (r2) {$\rho_1$};
				\node at (.3,.3) (r3) {$\rho_2$};
				\draw[dashed] (.25,0) arc[radius=.25, start angle = 0, end angle = 90];
				\draw[dashed] (0,1.75) arc[radius=.25, start angle = 270, end angle = 360];
				\draw[dashed] (1.75,0) arc[radius=.25, start angle = 180, end angle = 90];
				\draw[dashed] (1.75,2) arc[radius=.25, start angle = 180, end angle = 270];
				\draw[color = red,thick] (.25,0) to (1.75,0);
				\draw[color = red,thick] (0,0.25) to (0,1.75);
				\draw[color = red,thick] (.25,2) to (1.75,2);
				\draw[color = red,thick] (2,0.25) to (2,1.75);
		
				\draw[color = blue,thick] (0,1) to (2,1);
			\end{scope}
			\begin{scope}[shift={(0,0)}]
				\node at (.25,1.7) (r4) {$\rho_3$};
				\node at (1.7,1.7) (r1) {$\rho_4$};
				\node at (1.7,.3) (r2) {$\rho_1$};
				\node at (.3,.3) (r3) {$\rho_2$};
				\draw[dashed] (.25,0) arc[radius=.25, start angle = 0, end angle = 90];
				\draw[dashed] (0,1.75) arc[radius=.25, start angle = 270, end angle = 360];
				\draw[dashed] (1.75,0) arc[radius=.25, start angle = 180, end angle = 90];
				\draw[dashed] (1.75,2) arc[radius=.25, start angle = 180, end angle = 270];
				\draw[color = red,thick] (.25,0) to (1.75,0);
				\draw[color = red,thick] (0,0.25) to (0,1.75);
				\draw[color = red,thick] (.25,2) to (1.75,2);
				\draw[color = red,thick] (2,0.25) to (2,1.75);
		
				\draw[color = blue,thick] (0,1) to (2,1);

				\draw [line width = 1.5pt]
				(0,0.25) to (0,1.75);
			\end{scope}
			\begin{scope}[shift={(0,2)}]
				\node at (2.2,0.9) (a2) {$a$};
				\node at (1.7,.3) (r2) {$\rho_1$};
				\node at (.3,.3) (r3) {$\rho_2$};
				\draw[dashed] (.25,0) arc[radius=.25, start angle = 0, end angle = 90];
				\draw[dashed] (1.75,0) arc[radius=.25, start angle = 180, end angle = 90];
				\draw[color = red,thick] (.25,0) to (1.75,0);
				\draw[color = red,thick] (0,0.25) to (0,1);
				\draw[color = red,thick] (2,0.25) to (2,1);
		
				\draw[color = blue,thick] (0,1) to (2,1);
			\end{scope}
			\begin{scope}[shift={(-2,2)}]
				\node at (-0.2,0.9) (a2) {$a$};
				\node at (1.7,.3) (r2) {$\rho_1$};
				\node at (.3,.3) (r3) {$\rho_2$};
				\draw[dashed] (.25,0) arc[radius=.25, start angle = 0, end angle = 90];
				\draw[dashed] (1.75,0) arc[radius=.25, start angle = 180, end angle = 90];
				\draw[color = red,thick] (.25,0) to (1.75,0);
				\draw[color = red,thick] (0,0.25) to (0,1);
				\draw[color = red,thick] (2,0.25) to (2,1);
		
				\draw[color = blue,thick] (0,1) to (2,1);
			\end{scope}
		\end{tikzpicture}
		\caption [An illustration of a disk where move (3) can be performed.] {\textbf{An illustration of a disk where move (3) can be performed.} There is a cut in the disk drawn with black, and the edges on either side can be glued together.}
		\label{fig:move-3-illustration}
\end{figure}

	Move (3) can be performed on the disk illustrated in Figure~\ref{fig:move-3-illustration}, by gluing the cut drawn with black. The operation $m^0_8(a, \rho_{23}, \rho_2, \rho_1, \rho_{4123}, \rho_2, \rho_1, \rho_{41}) = U^2 a$ induces the operation $m^1_6(a, \rho_{23}, \rho_2, \rho_{12}, \rho_1, \rho_{41}) = U^2 a$.
	
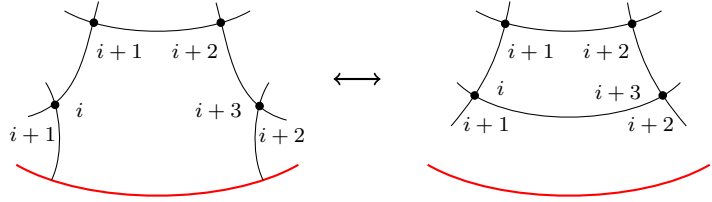
\begin{figure}
	\centering
	\begin{tikzpicture}[scale=1.25]
	\begin{pgfonlayer}{nodelayer}
		\node [style=none] (32) at (-0.5, 0) {};
		\node [style=none] (33) at (0.5, 0) {};
		\node [style=none] (34) at (-1.25, -1.825) {};
		\node [style=none] (35) at (-7.25, -1.825) {};
		\node [style=none] (64) at (-6.5, -2.15) {};
		\node [style=none] (65) at (-2, -2.15) {};
		\node [style=none] (66) at (-6.625, -0.075) {};
		\node [style=none] (67) at (-5.5, 1.65) {};
		\node [style=none] (68) at (-1.875, -0.075) {};
		\node [style=none] (69) at (-7, -0.8) {};
		\node [style=none] (70) at (-1.5, -0.875) {};
		\node [style=none] (71) at (-3, 1.625) {};
		\node [style=none] (72) at (-6.225, 1.45) {};
		\node [style=none] (73) at (-2.275, 1.45) {};
		\node [style=number-labels] (87) at (-5.05, 0.575) {$i+1$};
		\node [style=number-labels] (88) at (-3.45, 0.575) {$i+2$};
		\node [style=number-labels] (89) at (-2.95, -0.675) {$i+3$};
		\node [style=number-labels] (90) at (-5.9, -0.65) {$i$};
		\node [style=none] (91) at (7.5, -1.825) {};
		\node [style=none] (92) at (1.5, -1.825) {};
		\node [style=none] (95) at (2.125, -0.075) {};
		\node [style=none] (96) at (3.25, 1.65) {};
		\node [style=none] (97) at (6.875, -0.075) {};
		\node [style=none] (98) at (2, -1) {};
		\node [style=none] (99) at (7, -1) {};
		\node [style=none] (100) at (5.75, 1.625) {};
		\node [style=none] (101) at (2.525, 1.45) {};
		\node [style=none] (102) at (6.475, 1.45) {};
		\node [style=number-labels] (103) at (3.7, 0.575) {$i+1$};
		\node [style=number-labels] (104) at (5.3, 0.575) {$i+2$};
		\node [style=number-labels] (105) at (5.55, -0.3) {$i+3$};
		\node [style=number-labels] (106) at (3.05, -0.2) {$i$};
		\node [style=dot] (107) at (-5.6, 1.2) {};
		\node [style=dot] (108) at (-2.9, 1.2) {};
		\node [style=dot] (109) at (-6.425, -0.55) {};
		\node [style=dot] (110) at (-2.075, -0.575) {};
		\node [style=dot] (111) at (3.15, 1.175) {};
		\node [style=dot] (112) at (5.85, 1.175) {};
		\node [style=dot] (113) at (2.5, -0.35) {};
		\node [style=dot] (114) at (6.5, -0.35) {};
		\node [style=number-labels] (115) at (-6.9, -1.2) {$i+1$};
		\node [style=number-labels] (116) at (-1.6, -1.225) {$i+2$};
		\node [style=number-labels] (117) at (2.75, -0.975) {$i+1$};
		\node [style=number-labels] (118) at (6.25, -0.975) {$i+2$};
	\end{pgfonlayer}
	\begin{pgfonlayer}{edgelayer}
		\draw [style=double-sided] (32.center) to (33.center);
		\draw [style=red, bend right, looseness=0.75] (35.center) to (34.center);
		\draw [bend left, looseness=0.75] (66.center) to (64.center);
		\draw [in=255, out=15] (69.center) to (67.center);
		\draw [bend right, looseness=0.75] (68.center) to (65.center);
		\draw [in=165, out=-75] (71.center) to (70.center);
		\draw [bend right, looseness=0.75] (72.center) to (73.center);
		\draw [style=red, bend right, looseness=0.75] (92.center) to (91.center);
		\draw [bend right=15] (98.center) to (96.center);
		\draw [bend right=15] (100.center) to (99.center);
		\draw [bend right, looseness=0.75] (101.center) to (102.center);
		\draw [bend right=45, looseness=0.75] (95.center) to (97.center);
	\end{pgfonlayer}
\end{tikzpicture}
	\caption[Move (3) on module tiling patterns.]
	{\textbf{Move (3) on module tiling patterns.} The move creates a new internal face, and the labelings remain unchanged. The labelings around the newly formed edge remain consistent by virtue of the labeling conventions around vertices, as depicted in the figure.}
	\label{fig:move-3}
\end{figure}

	In terms of module tiling patterns, the move takes a five-edged red face and glues the two edges going out to the boundary to create a new internal face. This, along with a schema for the labelings, is illustrated in Figure~\ref{fig:move-3}.

	The three moves have evident inverses. An inverse of move (1) splits a module tiling pattern along an edge visible from both the red and blue boundaries, and yields two smaller module tiling patterns. An inverse of move (2) removes a red vertex. An inverse of move (3) punctures an internal face that shares an edge with the red boundary. We note that the three inverses result in simpler module tiling patterns, either in terms of the total number of vertices, or internal faces.
	
	The following proposition computes the index of an immersion resulting from performing any of the three moves.
	
\begin{proposition}
	Let $B \in \pi_2(\x, \y)$ and $B' \in \pi_2(\y, \z)$ be homology classes, and $\vec a = (a_1, \dots, a_n)$ and $\vec a' = (a'_1, \dots, a'_m)$ be sequences of Reeb chords compatible with $B$ and $B'$ respectively. Let $w$ and $w'$ be their respective number of simple Reeb orbits. 
	
	If $(\tilde B, \tilde {\vec a}, \tilde w)$ is obtained from performing move (1) on immersed disks corresponding to $(B, \vec a, w)$ and $(B', \vec a', w')$, then 
	\[\ind(\tilde B, \tilde {\vec a}, \tilde w) = \ind(B, \vec a, w) + \ind(B', \vec a', w') - 1.\]
	If $(\tilde B, \tilde {\vec a}, \tilde w)$ is obtained from some move (2) or move (3) on an immersed disk corresponding to $(B, \vec a, w)$, then
	\[\ind(\tilde B, \tilde {\vec a}, \tilde w) = \ind(B, \vec a, w).\]
	
	In particular, if $\ind(B, \vec a, w) = \ind(B', \vec a', w') = 1$,
	then $\ind(\tilde B, \tilde {\vec a}, \tilde w) = 1$.
\end{proposition}

\begin{proof}
	Move (1): The homology classes and weights are additive under move (1): $\tilde B$ is $B + B'$ and $\tilde w$ is $w + w'$. The chord sequence $\tilde {\vec a}$ is given by $(a_1, \dots, a_n a'_1, \dots, a'_m)$.
	
	The additivity of the index, Proposition~\ref{prop:index-additivity}, says that $\ind(B + B', (\vec a, \vec a'), w + w')$ equals $\ind(B, \vec a, w) + \ind(B, \vec a', w')$. It suffices to show that $\ind(B + B', \tilde {\vec a}, w + w')$ is one less than $\ind(B + B', (\vec a, \vec a'), w + w')$. This follows from the right hand side of the embedded index formula, Equation~(\ref{prop:embedded-index-formula}). The Euler measure, point measures, and weight remain the same. The terms $\iota((\vec a, \vec a'))$ and $\iota(\tilde{\vec a})$ can both be expressed as the Maslov component of the same product of gradings using Lemma~\ref{lem:iota-maslov}, and so are equal. The only term that changes is the length of the chord sequence, and it drops as desired by one.
	
	Move (2): Any move (2) adds a copy of the torus to the homology class, so $\tilde B = B + [\Sigma]$. The new chord sequence is $\tilde{\vec a} = (a_1, \dots, a_i \rho_{j + 1}, \rho_j, \rho_{j - 1}, \rho_{j - 2} a_{i + 1}, \dots, a_n)$ for some $1 \leq i < n$ with $a_i = b \rho_j$. The weight is left unchanged.
	
	The Euler measure $e(B + [\Sigma]) = e(B) - 1$, while the point measures of $B + [\Sigma]$ are one more than that of $B$.
	To compute the change in $\iota$, note that
	\begin{align*}
		\gr'(a_1) &\cdots \gr'(a_i \rho_{j + 1}) \gr'(\rho_j) \gr'(\rho_{j - 1}) \gr'(\rho_{j - 2} a_{i + 1}) \cdots \gr'(a_n) \\
		& = \gr'(a_1) \cdots \gr'(a_i) \left(\gr'(\rho_{j + 1}) \gr'(\rho_j) \gr'(\rho_{j - 1}) \gr'(\rho_{j - 2})\right) \gr'(a_{i + 1}) \cdots \gr'(a_n) \\
		& = \gr'(a_1) \cdots \gr'(a_i) (-3; 1, 1, 1, 1) \gr'(a_{i + 1}) \cdots \gr'(a_n)
	\end{align*}
	The element $(-3; 1, 1, 1, 1)$ is central, and in particular drops the Maslov component by three. So, by Lemma~\ref{lem:iota-maslov},
	$\iota(\tilde{\vec a}) = \iota(\vec a) - 3$. Putting it together,
	\begin{align*}
		\ind(B + [\Sigma], \tilde{\vec a}, \tilde w) &= e(B + [\Sigma]) + n_{\x}(B + [\Sigma]) + n_{\y}(B + [\Sigma]) + |\tilde{\vec a}| + \iota(\tilde{\vec a}) + \tilde w \\
		&= (e(B) - 1) + (n_{\x}(B) + 1) + (n_{\y}(B) + 1) \\
		&\qquad\qquad\qquad + (|\vec a| + 2) + (\iota(\vec a) - 3) + w \\
		&= e(B) + n_{\x}(B) + n_{\y}(B) + |\vec a| + \iota(\vec a) + w = \ind(B, \vec a, w).
	\end{align*}
	Move (3): The homology class is left unchanged as a result of move (3). The new chord sequence is $\tilde{\vec a} = (a_1, \dots, a_{i - 1} a_{i + 1}, \dots, a_n)$ for some $1 < i < n$ with $|a_i| = 4$; and $\tilde w = w + 1$.
	
	Observe that $\gr'(a_i) = (-1; 1, 1, 1, 1)$ for any Reeb chord $a_i$ with length four. Multiplication with $(-1; 1, 1, 1, 1)$ drops the Maslov component by one, so Lemma~\ref{lem:iota-maslov} implies that $\iota(\vec a) = \iota(\tilde{\vec a}) - 1$. Then,
		\begin{align*}
		\ind(\tilde B, \tilde{\vec a}, \tilde w) &= e(\tilde B) + n_{\x}(\tilde B) + n_{\y}(\tilde B) + |\tilde{\vec a}| + \iota(\tilde{\vec a}) + \tilde w \\
		&= e(B) + n_{\x}(B) + n_{\y}(B) + (|\vec a| - 2) + (\iota(\vec a) + 1) + (w + 1) \\
		&= e(B) + n_{\x}(B) + n_{\y}(B) + |\vec a| + \iota(\vec a) + w = \ind(B, \vec a, w). \qedhere
	\end{align*}
\end{proof}

\subsubsection{Characterization of all module tiling patterns} \label{sec:char-mod-tiling}

We have seen that immersed disks that form acute angles at the intersections of the $\alpha$-arcs and the $\beta$-circles contribute operations to the type A module. The following proposition shows that those with obtuse angles do not.

\begin{proposition}
	\label{prop:solid-torus-obtuse}
	Let $u: \bD \to T^2$ be a disk for $\HD_{st}$ from $a$ to $a$, and suppose further that it maps the corners at $-i$ and $i$ to \emph{obtuse} intersections of the $\alpha$-arcs and $\beta$-circles.
	
	As in Proposition~\ref{prop:acute-angled-disk}, let $B$ be the positive domain in $\pi_2(a, a)$, $\vec a$ the sequence of Reeb chords, and $w$ the number of simple Reeb orbits corresponding to $u$.	
	Then, $\ind(B, \vec a, w) > 1$.
\end{proposition}

\begin{proof}
	This too is a classical result in Heegaard Floer theory, c.f. \cite{Petkova2013, Hom2013}.

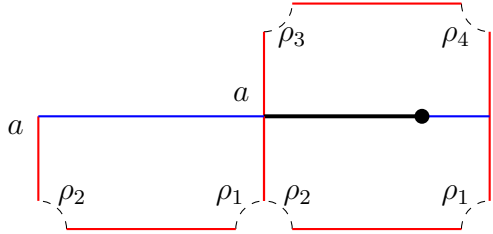
\begin{figure}[h!tbp]
	\centering
		\begin{tikzpicture}[scale=1.5]
			\begin{scope}[shift={(2,0)}]
				\node at (.25,1.7) (r4) {$\rho_3$};
				\node at (1.7,1.7) (r1) {$\rho_4$};
				\node at (1.7,.3) (r2) {$\rho_1$};
				\node at (.3,.3) (r3) {$\rho_2$};
				\draw[dashed] (.25,0) arc[radius=.25, start angle = 0, end angle = 90];
				\draw[dashed] (0,1.75) arc[radius=.25, start angle = 270, end angle = 360];
				\draw[dashed] (1.75,0) arc[radius=.25, start angle = 180, end angle = 90];
				\draw[dashed] (1.75,2) arc[radius=.25, start angle = 180, end angle = 270];
				\draw[color = red,thick] (.25,0) to (1.75,0);
				\draw[color = red,thick] (0,0.25) to (0,1.75);
				\draw[color = red,thick] (.25,2) to (1.75,2);
				\draw[color = red,thick] (2,0.25) to (2,1.75);
		
				\draw[color = blue,thick] (0,1) to (2,1);
				\draw [line width = 1.5pt] (0,1) to (1.4,1);
				\node at (1.4,1)[circle,fill,inner sep=2pt]{};
			\end{scope}
			\begin{scope}[shift={(0,0)}]
				\node at (1.8,1.2) (a2) {$a$};
				\node at (-0.2,0.9) (a2) {$a$};
				\node at (1.7,.3) (r2) {$\rho_1$};
				\node at (.3,.3) (r3) {$\rho_2$};
				\draw[dashed] (.25,0) arc[radius=.25, start angle = 0, end angle = 90];
				\draw[dashed] (1.75,0) arc[radius=.25, start angle = 180, end angle = 90];
				\draw[color = red,thick] (.25,0) to (1.75,0);
				\draw[color = red,thick] (0,0.25) to (0,1);
				\draw[color = red,thick] (2,0.25) to (2,1);
		
				\draw[color = blue,thick] (0,1) to (2,1);
			\end{scope}
		\end{tikzpicture}
		\caption [An illustration of a disk with an obtuse angle at the \texorpdfstring{$\alpha$}{alpha}-arc and \texorpdfstring{$\beta$}{beta}-circle boundary.] {\textbf{An illustration of a disk with an obtuse angle at the $\alpha$-arc and $\beta$-circle boundary.} There is a cut in the disk indicated by the black line, and a branch point at the end of the cut. When the boundary branch point goes out to the $\alpha$-arc, the disk decomposes into two disks. The chord sequences for the two disks in order are $(\rho_2, \rho_{12}, \rho_{1})$ and $(\rho_4, \rho_3)$.}
		\label{fig:solid-torus-index-2-disk}
\end{figure}
	A disk with an obtuse intersection along the $\alpha$-arc and $\beta$-circle admits a cut along the $\beta$-circle. Suppose for simplicity that there is one such obtuse intersection. When the boundary branch point goes out to the $\alpha$-arc, the disk decomposes into two disks with acute intersections along the $\alpha$-arcs and $\beta$-circles. These we have seen to have index one. That is, we can express $(B, \vec a, w)$ as $(B' + B'', (\vec a', \vec a''), w' + w'')$, where
	$\ind(B', \vec a', w') = \ind(B'', \vec a'', w'') = 1$. The additivity of the index tells us that $\ind(B, \vec a, w) > 1$.
\end{proof}

We are ready to characterize all operations on the type A module $\CFAm(\HD_{st})$ in terms of two simple module tiling patterns.

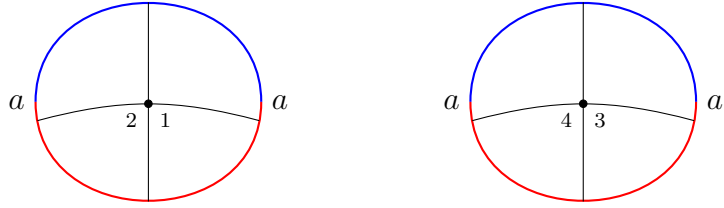
\begin{figure}
	\centering
	\begin{tikzpicture}
	\begin{pgfonlayer}{nodelayer}
		\node [style=none] (0) at (-3, 0) {};
		\node [style=none] (1) at (3, 0) {};
		\node [style=none] (2) at (-2.95, -0.5) {};
		\node [style=none] (3) at (2.95, -0.5) {};
		\node [style=none] (4) at (0, 2.65) {};
		\node [style=none] (5) at (0, -2.65) {};
		\node [style=none] (6) at (-3.5, 0) {$a$};
		\node [style=none] (7) at (3.5, 0) {$a$};
		\node [style=number-labels] (8) at (-0.45, -0.475) {$2$};
		\node [style=number-labels] (9) at (0.45, -0.475) {$1$};
		\node [style=dot] (10) at (0, -0.05) {};
	\end{pgfonlayer}
	\begin{pgfonlayer}{edgelayer}
		\draw [style=blue, bend left=90, looseness=1.50] (0.center) to (1.center);
		\draw [style=red, bend right=90, looseness=1.50] (0.center) to (1.center);
		\draw [bend left=15] (2.center) to (3.center);
		\draw (4.center) to (5.center);
	\end{pgfonlayer}
\end{tikzpicture}
	\qquad \qquad
	\begin{tikzpicture}
	\begin{pgfonlayer}{nodelayer}
		\node [style=none] (0) at (-3, 0) {};
		\node [style=none] (1) at (3, 0) {};
		\node [style=none] (2) at (-2.95, -0.5) {};
		\node [style=none] (3) at (2.95, -0.5) {};
		\node [style=none] (4) at (0, 2.65) {};
		\node [style=none] (5) at (0, -2.65) {};
		\node [style=none] (6) at (-3.5, 0) {$a$};
		\node [style=none] (7) at (3.5, 0) {$a$};
		\node [style=number-labels] (8) at (-0.45, -0.475) {$4$};
		\node [style=number-labels] (9) at (0.45, -0.475) {$3$};
		\node [style=dot] (10) at (0, -0.05) {};
	\end{pgfonlayer}
	\begin{pgfonlayer}{edgelayer}
		\draw [style=blue, bend left=90, looseness=1.50] (0.center) to (1.center);
		\draw [style=red, bend right=90, looseness=1.50] (0.center) to (1.center);
		\draw [bend left=15] (2.center) to (3.center);
		\draw (4.center) to (5.center);
	\end{pgfonlayer}
\end{tikzpicture}
	\caption[The building blocks for operations for the $0$-framed solid torus.]
	{\textbf{The building blocks for operations for the $0$-framed solid torus.}}
	\label{fig:solid-torus-building-blocks}
\end{figure}

\begin{proposition}
	Any operation on the type A module $\CFAm(\HD_{st})$ can be obtained as the operation associated to a module tiling pattern constructed by starting with the
	two building blocks in Figure~\ref{fig:solid-torus-building-blocks}, and
	performing a finite sequence of the three moves catalogued in
	Section~\ref{sec:three-moves}.
	\label{prop:solid-torus-building-blocks}
\end{proposition}

\begin{proof}
	Let $u: \bD \to T^2$ be an immersed disk from $a$ to $a$ for $\HD_{st}$. There is a neighborhood of the $\beta$-circle on the torus that does not contain the marked point $p$, but contain the two vertices in $\Gamma(\HD_{st})$. In this neighborhood, $u$ lifts to the universal cover
	of the torus. We can read off the graph $\Gamma(u)$ near the blue boundary from this lift. The two possibilities for the structure of this are depicted in Figure~\ref{fig:solid-torus-spine}. The description near the blue-red intersection is forced by virtue of the angle of intersection being acute. We call
	the top-most horizontal arc the spine of the tiling pattern.

\begin{figure}[h!tbp]
	\begin{tikzpicture}
	\begin{pgfonlayer}{nodelayer}
		\node [style=none] (0) at (-6, 0) {};
		\node [style=none] (1) at (6, 0) {};
		\node [style=none] (2) at (-5.95, -0.5) {};
		\node [style=none] (3) at (5.95, -0.5) {};
		\node [style=none] (4) at (-4.75, 1.7) {};
		\node [style=none] (5) at (-4.25, -0.65) {};
		\node [style=none] (6) at (-6.5, 0) {$a$};
		\node [style=none] (7) at (6.5, 0) {$a$};
		\node [style=number-labels] (8) at (-4.6, -0.55) {$2$};
		\node [style=number-labels] (9) at (-4.025, -0.5) {$1$};
		\node [style=none] (10) at (-5.55, -1.5) {};
		\node [style=none] (11) at (5.6, -1.5) {};
		\node [style=none] (12) at (-3.25, 2.275) {};
		\node [style=none] (13) at (-3, -0.5) {};
		\node [style=none] (14) at (-1.75, 2.55) {};
		\node [style=none] (15) at (-1.75, -0.5) {};
		\node [style=none] (16) at (1.75, 2.55) {};
		\node [style=none] (17) at (1.75, -0.5) {};
		\node [style=none] (18) at (4.5, 1.8) {};
		\node [style=none] (19) at (4.25, -0.75) {};
		\node [style=none] (20) at (3.25, 2.275) {};
		\node [style=none] (21) at (3, -0.5) {};
		\node [style=none] (22) at (0, 1.25) {};
		\node [style=none] (23) at (0, 1.25) {};
		\node [style=none] (24) at (0, 1.25) {$\cdots$};
		\node [style=number-labels] (25) at (-3.3, -0.425) {$2$};
		\node [style=number-labels] (26) at (-2.025, -0.35) {$2$};
		\node [style=number-labels] (27) at (1.5, -0.325) {$2$};
		\node [style=number-labels] (28) at (2.75, -0.4) {$2$};
		\node [style=number-labels] (29) at (4.025, -0.5) {$2$};
		\node [style=number-labels] (30) at (-2.75, -0.4) {$1$};
		\node [style=number-labels] (31) at (-1.5, -0.325) {$1$};
		\node [style=number-labels] (32) at (2.025, -0.35) {$1$};
		\node [style=number-labels] (33) at (3.3, -0.425) {$1$};
		\node [style=number-labels] (34) at (4.65, -0.55) {$1$};
		\node [style=dot] (35) at (4.4, -0.25) {};
		\node [style=dot] (36) at (3.1, -0.15) {};
		\node [style=dot] (37) at (1.8, -0.075) {};
		\node [style=dot] (38) at (-1.8, -0.075) {};
		\node [style=dot] (39) at (-3.1, -0.15) {};
		\node [style=dot] (40) at (-4.425, -0.275) {};
	\end{pgfonlayer}
	\begin{pgfonlayer}{edgelayer}
		\draw [bend left=15, looseness=0.50] (2.center) to (3.center);
		\draw [bend right=15, looseness=0.75] (4.center) to (5.center);
		\draw [style=red, bend right=15, looseness=0.75] (0.center) to (10.center);
		\draw [style=red, bend left=15, looseness=0.50] (1.center) to (11.center);
		\draw [bend right=15, looseness=0.50] (12.center) to (13.center);
		\draw [bend right=15, looseness=0.25] (14.center) to (15.center);
		\draw [bend left=15, looseness=0.25] (16.center) to (17.center);
		\draw [bend left=15, looseness=0.50] (20.center) to (21.center);
		\draw [bend left=15, looseness=0.75] (18.center) to (19.center);
		\draw [style=blue, bend left=90, looseness=0.75] (0.center) to (1.center);
	\end{pgfonlayer}
\end{tikzpicture} \qquad \qquad
	\begin{tikzpicture}
	\begin{pgfonlayer}{nodelayer}
		\node [style=none] (0) at (-6, 0) {};
		\node [style=none] (1) at (6, 0) {};
		\node [style=none] (2) at (-5.95, -0.5) {};
		\node [style=none] (3) at (5.95, -0.5) {};
		\node [style=none] (4) at (-4.75, 1.7) {};
		\node [style=none] (5) at (-4.25, -0.65) {};
		\node [style=none] (6) at (-6.5, 0) {$a$};
		\node [style=none] (7) at (6.5, 0) {$a$};
		\node [style=number-labels] (8) at (-4.6, -0.55) {$4$};
		\node [style=number-labels] (9) at (-4.025, -0.5) {$3$};
		\node [style=none] (10) at (-5.55, -1.5) {};
		\node [style=none] (11) at (5.6, -1.5) {};
		\node [style=none] (12) at (-3.25, 2.275) {};
		\node [style=none] (13) at (-3, -0.5) {};
		\node [style=none] (14) at (-1.75, 2.55) {};
		\node [style=none] (15) at (-1.75, -0.5) {};
		\node [style=none] (16) at (1.75, 2.55) {};
		\node [style=none] (17) at (1.75, -0.5) {};
		\node [style=none] (18) at (4.5, 1.8) {};
		\node [style=none] (19) at (4.25, -0.75) {};
		\node [style=none] (20) at (3.25, 2.275) {};
		\node [style=none] (21) at (3, -0.5) {};
		\node [style=none] (22) at (0, 1.25) {};
		\node [style=none] (23) at (0, 1.25) {};
		\node [style=none] (24) at (0, 1.25) {$\cdots$};
		\node [style=number-labels] (25) at (-3.3, -0.425) {$4$};
		\node [style=number-labels] (26) at (-2.025, -0.35) {$4$};
		\node [style=number-labels] (27) at (1.5, -0.325) {$4$};
		\node [style=number-labels] (28) at (2.75, -0.4) {$4$};
		\node [style=number-labels] (29) at (4.025, -0.5) {$4$};
		\node [style=number-labels] (30) at (-2.75, -0.4) {$3$};
		\node [style=number-labels] (31) at (-1.5, -0.325) {$3$};
		\node [style=number-labels] (32) at (2.025, -0.35) {$3$};
		\node [style=number-labels] (33) at (3.3, -0.425) {$3$};
		\node [style=number-labels] (34) at (4.65, -0.55) {$3$};
		\node [style=dot] (35) at (4.4, -0.25) {};
		\node [style=dot] (36) at (3.1, -0.15) {};
		\node [style=dot] (37) at (1.8, -0.075) {};
		\node [style=dot] (38) at (-1.8, -0.075) {};
		\node [style=dot] (39) at (-3.1, -0.15) {};
		\node [style=dot] (40) at (-4.425, -0.275) {};
	\end{pgfonlayer}
	\begin{pgfonlayer}{edgelayer}
		\draw [bend left=15, looseness=0.50] (2.center) to (3.center);
		\draw [bend right=15, looseness=0.75] (4.center) to (5.center);
		\draw [style=red, bend right=15, looseness=0.75] (0.center) to (10.center);
		\draw [style=red, bend left=15, looseness=0.50] (1.center) to (11.center);
		\draw [bend right=15, looseness=0.50] (12.center) to (13.center);
		\draw [bend right=15, looseness=0.25] (14.center) to (15.center);
		\draw [bend left=15, looseness=0.25] (16.center) to (17.center);
		\draw [bend left=15, looseness=0.50] (20.center) to (21.center);
		\draw [bend left=15, looseness=0.75] (18.center) to (19.center);
		\draw [style=blue, bend left=90, looseness=0.75] (0.center) to (1.center);
	\end{pgfonlayer}
\end{tikzpicture}
	\caption[The structure of a module tiling pattern for \texorpdfstring{$\HD_{st}$}{HDst} near the blue boundary.]
	{\textbf{The structure of a module tiling pattern for $\HD_{st}$ near the blue boundary.} The
	horizontal arc is the spine.}
		\label{fig:solid-torus-spine}
\end{figure}
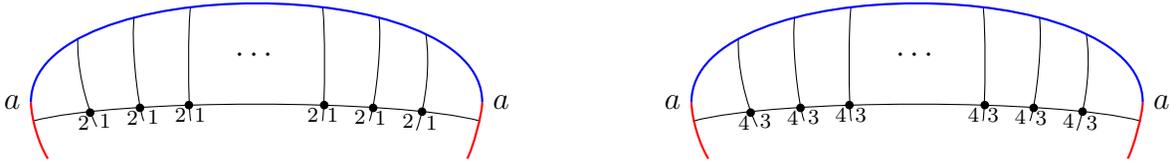

The graph $\Gamma(u)$ is connected and planar. Away from the blue boundary,
	ever vertex of $\Gamma(u)$ is $4$-valent, and every internal face is bound by four edges. This is entirely analogous to immersions of the disk from which the torus algebra arises.
	We prove the result by induction on the number of internal faces and then the number of
	vertices of $\Gamma(u)$.
	
	If there exists an internal face in $\Gamma(u)$, then there exists one
	that is adjacent to a red face. The inverse of move (3) allows us to puncture this
	internal face without modifying the rest of $\Gamma(u)$, reducing the total number of internal faces by one.
	
	Suppose now that there are no internal faces in $\Gamma(u)$. If the number of vertices along the spine of $\Gamma(u)$ is greater than
	one, then there exists an edge between two vertices
	that separates the blue and the red boundary. The
	inverse of move (1) allows us to decompose $\Gamma(u)$ into two smaller module tiling patterns, each with at least one fewer vertex.

	We are left with the case where there are no internal faces in $\Gamma(u)$, and there is exactly one vertex along the spine. Let us consider the graph obtained by deleting the spine. This graph is a tree, where every
	vertex is $4$-valent. It decomposes this into $l$ arcs intersecting each other
	transversely. Suppose we can find an arc that intersects the red
	boundary and exactly one other
	arc, and consider the face this arc would bound with the red boundary. If the face bounded
	contains a single edge, we observe that we can use the inverse of move (2) to get rid of the vertex of intersection in $\Gamma(u)$. Otherwise, the face bounded must contain a tree where every vertex is $4$-valent. We can recurse
	down to this tree to find a vertex that can be eliminated with the inverse of move (2).

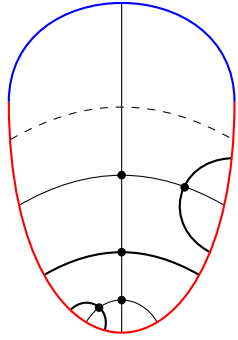
\begin{figure}[h!tbp]
	\centering
	\begin{tikzpicture}
	\begin{pgfonlayer}{nodelayer}
		\node [style=none] (0) at (-3, 1.5) {};
		\node [style=none] (1) at (3, 1.5) {};
		\node [style=none] (2) at (0, 4.125) {};
		\node [style=none] (3) at (0, -4.625) {};
		\node [style=none] (4) at (-2.725, -1.25) {};
		\node [style=none] (5) at (2.725, -1.25) {};
		\node [style=none] (6) at (2.35, -2.5) {};
		\node [style=none] (7) at (2.925, 0) {};
		\node [style=none] (8) at (-2.05, -3.1) {};
		\node [style=none] (9) at (2.05, -3.1) {};
		\node [style=none] (10) at (-0.95, -4.375) {};
		\node [style=none] (11) at (0.95, -4.375) {};
		\node [style=none] (12) at (-1.35, -4.025) {};
		\node [style=none] (13) at (-0.45, -4.575) {};
		\node [style=none] (14) at (-2.95, 0.5) {};
		\node [style=none] (15) at (2.95, 0.5) {};
		\node [style=dot] (17) at (0, -0.45) {};
		\node [style=dot] (18) at (1.675, -0.775) {};
		\node [style=dot] (19) at (0, -2.5) {};
		\node [style=dot] (20) at (0, -3.775) {};
		\node [style=dot] (21) at (-0.6, -3.975) {};
	\end{pgfonlayer}
	\begin{pgfonlayer}{edgelayer}
		\draw (2.center) to (3.center);
		\draw [bend left] (4.center) to (5.center);
		\draw [style=thickedge, bend right=75, looseness=1.50] (7.center) to (6.center);
		\draw [style=thick, bend right] (9.center) to (8.center);
		\draw [bend left=60, looseness=1.25] (10.center) to (11.center);
		\draw [style=blue, bend left=90, looseness=1.50] (0.center) to (1.center);
		\draw [style=thickedge, bend left=75, looseness=1.50] (12.center) to (13.center);
		\draw [style=red, bend right=90, looseness=3.50] (0.center) to (1.center);
		\draw [style=dashededged, bend left] (14.center) to (15.center);
	\end{pgfonlayer}
\end{tikzpicture}
	\vspace{-2em}
	\caption[A graph for \texorpdfstring{$\HD_{st}$}{HDst} illustrating possible simplifications with move (3).]
	{\textbf{A graph for $\HD_{st}$ illustrating possible simplifications with move (3).} In dashed is drawn the spine.
	The three arcs drawn in thick intersect exactly one other arc. Two of them can be
	removed by an inverse of move (2). The remaining arc us to recurse to a
	smaller graph.}
	\label{fig:simplify-torus}
\end{figure}
	
	It remains to find such an arc. Let the $l$ arcs intersect $i_1, \dots, i_l$
	other arcs, where each $i_k \geq 1$.
	The tree has $\sum_{k=1}^l i_k / 2$ vertices and
	$2l$ leaves. The number of edges is $\sum_{k=1}^l (i_k + 1)$. As a tree,
	the number of vertices is one more than the number of edges; so,
	\[\sum_{k=1}^l \frac{i_k}{2} + 2l = \sum_{k=1}^l (i_k + 1) + 1
		\implies \frac{1}{l} \sum_{k = 1}^l i_k = 2 - \frac{2}{l}.\]
	This implies that average of the $i_k$ is less than two. There must thus be a $m$
	with $i_m = 1$. This $m$th arc could intersect the blue boundary. If so, we can consider the tree obtained by deleting it. As there is only one arc that intersects the blue boundary, we are left with a tree where every vertex is 4-valent
	and every arc in the decomposition intersects only the red boundary. Repeating the
	same argument for this sub-tree finds the arc we are looking for.
\end{proof}

\begin{remark}
\label{rmk:solid-torus-simplify}
In the proof of Proposition~\ref{prop:solid-torus-building-blocks}, for the sake of simplicity of exposition we first do all move (1)s to reduce the length of the spine, and then do move (2)s to eliminate the red vertices. However, we could indeed first eliminate all red vertices, and then break up the tiling pattern into simpler pieces.
The enumeration of the \emph{basic} building blocks is a finite calculation---that is, the tiling patterns that do not have red vertices and cannot be broken apart with move (1) into simpler tiling patterns. 
\end{remark}

\begin{figure}[h!tbp]
\centering
\begin{tikzpicture}
	\node at (0,0) (x) {$\circ$};
	\node at (-1,0) (lphant) {};
	\node at (1,0) (rphant) {};
	\draw[->, bend left=80] (x) to node[inner sep=2ex,left]
	{\lab{U\rho_4\otimes\rho_3}} (lphant) to (x);
	\draw[->, bend right=80] (x) to node[inner sep=2ex,right]
	{\lab{\rho_2\otimes\rho_1}} (rphant) to (x);
\end{tikzpicture}
	\caption[A representation of the type A module \texorpdfstring{$\CFAm(\HD_{st})$}{CFA minus of the 0-framed solid torus}.]
	{\textbf{A representation of the type A module $\CFAm(\HD_{st})$.}}
	\label{fig:solid-torus-graph}
\end{figure}

The type A module $\CFAm(\HD_{st})$ can be represented by the graph in Figure~\ref{fig:solid-torus-graph}. The node corresponds to the generator $x$, and each edge represents an operation corresponding to the building blocks depicted in Figure~\ref{fig:solid-torus-building-blocks}; namely, $m^0_3(x, \rho_2, \rho_1) = x$ and $m^0_3(x, \rho_4, \rho_3) = Ux$. Every other operation can be generated from the edges using the three moves.

\subsubsection{The \texorpdfstring{$\Ainf$}{A infinity} relations}

The description of the operations for the type A module $\CFAm(\HD_{st})$ in terms of counting the module tiling patterns allows us to prove that the module satisfies the $\Ainf$ structure relations.

\begin{theorem}
	The type A module $\CFAm(\HD_{st})$, with operations $\{m^w_{1+n}\}$ given by counting its module tiling patterns of chord sequence length $n$ and weight $w$, satisfies the $\Ainf$ structure relations.
	\label{prop:solid-torus-ainf}
\end{theorem}

\begin{proof}
	We verify the $\Ainf$ relation for a fixed input $(x, a_1, \dots, a_n)$ and weight $w$. Recall that the $\Ainf$ relation states that the sum of all operations associated to weighted trees with two internal vertices is zero. We distinguish the operations corresponding to the two vertices as the inner and outer operations; or respectively as the first and second operations. We show the $\Ainf$ relation holds by pairing non-zero terms that appear in this $\Ainf$ relation: over $\FF_2$, then, the terms \emph{cancel}.

	It suffices to consider the case when all the $a_i$ are basic algebra elements, as the operations extend linearly.	 Let us consider the case when some $a_i$ is an idempotent. Recall that the torus algebra and our type A modules are strictly unital; that is, $\mu_m^w(a_1, \dots, a_m)$ and
	$m^w_m(x, a_1, \dots, a_{m-1})$ are both zero if any of the $a_i$ is an idempotent and $(m,w) \neq (2,0)$. The only possibility for an idempotent in a sequence of inputs $(x, a_1, \dots, a_n)$ corresponding to a non-zero term in the $\Ainf$ relation is then for just one $a_i$ to be an idempotent. If $a_1$ is an idempotent, the cancelling terms consist of doing a $\mu_2^0$ followed by a module operation and a $m_2^0$ followed by a module operation; if $a_n$ is an idempotent, we first do a module operation and then $m_2^0$ to cancel against $\mu_2^0$ followed by a module operation. For example, in the weight $0$ $\Ainf$ relation with input $(x, \rho_4, \rho_3, \iota_1)$, the following two non-zero terms cancel
	\[m^0_3(x, \rho_4, \mu_2^0(\rho_3, \iota_1)) = m^0_2(m^0_3(x, \rho_4, \rho_3), \iota_1).\]
		If $a_i$ is an idempotent for $1 < i < n$, the cancelling terms consist of a $\mu^0_2$ to the left and the right of the idempotent.
	
	We now consider the case where all
	the $a_i$ are Reeb elements. The following is how the non-zero
	terms cancel. Each term appears twice in the list to illustrate the pairing.
	
	First, we consider the cases where both internal vertices in our tree for the $\Ainf$ relation correspond to module operations.
	
	\begin{enumerate}[wide,label=(T-\arabic*)]
		\item \label{term:st1} Consider terms of the form $m_{1+n-j}^{w-w'}\left(m_{1+j}^{w'}\left(x,
		a_1,\dots,a_j\right),a_{j+1},\dots,a_n\right)$, where $a_j a_{j+1}\neq0$.
		These cancel against $m_{n}^w\left(
		x,a_1,\dots,\mu_2^0\left(a_j,a_{j+1}\right),\dots,a_n\right)$.
		Note that this is precisely the operation that move (1) produces, and is thus non-zero. In terms of the module tiling patterns it corresponds to concatenation as in Figure~\ref{fig:move-1}. For example, in the weight $1$ $\Ainf$ relation with inputs $x, \rho_{23}, \rho_2, \rho_{12}, \rho_{1}, \rho_{41}, \rho_2, \rho_1$, the following two non-zero terms cancel
		\[m^0_3(m^1_6(x, \rho_{23}, \rho_2, \rho_{12}, \rho_{1}, \rho_{41}), \rho_2, \rho_1) = m^1_7(x, \rho_{23}, \rho_2, \rho_{12}, \rho_{1}, \mu_2^0(\rho_{41}, \rho_2), \rho_1).\]
		
		\item \label{term:st2} Consider terms of the form $m_{1+n-j}^{w-w'}\left(m_{1+j}^{w'}(x,
		a_1,\dots,a_j),a_{j+1},\dots,a_n\right)$, where $a_j a_{j+1}=0$.
		
		Let $\Gamma_1$ and $\Gamma_2$ be the module tiling patterns for the first and second operations respectively. Because $a_ja_{j+1} = 0$, the two must have differing spines as in Figure~\ref{fig:solid-torus-spine}.
		We view the two tiling patterns as
		sharing a \emph{common} blue boundary. We can glue the two patterns along this common boundary, while identifying edges that touch this blue boundary. These edges can be located as the pre-image of the edge in $\Gamma^\beta(\HD_{st})$. We contract these edges as in our construction of the module tiling pattern associated to an immersion, and delete parallel edges. Contracting an edge identifies two blue vertices as one, and we say the new vertex is red. In other words, we glue the immersed disks for the two operations along their common blue boundary.

		For this case, suppose the length of the blue boundaries of $\Gamma_1$ and $\Gamma_2$ match. 
		In this case, the two patterns glue together perfectly. This yields a tiling pattern with no blue boundary; and it satisfies the properties of the torus algebra's tiling patterns. Placing the root at the left of the identified spines produces an algebra tiling pattern corresponding to a centered algebra operation,
		and the term cancels against
		$m_2^{0}\left(x,\mu_n^w\left(a_1,\dots,a_n\right)\right)$. The output of the centered algebra operation is a power of $U$ times the left idempotent corresponding to the first label after the root. The $m_2^0$ product of this idempotent with the generator $x$ is non-zero because $\Gamma_1$ is a disk from $x$. The powers of $U$ line up as the identification of the two spines produces an arc where every vertex carries precisely a labeling of $U$.
		A schema for this is
		illustrated in Figure~\ref{fig:solid-torus-cancellations-2a}.

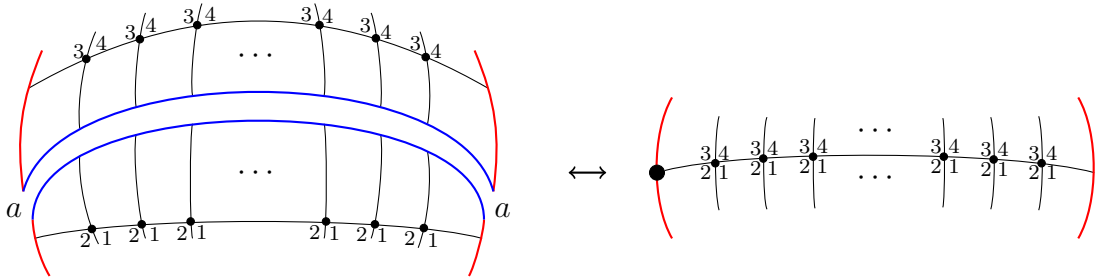
\begin{figure}[h!tbp]
	\centering
	\begin{tikzpicture}
	\begin{pgfonlayer}{nodelayer}
		\node [style=none] (0) at (-14.75, -1.25) {};
		\node [style=none] (1) at (-2.75, -1.25) {};
		\node [style=none] (2) at (-14.7, -1.75) {};
		\node [style=none] (3) at (-2.8, -1.75) {};
		\node [style=none] (4) at (-13.5, 0.45) {};
		\node [style=none] (5) at (-13, -1.9) {};
		\node [style=none] (6) at (-15.25, -1) {$a$};
		\node [style=none] (7) at (-2.25, -1) {$a$};
		\node [style=number-labels] (8) at (-13.35, -1.8) {$2$};
		\node [style=number-labels] (9) at (-12.775, -1.75) {$1$};
		\node [style=none] (10) at (-14.3, -2.75) {};
		\node [style=none] (11) at (-3.15, -2.75) {};
		\node [style=none] (12) at (-12, 1.025) {};
		\node [style=none] (13) at (-11.75, -1.75) {};
		\node [style=none] (14) at (-10.5, 1.3) {};
		\node [style=none] (15) at (-10.5, -1.75) {};
		\node [style=none] (16) at (-7, 1.3) {};
		\node [style=none] (17) at (-7, -1.75) {};
		\node [style=none] (18) at (-4.25, 0.55) {};
		\node [style=none] (19) at (-4.5, -2) {};
		\node [style=none] (20) at (-5.5, 1.025) {};
		\node [style=none] (21) at (-5.75, -1.75) {};
		\node [style=none] (22) at (-8.75, 0) {};
		\node [style=none] (23) at (-8.75, 0) {};
		\node [style=none] (24) at (-8.75, 0) {$\cdots$};
		\node [style=number-labels] (25) at (-12.05, -1.675) {$2$};
		\node [style=number-labels] (26) at (-10.775, -1.6) {$2$};
		\node [style=number-labels] (27) at (-7.25, -1.575) {$2$};
		\node [style=number-labels] (28) at (-6, -1.65) {$2$};
		\node [style=number-labels] (29) at (-4.725, -1.75) {$2$};
		\node [style=number-labels] (30) at (-11.5, -1.65) {$1$};
		\node [style=number-labels] (31) at (-10.25, -1.575) {$1$};
		\node [style=number-labels] (32) at (-6.725, -1.6) {$1$};
		\node [style=number-labels] (33) at (-5.45, -1.675) {$1$};
		\node [style=number-labels] (34) at (-4.1, -1.8) {$1$};
		\node [style=none] (35) at (-15, -0.5) {};
		\node [style=none] (36) at (-2.5, -0.5) {};
		\node [style=none] (37) at (-14.5, 3.25) {};
		\node [style=none] (38) at (-3, 3.25) {};
		\node [style=none] (39) at (-13.5, 1.175) {};
		\node [style=none] (40) at (-12, 1.75) {};
		\node [style=none] (41) at (-10.5, 2.05) {};
		\node [style=none] (42) at (-7, 2.05) {};
		\node [style=none] (43) at (-5.5, 1.75) {};
		\node [style=none] (44) at (-4.25, 1.3) {};
		\node [style=none] (45) at (-14.85, 2.25) {};
		\node [style=none] (46) at (-2.65, 2.25) {};
		\node [style=none] (47) at (-13.075, 3.75) {};
		\node [style=none] (48) at (-11.75, 4.25) {};
		\node [style=none] (49) at (-10.25, 4.5) {};
		\node [style=none] (50) at (-7.25, 4.5) {};
		\node [style=none] (51) at (-5.75, 4.25) {};
		\node [style=none] (52) at (-4.45, 3.75) {};
		\node [style=none] (53) at (-8.75, 3.1) {$\cdots$};
		\node [style=number-labels] (54) at (-4.025, 3.225) {$4$};
		\node [style=number-labels] (55) at (-4.575, 3.45) {$3$};
		\node [style=number-labels] (56) at (-5.35, 3.725) {$4$};
		\node [style=number-labels] (57) at (-5.9, 3.925) {$3$};
		\node [style=number-labels] (58) at (-6.9, 4.15) {$4$};
		\node [style=number-labels] (59) at (-7.475, 4.225) {$3$};
		\node [style=number-labels] (60) at (-10.025, 4.225) {$4$};
		\node [style=number-labels] (61) at (-10.6, 4.15) {$3$};
		\node [style=number-labels] (62) at (-11.6, 3.925) {$4$};
		\node [style=number-labels] (63) at (-12.15, 3.725) {$3$};
		\node [style=number-labels] (64) at (-12.925, 3.45) {$4$};
		\node [style=number-labels] (65) at (-13.525, 3.225) {$3$};
		\node [style=none] (66) at (-0.5, 0) {};
		\node [style=none] (67) at (0.5, 0) {};
		\node [style=none] (68) at (2.25, 2) {};
		\node [style=none] (69) at (2.25, -1.75) {};
		\node [style=bigdot] (74) at (1.85, 0) {};
		\node [style=none] (75) at (13.475, 0) {};
		\node [style=number-labels] (78) at (3.15, 0) {$2$};
		\node [style=number-labels] (79) at (3.625, 0) {$1$};
		\node [style=number-labels] (90) at (4.475, 0.125) {$2$};
		\node [style=number-labels] (91) at (5.775, 0.175) {$2$};
		\node [style=number-labels] (92) at (9.275, 0.175) {$2$};
		\node [style=number-labels] (93) at (10.6, 0.1) {$2$};
		\node [style=number-labels] (94) at (11.875, 0) {$2$};
		\node [style=number-labels] (95) at (4.925, 0.125) {$1$};
		\node [style=number-labels] (96) at (6.25, 0.175) {$1$};
		\node [style=number-labels] (97) at (9.75, 0.175) {$1$};
		\node [style=number-labels] (98) at (11.05, 0.1) {$1$};
		\node [style=number-labels] (99) at (12.35, 0) {$1$};
		\node [style=none] (100) at (13.05, 2) {};
		\node [style=none] (101) at (13.05, -1.75) {};
		\node [style=none] (102) at (3.5, 1.5) {};
		\node [style=none] (103) at (3.5, -1) {};
		\node [style=none] (104) at (4.775, 1.475) {};
		\node [style=none] (105) at (4.775, -0.975) {};
		\node [style=none] (106) at (6.05, 1.45) {};
		\node [style=none] (107) at (6.05, -0.95) {};
		\node [style=none] (108) at (9.45, 1.45) {};
		\node [style=none] (109) at (9.45, -0.925) {};
		\node [style=none] (110) at (10.725, 1.475) {};
		\node [style=none] (111) at (10.725, -0.975) {};
		\node [style=none] (112) at (12, 1.5) {};
		\node [style=none] (113) at (12, -1) {};
		\node [style=none] (114) at (7.725, 1.1) {$\cdots$};
		\node [style=none] (115) at (7.725, -0.15) {$\cdots$};
		\node [style=number-labels] (116) at (3.15, 0.525) {$3$};
		\node [style=number-labels] (117) at (3.625, 0.525) {$4$};
		\node [style=number-labels] (118) at (4.475, 0.65) {$3$};
		\node [style=number-labels] (119) at (5.775, 0.7) {$3$};
		\node [style=number-labels] (120) at (9.275, 0.7) {$3$};
		\node [style=number-labels] (121) at (10.6, 0.625) {$3$};
		\node [style=number-labels] (122) at (11.875, 0.525) {$3$};
		\node [style=number-labels] (123) at (4.925, 0.65) {$4$};
		\node [style=number-labels] (124) at (6.25, 0.7) {$4$};
		\node [style=number-labels] (125) at (9.75, 0.7) {$4$};
		\node [style=number-labels] (126) at (11.05, 0.625) {$4$};
		\node [style=number-labels] (127) at (12.35, 0.525) {$4$};
		\node [style=dot] (128) at (3.4, 0.25) {};
		\node [style=dot] (129) at (4.675, 0.375) {};
		\node [style=dot] (130) at (6, 0.425) {};
		\node [style=dot] (131) at (9.475, 0.425) {};
		\node [style=dot] (132) at (10.8, 0.35) {};
		\node [style=dot] (133) at (12.075, 0.25) {};
		\node [style=dot] (134) at (-13.325, 3.025) {};
		\node [style=dot] (135) at (-11.9, 3.55) {};
		\node [style=dot] (136) at (-10.375, 3.925) {};
		\node [style=dot] (137) at (-7.125, 3.925) {};
		\node [style=dot] (138) at (-5.625, 3.575) {};
		\node [style=dot] (139) at (-4.3, 3.075) {};
		\node [style=dot] (140) at (-4.35, -1.475) {};
		\node [style=dot] (141) at (-5.65, -1.375) {};
		\node [style=dot] (142) at (-6.95, -1.3) {};
		\node [style=dot] (143) at (-10.55, -1.3) {};
		\node [style=dot] (144) at (-11.85, -1.375) {};
		\node [style=dot] (145) at (-13.175, -1.5) {};
	\end{pgfonlayer}
	\begin{pgfonlayer}{edgelayer}
		\draw [bend left=15, looseness=0.50] (2.center) to (3.center);
		\draw [bend right=15, looseness=0.75] (4.center) to (5.center);
		\draw [style=red, bend right=15, looseness=0.75] (0.center) to (10.center);
		\draw [style=red, bend left=15, looseness=0.50] (1.center) to (11.center);
		\draw [bend right=15, looseness=0.50] (12.center) to (13.center);
		\draw [bend right=15, looseness=0.25] (14.center) to (15.center);
		\draw [bend left=15, looseness=0.25] (16.center) to (17.center);
		\draw [bend left=15, looseness=0.50] (20.center) to (21.center);
		\draw [bend left=15, looseness=0.75] (18.center) to (19.center);
		\draw [style=blue, bend left=90, looseness=0.75] (0.center) to (1.center);
		\draw [bend left] (45.center) to (46.center);
		\draw [bend right=15, looseness=0.75] (47.center) to (39.center);
		\draw [bend right=15, looseness=0.50] (48.center) to (40.center);
		\draw [bend right=15, looseness=0.50] (49.center) to (41.center);
		\draw [bend left=15, looseness=0.50] (50.center) to (42.center);
		\draw [bend left=15, looseness=0.50] (51.center) to (43.center);
		\draw [bend left=15, looseness=0.75] (52.center) to (44.center);
		\draw [style=red, bend left=15] (35.center) to (37.center);
		\draw [style=red, bend left=15] (38.center) to (36.center);
		\draw [style=blue, bend left=75, looseness=0.75] (35.center) to (36.center);
		\draw [style=double-sided] (67.center) to (66.center);
		\draw [style=red, in=120, out=-120, looseness=0.75] (68.center) to (69.center);
		\draw [bend left=15, looseness=0.50] (74) to (75.center);
		\draw [style=red, bend left, looseness=0.75] (100.center) to (101.center);
		\draw [bend left=15, looseness=0.50] (103.center) to (102.center);
		\draw [bend right, looseness=0.25] (104.center) to (105.center);
		\draw [bend right=15, looseness=0.25] (106.center) to (107.center);
		\draw [bend left=15, looseness=0.25] (108.center) to (109.center);
		\draw [bend left=15, looseness=0.50] (110.center) to (111.center);
		\draw [bend left=15, looseness=0.50] (112.center) to (113.center);
	\end{pgfonlayer}
\end{tikzpicture}
	\caption[The term cancellations in \ref{term:st2} in Theorem~\ref{prop:solid-torus-ainf}.]
	{\textbf{The term cancellations in \ref{term:st2} in Theorem~\ref{prop:solid-torus-ainf}.}
	On the left, the first module tiling pattern is drawn on the bottom, and the second on
	the top.
	The tiling pattern for the centered algebra operation is on the right.}
	\label{fig:solid-torus-cancellations-2a}
\end{figure}

		For example, in the weight $0$ $\Ainf$ relation with inputs $x, \rho_2, \rho_1, \rho_4, \rho_3$, the following two non-zero terms cancel:
		\[m^0_2(m^0_2(x, \rho_2, \rho_1), \rho_4, \rho_3) = m^0_2(x, \mu_4^0(\rho_2, \rho_1, \rho_4, \rho_3)).\]

	\item \label{term:st3}
	  Consider terms of the form $m_{1+n-j}^{w-w'}\left(m_{1+j}^{w'}(x,
		a_1,\dots,a_j),a_{j+1},\dots,a_n\right)$, where $a_j a_{j+1}=0$.
		Let $\Gamma_1$ and $\Gamma_2$ be the module tiling patterns for the first and second operations respectively.

		For this case, suppose that the length of the blue boundaries do not match. Say the length of the blue boundary of $\Gamma_1$ is longer than that of $\Gamma_2$. The other case is entirely analogous; for completeness, we will run through the differences after this discussion. When we perform the gluing of $\Gamma_1$ with $\Gamma_2$, the red boundary of $\Gamma_2$ juts
		against the blue boundary of $\Gamma_1$, along a face $f'$ of $\Gamma_1$. Consider the penultimate edge $e$ of this
		face---that is, the edge preceding an edge $e'$ to
		the blue boundary. We will perform a maneuver we call pushing out the edge $e$ to the red boundary just before
		the blue-red intersection, breaking it into two edges.

		Consider the face $f$ of $\Gamma_1$ that lay on the other side of the edge $e$.
		There are two possibilities for the $f$: it is either an internal face and is bound by four edges, or it
		touches the boundary. If it is an internal face, pushing $e$ out
		to the boundary punctures this face. This produces a module tiling pattern with
		$\rho_i \rho_{i+1} \rho_{i+2} \rho_{i+3}$ as the last term of the chord sequence of the tiling pattern, while the other terms remain the same. The weight of the tiling pattern is one less than the sum of the two, as we puncture an internal face to construct it. It is thus a pattern for the operation $m^{w-1}_{2+n}(x, a_1, \dots, a_n, \mu_0^1)$ with the same output, and our term cancels against this. We observe that there is exactly one length four chord in $\mu_0^1$ for which there could be a tiling pattern because of the labeling restrictions. In particular, this is the length four chord $a$ such that $a_n \tensor a \neq 0$ so the idempotents match; but $a_n a = 0$.
		
		A schema for this cancellation is illustrated at the top of Figure~\ref{fig:solid-torus-cancellations-2}. As an example, for the weight $1$ $\Ainf$ relation with inputs $x, \rho_{23}, \rho_2, \rho_{12}, \rho_1, \rho_{41}, \rho_4, \rho_3$, the following two non-zero terms cancel
		\begin{align*}
			m^0_3(m^1_6(x, \rho_{23}, \rho_2, \rho_{12}, \rho_1, \rho_{41}), \rho_4, \rho_3) &= m^0_9(x, \rho_{23}, \rho_2, \rho_{12}, \rho_1, \rho_{41}, \rho_4, \rho_3, \rho_{2341}) \\
			&= m^0_9(x, \rho_{23}, \rho_2, \rho_{12}, \rho_1, \rho_{41}, \rho_4, \rho_3, \mu^1_0).
		\end{align*}

		If the face $f$ touches the boundary, it must touch the red boundary. This is because there are no internal edges on which either side lies the blue boundary.
		Pushing $e$ out disconnects the module tiling pattern, producing two patterns
		$\Gamma_m$ and $\Gamma_a$. We call this the disconnected tiling pattern a composite tiling pattern, and represent it as $\Gamma_m~\#~ \Gamma_a$. The pattern $\Gamma_m$ is a valid module tiling pattern, while the pattern $\Gamma_a$ has only red vertices. We claim that there is a right-extended algebra operation $\mu^{w_1}_k\left(a_k,
		\dots,a_n\right)$ such that $\Gamma_m$ is a module tiling pattern for the operation $m_{1+n-k}^{w-w_1}\left(x,a_{1},\dots,a_{n-k},\mu^{w_1}_k\left(a_{n-k+1},
		\dots,a_n\right)\right)$. This is then a non-zero term which cancels against our given term.
		
		Let the last term of the chord sequence of $\Gamma_m$ be $\rho_i \rho_{i+1} \dots \rho_{i+l}$. Along the edge pushed out to the boundary in $\Gamma_a$, add $l+1$ $2$-valent vertices. Label these $i, i+1, \dots, i+l$ on the right. The labelings remain consistent for $\Gamma_a$ as these were the labelings for an erstwhile face of $\Gamma_1$. Mark as root the leaf where the pushed out edge touches the boundary. This makes $\Gamma_a$ an algebra tiling pattern for the right-extended algebra option $\mu^{w_1}_k\left(a_{n-k+1},
		\dots,a_n\right)$, as desired.
		
		A schema for this cancellation is illustrated at the bottom of Figure~\ref{fig:solid-torus-cancellations-2}. As an example, for the weight 0 $\Ainf$ relation with inputs $x, \rho_4, \rho_{34}, \rho_3, \rho_2, \rho_1$, the following two non-zero terms cancel
		\begin{align*}m^3_0(m^4_0(x, \rho_4, \rho_{34}, \rho_3), \rho_2, \rho_1)
		&= m^3_0(x, \rho_4, \rho_3) \\
		&= m^3_0(x, \rho_4, \mu^0_4(\rho_{34}, \rho_3, \rho_2, \rho_1)).\end{align*}
		
		As another example, in the weight $0$ $\Ainf$ relation with inputs $x, \rho_4, \rho_{341}, \rho_4, \rho_3, \rho_{23}, \rho_2, \rho_1$, the following two non-zero terms cancel
		\[m^0_3(m^0_6(x, \rho_4, \rho_{341}, \rho_4, \rho_3, \rho_{23}), \rho_2, \rho_1) = m^0_3(x, \rho_4, \mu^0_6(\rho_{341}, \rho_4, \rho_3, \rho_{23}, \rho_2, \rho_1)).\]

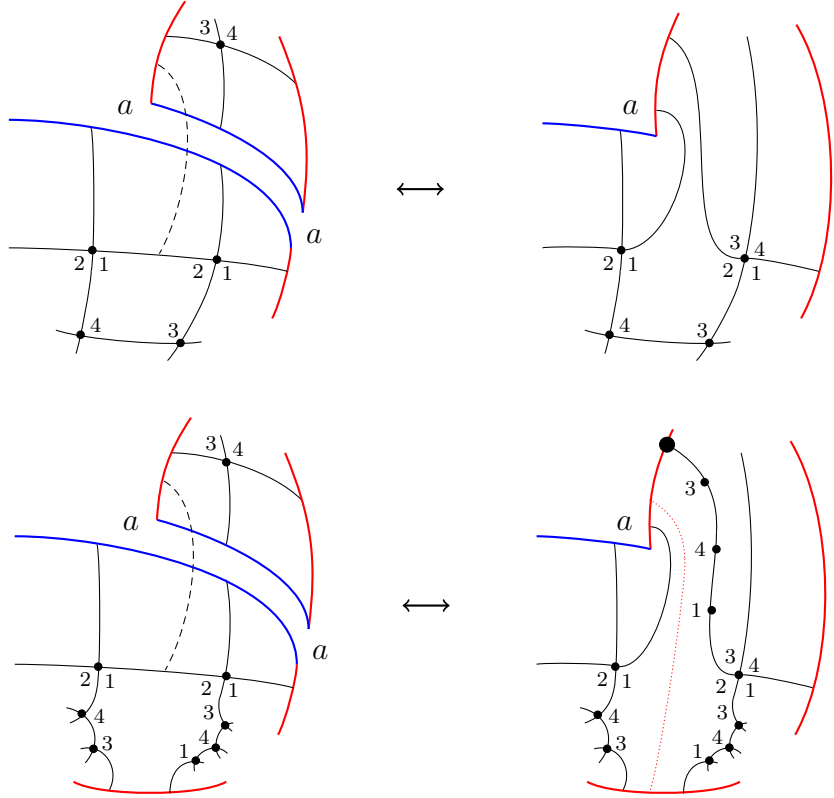
\begin{figure}[h!tbp]
	\centering
	\begin{tikzpicture}[scale=1.25]
	\begin{pgfonlayer}{nodelayer}
		\node [style=none] (1) at (-2.75, -1.25) {};
		\node [style=none] (3) at (-2.825, -1.75) {};
		\node [style=none] (7) at (-2.25, -1) {$a$};
		\node [style=none] (11) at (-3.15, -2.75) {};
		\node [style=none] (16) at (-7, 1.3) {};
		\node [style=none] (17) at (-7, -1.75) {};
		\node [style=none] (18) at (-4.25, 0.525) {};
		\node [style=none] (19) at (-4.475, -2) {};
		\node [style=number-labels] (27) at (-7.25, -1.575) {$2$};
		\node [style=number-labels] (29) at (-4.725, -1.75) {$2$};
		\node [style=number-labels] (32) at (-6.725, -1.6) {$1$};
		\node [style=number-labels] (34) at (-4.1, -1.8) {$1$};
		\node [style=none] (35) at (-5.725, 1.825) {};
		\node [style=none] (36) at (-2.5, -0.5) {};
		\node [style=none] (38) at (-3, 3.25) {};
		\node [style=none] (44) at (-4.25, 1.3) {};
		\node [style=none] (45) at (-5.425, 3.25) {};
		\node [style=none] (46) at (-2.65, 2.25) {};
		\node [style=none] (52) at (-4.375, 3.625) {};
		\node [style=number-labels] (54) at (-4, 3.325) {$4$};
		\node [style=number-labels] (55) at (-4.575, 3.45) {$3$};
		\node [style=none] (66) at (-0.5, 0) {};
		\node [style=none] (67) at (0.5, 0) {};
		\node [style=none] (128) at (-5, 4) {};
		\node [style=none] (129) at (-6.275, 1.75) {$a$};
		\node [style=none] (130) at (-7.325, -3.5) {};
		\node [style=none] (131) at (-5.375, -3.65) {};
		\node [style=none] (132) at (-7.75, -3) {};
		\node [style=none] (133) at (-4.65, -3.25) {};
		\node [style=number-labels] (134) at (-5.25, -3) {$3$};
		\node [style=number-labels] (135) at (-6.9, -2.9) {$4$};
		\node [style=none] (136) at (-5.575, 2.65) {};
		\node [style=none] (138) at (-5.55, -1.375) {};
		\node [style=none] (142) at (8.475, -1.75) {};
		\node [style=none] (150) at (8.1, -2.75) {};
		\node [style=none] (155) at (4.25, 1.25) {};
		\node [style=none] (156) at (4.25, -1.75) {};
		\node [style=none] (157) at (6.95, 3.25) {};
		\node [style=none] (158) at (6.775, -2) {};
		\node [style=number-labels] (164) at (4, -1.575) {$2$};
		\node [style=number-labels] (165) at (6.55, -1.725) {$2$};
		\node [style=number-labels] (168) at (4.525, -1.6) {$1$};
		\node [style=number-labels] (169) at (7.15, -1.8) {$1$};
		\node [style=none] (170) at (5.025, 1.125) {};
		\node [style=none] (172) at (8, 3.5) {};
		\node [style=number-labels] (177) at (7.225, -1.275) {$4$};
		\node [style=number-labels] (178) at (6.7, -1.125) {$3$};
		\node [style=none] (179) at (5.5, 3.75) {};
		\node [style=none] (180) at (4.475, 1.75) {$a$};
		\node [style=none] (181) at (3.925, -3.5) {};
		\node [style=none] (182) at (5.875, -3.65) {};
		\node [style=none] (183) at (3.5, -3) {};
		\node [style=none] (184) at (6.6, -3.25) {};
		\node [style=number-labels] (185) at (6, -3) {$3$};
		\node [style=number-labels] (186) at (4.35, -2.9) {$4$};
		\node [style=none] (187) at (5.025, 1.675) {};
		\node [style=none] (188) at (4.275, -1.3) {};
		\node [style=dot] (189) at (6.9, -1.475) {};
		\node [style=none] (190) at (5.3, 3.225) {};
		\node [style=none] (191) at (-2.75, -1.25) {};
		\node [style=none] (192) at (-8.75, 1.475) {};
		\node [style=none] (193) at (-2.825, -1.75) {};
		\node [style=none] (194) at (-8.75, -1.25) {};
		\node [style=none] (195) at (5.025, 1.125) {};
		\node [style=dot] (196) at (4.275, -1.3) {};
		\node [style=none] (197) at (2.6, 1.4) {};
		\node [style=none] (198) at (2.6, -1.25) {};
		\node [style=dot] (199) at (6.15, -3.275) {};
		\node [style=dot] (200) at (4.025, -3.1) {};
		\node [style=dot] (201) at (-4.325, -1.5) {};
		\node [style=dot] (202) at (-5.1, -3.275) {};
		\node [style=dot] (203) at (-7.225, -3.1) {};
		\node [style=dot] (204) at (-6.975, -1.3) {};
		\node [style=dot] (205) at (-4.25, 3.075) {};
	\end{pgfonlayer}
	\begin{pgfonlayer}{edgelayer}
		\draw [bend left=15, looseness=0.25] (16.center) to (17.center);
		\draw [bend left=15, looseness=0.75] (18.center) to (19.center);
		\draw [in=135, out=0, looseness=0.75] (45.center) to (46.center);
		\draw [bend left=15, looseness=0.75] (52.center) to (44.center);
		\draw [style=red, bend left=15] (38.center) to (36.center);
		\draw [style=blue, in=90, out=-15, looseness=0.75] (35.center) to (36.center);
		\draw [style=double-sided] (67.center) to (66.center);
		\draw [style=red, bend left=15] (35.center) to (128.center);
		\draw [bend left=15, looseness=0.25] (17.center) to (130.center);
		\draw [bend left=15, looseness=0.50] (19.center) to (131.center);
		\draw [bend right=15, looseness=0.50] (132.center) to (133.center);
		\draw [style=red, bend left=15, looseness=0.50] (1.center) to (11.center);
		\draw [style=densedashedge, in=60, out=-30, looseness=0.75] (136.center) to (138.center);
		\draw [bend left=15, looseness=0.25] (155.center) to (156.center);
		\draw [bend left=15, looseness=0.75] (157.center) to (158.center);
		\draw [style=red, bend left=15] (170.center) to (179.center);
		\draw [bend left=15, looseness=0.25] (156.center) to (181.center);
		\draw [bend left=15, looseness=0.50] (158.center) to (182.center);
		\draw [bend right=15, looseness=0.50] (183.center) to (184.center);
		\draw [style=red, bend left, looseness=0.75] (172.center) to (150.center);
		\draw [in=165, out=0, looseness=0.50] (189) to (142.center);
		\draw [in=0, out=0] (188.center) to (187.center);
		\draw [in=180, out=-30, looseness=0.75] (190.center) to (189);
		\draw [style=blue, in=90, out=0, looseness=0.75] (192.center) to (191.center);
		\draw [in=165, out=0, looseness=0.50] (194.center) to (193.center);
		\draw [bend right=15, looseness=0.25] (196) to (198.center);
		\draw [style=blue, in=165, out=0, looseness=0.50] (197.center) to (195.center);
	\end{pgfonlayer}
\end{tikzpicture}

	\bigskip \bigskip
	\begin{tikzpicture}[scale=1.25]
	\begin{pgfonlayer}{nodelayer}
		\node [style=none] (1) at (-2.75, -1.25) {};
		\node [style=none] (3) at (-2.825, -1.75) {};
		\node [style=none] (7) at (-2.25, -1) {$a$};
		\node [style=none] (11) at (-3.15, -2.75) {};
		\node [style=none] (16) at (-7, 1.3) {};
		\node [style=none] (17) at (-7, -1.75) {};
		\node [style=none] (18) at (-4.25, 0.55) {};
		\node [style=none] (19) at (-4.4, -2) {};
		\node [style=number-labels] (27) at (-7.25, -1.575) {$2$};
		\node [style=number-labels] (29) at (-4.675, -1.75) {$2$};
		\node [style=number-labels] (32) at (-6.725, -1.6) {$1$};
		\node [style=number-labels] (34) at (-4.1, -1.8) {$1$};
		\node [style=none] (35) at (-5.725, 1.825) {};
		\node [style=none] (36) at (-2.5, -0.5) {};
		\node [style=none] (38) at (-3, 3.25) {};
		\node [style=none] (44) at (-4.25, 1.3) {};
		\node [style=none] (45) at (-5.425, 3.25) {};
		\node [style=none] (46) at (-2.65, 2.25) {};
		\node [style=none] (52) at (-4.375, 3.625) {};
		\node [style=number-labels] (54) at (-4, 3.325) {$4$};
		\node [style=number-labels] (55) at (-4.575, 3.45) {$3$};
		\node [style=none] (66) at (-0.5, 0) {};
		\node [style=none] (67) at (0.5, 0) {};
		\node [style=none] (128) at (-5, 4) {};
		\node [style=none] (129) at (-6.275, 1.75) {$a$};
		\node [style=none] (130) at (-7.575, -2.5) {};
		\node [style=none] (131) at (-4.125, -2.675) {};
		\node [style=none] (136) at (-5.575, 2.65) {};
		\node [style=none] (138) at (-5.55, -1.375) {};
		\node [style=none] (142) at (8.225, -1.75) {};
		\node [style=none] (150) at (7.85, -2.75) {};
		\node [style=none] (155) at (4, 1.325) {};
		\node [style=none] (156) at (4, -1.75) {};
		\node [style=none] (157) at (6.7, 3.25) {};
		\node [style=none] (158) at (6.525, -2) {};
		\node [style=number-labels] (164) at (3.75, -1.575) {$2$};
		\node [style=number-labels] (165) at (6.275, -1.725) {$2$};
		\node [style=number-labels] (168) at (4.275, -1.6) {$1$};
		\node [style=number-labels] (169) at (6.9, -1.8) {$1$};
		\node [style=none] (170) at (4.775, 1.2) {};
		\node [style=none] (172) at (7.75, 3.5) {};
		\node [style=number-labels] (177) at (6.975, -1.275) {$4$};
		\node [style=number-labels] (178) at (6.475, -1.075) {$3$};
		\node [style=none] (179) at (5.25, 3.75) {};
		\node [style=none] (180) at (4.225, 1.75) {$a$};
		\node [style=none] (187) at (4.775, 1.675) {};
		\node [style=dot] (188) at (4.025, -1.3) {};
		\node [style=dot] (189) at (6.65, -1.475) {};
		\node [style=none] (190) at (6.175, 1.55) {};
		\node [style=none] (191) at (-7.65, -2.175) {};
		\node [style=none] (192) at (-7.25, -3.35) {};
		\node [style=none] (193) at (-7.35, -3.05) {};
		\node [style=none] (194) at (-6.75, -3.95) {};
		\node [style=none] (195) at (-7.5, -3.75) {};
		\node [style=none] (196) at (-4.25, -3.75) {};
		\node [style=none] (197) at (-4.1, -2.5) {};
		\node [style=none] (198) at (-4.375, -3.25) {};
		\node [style=none] (199) at (-4.225, -3.15) {};
		\node [style=none] (200) at (-4.925, -3.525) {};
		\node [style=none] (201) at (-4.725, -3.35) {};
		\node [style=none] (202) at (-5.45, -4) {};
		\node [style=none] (203) at (6.525, -2) {};
		\node [style=none] (205) at (6.8, -2.675) {};
		\node [style=none] (206) at (3.275, -2.175) {};
		\node [style=none] (207) at (3.675, -3.35) {};
		\node [style=none] (208) at (3.575, -3.05) {};
		\node [style=none] (209) at (4.175, -3.95) {};
		\node [style=none] (210) at (3.425, -3.75) {};
		\node [style=none] (211) at (6.675, -3.75) {};
		\node [style=none] (212) at (6.825, -2.5) {};
		\node [style=none] (213) at (6.55, -3.25) {};
		\node [style=none] (214) at (6.7, -3.15) {};
		\node [style=none] (215) at (6, -3.525) {};
		\node [style=none] (216) at (6.2, -3.35) {};
		\node [style=none] (217) at (5.475, -4) {};
		\node [style=none] (218) at (4, -1.75) {};
		\node [style=none] (219) at (3.375, -2.5) {};
		\node [style=number-labels] (220) at (4, -2.35) {$4$};
		\node [style=number-labels] (221) at (4.15, -2.925) {$3$};
		\node [style=number-labels] (222) at (6.3, -2.275) {$3$};
		\node [style=number-labels] (223) at (6.2, -2.8) {$4$};
		\node [style=number-labels] (224) at (5.775, -3.1) {$1$};
		\node [style=dot] (225) at (5.925, 2.625) {};
		\node [style=dot] (226) at (6.075, -0.1) {};
		\node [style=dot] (227) at (6.175, 1.2) {};
		\node [style=none] (228) at (4.75, 2.25) {};
		\node [style=none] (229) at (5.5, 1) {};
		\node [style=number-labels] (230) at (-6.95, -2.35) {$4$};
		\node [style=number-labels] (231) at (-6.775, -2.925) {$3$};
		\node [style=number-labels] (232) at (-4.625, -2.275) {$3$};
		\node [style=number-labels] (233) at (-4.725, -2.8) {$4$};
		\node [style=number-labels] (234) at (-5.15, -3.1) {$1$};
		\node [style=bigdot] (235) at (5.125, 3.425) {};
		\node [style=none] (236) at (4.75, -3.975) {};
		\node [style=number-labels] (237) at (5.625, 2.5) {$3$};
		\node [style=number-labels] (238) at (5.825, 1.2) {$4$};
		\node [style=number-labels] (239) at (5.75, -0.1) {$1$};
		\node [style=none] (240) at (-8.75, 1.475) {};
		\node [style=none] (241) at (-8.75, -1.25) {};
		\node [style=none] (242) at (8.35, -1.25) {};
		\node [style=none] (243) at (8.275, -1.75) {};
		\node [style=none] (244) at (5.55, -1.375) {};
		\node [style=none] (245) at (2.35, 1.475) {};
		\node [style=none] (246) at (2.35, -1.25) {};
		\node [style=dot] (247) at (-4.25, 3.05) {};
		\node [style=dot] (248) at (-6.975, -1.3) {};
		\node [style=dot] (249) at (-4.25, -1.5) {};
		\node [style=dot] (250) at (-7.325, -2.3) {};
		\node [style=dot] (251) at (-4.275, -2.55) {};
		\node [style=dot] (252) at (-7.075, -3.05) {};
		\node [style=dot] (253) at (-4.475, -3.025) {};
		\node [style=dot] (254) at (-4.9, -3.3) {};
		\node [style=dot] (256) at (3.65, -2.325) {};
		\node [style=dot] (257) at (6.65, -2.55) {};
		\node [style=dot] (258) at (3.85, -3.05) {};
		\node [style=dot] (259) at (6.45, -3.025) {};
		\node [style=dot] (260) at (6.025, -3.3) {};
	\end{pgfonlayer}
	\begin{pgfonlayer}{edgelayer}
		\draw [bend left=15, looseness=0.25] (16.center) to (17.center);
		\draw [bend left=15, looseness=0.75] (18.center) to (19.center);
		\draw [in=135, out=0, looseness=0.75] (45.center) to (46.center);
		\draw [bend left=15, looseness=0.75] (52.center) to (44.center);
		\draw [style=red, bend left=15] (38.center) to (36.center);
		\draw [style=blue, in=90, out=-15, looseness=0.75] (35.center) to (36.center);
		\draw [style=double-sided] (67.center) to (66.center);
		\draw [style=red, bend left=15] (35.center) to (128.center);
		\draw [bend left] (17.center) to (130.center);
		\draw [bend right] (19.center) to (131.center);
		\draw [style=red, bend left=15, looseness=0.50] (1.center) to (11.center);
		\draw [style=densedashedge, in=60, out=-30, looseness=0.75] (136.center) to (138.center);
		\draw [bend left=15, looseness=0.25] (155.center) to (156.center);
		\draw [bend left=15, looseness=0.75] (157.center) to (158.center);
		\draw [style=red, bend left=15] (170.center) to (179.center);
		\draw [in=165, out=0, looseness=0.50] (189) to (142.center);
		\draw [in=0, out=0, looseness=0.75] (188) to (187.center);
		\draw [in=180, out=-90, looseness=0.75] (190.center) to (189);
		\draw [bend left=60, looseness=1.25] (191.center) to (192.center);
		\draw [in=60, out=15, looseness=1.25] (193.center) to (194.center);
		\draw [style=red, bend right, looseness=0.50] (195.center) to (196.center);
		\draw [in=135, out=-180] (197.center) to (198.center);
		\draw [bend right=75, looseness=1.25] (199.center) to (200.center);
		\draw [in=90, out=165, looseness=1.25] (201.center) to (202.center);
		\draw [bend right] (203.center) to (205.center);
		\draw [bend left=60, looseness=1.25] (206.center) to (207.center);
		\draw [in=60, out=15, looseness=1.25] (208.center) to (209.center);
		\draw [style=red, bend right, looseness=0.50] (210.center) to (211.center);
		\draw [in=135, out=-180] (212.center) to (213.center);
		\draw [bend right=75, looseness=1.25] (214.center) to (215.center);
		\draw [in=90, out=165, looseness=1.25] (216.center) to (217.center);
		\draw [bend left] (218.center) to (219.center);
		\draw [style=denselydottedred, bend left] (228.center) to (229.center);
		\draw [bend left] (235) to (190.center);
		\draw [style=denselydottedred, bend left=15, looseness=0.25] (229.center) to (236.center);
		\draw [in=165, out=0, looseness=0.50] (241.center) to (3.center);
		\draw [style=blue, in=90, out=0, looseness=0.75] (240.center) to (1.center);
		\draw [style=red, bend left, looseness=0.75] (172.center) to (150.center);
		\draw [bend right=15, looseness=0.25] (188) to (246.center);
		\draw [style=blue, in=165, out=0, looseness=0.50] (245.center) to (170.center);
	\end{pgfonlayer}
\end{tikzpicture}
	\caption[The term cancellations in \ref{term:st3} in Theorem~\ref{prop:solid-torus-ainf}.]
	{\textbf{The term cancellations in \ref{term:st3} in Theorem~\ref{prop:solid-torus-ainf}.}
	On the left, the first module operation is drawn on the bottom, and the second on
	the top. The labels are arbitrarily chosen.

	Top: pushing out the edge
	$e$ to the red boundary punctures an internal face, producing a module operation
	on the right.
	
	Bottom: pushing out the edge $e$ to the red boundary disconnects
	the tiling pattern. Along the red dotted arc, it splits into a composite pattern $\Gamma_m~\#~\Gamma_a$. On the left of
	the dotted red arc is $\Gamma_a$. On the right of the arc
	is $\Gamma_a$ with $2$-valent vertices added.}
	\label{fig:solid-torus-cancellations-2}
\end{figure}

We briefly describe the analogous case where the length of the blue boundary of $\Gamma_1$ is shorter than that of $\Gamma_2$. When we perform the gluing of $\Gamma_1$ with $\Gamma_2$, the red
		boundary of the $\Gamma_1$ juts
		against the blue boundary of $\Gamma_2$, along a face $f'$ of $\Gamma_2$. Consider the second edge $e$ of this
		face, and push $e$ out to the red boundary just after
		the blue-red intersection.

Consider the face $f$ of $\Gamma_2$ that lay on the other side of the edge $e$.
		If $f$ is an internal face, we produce a module tiling pattern with
		$\rho_i \rho_{i+1} \rho_{i+2} \rho_{i+3}$ as the first term of the chord sequence. Our term cancels against $m^{w-1}_{2+n}(x, \mu_0^1, a_1, \dots, a_n)$.
		If the face $f$ touches the boundary, pushing $e$ out produces a composite pattern
		$\Gamma_a~\#~\Gamma_m$. We can label $2$-valent vertices along the left of the pushed out edge in $\Gamma_a$ such that we get an algebra tiling pattern for a left-extended algebra operation $\mu^{w_1}_k\left(a_1,\dots,
		a_k\right)$. Our term
		cancels against $m_{1+n-k}^{w-w_1}\left(x,\mu^{w_1}_k\left(a_1,\dots,
		a_k\right),a_{k+1},\dots,a_n\right)$.
	\end{enumerate}
	
Next, we consider the cases where both the first internal vertex in our tree for the $\Ainf$ relation corresponds to an algebra operation. In order, we consider the cases where: a non-centered algebra operation occurs non-\emph{extremally}; a non-centered algebra operation occurs extremally; $\mu_0^1$ occurs extremally; $\mu_0^1$ occurs non-extremally; there is a centered algebra operation; and lastly, when the algebra operation is $\mu_2^0$. These exhaust the cases.

	\begin{enumerate}[wide,label=(T-\arabic*),resume]
		\item \label{term:st4} \sloppy Consider terms of the form 
		$m_{1+n-k}^{w-w'}(x,a_1,\dots,a_j,\mu^{w'}_k(a_{j+1},\dots,
		a_{j+k}),\allowbreak a_{j+k+1},\dots,a_n)$ for a left-extended algebra operation $\mu^{w'}_k$. In this case, the algebra operation is not permitted to be the first input to the module operation.
		We cancel these terms against
		$m_{1+n}^w(x,a_1,\dots,\mu^0_2(a_j,a_{j+1}),a_{j+k+1}\dots,a_n)$.

\begin{figure}[h!tbp]
	\centering
	\begin{tikzpicture}
	\begin{pgfonlayer}{nodelayer}
		\node [style=none] (0) at (-14.925, -0.75) {};
		\node [style=none] (1) at (-13.75, 0.7) {};
		\node [style=none] (2) at (-14, 0.5) {};
		\node [style=none] (3) at (-12.1, 1.275) {};
		\node [style=none] (4) at (-12.5, 1.075) {};
		\node [style=none] (5) at (-9, 4) {};
		\node [style=none] (6) at (-10.45, 1.85) {};
		\node [style=none] (7) at (-7.3, 1.475) {};
		\node [style=none] (8) at (-7.75, 1.6) {};
		\node [style=none] (9) at (-6, 0) {};
		\node [style=none] (12) at (-13.575, -0.35) {};
		\node [style=none] (13) at (-9.65, 1) {};
		\node [style=none] (14) at (-6.25, 0.35) {};
		\node [style=none] (15) at (-5.375, -1.1) {};
		\node [style=none] (16) at (-8.275, 0.925) {};
		\node [style=none] (17) at (-6.35, -0.525) {};
		\node [style=none] (18) at (-9.25, 3.775) {};
		\node [style=none] (19) at (-6.25, 4.5) {};
		\node [style=none] (20) at (-6.6, 4.5) {};
		\node [style=none] (21) at (-3.75, 3.25) {};
		\node [style=none] (22) at (-4.075, 3.5) {};
		\node [style=none] (23) at (-2.05, 0.95) {};
		\node [style=none] (24) at (-15.5, 1) {};
		\node [style=none] (25) at (-13, -2.5) {};
		\node [style=none] (26) at (-5.25, -1) {};
		\node [style=none] (29) at (-1.5, 1.75) {};
		\node [style=number-labels] (30) at (-11.25, -0.25) {$a_j$};
		\node [style=number-labels] (31) at (-7.8, -0.325) {$a_{j+1}$};
		\node [style=none] (32) at (-9, 1.15) {};
		\node [style=none] (33) at (-9.75, -2.725) {};
		\node [style=none] (34) at (2.075, -0.75) {};
		\node [style=none] (35) at (3.25, 0.7) {};
		\node [style=none] (36) at (3, 0.5) {};
		\node [style=none] (37) at (4.9, 1.275) {};
		\node [style=none] (38) at (4.5, 1.075) {};
		\node [style=none] (39) at (8, 4) {};
		\node [style=none] (40) at (6.55, 1.85) {};
		\node [style=none] (41) at (9.7, 1.475) {};
		\node [style=none] (42) at (9.25, 1.6) {};
		\node [style=none] (43) at (11, 0) {};
		\node [style=none] (44) at (3.425, -0.35) {};
		\node [style=none] (45) at (6.775, 1) {};
		\node [style=none] (46) at (10.75, 0.35) {};
		\node [style=none] (47) at (11.625, -1.1) {};
		\node [style=none] (48) at (9.25, 0.5) {};
		\node [style=none] (49) at (10.65, -0.525) {};
		\node [style=none] (50) at (7.75, 3.775) {};
		\node [style=none] (51) at (10.75, 4.5) {};
		\node [style=none] (52) at (10.4, 4.5) {};
		\node [style=none] (53) at (13.25, 3.25) {};
		\node [style=none] (54) at (12.925, 3.5) {};
		\node [style=none] (55) at (14.95, 0.95) {};
		\node [style=none] (56) at (1.5, 1) {};
		\node [style=none] (57) at (4, -2.5) {};
		\node [style=none] (58) at (11.75, -1) {};
		\node [style=none] (59) at (15.5, 1.75) {};
		\node [style=number-labels] (60) at (5.425, -0.275) {$a_j$};
		\node [style=number-labels] (61) at (9.65, -0.6) {$a_{j+1}$};
		\node [style=none] (62) at (7.9, 1.15) {};
		\node [style=none] (63) at (7.25, -2.725) {};
		\node [style=none] (64) at (-0.5, 0) {};
		\node [style=none] (65) at (0.5, 0) {};
		\node [style=none] (66) at (7.25, 0.75) {};
		\node [style=none] (67) at (7, -2.75) {};
		\node [style=dot] (68) at (8.8, 0.575) {};
		\node [style=bigdot] (69) at (8.525, -2.5) {};
		\node [style=none] (70) at (9.25, 3) {};
		\node [style=none] (71) at (14.4, 0.575) {};
		\node [style=dot] (72) at (-13.725, 0.35) {};
		\node [style=dot] (73) at (-12.175, 1) {};
		\node [style=dot] (74) at (-10.25, 1.575) {};
		\node [style=dot] (75) at (-9.025, 3.65) {};
		\node [style=dot] (76) at (-6.45, 4.2) {};
		\node [style=dot] (77) at (-4.05, 3.15) {};
		\node [style=dot] (78) at (-7.65, 1.3) {};
		\node [style=dot] (79) at (-6.25, 0.05) {};
		\node [style=dot] (80) at (3.25, 0.35) {};
		\node [style=dot] (81) at (4.8, 1) {};
		\node [style=dot] (82) at (6.675, 1.55) {};
		\node [style=dot] (83) at (7.975, 3.65) {};
		\node [style=dot] (84) at (10.525, 4.2) {};
		\node [style=dot] (85) at (12.95, 3.15) {};
		\node [style=dot] (86) at (9.375, 1.2) {};
		\node [style=dot] (87) at (10.75, 0.05) {};
		\node [style=dot] (88) at (8.65, -0.175) {};
		\node [style=dot] (89) at (8.55, -0.95) {};
		\node [style=dot] (90) at (8.525, -1.725) {};
		\node [style=number-labels] (91) at (8.275, 3.15) {$2$};
		\node [style=number-labels] (92) at (10.5, 3.575) {$3$};
		\node [style=number-labels] (93) at (12.675, 2.675) {$4$};
		\node [style=number-labels] (94) at (7.25, 1.5) {$1$};
		\node [style=number-labels] (95) at (8.425, 0.625) {$1$};
		\node [style=number-labels] (96) at (8.275, -0.15) {$2$};
		\node [style=number-labels] (97) at (8.175, -0.925) {$3$};
		\node [style=number-labels] (98) at (8.15, -1.7) {$4$};
		\node [style=number-labels] (99) at (-8.725, 3.15) {$2$};
		\node [style=number-labels] (100) at (-6.5, 3.575) {$3$};
		\node [style=number-labels] (101) at (-4.325, 2.675) {$4$};
		\node [style=number-labels] (102) at (-9.575, 1.65) {$1$};
	\end{pgfonlayer}
	\begin{pgfonlayer}{edgelayer}
		\draw [in=285, out=15] (0.center) to (1.center);
		\draw [bend right=60] (2.center) to (3.center);
		\draw [in=-90, out=-15] (4.center) to (5.center);
		\draw [in=-150, out=-60, looseness=0.75] (6.center) to (7.center);
		\draw [bend right, looseness=1.25] (8.center) to (9.center);
		\draw [style=braceedge] (13.center) to (12.center);
		\draw [in=135, out=-105, looseness=0.75] (14.center) to (15.center);
		\draw [style=braceedge] (17.center) to (16.center);
		\draw [bend right=45] (18.center) to (19.center);
		\draw [bend right=45] (20.center) to (21.center);
		\draw [in=150, out=-90] (22.center) to (23.center);
		\draw [style=red, in=165, out=-75] (24.center) to (25.center);
		\draw [style=red, bend right, looseness=0.75] (25.center) to (26.center);
		\draw [style=red, in=-105, out=45, looseness=0.75] (26.center) to (29.center);
		\draw [style=dashededged, bend right=15, looseness=0.50] (32.center) to (33.center);
		\draw [in=285, out=15] (34.center) to (35.center);
		\draw [bend right=60] (36.center) to (37.center);
		\draw [in=-90, out=-15] (38.center) to (39.center);
		\draw [bend right, looseness=1.25] (42.center) to (43.center);
		\draw [style=braceedge] (45.center) to (44.center);
		\draw [in=135, out=-105, looseness=0.75] (46.center) to (47.center);
		\draw [style=braceedge] (49.center) to (48.center);
		\draw [bend right=45] (50.center) to (51.center);
		\draw [bend right=45] (52.center) to (53.center);
		\draw [in=150, out=-90] (54.center) to (55.center);
		\draw [style=red, in=165, out=-75] (56.center) to (57.center);
		\draw [style=red, bend right, looseness=0.75] (57.center) to (58.center);
		\draw [style=red, in=-105, out=45, looseness=0.75] (58.center) to (59.center);
		\draw [style=denselydottedred, in=90, out=-105, looseness=0.50] (62.center) to (63.center);
		\draw [style=double-sided] (65.center) to (64.center);
		\draw [in=90, out=-75] (40.center) to (66.center);
		\draw [in=105, out=-90, looseness=0.50] (66.center) to (67.center);
		\draw [in=75, out=-135] (41.center) to (68);
		\draw [in=90, out=-105, looseness=0.75] (68) to (69);
		\draw [style=denselydottedred, in=-165, out=75, looseness=0.75] (62.center) to (70.center);
		\draw [style=denselydottedred, in=150, out=15, looseness=0.75] (70.center) to (71.center);
	\end{pgfonlayer}
\end{tikzpicture}
	\caption[The term cancellations in \ref{term:st4} in Theorem~\ref{prop:solid-torus-ainf}.]
	{\textbf{The term cancellations in \ref{term:st4} in Theorem~\ref{prop:solid-torus-ainf}.}
	The labels are arbitrarily chosen. 
	Left: a module tiling pattern with a chord factorized as the product of
	$a_j$ and $a_{j+1}$. Right:
	pushing out the edge corresponding to the factorization to the red boundary
	produces a composite pattern $\Gamma_m~\#~\Gamma_a$. Bounded by the dotted red arc
	and the red boundary is $\Gamma_a$, with $2$-valent vertices added in,
	corresponding to a left-extended algebra operation. $\Gamma_m$ is a module tiling pattern.}
	\label{fig:solid-torus-cancellations-ext}
\end{figure}
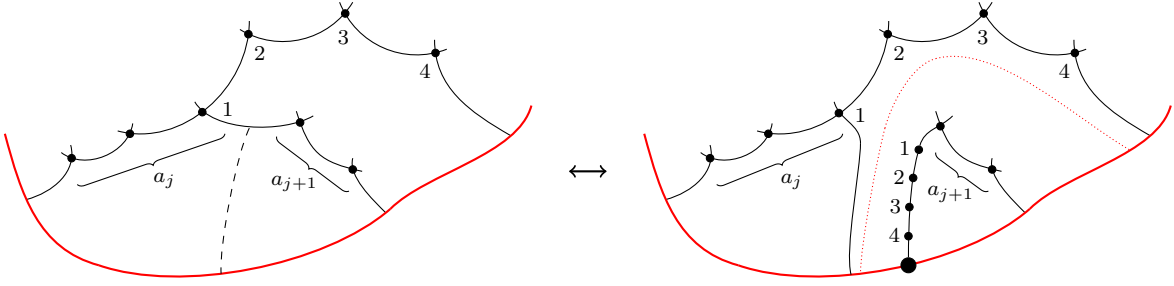

		Consider a module tiling pattern $\Gamma_m$ for the module operation and an algebra tiling pattern $\Gamma_a$ for the algebra operation. Along the red boundary in $\Gamma_m$, there is a face $f$ the labels of which match with the output of the algebra operation. In other words, the labels of this face match with the labels along the $2$-valent vertices of $\Gamma_a$. Decompose $\Gamma_a$ into arcs intersecting transversely. We glue in $\Gamma_a$ into $\Gamma_m$, identifying the arc containing the $2$-valent vertices with the edges of the face $f$. The orientation respects the fact that the $2$-valent vertices are labelled along the left. This is illustrated in Figure~\ref{fig:solid-torus-cancellations-ext}.

		This procedure of gluing in produces a bona fide module tiling pattern. In particular, the fact that the algebra operation is not the first input to the module operation, along with the fact that we glue in to the left, ensures that we do not modify the module tiling pattern near the blue-red boundary. With this new module tiling pattern, we thus demonstrate that our term cancels against the non-zero term $m_{1+n}^w(x,a_1,\dots,\mu^0_2(a_j,a_{j+1}),a_{j+k+1}\dots,a_n)$.
		
		As an example, in the weight $0$ $\Ainf$ relation with inputs $x, \rho_4, \rho_1, \rho_4, \rho_3, \rho_{234}, \rho_3$, the following non-zero terms cancel
		\[ m^0_4(x, \rho_4, \mu_4^0(\rho_1, \rho_4, \rho_3, \rho_{234}), \rho_3) = m^0_4(x, \mu^0_2(\rho_4, \rho_1), \rho_3, \rho_{234}, \rho_3). \]
		As another example where the left-extended algebra operation is the last input, in the weight $0$ $\Ainf$ relation with inputs $x, \rho_2, \rho_3, \rho_2, \rho_1, \rho_{41}$, the following non-zero terms cancel
		\[ m^0_3(x, \rho_2, \mu_4^0(\rho_3, \rho_2, \rho_1, \rho_{41})) = m^0_4(x, \mu^0_2(\rho_2, \rho_3), \rho_2, \rho_1, \rho_{41}). \]
		
		\item \label{term:st5} \sloppy Consider terms of the form 
		$m_{1+n-k}^{w-w'}(x,a_1,\dots,a_j,\mu^{w'}_k(a_{j+1},\dots,
		a_{j+k}),\allowbreak a_{j+k+1},\dots,a_n)$ for a right-extended algebra operation $\mu^{w'}_k$. In this case, the algebra operation is not permitted to be the last input to the module operation.
		We cancel these terms against
		$m_{1+n}^w\left(x,a_1,\dots,a_{j+k-1},\mu^0_2(a_{j+k},a_{j+k+1})\dots,a_n\right)$.
		
		We produce a module tiling pattern corresponding to $m_{1+n}^w\left(x,a_1,\dots,a_{j+k-1},\mu^0_2(a_{j+k},a_{j+k+1})\dots,a_n\right)$ entirely analogously to \ref{term:st4}. In this case, the gluing of the tiling algebra pattern into the module tiling pattern respects the labeling of the $2$-valent vertices on the right. The fact that the algebra operation is not the last input ensures we do not modify the module tiling pattern near the blue-red boundary, and that we produce a bona fide module tiling pattern for the cancelling term.
	
		\item \label{term:st6}
		Consider terms of the form $m_{1+n-k}^{w-w'}\left(x,\mu^{w'}_k\left(a_1,\dots,
		a_k\right),a_{k+1},\dots,a_n\right)$ for a left-extended algebra operation $\mu^{w'}_k$, or terms of the form
		$m_{1+n-k}^{w-w'}\left(x,a_{1},\dots,a_{n-k},\mu^{w'}_k\left(a_{n-k+1},
		\dots,a_n\right)\right)$ for a right-extended algebra operation $\mu^{w'}_k$. We describe the cancellation for the right-extended case; the left-extended case follows analogously as in \ref{term:st3}.
		
		Consider a module tiling pattern $\Gamma_m$ for the module operation and an algebra tiling pattern $\Gamma_a$ for the algebra operation. Along the red boundary in $\Gamma_m$, the labels of the last face match with the output of the algebra operation. In other words, the labels of the face match with the labels along the $2$-valent vertices of $\Gamma_a$. Decompose $\Gamma_a$ into arcs intersecting transversely. Just as in \ref{term:st4}, we wish to glue in $\Gamma_a$ into $\Gamma_m$, identifying the arc $s$ containing the $2$-valent vertices with the edges of the face $f$.
		
		 Such a gluing, however, would not produce a legitimate module tiling pattern. The corresponding immersed disk, in particular, would have an obtuse angle at the intersection of the $\alpha$-arc and $\beta$-circle. Instead, after gluing in $\Gamma_a$, we \emph{split} the red vertices on the arc $s$ using the edges in $\Gamma^\beta(\HD_{st})$. In other words, we perform a (partial) inverse of the operation depicted in Figure~\ref{fig:solid-torus-cancellations-2a}. This decomposes our tiling pattern into two module tiling patterns, $\Gamma_1$ and $\Gamma_2$, that can be performed in order. This describes an inverse to the cancellation of this term in \ref{term:st3}; the inverse is illustrated in
		Figure~\ref{fig:solid-torus-cancellations-2}.
		 
		 The inverse admits a clean description in terms of the immersed disks. As we observed, gluing the immersed disks corresponding to the two module tiling operations along their common blue boundary---aligning the \emph{target} of the first with the \emph{source} of the second---yields an immersed disk with an obtuse angle at the intersection of the $\alpha$-arc and $\beta$-circle. This admits a cut along the $\beta$-circle, and when the boundary branch point goes out to the $\alpha$-arc, the disk decomposes into two disks with acute angles. $\Gamma_1$ and $\Gamma_2$ are the module tiling patterns associated to these two disks.

		\item \label{term:st7} Consider terms of the form $m^{w-1}_{2+n}(x, \mu_0^1, a_1, \dots, a_n)$ or $m^{w-1}_{2+n}(x, a_1, \dots, a_n, \mu_0^1)$. We describe the cancellation for the case where $\mu_0^1$ is the first input; the other case follows analogously. As in \ref{term:st3}, observe that there is exactly one length four chord in $\mu_0^1$ for which the labeling restrictions allow a module tiling pattern.
		
		Let $\Gamma$ be a module tiling pattern for the operation. The first face along the red boundary has four labels, and so five edges; two of these edges go out to the red boundary. We wish to push \emph{in} these two edges and glue them together, increasing the number of short cycles (and thus weight) by one. However, akin to \ref{term:st6}, such a gluing modifies the tiling pattern near the blue-red boundary. Instead, we follow the same procedure as in \ref{term:st6}. After gluing the two edges, we follow the edge and \emph{split} the vertices along it like in Figure~\ref{fig:solid-torus-cancellations-2a}. This decomposes our tiling pattern into two module tiling patterns, $\Gamma_1$ and $\Gamma_2$, that can be performed in order. This describes an inverse to the cancellation of this term in \ref{term:st3}; the inverse is illustrated in Figure~\ref{fig:solid-torus-cancellations-2}.
	
		 The inverse admits a clean description in terms of the immersed disks. Consider the $\alpha$-arcs on either side of the length four Reeb chord that occurs as the first or last Reeb chord for an immersed disk. If we glue these $\alpha$-arcs together, we obtain an immersed disk with an obtuse angle at the intersection of the $\alpha$-arc and $\beta$-circle. This admits a cut along the $\beta$-circle, and when the boundary branch point goes out to the $\alpha$-arc, the disk decomposes into two disks with acute angles. $\Gamma_1$ and $\Gamma_2$ are the module tiling patterns associated to these two disks.

		As an example, in the weight $1$ $\Ainf$ relation with inputs $x, \rho_2, \rho_1, \rho_{41}, \rho_4, \rho_{34}, \rho_3, \rho_{23}$, the following two non-zero terms cancel
		\begin{align*}
		m^0_9(x, \mu^1_0, \rho_2, \rho_1, \rho_{41}, \rho_4, \rho_{34}, \rho_3, \rho_{23}) &= m^0_9(x, \rho_{4123}, \rho_2, \rho_1, \rho_{41}, \rho_4, \rho_{34}, \rho_3, \rho_{23}) \\
		&= m^1_6(m_3^0(x, \rho_2, \rho_1), \rho_{41}, \rho_4, \rho_{34}, \rho_3, \rho_{23}).
		\end{align*}
		
		The immersed disks corresponding to this cancellation are depicted in Figure~\ref{fig:solid-torus-cancellation-7}.
		
\begin{figure}[h!tbp]
	\centering
		\begin{tikzpicture}[scale=1.5]
					\begin{scope}[shift={(0,2)}]
				\node at (.25,1.7) (r4) {$\rho_3$};
				\node at (1.7,1.7) (r1) {$\rho_4$};
				\node at (1.7,.3) (r2) {$\rho_1$};
				\node at (.3,.3) (r3) {$\rho_2$};
				\draw[dashed] (.25,0) arc[radius=.25, start angle = 0, end angle = 90];
				\draw[dashed] (0,1.75) arc[radius=.25, start angle = 270, end angle = 360];
				\draw[dashed] (1.75,0) arc[radius=.25, start angle = 180, end angle = 90];
				\draw[dashed] (1.75,2) arc[radius=.25, start angle = 180, end angle = 270];
				\draw[color = red,thick] (.25,0) to (1.75,0);
				\draw[color = red,thick] (0,0.25) to (0,1.75);
				\draw[color = red,thick] (.25,2) to (1.75,2);
				\draw[color = red,thick] (2,0.25) to (2,1.75);
		
				\draw[color = blue,thick] (0,1) to (2,1);
			\end{scope}
			\begin{scope}[shift={(2,2)}]
				\node at (.25,1.7) (r4) {$\rho_3$};
				\node at (1.7,1.7) (r1) {$\rho_4$};
				\node at (1.7,.3) (r2) {$\rho_1$};
				\node at (.3,.3) (r3) {$\rho_2$};
				\draw[dashed] (.25,0) arc[radius=.25, start angle = 0, end angle = 90];
				\draw[dashed] (0,1.75) arc[radius=.25, start angle = 270, end angle = 360];
				\draw[dashed] (1.75,0) arc[radius=.25, start angle = 180, end angle = 90];
				\draw[dashed] (1.75,2) arc[radius=.25, start angle = 180, end angle = 270];
				\draw[color = red,thick] (.25,0) to (1.75,0);
				\draw[color = red,thick] (0,0.25) to (0,1.75);
				\draw[color = red,thick] (.25,2) to (1.75,2);
				\draw[color = red,thick] (2,0.25) to (2,1.75);
		
				\draw[color = blue,thick] (0,1) to (2,1);
			\end{scope}
			\begin{scope}[shift={(2,0)}]
				\node at (.25,1.7) (r4) {$\rho_3$};
				\node at (1.7,1.7) (r1) {$\rho_4$};
				\node at (1.7,.3) (r2) {$\rho_1$};
				\node at (.3,.3) (r3) {$\rho_2$};
				\draw[dashed] (.25,0) arc[radius=.25, start angle = 0, end angle = 90];
				\draw[dashed] (0,1.75) arc[radius=.25, start angle = 270, end angle = 360];
				\draw[dashed] (1.75,0) arc[radius=.25, start angle = 180, end angle = 90];
				\draw[dashed] (1.75,2) arc[radius=.25, start angle = 180, end angle = 270];
				\draw[color = red,thick] (.25,0) to (1.75,0);
				\draw[color = red,thick] (0,0.25) to (0,1.75);
				\draw[color = red,thick] (.25,2) to (1.75,2);
				\draw[color = red,thick] (2,0.25) to (2,1.75);
		
				\draw[color = blue,thick] (0,1) to (2,1);
				\draw [line width = 1.5pt] (0,1) to (1.4,1);
				\node at (1.4,1)[circle,fill,inner sep=2pt]{};
			\end{scope}
			\begin{scope}[shift={(0,0)}]
				\node at (.25,1.7) (r4) {$\rho_3$};
				\node at (1.7,1.7) (r1) {$\rho_4$};
				\draw[dashed] (0,1.75) arc[radius=.25, start angle = 270, end angle = 360];
				\draw[dashed] (1.75,2) arc[radius=.25, start angle = 180, end angle = 270];
				\draw[color = red,thick] (0,1) to (0,1.75);
				\draw[color = red,thick] (.25,2) to (1.75,2);
				\draw[color = red,thick] (2,1) to (2,1.75);
		
				\draw[color = blue,thick] (0,1) to (2,1);
				
				\draw [line width = 4pt, draw=purple, opacity=0.4] (2,1) to (2,1.75);
			\end{scope}
		\end{tikzpicture}
		\caption [An immersed disk corresponding to the term cancellation in \ref{term:st7} in Theorem~\ref{prop:solid-torus-ainf}.] {\textbf{An immersed disk corresponding to the term cancellation in \ref{term:st7} in Theorem~\ref{prop:solid-torus-ainf}.} There is a cut in the disk along the shaded edge, such that the first Reeb chord is $\rho_{4123}$. Gluing this and cutting along the black line decomposes the disk into two disks, the module tiling operations for which are $\Gamma_1$ and $\Gamma_2$.}
		\label{fig:solid-torus-cancellation-7}
\end{figure}
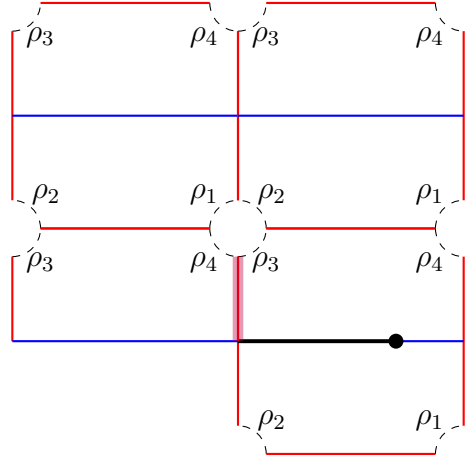

		\item \label{term:st8} Consider of the form $m_{2+n}^{w-1}\left(x,a_1,\dots,a_i,\mu^1_0,a_{i+1},
		\dots,a_n\right)$. In this case, the $\mu_0^1$ does not appear extremally. These cancel against terms of the form
		$m_{n}^w\left(x,a_1,\dots,\mu^0_2(a_i,a_{i+1}),\dots,a_n\right)$. Note that this is precisely the operation that move (3) produces, and is thus non-zero. In terms of module tiling patterns it corresponds to the operation depicted in Figure~\ref{fig:move-3}. That is, let $\Gamma$ be a module tiling pattern for the operation. The face corresponding to $\mu_0^1$ has five edges; we push \emph{in} and glue together the two edges that go to the red boundary.
		
		\item \label{term:st9} Terms with any centered algebra operation have an idempotent output. As our type A module is strictly unital, for the term to be non-zero the module operation must be a $m_2^0$. These
		were cancelled against in \ref{term:st2}. That is, given a tiling pattern for a centered algebra operation, we \emph{split} the vertices along the arc that follows from the root. This decomposes the tiling pattern into two module tiling patterns $\Gamma_1$ and $\Gamma_2$ that can be performed in order. The inverse is illustrated in
		Figure~\ref{fig:solid-torus-cancellations-2a}.
		
		\item \label{term:st10} As our last case, consider any term where the algebra operation is a $\mu_2^0$. We claim that all such terms have been cancelled against in \ref{term:st1}, \ref{term:st4}, \ref{term:st5}, and \ref{term:st8}.
		
		Let $\Gamma$ be a module tiling pattern for the module operation. The algebra operation $\mu_0^2$ factorizes a term in its chord sequence. Along the face corresponding to this term, there is an edge $e$ that \emph{exhibits} this factorization. We push out the edge to the red boundary.
		
		Consider the face $f$ on the other side of the edge $e$. 
		If $f$ is an internal face, then pushing out the edge $e$ punctures this face. This is the inverse to the cancellation in \ref{term:st8}.
		Otherwise,
		$f$ touches the boundary.
		If it touches the blue boundary,
		pushing out the edge breaks $\Gamma$ into two module tiling patterns. This is the inverse to the cancellation in \ref{term:st1}.
		Otherwise, $f$ touches the red boundary. If the other side of $e$
		is visible from the red boundary after the edge $e$, pushing out the edge is the inverse to the cancellation in \ref{term:st4}. Otherwise, it is the inverse to the cancellation in \ref{term:st5}. 
	\end{enumerate}
	
	We can now breathe. We have considered the contribution to the $\Ainf$ relations of every weighted tree with two internal vertices, and shown how to pair these contributions up so they cancel. This verifies the relations.
\end{proof}

\subsection{Cables}

With our motivating example out of the way, we are ready to tackle the case of cables. We do this piecemeal as well, starting with the $(2,1)$-cable. Our constructions and arguments will follow closely those for the solid torus.

To have utility in the computation of knot Floer invariants, our bordered Heegaard diagrams need to be doubly-pointed, with basepoints $w$ and $z$. We use a variable $U$ to count the intersection of a pseudoholomorphic disc with $w$, and a variable $V$ to count the intersection of a disc with $z$. We consider the enriched torus algebra $\Algm^{U,V} = \Algm \tensor_{\FF_2[U]} \FF_2[U,V]$. An operation $\mu^w_n(\vec a) = U^k y$ for $\Algm$ corresponds to an operation $\mu^w_n(\vec a) = U^kV^k y$ for $\Algm^{U,V}$; that is, each copy of $U$ is replaced by a copy of $UV$.

The weighted type A module for the doubly-pointed Heegaard diagrams are modules over $\Algm^{U,V}$. The operations for these modules correspond to module tiling patterns, as defined in Section~\ref{sec:tiling-patterns}. The sole modification we make is that each vertex records a power of both $U$ and $V$ to keep track of the multiplicity of intersection with $w$ and $z$. The red vertices have the property that $U_v = U$ and $V_v = V$.

\subsubsection{The \texorpdfstring{$(2,1)$}{(2,1)}-cable}

\begin{figure}
	\centering
	\begin{tikzpicture}[scale=2]
		\draw[dashed] (.25,0) arc[radius=.25, start angle = 0, end angle = 90];
		\draw[dashed] (0,1.75) arc[radius=.25, start angle = 270, end angle = 360];
		\draw[dashed] (1.75,0) arc[radius=.25, start angle = 180, end angle = 90];
		\draw[dashed] (1.75,2) arc[radius=.25, start angle = 180, end angle = 270];
		\draw[color = red,thick] (.25,0) to (1.75,0);
		\draw[color = red,thick] (0,0.25) to (0,1.75);
		\draw[color = red,thick] (.25,2) to (1.75,2);
		\draw[color = red,thick] (2,0.25) to (2,1.75);
		\node at (.25,1.7) (r4) {$\rho_4$};
		\node at (1.7,1.7) (r1) {$\rho_1$};
		\node at (1.7,.3) (r2) {$\rho_2$};
		\node at (.3,.3) (r3) {$\rho_3$};
		\node at (.4,1.9) (z) {$z$};
		\node at (1,0.3) (w) {$w$};
		\draw[color = blue,thick] (0,1) to[in=270,out=0] (.667,2);
		\draw[color = blue,thick] (2,1) to[in=270,out=180] (1.333,2);
		\draw[color = blue,thick] (.667,0) to[in=180,out=90] (1,0.75) to[in=90,out=0] (1.333,0);
		\node at (-.2,1) (x1) {$x$};
		\node at (2.2,1) (x2) {\phantom{$x$}};
		\node at (0.667,-.2) (b1) {$b$};
		\node at (1.333,-.225) (c1) {$c$};
		
		\draw[color = gray] (1,1) to[in=0,out=210] (0,.7);
		\draw[color = gray] (1,1) to[in=180,out=330] (2,.7);
		\draw[color = purple,densely dotted] (1,1) to (0.333,1.5);
		\draw[color = purple,densely dotted] (1,1) to (1.667,1.5);
		\draw[color = purple,densely dotted] (1,1) to (1,.5);
		\draw[color = gray] (1,.5) to (1,0);
		\draw[color = gray] (1,1) to (1,2);
		\draw[color = gray] [in=270,out=45] (0.333,1.5) to (0.5,2);
		\draw[color = gray] [in=0,out=225] (0.333,1.5) to (0,1.333);
		\draw[color = gray] [in=270,out=135] (1.667,1.5) to (1.5,2);
		\draw[color = gray] [in=180,out=315] (1.667,1.5) to (2,1.333);
		\draw[color = gray] [in=90,out=210] (1,1) to (0.5,0);
		\draw[color = gray] [in=90,out=330] (1,1) to (1.5,0);
	
		\node at (1,1) [circle,fill,inner sep=1pt]{};
		\node at (1,.5) [circle,fill,inner sep=1pt]{};
		\node at (0.333,1.5) [circle,fill,inner sep=1pt]{};
		\node at (1.667,1.5) [circle,fill,inner sep=1pt]{};
	\end{tikzpicture}  
	\caption[A bordered Heegaard diagram \texorpdfstring{$\CD_2$}{CD2} for the \texorpdfstring{$(2,1)$}{(2,1)}-cable, and its dual graph.]
	{\textbf{A bordered Heegaard diagram $\CD_2$ for the $(2,1)$-cable, and its dual graph.}}
	\label{fig:cable-bhd}
\end{figure}
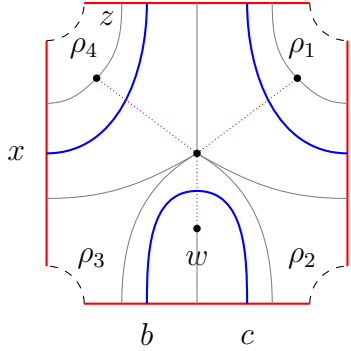

The bordered Heegaard diagram $\CD_2$ for the $(2,1)$-cable is pictured in Figure~\ref{fig:cable-bhd}. As a left $\FF_2[U,V]\langle \iota_0, \iota_1 \rangle$-module, $\CFAm(\CD_2)$ is freely generated by $x$, $b$, and $c$. Figure~\ref{fig:cable-bhd} also depicts the dual graph $\Gamma(\CD_2)$.

An immersed disk for $\CD_2$ from any $\x$ to $\y$ with an obtuse intersection along the $\alpha$-arc and $\beta$-circle admits a cut along the $\beta$-circle. This decomposes the disk into two disks of index one, and so the proof Proposition~\ref{prop:solid-torus-obtuse} repeats.

\begin{proposition}
	\label{prop:c2-obtuse}
	Let $u: \bD \to T^2$ be a disk for $\CD_{2}$ from $\x$ to $\y$, and suppose further that it maps the corners at $-i$ and $i$ to \emph{obtuse} intersections of the $\alpha$-arcs and $\beta$-circles.
	
	As in Proposition~\ref{prop:acute-angled-disk}, let $B$ be the positive domain in $\pi_2(\x, \y)$, $\vec a$ the sequence of Reeb chords, and $w$ the number of simple Reeb orbits corresponding to $u$.	
	Then, $\ind(B, \vec a, w) > 1$.
\end{proposition}

\begin{figure}
	\centering
	\begin{tikzpicture}
	\begin{pgfonlayer}{nodelayer}
		\node [style=none] (0) at (-3, 0) {};
		\node [style=none] (1) at (3, 0) {};
		\node [style=none] (2) at (-2.95, -0.5) {};
		\node [style=none] (3) at (2.95, -0.5) {};
		\node [style=none] (4) at (0, 2.65) {};
		\node [style=none] (5) at (0, -2.65) {};
		\node [style=none] (6) at (-3.5, 0) {$b$};
		\node [style=none] (7) at (3.5, 0) {$c$};
		\node [style=dot] (8) at (0, 0) {};
	\end{pgfonlayer}
	\begin{pgfonlayer}{edgelayer}
		\draw [style=blue, bend left=90, looseness=1.50] (0.center) to (1.center);
		\draw [style=red, bend right=90, looseness=1.50] (0.center) to (1.center);
		\draw (4.center) to (8);
		\draw (8) to (5.center);
	\end{pgfonlayer}
\end{tikzpicture} \qquad \qquad
	\begin{tikzpicture}
	\begin{pgfonlayer}{nodelayer}
		\node [style=none] (0) at (-3, 0) {};
		\node [style=none] (1) at (3, 0) {};
		\node [style=none] (2) at (-2.95, -0.5) {};
		\node [style=none] (3) at (2.95, -0.5) {};
		\node [style=none] (4) at (0, 2.65) {};
		\node [style=none] (6) at (-3.5, 0) {$x$};
		\node [style=none] (7) at (3.5, 0) {$c$};
		\node [style=dot] (10) at (0, -0.05) {};
		\node [style=number-labels] (11) at (0, -0.4) {$1$};
	\end{pgfonlayer}
	\begin{pgfonlayer}{edgelayer}
		\draw [style=blue, bend left=90, looseness=1.50] (0.center) to (1.center);
		\draw [style=red, bend right=90, looseness=1.50] (0.center) to (1.center);
		\draw [bend left=15] (2.center) to (3.center);
		\draw (10) to (4.center);
	\end{pgfonlayer}
\end{tikzpicture} \qquad \qquad
	\begin{tikzpicture}
	\begin{pgfonlayer}{nodelayer}
		\node [style=none] (0) at (-3, 0) {};
		\node [style=none] (1) at (3, 0) {};
		\node [style=none] (2) at (-2.95, -0.5) {};
		\node [style=none] (3) at (2.95, -0.5) {};
		\node [style=none] (4) at (0, 2.65) {};
		\node [style=none] (6) at (-3.5, 0) {$b$};
		\node [style=none] (7) at (3.5, 0) {$x$};
		\node [style=dot] (10) at (0, -0.05) {};
		\node [style=number-labels] (11) at (0, -0.4) {$4$};
	\end{pgfonlayer}
	\begin{pgfonlayer}{edgelayer}
		\draw [style=blue, bend left=90, looseness=1.50] (0.center) to (1.center);
		\draw [style=red, bend right=90, looseness=1.50] (0.center) to (1.center);
		\draw [bend left=15] (2.center) to (3.center);
		\draw (10) to (4.center);
	\end{pgfonlayer}
\end{tikzpicture}
	\begin{tikzpicture}
	\begin{pgfonlayer}{nodelayer}
		\node [style=none] (0) at (-3, 0) {};
		\node [style=none] (1) at (3, 0) {};
		\node [style=none] (2) at (-2.95, -0.5) {};
		\node [style=none] (3) at (2.95, -0.5) {};
		\node [style=none] (4) at (0, 2.65) {};
		\node [style=none] (6) at (-3.5, 0) {$x$};
		\node [style=none] (7) at (3.5, 0) {$b$};
		\node [style=dot] (10) at (0, -0.05) {};
		\node [style=number-labels] (11) at (-0.5, -0.4) {$3$};
		\node [style=number-labels] (12) at (0.5, -0.4) {$1$};
		\node [style=none] (13) at (-1, -2.5) {};
		\node [style=none] (14) at (1, -2.5) {};
		\node [style=number-labels] (15) at (0, -0.65) {$2$};
	\end{pgfonlayer}
	\begin{pgfonlayer}{edgelayer}
		\draw [style=blue, bend left=90, looseness=1.50] (0.center) to (1.center);
		\draw [bend left=15] (2.center) to (3.center);
		\draw (10) to (4.center);
		\draw [bend right=15, looseness=0.25] (10) to (13.center);
		\draw [bend left=15, looseness=0.25] (10) to (14.center);
		\draw [style=red, bend right=90, looseness=1.50] (0.center) to (1.center);
	\end{pgfonlayer}
\end{tikzpicture} \qquad \qquad
	\begin{tikzpicture}
	\begin{pgfonlayer}{nodelayer}
		\node [style=none] (0) at (-3, 0) {};
		\node [style=none] (1) at (3, 0) {};
		\node [style=none] (2) at (-2.95, -0.5) {};
		\node [style=none] (3) at (2.95, -0.5) {};
		\node [style=none] (4) at (0, 2.65) {};
		\node [style=none] (6) at (-3.5, 0) {$c$};
		\node [style=none] (7) at (3.5, 0) {$x$};
		\node [style=dot] (10) at (0, -0.05) {};
		\node [style=number-labels] (11) at (-0.5, -0.4) {$4$};
		\node [style=number-labels] (12) at (0.5, -0.4) {$2$};
		\node [style=none] (13) at (-1, -2.5) {};
		\node [style=none] (14) at (1, -2.5) {};
		\node [style=number-labels] (15) at (0, -0.65) {$3$};
	\end{pgfonlayer}
	\begin{pgfonlayer}{edgelayer}
		\draw [style=blue, bend left=90, looseness=1.50] (0.center) to (1.center);
		\draw [bend left=15] (2.center) to (3.center);
		\draw (10) to (4.center);
		\draw [bend right=15, looseness=0.25] (10) to (13.center);
		\draw [bend left=15, looseness=0.25] (10) to (14.center);
		\draw [style=red, bend right=90, looseness=1.50] (0.center) to (1.center);
	\end{pgfonlayer}
\end{tikzpicture} \qquad \qquad
	\begin{tikzpicture}
	\begin{pgfonlayer}{nodelayer}
		\node [style=none] (0) at (-3, 0) {};
		\node [style=none] (1) at (3, 0) {};
		\node [style=none] (2) at (-2.95, -0.5) {};
		\node [style=none] (3) at (2.95, -0.5) {};
		\node [style=none] (4) at (0, 2.65) {};
		\node [style=none] (6) at (-3.5, 0) {$c$};
		\node [style=none] (7) at (3.5, 0) {$b$};
		\node [style=dot] (10) at (0, -0.05) {};
		\node [style=number-labels] (11) at (-1, -0.45) {$2$};
		\node [style=number-labels] (12) at (1, -0.45) {$3$};
		\node [style=none] (13) at (-1.6, -2.275) {};
		\node [style=none] (14) at (1.6, -2.275) {};
		\node [style=none] (15) at (0, -2.625) {};
		\node [style=number-labels] (16) at (-0.425, -0.925) {$1$};
		\node [style=number-labels] (18) at (0.425, -0.925) {$4$};
	\end{pgfonlayer}
	\begin{pgfonlayer}{edgelayer}
		\draw [style=blue, bend left=90, looseness=1.50] (0.center) to (1.center);
		\draw [bend left=15] (2.center) to (3.center);
		\draw (10) to (4.center);
		\draw [bend right, looseness=0.50] (10) to (13.center);
		\draw [bend left=15, looseness=0.75] (10) to (14.center);
		\draw [style=red, bend right=90, looseness=1.50] (0.center) to (1.center);
		\draw (10) to (15.center);
	\end{pgfonlayer}
\end{tikzpicture}
	\begin{tikzpicture}
	\begin{pgfonlayer}{nodelayer}
		\node [style=none] (0) at (-3, 0) {};
		\node [style=none] (1) at (3, 0) {};
		\node [style=none] (2) at (-2.95, -0.5) {};
		\node [style=none] (3) at (2.95, -0.5) {};
		\node [style=none] (4) at (0, 2.65) {};
		\node [style=none] (6) at (-3.5, 0) {$x$};
		\node [style=none] (7) at (3.5, 0) {$x$};
		\node [style=dot] (10) at (0, -0.05) {};
		\node [style=none] (11) at (0, -2.65) {};
		\node [style=dot] (12) at (0, 1.25) {};
		\node [style=none] (13) at (-2.7, 1.225) {};
		\node [style=none] (14) at (2.7, 1.225) {};
		\node [style=number-labels] (15) at (-0.45, -0.5) {$3$};
		\node [style=number-labels] (16) at (0.45, -0.5) {$2$};
	\end{pgfonlayer}
	\begin{pgfonlayer}{edgelayer}
		\draw [bend left=15] (2.center) to (3.center);
		\draw (10) to (11.center);
		\draw (4.center) to (12);
		\draw (12) to (10);
		\draw [bend right=15, looseness=0.50] (13.center) to (10);
		\draw [bend right=15, looseness=0.50] (10) to (14.center);
		\draw [style=blue, bend left=90, looseness=1.50] (0.center) to (1.center);
		\draw [style=red, bend right=90, looseness=1.50] (0.center) to (1.center);
	\end{pgfonlayer}
\end{tikzpicture}
	\caption[The building blocks for operations for the \texorpdfstring{$(2,1)$}{(2,1)}-cable.]
	{\textbf{The building blocks for operations for the $(2,1)$-cable.}}
	\label{fig:cable-2-building-blocks}
\end{figure}
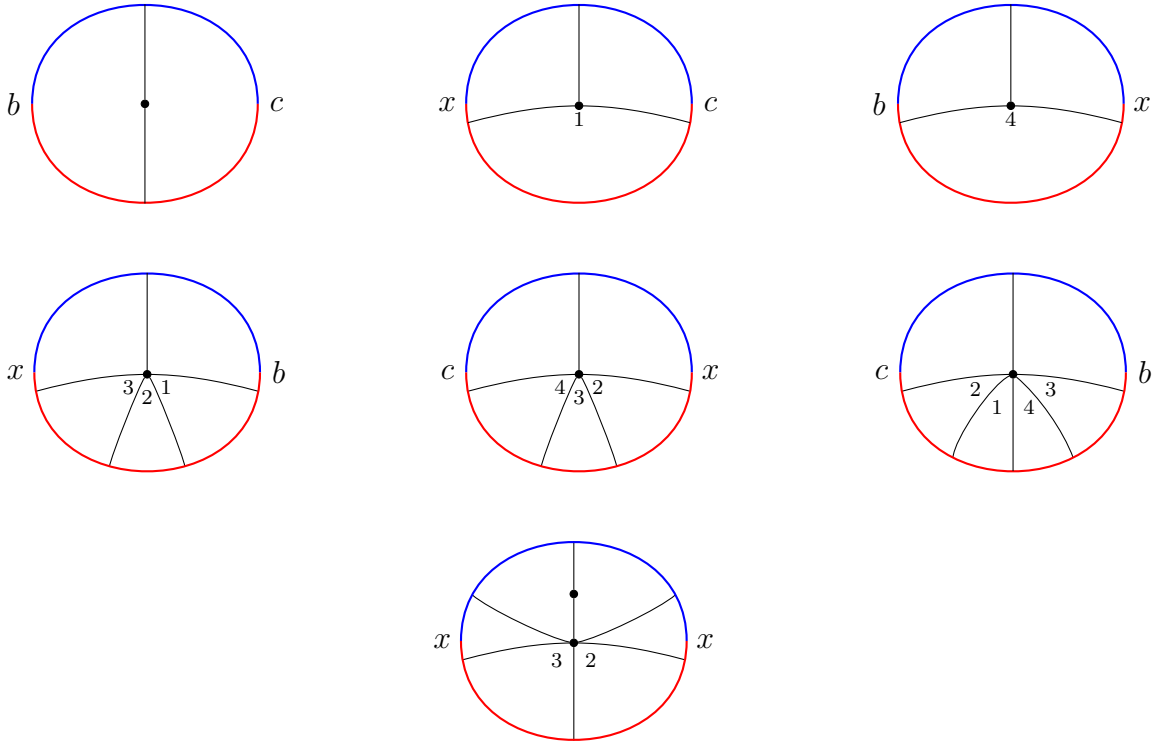

We can characterize all operations on the type A module $\CFAm(\CD_2)$ in terms of the simple module tiling patterns of Figure~\ref{fig:cable-2-building-blocks}. These are precisely the module tiling patterns with no red vertices, and where move (1) cannot be used to simplify them.

\begin{proposition}
	\label{prop:c2-building-blocks}
	Any operation on the type A module $\CFAm(\CD_2)$ can be obtained as the operation associated to a module tiling pattern constructed by starting with the building blocks in Figure~\ref{fig:cable-2-building-blocks}, and performing a finite sequence of the moves catalogued in Section~\ref{sec:three-moves}.
\end{proposition}

\begin{proof}
	This follows closely the proof of Proposition~\ref{prop:solid-torus-building-blocks}. Let $u: \bD \to T^2$ be an immersed disk from $\x$ to $\y$ for $\CD_2$. There is a neighborhood of the $\beta$-circle on the torus that does not contain the marked point $p$, but contains the vertices in $\Gamma(\CD_2)$. In this neighborhood, $u$ lifts to the universal cover of the torus. We can read off the graph $\Gamma(u)$ near the blue boundary from this blue. In particular, near the blue boundary $\Gamma(u)$ is a composition of the graphs in Figure~\ref{fig:cable-2-building-blocks}.

	The graph $\Gamma(u)$ retains its nice properties from those for $\HD_{st}$. It is connected and planar. Away from the blue boundary,
	every vertex of $\Gamma(u)$ is $4$-valent, and every internal face is bounded by four edges.
	
	Using the inverse of moves (1) and (3), and repeating the arguments from the
	proof of Proposition~\ref{prop:solid-torus-building-blocks}, we reduce
	to the case where $\Gamma(u)$ has no internal faces and no edge that separates the blue and the red boundary.
	To find an arc using which we can eliminate a red vertex with the inverse of move (2) in Proposition~\ref{prop:solid-torus-building-blocks}, we started by deleting the spine of the underlying graph. For $\CD_2$, a description of the spine can get hairy. Equivalently, we delete the vertices and edges visible from the blue boundary. We are left with a collection of
	trees, and all red vertices from these can
	be eliminated as in the proof of
	Proposition~\ref{prop:solid-torus-building-blocks}, using inverses of
	move (2).
\end{proof}

\begin{figure}
\centering
\begin{tikzpicture}[scale=2]
	\node at (0,0) (x) {$x$};
	\node at (2,1) (c) {$c$};
	\node at (2,-1) (b) {$b$};
	\draw[->, bend right=15] (x) to node[below]{\lab{\rho_1}} (c);
	\draw[->, bend right=15] (c) to node[above,sloped]{\lab{UV\rho_4\otimes\rho_3\otimes\rho_2}} (x);
	\draw[->, bend left=80] (x) to node[left]{\lab{U^2\rho_3\otimes\rho_2}} (-.5,0) to (x);
	\draw[->, bend right=15] (b) to node[right]{\lab{U}} (c);
	\draw[->, bend right=15] (c) to node[above, sloped]{\lab{V\rho_2\otimes\rho_1\otimes\rho_4\otimes\rho_3}} (b);
	\draw[->, bend right=15] (x) to node[below, sloped]{\lab{U\rho_3\otimes\rho_2\otimes\rho_1}} (b);
	\draw[->, bend right=15] (b) to node[above]{\lab{V\rho_4}} (x);
\end{tikzpicture}
	\caption[A representation of the type A module \texorpdfstring{$\CFAm(\CD_2)$}{CFA minus of the (2,1)-cable}.]
	{\textbf{A representation of the type A module $\CFAm(\CD_2)$.}}
	\label{fig:cable-2-graph}
\end{figure}
The type A module $\CFAm(\CD_2)$ can be represented by the graph in Figure~\ref{fig:cable-2-graph}. The nodes correspond to the generators, and the edges to the building blocks depicted in Figure~\ref{fig:cable-2-building-blocks}. Every other operation can be generated from the edges using the three moves.

A similar calculus to that of Theorem~\ref{prop:solid-torus-ainf} for $\CFAm(\HD_{st})$ allows us to prove that the module $\CFAm(\CD_2)$ satisfies the $\Ainf$ relations.

\begin{theorem}
	The type A module $\CFAm(\CD_2)$, with operations $\{m^w_{1+n}\}$ given by counting its module tiling patterns of chord sequence length $n$ and weight $w$, satisfies the $\Ainf$ structure relations.
	\label{prop:cable-2-ainf}
\end{theorem}

\begin{proof}
	We verify the $\Ainf$ relation for a fixed input $(\x, a_1, \dots, a_n)$ and weight $w$. This proof follows the same schema as that of Theorem~\ref{prop:solid-torus-ainf}. In particular, we show the $\Ainf$ relation holds by pairing non-zero terms that appear in it.
	
	It suffices to consider the case where all the $a_i$ are Reeb elements.
	Most pairings for the cancellations sail through in analogy with Theorem~\ref{prop:solid-torus-ainf}. A cause of trouble is the $1$-input operation $m_1^0(b) = Uc$. In terms of its module tiling pattern, we observe a $2$-valent vertex and an internal edge on either side of which lies the blue boundary. This we claimed in \ref{term:st3} does not happen. There is another $2$-valent vertex and such an internal edge in the operation $m_3^0(x, \rho_3, \rho_2) = U^2 x$, and we shall see that this leads to cancellations.

	The following is how the terms cancel.
	
	\begin{enumerate}[wide,label=(TC-\arabic*)]
		\item \label{term:c1} Consider terms of the form $m_{1+n-j}^{w-w'}\left(m_{1+j}^{w'}\left(\x,
		a_1,\dots,a_j\right),a_{j+1},\dots,a_n\right)$, where $a_j a_{j+1}\neq0$.
		These cancel against $m_{n}^w\left(
		\x,a_1,\dots,\mu_2^0\left(a_j,a_{j+1}\right),\dots,a_n\right)$.
		This remains precisely the operation that move (1) produces. For example, in the weight $0$ $\Ainf$ relation with inputs $x, \rho_{3}, \rho_2, \rho_{3}, \rho_{2}, \rho_{1}$, the following two non-zero terms cancel
		\[m^0_4(m^0_3(x, \rho_{3}, \rho_2), \rho_{3}, \rho_{2}, \rho_{1}) = U^3 b = m^0_6(x, \rho_{3}, \mu^0_2(\rho_2, \rho_{3}), \rho_2, \rho_1).\]
		
		\item \label{term:c2} Consider terms of the form $m_{1+n-j}^{w-w'}\left(m_{1+j}^{w'}(\x,
		a_1,\dots,a_j),a_{j+1},\dots,a_n\right)$, where $a_j a_{j+1}=0$.
		
		Let $\Gamma_1$ and $\Gamma_2$ be the module tiling patterns for the first and second operations respectively. Because $a_ja_{j+1} = 0$, the two tiling patterns
		share a \emph{common} blue boundary. We can glue the two patterns along this common boundary, while identifying edges that touch this blue boundary. These edges can be located as the pre-image of the edge in $\Gamma^\beta(\HD_{st})$. We contract these edges and delete parallel edges. In other words, we glue the immersed disks for the two operations along their common blue boundary and consider the associated tiling pattern.

		For this case, suppose the length of the blue boundaries of $\Gamma_1$ and $\Gamma_2$ match. 
		In this case, the output of the second tiling operation is necessarily a scalar multiple of $\x$. The two patterns glue together perfectly and yield a tiling pattern with no blue boundary. Placing the root at the first edge to the boundary produces an algebra tiling pattern corresponding to a centered algebra operation.
		Our term cancels against
		$m_2^{0}\left(\x,\mu_n^w\left(a_1,\dots,a_n\right)\right)$. The idempotents align as $\Gamma_1$ is a disk from $\x$, and so the first label must have compatible idempotents.
		An example of this is
		illustrated in Figure~\ref{fig:cable-2-ainf-2}.

\begin{figure}
	\centering
	\begin{tikzpicture}
	\begin{pgfonlayer}{nodelayer}
		\node [style=none] (0) at (-14.75, -2) {};
		\node [style=none] (1) at (-2.75, -2) {};
		\node [style=none] (2) at (-14.7, -2.5) {};
		\node [style=none] (3) at (-2.8, -2.5) {};
		\node [style=none] (4) at (-13.35, -0.25) {};
		\node [style=none] (6) at (-15.25, -1.75) {$x$};
		\node [style=none] (7) at (-2.25, -1.75) {$x$};
		\node [style=none] (10) at (-14.3, -3.5) {};
		\node [style=none] (11) at (-3.15, -3.5) {};
		\node [style=none] (12) at (-12.425, 0.125) {};
		\node [style=none] (14) at (-11.425, 0.4) {};
		\node [style=none] (16) at (-5.95, 0.375) {};
		\node [style=none] (18) at (-4.125, -0.25) {};
		\node [style=none] (20) at (-4.975, 0.1) {};
		\node [style=none] (35) at (-15, -1.25) {};
		\node [style=none] (36) at (-2.5, -1.25) {};
		\node [style=none] (37) at (-14.25, 3.25) {};
		\node [style=none] (38) at (-3.25, 3.25) {};
		\node [style=none] (46) at (-2.65, 1.5) {};
		\node [style=none] (66) at (-0.5, 0) {};
		\node [style=none] (67) at (0.5, 0) {};
		\node [style=none] (68) at (2.75, 2.75) {};
		\node [style=none] (69) at (2.75, -2.25) {};
		\node [style=bigdot] (74) at (2.2, 0) {};
		\node [style=none] (75) at (14.1, 0) {};
		\node [style=number-labels] (90) at (4.925, 0.1) {$3$};
		\node [style=number-labels] (92) at (7.975, 0.175) {$3$};
		\node [style=number-labels] (93) at (11.15, 0.075) {$3$};
		\node [style=number-labels] (95) at (5.375, 0.1) {$2$};
		\node [style=number-labels] (97) at (8.45, 0.175) {$2$};
		\node [style=number-labels] (98) at (11.6, 0.075) {$2$};
		\node [style=none] (100) at (13.55, 2.75) {};
		\node [style=none] (101) at (13.55, -2.25) {};
		\node [style=none] (104) at (5.275, 2.475) {};
		\node [style=none] (105) at (5.275, -1.225) {};
		\node [style=none] (108) at (8.2, 2.45) {};
		\node [style=none] (109) at (8.2, -1.175) {};
		\node [style=none] (110) at (11.225, 2.475) {};
		\node [style=none] (111) at (11.225, -1.225) {};
		\node [style=number-labels] (118) at (4.925, 0.675) {$4$};
		\node [style=number-labels] (120) at (7.975, 0.75) {$4$};
		\node [style=number-labels] (121) at (11.15, 0.65) {$4$};
		\node [style=number-labels] (123) at (5.375, 0.675) {$1$};
		\node [style=number-labels] (125) at (8.45, 0.75) {$1$};
		\node [style=number-labels] (126) at (11.6, 0.65) {$1$};
		\node [style=dot] (129) at (5.125, 0.375) {};
		\node [style=dot] (131) at (8.2, 0.45) {};
		\node [style=dot] (132) at (11.35, 0.35) {};
		\node [style=dot] (140) at (-12, -2.125) {};
		\node [style=dot] (141) at (-8.75, -2.05) {};
		\node [style=dot] (142) at (-5.5, -2.15) {};
		\node [style=dot] (143) at (-12.225, -0.925) {};
		\node [style=dot] (144) at (-8.75, -0.75) {};
		\node [style=dot] (145) at (-5.2, -1) {};
		\node [style=none] (146) at (-11.75, -3.25) {};
		\node [style=none] (147) at (-8.75, -3.25) {};
		\node [style=none] (148) at (-5.75, -3.25) {};
		\node [style=none] (149) at (-9.75, 0.625) {};
		\node [style=none] (150) at (-8.75, 0.65) {};
		\node [style=none] (151) at (-7.75, 0.6) {};
		\node [style=none] (152) at (-13.375, 0.5) {};
		\node [style=none] (153) at (-12.45, 0.875) {};
		\node [style=none] (154) at (-11.35, 1.175) {};
		\node [style=none] (155) at (-5.9, 1.1) {};
		\node [style=none] (156) at (-3.925, 0.425) {};
		\node [style=none] (157) at (-4.85, 0.8) {};
		\node [style=none] (158) at (-9.775, 1.375) {};
		\node [style=none] (159) at (-8.75, 1.4) {};
		\node [style=none] (160) at (-7.675, 1.375) {};
		\node [style=dot] (161) at (-11.45, 2.25) {};
		\node [style=dot] (162) at (-9.75, 2.5) {};
		\node [style=dot] (163) at (-7.65, 2.5) {};
		\node [style=dot] (164) at (-5.75, 2.25) {};
		\node [style=dot] (165) at (-3.85, 1.975) {};
		\node [style=dot] (166) at (-4.75, 3) {};
		\node [style=dot] (167) at (-8.75, 3.5) {};
		\node [style=dot] (168) at (-12.25, 3) {};
		\node [style=dot] (169) at (-13.25, 1.75) {};
		\node [style=none] (170) at (-14.875, 1.5) {};
		\node [style=none] (171) at (-13.075, 3.4) {};
		\node [style=none] (172) at (-11.45, 3.45) {};
		\node [style=none] (173) at (-12, 3.9) {};
		\node [style=none] (174) at (-9.575, 3.875) {};
		\node [style=none] (175) at (-7.975, 3.875) {};
		\node [style=none] (176) at (-8.75, 4.225) {};
		\node [style=none] (177) at (-5.525, 3.45) {};
		\node [style=none] (178) at (-3.9, 3.3) {};
		\node [style=none] (179) at (-4.775, 3.875) {};
		\node [style=number-labels] (180) at (-12.325, -2.525) {$3$};
		\node [style=number-labels] (181) at (-11.6, -2.475) {$2$};
		\node [style=number-labels] (182) at (-9.1, -2.4) {$3$};
		\node [style=number-labels] (183) at (-8.375, -2.425) {$2$};
		\node [style=number-labels] (184) at (-5.925, -2.45) {$3$};
		\node [style=number-labels] (185) at (-5.175, -2.5) {$2$};
		\node [style=number-labels] (186) at (-13.475, 2.1) {$4$};
		\node [style=number-labels] (187) at (-12.775, 2.975) {$3$};
		\node [style=number-labels] (188) at (-11.75, 3) {$2$};
		\node [style=number-labels] (189) at (-11.225, 2.525) {$1$};
		\node [style=number-labels] (190) at (-10, 2.775) {$4$};
		\node [style=number-labels] (191) at (-9.3, 3.475) {$3$};
		\node [style=number-labels] (192) at (-8.225, 3.5) {$2$};
		\node [style=number-labels] (193) at (-7.475, 2.8) {$1$};
		\node [style=number-labels] (194) at (-6.025, 2.55) {$4$};
		\node [style=number-labels] (195) at (-5.275, 3.05) {$3$};
		\node [style=number-labels] (196) at (-4.225, 2.925) {$2$};
		\node [style=number-labels] (197) at (-3.65, 2.275) {$1$};
		\node [style=none] (207) at (10.625, 1.85) {};
		\node [style=none] (208) at (13.95, 1.55) {};
		\node [style=none] (209) at (7.5, 1.875) {};
		\node [style=none] (210) at (8.875, 1.875) {};
		\node [style=none] (211) at (2.325, 1.5) {};
		\node [style=none] (212) at (5.825, 1.85) {};
		\node [style=number-labels] (213) at (11.075, 1.575) {$3$};
		\node [style=number-labels] (214) at (11.575, 1.575) {$2$};
		\node [style=dot] (217) at (11.325, 1.875) {};
		\node [style=number-labels] (218) at (7.95, 1.65) {$3$};
		\node [style=number-labels] (219) at (8.425, 1.65) {$2$};
		\node [style=dot] (222) at (8.2, 1.9) {};
		\node [style=number-labels] (223) at (4.9, 1.575) {$3$};
		\node [style=number-labels] (224) at (5.4, 1.575) {$2$};
		\node [style=dot] (227) at (5.175, 1.875) {};
	\end{pgfonlayer}
	\begin{pgfonlayer}{edgelayer}
		\draw [bend left=15, looseness=0.50] (2.center) to (3.center);
		\draw [style=red, bend right=15, looseness=0.75] (0.center) to (10.center);
		\draw [style=red, bend left=15, looseness=0.50] (1.center) to (11.center);
		\draw [style=blue, bend left=90, looseness=0.75] (0.center) to (1.center);
		\draw [style=red, bend left=15] (35.center) to (37.center);
		\draw [style=red, bend left=15] (38.center) to (36.center);
		\draw [style=blue, bend left=75, looseness=0.75] (35.center) to (36.center);
		\draw [style=double-sided] (67.center) to (66.center);
		\draw [style=red, in=120, out=-120, looseness=0.75] (68.center) to (69.center);
		\draw [bend left=15, looseness=0.50] (74) to (75.center);
		\draw [style=red, bend left, looseness=0.75] (100.center) to (101.center);
		\draw [bend right, looseness=0.25] (104.center) to (105.center);
		\draw (108.center) to (109.center);
		\draw [bend left=15, looseness=0.50] (110.center) to (111.center);
		\draw (140) to (143);
		\draw (144) to (141);
		\draw [in=255, out=75, looseness=0.50] (142) to (145);
		\draw (12.center) to (143);
		\draw [in=150, out=-75, looseness=0.75] (4.center) to (140);
		\draw [bend left=15, looseness=0.75] (14.center) to (140);
		\draw [in=135, out=-90, looseness=0.25] (140) to (146.center);
		\draw (141) to (147.center);
		\draw [in=60, out=-105, looseness=0.25] (142) to (148.center);
		\draw (150.center) to (144);
		\draw [in=30, out=-90, looseness=0.50] (151.center) to (141);
		\draw [in=-90, out=135, looseness=0.50] (141) to (149.center);
		\draw [bend right=15, looseness=0.75] (16.center) to (142);
		\draw (145) to (20.center);
		\draw [in=30, out=-105, looseness=0.75] (18.center) to (142);
		\draw [bend left=15, looseness=0.75] (170.center) to (169);
		\draw [bend right=15, looseness=0.50] (169) to (152.center);
		\draw [bend right, looseness=0.50] (168) to (169);
		\draw [bend right=15, looseness=0.50] (168) to (153.center);
		\draw [bend left=15, looseness=0.75] (168) to (161);
		\draw [bend left=15, looseness=0.50] (161) to (154.center);
		\draw [bend right=15, looseness=0.50] (161) to (162);
		\draw [bend right=15, looseness=0.75] (162) to (158.center);
		\draw [bend left=15, looseness=0.75] (162) to (167);
		\draw (167) to (159.center);
		\draw [bend left=15] (167) to (163);
		\draw [bend left=15, looseness=0.50] (163) to (160.center);
		\draw [bend right=15, looseness=0.50] (163) to (164);
		\draw [bend left=15, looseness=0.75] (164) to (166);
		\draw [bend left=15, looseness=0.75] (166) to (165);
		\draw [bend left=15, looseness=0.50] (165) to (156.center);
		\draw [bend left=15, looseness=0.75] (165) to (46.center);
		\draw [bend left=15, looseness=0.50] (166) to (157.center);
		\draw [bend right=15, looseness=0.50] (164) to (155.center);
		\draw [bend left=15, looseness=0.75] (171.center) to (168);
		\draw [bend left=15, looseness=0.50] (168) to (173.center);
		\draw [bend left=15, looseness=0.75] (168) to (172.center);
		\draw [bend left=15, looseness=0.75] (174.center) to (167);
		\draw (167) to (176.center);
		\draw [bend left=15] (167) to (175.center);
		\draw [bend left=15] (177.center) to (166);
		\draw [bend right=15, looseness=0.50] (166) to (179.center);
		\draw [bend left=15] (166) to (178.center);
		\draw [bend left=15, looseness=0.50] (207.center) to (208.center);
		\draw [bend left=15, looseness=0.50] (209.center) to (210.center);
		\draw [bend left=15, looseness=0.50] (211.center) to (212.center);
	\end{pgfonlayer}
\end{tikzpicture}
	\caption[The term cancellations in \ref{term:c2} of Theorem~\ref{prop:cable-2-ainf}.]
	{\textbf{The term cancellations in \ref{term:c2} of Theorem~\ref{prop:cable-2-ainf}.}
	On the left, the first module operation is drawn on the bottom, and the second on
	the top. The lengths of the blue boundaries match, and the two tiling patterns glue together
	to give a centered algebra operation on the right.}
	\label{fig:cable-2-ainf-2}
\end{figure}
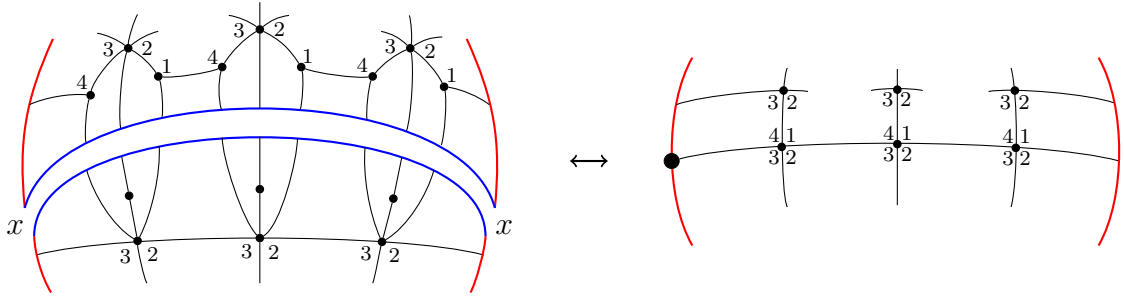

		For example, in the weight $0$ $\Ainf$ relation with inputs $x, \rho_3, \rho_2, \rho_{12}, \rho_1, \rho_4, \rho_{34}$, the following two non-zero terms cancel:
		\[m^0_5(m^0_3(x, \rho_3, \rho_2), \rho_{12}, \rho_1, \rho_4, \rho_{34}))
		= m^0_2(x, \mu^0_6(\rho_3, \rho_2, \rho_{12}, \rho_1, \rho_4, \rho_{34})).\]

	\item \label{term:c3}
	  Consider terms of the form $m_{1+n-j}^{w-w'}\left(m_{1+j}^{w'}(\x,
		a_1,\dots,a_j),a_{j+1},\dots,a_n\right)$, where $a_j a_{j+1}=0$. In this case, we admit the possibility that either of the two module operations are $m_1^0(b) = Uc$.
		Let $\Gamma_1$ and $\Gamma_2$ be the module tiling patterns for the first and second operations respectively.

		The $1$-input module operation results in cases where the two module tiling patterns do not share a common blue boundary. An analysis of the building blocks in Figure~\ref{fig:cable-2-building-blocks} shows that it is terms of the form $m^0_1(m^0_4(x, \rho_3, \rho_2, \rho_1))$ and $m^0_4(m^0_1(b), \rho_4, \rho_3, \rho_2)$ where the two patterns do not share a blue boundary. Here, we choose the simplest $4$-input module operations to illustrate the terms; they might have more inputs and weight.
		
		Let us first consider the case where the two patterns share a blue boundary.
		In this case we were dealing with the supposition that the length of the blue boundaries do not match. Say the length of the blue boundary of $\Gamma_1$ is longer than that of $\Gamma_2$. The other case is entirely analogous like in \ref{term:st3}, and we will omit its description. When we perform the gluing of $\Gamma_1$ with $\Gamma_2$, the red boundary of $\Gamma_2$ juts
		against the blue boundary of $\Gamma_1$, along a face $f'$ of $\Gamma_1$. Consider the penultimate edge $e$ of this
		face, and push $e$ out to the red boundary just before
		the blue-red intersection.

Consider the face $f$ of $\Gamma_1$ that lay on the other side of the edge $e$. If $f$ is an internal face, pushing $e$ out
		to the boundary punctures this face. This produces a module tiling pattern with
		$\rho_i \rho_{i+1} \rho_{i+2} \rho_{i+3}$ as the last term of the chord sequence of the tiling pattern, and our term cancels against a non-zero term $m^{w-1}_{2+n}(\x, a_1, \dots, a_n, \mu_0^1)$ for which we have produced a pattern. As an example, for the weight $1$ $\Ainf$ relation with inputs $x, \rho_{34}, \rho_3, \rho_{23}, \rho_2, \rho_{12}, \rho_{12}, \rho_1, \rho_4, \rho_{34}$, the following two non-zero terms cancel
		\[ m^1_{10}(x, \rho_{34}, \rho_3, \rho_{23}, \rho_2, \rho_{12}, \rho_{12}, \rho_1, \rho_4, \rho_{34}) = m^0_{11}(x, \rho_{34}, \rho_3, \rho_{23}, \rho_2, \rho_{12}, \rho_{12}, \rho_1, \rho_4, \rho_{34}, \mu^1_0). \]
		
		If the face $f$ touches the red boundary, pushing $e$ out disconnects the module tiling pattern, producing a composite pattern $\Gamma_m~\#~ \Gamma_a$. The pattern $\Gamma_m$ is a valid module tiling pattern, while the pattern $\Gamma_a$ has only red vertices. We claim that there is a right-extended algebra operation $\mu^{w_1}_k\left(a_k,
		\dots,a_n\right)$ such that $\Gamma_m$ is a module tiling pattern for the operation $m_{1+n-k}^{w-w_1}\left(\x,a_{1},\dots,a_{n-k},\mu^{w_1}_k\left(a_{n-k+1},
		\dots,a_n\right)\right)$. This is then a non-zero term which cancels against our given term. The construction of the algebra tiling pattern for $\mu^{w_1}_k\left(a_{n-k+1},
		\dots,a_n\right)$ from $\Gamma_a$ follows precisely as in \ref{term:st3}.
		
		As an example, for the weight 0 $\Ainf$ relation with inputs $x, \rho_3, \rho_{234}, \rho_3, \rho_2, \rho_{12}, \rho_{12}, \rho_1, \rho_4, \rho_{34}$, the following two non-zero terms cancel
		\[ m^0_3(x, \rho_3, \mu^0_8(\rho_{234}, \rho_3, \rho_2, \rho_{12}, \rho_{12}, \rho_1, \rho_4, \rho_{34})) = m^0_5((m^0_6(x, \rho_3, \rho_{234}, \rho_3, \rho_2, \rho_{12}), \rho_{12}, \rho_1, \rho_4, \rho_{34}). \]

		The face $f$ could also touch the blue boundary. If on the other side of $e$ lies the blue boundary, then the vertex $v = e \cap e'$ must necessarily be $2$-valent. Pushing $e$ out to the red boundary disconnects $v$ from the rest of the graph and produces a composite pattern $\Gamma_1~\#~\Gamma_2$. $\Gamma_2$ connects exactly one vertex. This vertex is $2$-valent, and $\Gamma_2$ is a module tiling pattern for the module operation $m^0_1(b) = Uc$. Our term then cancels against the term $m^0_1(m^w_{1+n}(\x, a_1, \dots, a_n))$ where the inner operation yields a scalar multiple of $b$. A schema for this is illustrated in Figure~\ref{fig:cable-2-ainf-3}. For completeness, we note that in the case where the blue boundary of the first module operation is longer than the second, this procedure yields a composite pattern $\Gamma_1~\#~\Gamma_2$, where $\Gamma_1$ is a module tiling pattern for the module operation $m^0_1(b) = Uc$. The cancelling term is $m^w_{1+n}(m^0_1(b), a_1, \dots, a_n)$
		
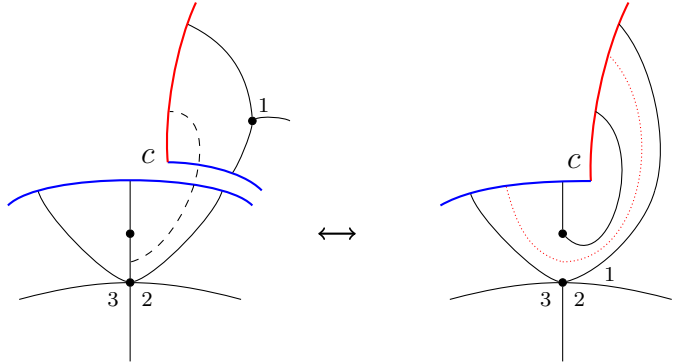
\begin{figure}[h!tbp]
	\centering
	\begin{tikzpicture}
	\begin{pgfonlayer}{nodelayer}
		\node [style=none] (0) at (-8.75, 0.75) {};
		\node [style=none] (1) at (-2.25, 0.75) {};
		\node [style=none] (2) at (-8.45, -1.75) {};
		\node [style=none] (3) at (-2.55, -1.75) {};
		\node [style=none] (4) at (-5.5, 1.425) {};
		\node [style=none] (7) at (-5.025, 2.05) {$c$};
		\node [style=dot] (10) at (-5.5, -1.3) {};
		\node [style=none] (11) at (-5.5, -3.4) {};
		\node [style=dot] (12) at (-5.5, 0) {};
		\node [style=none] (13) at (-7.95, 1.15) {};
		\node [style=none] (14) at (-3.05, 1.15) {};
		\node [style=number-labels] (15) at (-5.95, -1.75) {$3$};
		\node [style=number-labels] (16) at (-5.05, -1.75) {$2$};
		\node [style=none] (17) at (-4.5, 1.9) {};
		\node [style=none] (18) at (-2, 1.15) {};
		\node [style=none] (19) at (-3.75, 6.15) {};
		\node [style=none] (36) at (-0.5, 0) {};
		\node [style=none] (37) at (0.5, 0) {};
		\node [style=dot] (38) at (-2.25, 3) {};
		\node [style=none] (39) at (-2.775, 1.625) {};
		\node [style=none] (40) at (-3.975, 5.575) {};
		\node [style=none] (41) at (-1.25, 3) {};
		\node [style=number-labels] (42) at (-1.95, 3.425) {$1$};
		\node [style=none] (43) at (-4.45, 3.25) {};
		\node [style=none] (44) at (-5.475, -0.75) {};
		\node [style=none] (45) at (2.75, 0.75) {};
		\node [style=none] (49) at (6, 1.4) {};
		\node [style=none] (50) at (6.3, 1.925) {$c$};
		\node [style=dot] (51) at (6, -1.3) {};
		\node [style=none] (52) at (6, -3.4) {};
		\node [style=dot] (53) at (6, 0) {};
		\node [style=none] (54) at (3.55, 1.075) {};
		\node [style=none] (55) at (8.45, 1.15) {};
		\node [style=number-labels] (56) at (5.55, -1.75) {$3$};
		\node [style=number-labels] (57) at (6.45, -1.75) {$2$};
		\node [style=none] (58) at (6.75, 1.4) {};
		\node [style=none] (60) at (7.75, 6.15) {};
		\node [style=none] (66) at (6.875, 3.25) {};
		\node [style=none] (67) at (3, -1.75) {};
		\node [style=none] (68) at (8.9, -1.75) {};
		\node [style=none] (69) at (7.5, 5.575) {};
		\node [style=number-labels] (70) at (7.25, -1.075) {$1$};
		\node [style=none] (71) at (4.5, 1.25) {};
		\node [style=none] (72) at (6, -0.75) {};
		\node [style=none] (73) at (7.225, 4.75) {};
	\end{pgfonlayer}
	\begin{pgfonlayer}{edgelayer}
		\draw [bend left=15] (2.center) to (3.center);
		\draw (10) to (11.center);
		\draw (4.center) to (12);
		\draw (12) to (10);
		\draw [bend right, looseness=0.50] (13.center) to (10);
		\draw [bend right, looseness=0.50] (10) to (14.center);
		\draw [style=blue, bend left=45, looseness=0.50] (0.center) to (1.center);
		\draw [style=blue, in=135, out=0, looseness=0.75] (17.center) to (18.center);
		\draw [style=double-sided] (37.center) to (36.center);
		\draw [bend left=15, looseness=0.50] (38) to (39.center);
		\draw [bend left] (40.center) to (38);
		\draw [bend left, looseness=0.50] (38) to (41.center);
		\draw [style=dashededged, in=0, out=15] (44.center) to (43.center);
		\draw (51) to (52.center);
		\draw (49.center) to (53);
		\draw [bend right, looseness=0.50] (54.center) to (51);
		\draw [in=-105, out=15, looseness=0.75] (51) to (55.center);
		\draw [style=red, bend right=15, looseness=0.75] (60.center) to (58.center);
		\draw [style=blue, in=180, out=30, looseness=0.75] (45.center) to (58.center);
		\draw [in=-30, out=-45, looseness=1.25] (53) to (66.center);
		\draw [bend left=15] (67.center) to (68.center);
		\draw [style=red, bend left=15, looseness=0.75] (17.center) to (19.center);
		\draw [bend right, looseness=0.75] (55.center) to (69.center);
		\draw [style=denselydottedred, in=165, out=-75] (71.center) to (72.center);
		\draw [style=denselydottedred, in=-45, out=0] (72.center) to (73.center);
	\end{pgfonlayer}
\end{tikzpicture}
	\caption[The term cancellations in \ref{term:c3} of Theorem~\ref{prop:cable-2-ainf}.]
	{\textbf{The term cancellations in \ref{term:c3} of Theorem~\ref{prop:cable-2-ainf}.}
	On the left, the first module operation is drawn on the bottom, and the second on
	the top.
	The vertex
	$e \cap e'$ is 2-valent, and pushing out the edge
	$e$ to the red boundary produces a composite pattern $\Gamma_1~\#~\Gamma_2$ where
	$\Gamma_1$ and $\Gamma_2$ are both module operations. Bounded between
	the dotted red arc and the boundary is the module tiling pattern for the operation $m^0_1(b) = Uc$.}
	\label{fig:cable-2-ainf-3}
\end{figure}

		The case where the two module tiling patterns do not share a common blue boundary remains. We had observed that one of the two module operations must correspond to the operation $m^1_0(b) = Uc$; let us say the second. We view the two module operations as instead sharing a common \emph{red} boundary. Glue the two patterns along this common red boundary to obtain a composite pattern. We wish to push \emph{in} the last two edges of this composite pattern and glue them together. This must be followed by a procedure similar to \ref{term:st6} and \ref{term:st7} where we follow the edge and \emph{split} the vertices along it. This decomposes our tiling pattern into two module tiling patterns, $\Gamma_1$ and $\Gamma_2$, and describes an inverse to the operation described in the previous paragraph.
		
		The pairing of the last two paragraphs has the following description in terms of the immersed disks. In the first paragraph, we glue the two disks along the blue boundary. The cut at the obtuse angle is along the $\alpha$-arc. In the second paragraph, we glue the two disks along the red boundary. The cut at the obtuse angle is along the $\beta$-circle.

		\item \label{term:c4} \sloppy Consider terms of the form 
		$m_{1+n-k}^{w-w'}(\x,a_1,\dots,a_j,\mu^{w'}_k(a_{j+1},\dots,
		a_{j+k}),\allowbreak a_{j+k+1},\dots,a_n)$ for a left-extended algebra operation $\mu^{w'}_k$. In this case, the algebra operation is not permitted to be the first input to the module operation.
		We cancel these terms against
		$m_{1+n}^w(\x,a_1,\dots,\mu^0_2(a_j,a_{j+1}),a_{j+k+1}\dots,a_n)$.

		Consider a module tiling pattern $\Gamma_m$ for the module operation and an algebra tiling pattern $\Gamma_a$ for the algebra operation. The pattern $\Gamma_a$ glues into $\Gamma_m$ entirely analogously to \ref{term:st4}, producing a tiling pattern for $m_{1+n}^w(\x,a_1,\dots,\mu^0_2(a_j,a_{j+1}),a_{j+k+1}\dots,a_n)$.

		As an example, in the weight $0$ $\Ainf$ relation with inputs $x, \rho_3, \rho_4, \rho_3, \rho_2, \rho_{123}, \rho_2$, the following non-zero terms cancel
		\[ m^0_4(x, \rho_3, \mu^0_4(\rho_4, \rho_3, \rho_2, \rho_{123}), \rho_2) = m^0_6(x, \mu^0_2(\rho_3, \rho_4), \rho_3, \rho_2, \rho_{123}, \rho_2). \]
		
		\item \label{term:c5} \sloppy Consider terms of the form 
		$m_{1+n-k}^{w-w'}(\x,a_1,\dots,a_j,\mu^{w'}_k(a_{j+1},\dots,
		a_{j+k}),\allowbreak a_{j+k+1},\dots,a_n)$ for a right-extended algebra operation $\mu^{w'}_k$. In this case, the algebra operation is not permitted to be the last input to the module operation.
		We cancel these terms against
		$m_{1+n}^w\left(\x,a_1,\dots,a_{j+k-1},\mu^0_2(a_{j+k},a_{j+k+1})\dots,a_n\right)$.

		We produce a module tiling pattern corresponding to $m_{1+n}^w\left(\x,a_1,\dots,a_{j+k-1},\mu^0_2(a_{j+k},a_{j+k+1})\dots,a_n\right)$ entirely analogously to \ref{term:c4} (and \ref{term:st4}).

		\item \label{term:c6}
		Consider terms of the form $m_{1+n-k}^{w-w'}\left(\x,\mu^{w'}_k\left(a_1,\dots,
		a_k\right),a_{k+1},\dots,a_n\right)$ for a left-extended algebra operation $\mu^{w'}_k$, or terms of the form
		$m_{1+n-k}^{w-w'}\left(\x,a_{1},\dots,a_{n-k},\mu^{w'}_k\left(a_{n-k+1},
		\dots,a_n\right)\right)$ for a right-extended algebra operation $\mu^{w'}_k$. These terms were cancelled against in \ref{term:c3}. We describe the inverse of the cancellation for the right-extended case; the left-extended case follows analogously as in \ref{term:c3}.
		
		As in \ref{term:st3}, the inverse admits a clean description in terms of the immersed disks. Gluing the immersed disks corresponding to the two module tiling operations along their common blue boundary---aligning the \emph{target} of the first with the \emph{source} of the second---yields an immersed disk with an obtuse angle at the intersection of the $\alpha$-arc and $\beta$-circle. This admits a cut along the $\beta$-circle, and when the boundary branch point goes out to the $\alpha$-arc, the disk decomposes into two disks with acute angles. $\Gamma_1$ and $\Gamma_2$ are the module tiling patterns associated to these two disks.
				 
		We can describe this inverse in terms of a module tiling pattern $\Gamma_m$ for the module operation and an algebra tiling pattern $\Gamma_a$ for the algebra operation, as in \ref{term:st3}. We start by gluing in $\Gamma_a$ into $\Gamma_m$, identifying the arc $s$containing the $2$-valent vertices with the edges of the face $f$ that has as labels the output of the algebra operation.
		 After the gluing, we \emph{split} the red vertices on the arc $s$ using the edges in $\Gamma^\beta(\CD_2)$. In other words, we perform a (partial) inverse of the operation depicted in Figure~\ref{fig:cable-2-ainf-2}. This decomposes our tiling pattern into two module tiling patterns, $\Gamma_1$ and $\Gamma_2$.
		
		\item \label{term:c7} Consider terms of the form $m^{w-1}_{2+n}(\x, \mu_0^1, a_1, \dots, a_n)$ or $m^{w-1}_{2+n}(\x, a_1, \dots, a_n, \mu_0^1)$.  These terms were cancelled against in \ref{term:c3}.
			The inverse admits a clean description in terms of the immersed disks that is entirely analogous to \ref{term:st7}. We can further describe it in terms of the tiling patterns, like we did in \ref{term:c6}.

		\item \label{term:c8} Consider of the form $m_{2+n}^{w-1}\left(\x,a_1,\dots,a_i,\mu^1_0,a_{i+1},
		\dots,a_n\right)$. In this case, the $\mu_0^1$ does not appear extremally. These cancel against terms of the form
		$m_{n}^w\left(\x,a_1,\dots,\mu^0_2(a_i,a_{i+1}),\dots,a_n\right)$. This remains precisely the operation that move (3) produces.
		
		\item \label{term:c9} Terms with any centered algebra operation have an idempotent output. As our type A module is strictly unital, for the term to be non-zero the module operation must be a $m_2^0$. These
		were cancelled against in \ref{term:c2}. That is, given a tiling pattern for a centered algebra operation, we \emph{split} the vertices along the arc that follows from the root. This decomposes the tiling pattern into two module tiling patterns $\Gamma_1$ and $\Gamma_2$ that can be performed in order. The inverse is illustrated in
		Figure~\ref{fig:cable-2-ainf-2}. In terms of immersed disks, we cut along an immersed disk for the centered torus algebra along the $\beta$-circle.
		
		\item \label{term:c10} As our last case, consider any term where the algebra operation is a $\mu_2^0$. We claim that all such terms have been cancelled against in \ref{term:c1}, \ref{term:c4}, \ref{term:c5}, and \ref{term:c8}.
		
		Let $\Gamma$ be a module tiling pattern for the module operation. The algebra operation $\mu_0^2$ factorizes a term in its chord sequence. Along the face corresponding to this term, there is an edge $e$ that \emph{exhibits} this factorization. We push out the edge to the red boundary.
		
		Consider the face $f$ on the other side of the edge $e$. 
		If $f$ is an internal face, then pushing out the edge $e$ punctures this face. This is the inverse to the cancellation in \ref{term:c8}.
		Otherwise,
		$f$ touches the boundary.
		If it touches the blue boundary,
		pushing out the edge breaks $\Gamma$ into two module tiling patterns. This is the inverse to the cancellation in \ref{term:c1}.
		Otherwise, $f$ touches the red boundary. If the other side of $e$
		is visible from the red boundary after the edge $e$, pushing out the edge is the inverse to the cancellation in \ref{term:c4}. Otherwise, it is the inverse to the cancellation in \ref{term:c5}. \qedhere 
	\end{enumerate}
\end{proof}

\subsubsection{The \texorpdfstring{$(p,1)$}{(p,1)}-cable}

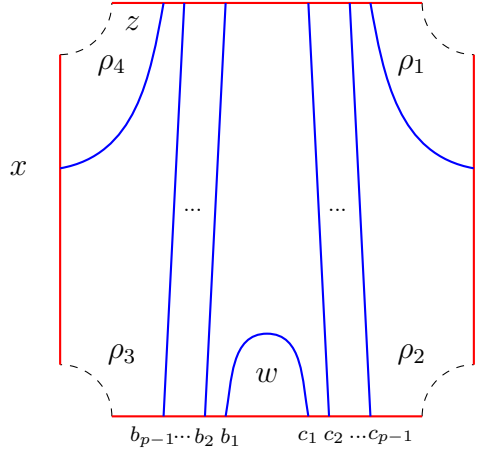
\begin{figure}
	\centering
	\begin{tikzpicture}[scale=2.75]
		\draw[dashed] (.25,0) arc[radius=.25, start angle = 0, end angle = 90];
		\draw[dashed] (0,1.75) arc[radius=.25, start angle = 270, end angle = 360];
		\draw[dashed] (1.75,0) arc[radius=.25, start angle = 180, end angle = 90];
		\draw[dashed] (1.75,2) arc[radius=.25, start angle = 180, end angle = 270];
		\draw[color = red,thick] (.25,0) to (1.75,0);
		\draw[color = red,thick] (0,0.25) to (0,1.75);
		\draw[color = red,thick] (.25,2) to (1.75,2);
		\draw[color = red,thick] (2,0.25) to (2,1.75);
		\node at (.25,1.7) (r4) {$\rho_4$};
		\node at (1.7,1.7) (r1) {$\rho_1$};
		\node at (1.7,.3) (r2) {$\rho_2$};
		\node at (.3,.3) (r3) {$\rho_3$};
		\node at (.35,1.9) (z) {$z$};
		\node at (1,0.2) (w) {$w$};
		\draw[color = blue,thick] (0,1.2) to[in=260,out=10] (.5,2);
		\draw[color = blue,thick] (0.5,0) to (.6,2);
		\draw[color = blue,thick] (0.7,0) to (0.8,2);
		\node at (0.65,1) {\lab{\dots}};
		\draw[color = blue,thick] (2,1.2) to[in=280,out=170] (1.5,2);
		\draw[color = blue,thick] (1.5,0) to (1.4,2);
		\draw[color = blue,thick] (1.3,0) to (1.2,2);
		\node at (1.35,1) {\lab{\dots}};
		
		\draw[color = blue,thick] (.8,0) to[in=180,out=80] (1,0.4) to[in=100,out=0] (1.2,0);
		\node at (-.2,1.2) (x1) {$x$};
		\node at (2.2,1.2) (x2) {\phantom{$x$}};

		\node at (0.45,-.1) (bp1) {\lab{b_{p-1}}};
		\node at (0.6,-.1) (bd) {\lab{\dots}};
		\node at (0.7,-.1) (b2) {\lab{b_{2}}};
		\node at (0.825,-.1) (b1) {\lab{b_{1}}};

		\node at (1.6,-.1) (cp1) {\lab{c_{p-1}}};
		\node at (1.45,-.1) (cd) {\lab{\dots}};
		\node at (1.325,-.1) (c2) {\lab{c_{2}}};
		\node at (1.2,-.1) (c1) {\lab{c_{1}}};
	\end{tikzpicture}
	\caption [A bordered Heegaard diagram \texorpdfstring{$\CD_p$}{CDp} for the $(p,1)$-cable.] {\textbf{A bordered Heegaard diagram $\CD_p$ for the $(p,1)$-cable.}}
	\label{fig:cable-p-bhd}
\end{figure}

A bordered Heegaard diagram $\CD_p$ for the $(p,1)$-cable is pictured in Figure~\ref{fig:cable-p-bhd}. As a left $\FF_2[U,V]\langle \iota_0, \iota_1 \rangle$-module, $\CFAm(\CD_p)$ is freely generated by $x, b_1, \dots, b_{p-1}, c_1, \dots, c_{p-1}$.

We can repeat, word for word, the analysis of Proposition~\ref{prop:c2-obtuse} and Proposition~\ref{prop:c2-building-blocks} to characterize all operations on the type A module $\CFAm(\CD_p)$ in terms of certain simple module tiling patterns. Recall from Remark~\ref{rmk:solid-torus-simplify} that the enumeration of these \emph{basic} building blocks is a finite calculation of tiling patterns that do not have red vertices, and cannot be broken apart with move (1). We do not draw the tiling patterns corresponding to these building blocks, but record the corresponding operations in our graphical format in Figure~\ref{fig:cable-p-bhd}.

\begin{proposition}
	\begin{figure}[h!]
	\centering
	\begin{tikzpicture}[scale=2]
		\node at (0,0) (x) {$x$};
		\node at (2,1) (cpm1) {$c_{p-1}$};
		\node at (3,1) (cpm2) {$c_{p-2}$};
		\node at (4,1) (dotsup1) {$\dots$};
		\node at (5,1) (ci) {$c_i$};
		\node at (6,1) (dotsup2) {$\dots$};
		\node at (7,1) (c1) {$c_1$};
		
		\node at (2,-1) (bpm1) {$b_{p-1}$};
		\node at (3,-1) (bpm2) {$b_{p-2}$};
		\node at (4,-1) (dotsdown1) {$\dots$};
		\node at (5,-1) (bi) {$b_i$};
		\node at (6,-1) (dotsdown2) {$\dots$};
		\node at (7,-1) (b1) {$b_1$};
		
		\draw[->, bend right=15] (x) to node[below]{\lab{\rho_1}} (cpm1);
		\draw[->, bend right=15] (cpm1) to node[above,sloped]{\lab{UV\rho_4\otimes\rho_3\otimes\rho_2}} (x);
		\draw[->, bend left=80] (x) to node[left]{\lab{U^p\rho_3\otimes\rho_2}} (-.5,0) to (x);
		\draw[->] (bpm1) to node[right]{\lab{U^{p-1}}} (cpm1);
		\draw[->, bend right=15] (x) to node[below, sloped]{\lab{U\rho_3\otimes\rho_2\otimes\rho_1}} (bpm1);
		\draw[->, bend right=15] (bpm1) to node[above]{\lab{V\rho_4}} (x);

		\draw[->, bend left=15](cpm1) to node[above]{\lab{\rho_2\otimes \rho_1}} (cpm2);
		\draw[->, bend left=15](cpm2) to node[below]{\lab{UV\rho_4 \otimes \rho_3}} (cpm1);
		\draw[->, bend right=15](bpm1) to node[below]{\lab{U\rho_2\otimes \rho_1}} (bpm2);
		\draw[->, bend right=15](bpm2) to node[above]{\lab{V\rho_4 \otimes \rho_3}} (bpm1);

		\draw[->] (bpm2) to node[right]{\lab{U^{p-2}}} (cpm2);

		\draw[->, bend left=15](cpm2) to node[above]{\lab{\rho_2\otimes \rho_1}} (dotsup1);
		\draw[->, bend left=15](dotsup1) to node[below]{\lab{UV\rho_4 \otimes \rho_3}} (cpm2);
		\draw[->, bend right=15](bpm2) to node[below]{\lab{U\rho_2\otimes \rho_1}} (dotsdown1);
		\draw[->, bend right=15](dotsdown1) to node[above]{\lab{V\rho_4 \otimes \rho_3}} (bpm2);

		\draw[->] (bi) to node[right]{\lab{U^{i}}} (ci);

		\draw[->, bend left=15](dotsup1) to node[above]{\lab{\rho_2\otimes \rho_1}} (ci);
		\draw[->, bend left=15](ci) to node[below]{\lab{UV\rho_4 \otimes \rho_3}} (dotsup1);
		\draw[->, bend right=15](dotsdown1) to node[below]{\lab{U\rho_2\otimes \rho_1}} (bi);
		\draw[->, bend right=15](bi) to node[above]{\lab{V\rho_4 \otimes \rho_3}} (dotsdown1);

		\draw[->, bend left=15](ci) to node[above]{\lab{\rho_2\otimes \rho_1}} (dotsup2);
		\draw[->, bend left=15](dotsup2) to node[below]{\lab{UV\rho_4 \otimes \rho_3}} (ci);
		\draw[->, bend right=15](bi) to node[below]{\lab{U\rho_2\otimes \rho_1}} (dotsdown2);
		\draw[->, bend right=15](dotsdown2) to node[above]{\lab{V\rho_4 \otimes \rho_3}} (bi);

		\draw[->, bend left=15](dotsup2) to node[above]{\lab{\rho_2\otimes \rho_1}} (c1);
		\draw[->, bend left=15](c1) to node[below]{\lab{UV\rho_4 \otimes \rho_3}} (dotsup2);
		\draw[->, bend right=15](dotsdown2) to node[below]{\lab{U\rho_2\otimes \rho_1}} (b1);
		\draw[->, bend right=15](b1) to node[above]{\lab{V\rho_4 \otimes \rho_3}} (dotsdown2);

		\draw[->, bend right=15] (b1) to node[left]{\lab{U}} (c1);
		\draw[->, bend right=15] (c1) to node[above,sloped]{\lab{V\rho_2\otimes\rho_1\otimes\rho_4\otimes\rho_3}} (b1);
	\end{tikzpicture}
		\caption [A representation of the type A module \texorpdfstring{$\CFAm(\CD_p)$}{CFA minus of CDp}.] {\textbf{A representation of the type A module $\CFAm(\CD_p)$.}}
	\label{fig:cable-p-graph}
	\end{figure}
	Any operation on the type A module $\CFAm(\CD_p)$ can be obtained as the operation constructed by starting with the operations recorded in the graph in Figure~\ref{fig:cable-p-graph}, and performing a finite sequence of the moves catalogued in Section~\ref{sec:three-moves}.
\end{proposition}

A calculation of the (non-weighted) type A module $\CFAa(\CD_p)$ was performed by Petkova \cite{Petkova2013}. In \cite{O17}, Ozsv\'ath, Stipsicz, and Szab\'o record this calculation with a similar graphical presentation. We remark that it agrees precisely with the subgraph of Figure~\ref{fig:cable-p-graph} that does not feature any of the edges with a  power of $V$.  

\begin{theorem}
	\label{thm:cable-p-ainf}
	The type A module $\CFAm(\CD_p)$, with operations $\{m^w_{1+n}\}$ given by counting its module tiling patterns of chord sequence length $n$ and weight $w$, satisfies the $\Ainf$ structure relations.
\end{theorem}

\begin{proof}
	We verify the $\Ainf$ relation for a fixed input $(\x, a_1, \dots, a_n)$ and weight $w$. This proof follows the same schema as that of Theorem~\ref{prop:solid-torus-ainf} and Theorem~\ref{prop:cable-2-ainf}. In particular, we show the $\Ainf$ relation holds by pairing non-zero terms that appear in it.

	It suffices to consider the case where all the $a_i$ are Reeb elements. The analysis in Theorem~\ref{prop:cable-2-ainf} of the terms considered in \ref{term:c1},\ref{term:c2}, \ref{term:c4}, \ref{term:c5}, \ref{term:c6}, \ref{term:c7}, \ref{term:c8}, \ref{term:c9}, \ref{term:c10}---that is, all terms but \ref{term:c3}---extends precisely here; just as the analysis extended from $\CFAm(\HD_{st})$ to $\CFAm(\CD_2)$.
	
	It remains to extend the cancellations considered in \ref{term:c3} for the $(p,1)$-cable.

	\begin{enumerate}[wide,label=(TC$'$-\arabic*)]
		\setcounter{enumi}{2}
		\item \label{term:cp3} Recall that in this case we consider terms of the form $m^{w-w'}_{1+n-j}(m^{w'}_{1+j}(\x, a_1, \dots, a_j), a_{j+1}, \dots, a_n)$, where $a_j a_{j+1} = 0$. We admit the possibility either of the two module operations are $m_1^0(b_i) = U^i c_i$. Let $\Gamma_1$ and $\Gamma_2$ be the module tiling patterns for the first and second operations respectively.
		
		As in \ref{term:c3}, the $1$-input module operations results in cases where the two module tiling patterns do not share a common blue boundary; these involve the operation $m_1^0(b_i) = U^i c_i$. Let us first consider the case where the two patterns share a blue boundary. We deal with the supposition that the length of the blue boundaries do not match, as the case we do is \ref{term:c2}. Say the length of the blue boundary of $\Gamma_1$ is longer than that of $\Gamma_2$. When we perform the gluing of $\Gamma_1$ with $\Gamma_2$, the red boundary of $\Gamma_2$ juts against the blue boundary of $\Gamma_1$, along a face $f'$ of $\Gamma_1$. Consider the penultimate edge $e$ of this face, and push $e$ out to the red boundary just before the blue-red intersection.
		
		Consider the face $f$ of $\Gamma_1$ that lay on the other side of this edge $e$. If $f$ is an internal face, pushing out to the boundary punctures this face. If the face $f$ touches the red boundary, pushing $e$ out disconnects the module tiling pattern, producing a composite pattern $\Gamma_m~\#~\Gamma_a$. In the former case, our term cancels against a non-zero term $m^{w-1}_{2+n}(\x, a_1, \dots, a_n), \mu^1_0)$; in the latter, against a term $m^{w-w_1}_{1+n-k}(\x, a_1, \dots, a_{n-k}, \mu^{w_1}_k(a_{n-k+1}, \dots, a_n))$. Both these cases are analogous to \ref{term:c3}.

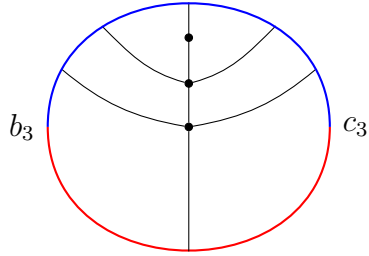
\begin{figure}[h!tbp]
	\centering
	\begin{tikzpicture}[scale=1.25]
	\begin{pgfonlayer}{nodelayer}
		\node [style=none] (0) at (-3, 0) {};
		\node [style=none] (1) at (3, 0) {};
		\node [style=none] (2) at (-2.95, -0.5) {};
		\node [style=none] (3) at (2.95, -0.5) {};
		\node [style=none] (4) at (0, 2.65) {};
		\node [style=none] (5) at (0, -2.65) {};
		\node [style=none] (6) at (-3.55, 0) {$b_3$};
		\node [style=none] (7) at (3.55, 0) {$c_3$};
		\node [style=dot] (8) at (0, 0) {};
		\node [style=dot] (9) at (0, 0.925) {};
		\node [style=dot] (10) at (0, 1.9) {};
		\node [style=none] (11) at (2.7, 1.225) {};
		\node [style=none] (12) at (-2.7, 1.225) {};
		\node [style=none] (13) at (1.85, 2.15) {};
		\node [style=none] (14) at (-1.85, 2.15) {};
	\end{pgfonlayer}
	\begin{pgfonlayer}{edgelayer}
		\draw [style=red, bend right=90, looseness=1.50] (0.center) to (1.center);
		\draw (4.center) to (8);
		\draw (8) to (5.center);
		\draw [bend right=15] (8) to (11.center);
		\draw [in=225, out=15] (9) to (13.center);
		\draw [in=315, out=165] (9) to (14.center);
		\draw [bend left=15] (8) to (12.center);
		\draw [style=blue, bend left=90, looseness=1.50] (0.center) to (1.center);
	\end{pgfonlayer}
\end{tikzpicture}
	\caption[A module tiling pattern for \texorpdfstring{$\CD_4$}{CD4}.]
	{\textbf{A module tiling pattern for $\CD_4$.} The operation is $m_1^0(b_3) = U^3 c_3$. The corresponding immersed disk is drawn on the left in solid purple in Figure~\ref{fig:cable-4-disks}.}
	\label{fig:cable-4-pattern}
\end{figure}

		The face $f$ could also touch the boundary. In \ref{term:c3}, there was precisely one such edge on either side of which lay the blue boundary; we located it in the module tiling pattern for $m^0_1(b) = Uc$. For a $(p,1)$-cable, the module tiling patterns for all operations $m^0_1(b_i) = U^i c_i$ feature such edges. For example, Figure~\ref{fig:cable-4-pattern} has two such edges in a module tiling pattern for the operation $m_1^0(b_3) = U^3 c_3$ for $\CFAm(\CD_4)$. Pushing such an $e$ out produces a composite pattern $\Gamma_1~\#~\Gamma_2$, where $\Gamma_2$ is a module tiling pattern for a module operation $m^0_1(b_i) = U^i c_i$ for some $i$. Our term then cancels against the term $m^0_1(m^w_{1+n}(\x, a_1, \dots, a_n))$ where the inner operation yields a scalar multiple of $b_i$.

		The case where the two module tiling patterns share a common red boundary were cancelled exactly against in the previous paragraph. In terms of the immersed disks corresponding to $\Gamma_1$ and $\Gamma_2$, the inverse is described by gluing the disks along the red boundary and cutting along the $\beta$-circle. This is illustrated in Figure~\ref{fig:cable-4-disks}. In the weight $0$ $\Ainf$ relation for $\CFAm(\CD_4)$ with inputs $b_3, \rho_{23}, \rho_2, \rho_1, \rho_{41}$, the following two non-zero terms are illustrated to cancel in Figure~\ref{fig:cable-4-disks}.
		\[ m^0_5(m^0_1(b_3), \rho_{23}, \rho_2, \rho_1, \rho_{41}) = 
		m^0_1(m^0_5(b_3, \rho_{23}, \rho_2, \rho_1, \rho_{41})). \]
		The disks for the term on the left and right are illustrated on the left and right respectively. It is the edge between the two $4$-valent vertices in Figure~\ref{fig:cable-4-pattern} that is pushed out in the corresponding cancellation for the module tiling patterns. \qedhere

\begin{figure}[h!]
	\centering
	\begin{tikzpicture}[scale=1.25]
		\begin{scope}[shift={(0,4)}]
		\draw[dashed] (.25,0) arc[radius=.25, start angle = 0, end angle = 90];
		\draw[dashed] (0,1.75) arc[radius=.25, start angle = 270, end angle = 360];
		\draw[dashed] (1.75,0) arc[radius=.25, start angle = 180, end angle = 90];
		\draw[dashed] (1.75,2) arc[radius=.25, start angle = 180, end angle = 270];
		\draw[color = red,thick] (.25,0) to (1.75,0);
		\draw[color = red,thick] (0,0.25) to (0,1.75);
		\draw[color = red,thick] (.25,2) to (1.75,2);
		\draw[color = red,thick] (2,0.25) to (2,1.75);

		\draw[color = blue,thick] (0,1.2) to[in=260,out=10] (.5,2);
		\draw[color = blue,thick] (0.5,0) to (.65,2);
		\draw[color = blue,thick] (0.65,0) to (0.8,2);
		
		\draw[color = blue,thick] (2,1.2) to[in=280,out=170] (1.5,2);
		\draw[color = blue,thick] (1.5,0) to (1.35,2);
		\draw[color = blue,thick] (1.35,0) to (1.2,2);
		
		\draw[color = blue,thick] (.8,0) to[in=180,out=80] (1,0.4) to[in=100,out=0] (1.2,0);
		\fill[fill=purple,opacity=0.2] (0.8,0) to[in=180,out=80] (1,0.4) to[in=100,out=0] (1.2,0) to cycle;
		\end{scope}
		\begin{scope}[shift={(0,2)}]
		\draw[dashed] (.25,0) arc[radius=.25, start angle = 0, end angle = 90];
		\draw[dashed] (0,1.75) arc[radius=.25, start angle = 270, end angle = 360];
		\draw[dashed] (1.75,0) arc[radius=.25, start angle = 180, end angle = 90];
		\draw[dashed] (1.75,2) arc[radius=.25, start angle = 180, end angle = 270];
		\draw[color = red,thick] (.25,0) to (1.75,0);
		\draw[color = red,thick] (0,0.25) to (0,1.75);
		\draw[color = red,thick] (.25,2) to (1.75,2);
		\draw[color = red,thick] (2,0.25) to (2,1.75);
		\draw[color = blue,thick] (0,1.2) to[in=260,out=10] (.5,2);
		\draw[color = blue,thick] (0.5,0) to (.65,2);
		\draw[color = blue,thick] (0.65,0) to (0.8,2);
		
		\draw[color = blue,thick] (2,1.2) to[in=280,out=170] (1.5,2);
		\draw[color = blue,thick] (1.5,0) to (1.35,2);
		\draw[color = blue,thick] (1.35,0) to (1.2,2);
		
		\draw[color = blue,thick] (.8,0) to[in=180,out=80] (1,0.4) to[in=100,out=0] (1.2,0);
		\fill[fill=purple,opacity=0.2] (0.65,0) to (0.8, 2) to (1.2,2) to (1.35,0) to (0.65,0);
		\end{scope}
		\begin{scope}[shift={(0,0)}]
		\draw[dashed] (.25,0) arc[radius=.25, start angle = 0, end angle = 90];
		\draw[dashed] (0,1.75) arc[radius=.25, start angle = 270, end angle = 360];
		\draw[dashed] (1.75,0) arc[radius=.25, start angle = 180, end angle = 90];
		\draw[dashed] (1.75,2) arc[radius=.25, start angle = 180, end angle = 270];
		\draw[color = red,thick] (.25,0) to (1.75,0);
		\draw[color = red,thick] (0,0.25) to (0,1.75);
		\draw[color = red,thick] (.25,2) to (1.75,2);
		\draw[color = red,thick] (2,0.25) to (2,1.75);
		\draw[color = blue,thick] (0,1.2) to[in=260,out=10] (.5,2);
		\draw[color = blue,thick] (0.5,0) to (.65,2);
		\draw[color = blue,thick] (0.65,0) to (0.8,2);
		
		\draw[color = blue,thick] (2,1.2) to[in=280,out=170] (1.5,2);
		\draw[color = blue,thick] (1.5,0) to (1.35,2);
		\draw[color = blue,thick] (1.35,0) to (1.2,2);
		
		\draw[color = blue,thick] (.8,0) to[in=180,out=80] (1,0.4) to[in=100,out=0] (1.2,0);
		
		\fill[fill=purple,opacity=0.2] (0.5,0) to (0.65, 2) to (1.35,2) to (1.5,0) to (0.5,0);
		\fill[pattern=north west lines, pattern color=green,opacity=0.8] (1.5, 0) to (1.75, 0) arc[radius=.25, start angle = 180, end angle = 90] to (2, 1.75) arc[radius=.25, start angle = 270, end angle=180] to (1.35, 2);
		\end{scope}
		\begin{scope}[shift={(2,0)}]
		\draw[dashed] (.25,0) arc[radius=.25, start angle = 0, end angle = 90];
		\draw[dashed] (0,1.75) arc[radius=.25, start angle = 270, end angle = 360];
		\draw[dashed] (1.75,0) arc[radius=.25, start angle = 180, end angle = 90];
		\draw[dashed] (1.75,2) arc[radius=.25, start angle = 180, end angle = 270];
		\draw[color = red,thick] (.25,0) to (1.75,0);
		\draw[color = red,thick] (0,0.25) to (0,1.75);
		\draw[color = red,thick] (.25,2) to (1.75,2);
		\draw[color = red,thick] (2,0.25) to (2,1.75);
		\draw[color = blue,thick] (0,1.2) to[in=260,out=10] (.5,2);
		\draw[color = blue,thick] (0.5,0) to (.65,2);
		\draw[color = blue,thick] (0.65,0) to (0.8,2);
		
		\draw[color = blue,thick] (2,1.2) to[in=280,out=170] (1.5,2);
		\draw[color = blue,thick] (1.5,0) to (1.35,2);
		\draw[color = blue,thick] (1.35,0) to (1.2,2);
		
		\draw[color = blue,thick] (.8,0) to[in=180,out=80] (1,0.4) to[in=100,out=0] (1.2,0);

		\fill[pattern=north west lines, pattern color=green,opacity=0.8] (0.25,0) arc[radius=.25, start angle = 0, end angle = 90] to (0, 1.75) arc[radius=.25, start angle = 270, end angle = 360] to (1.75, 2) arc[radius=.25, start angle = 180, end angle = 270] to (2, 0.25) 
		arc[radius=.25, start angle = 90, end angle = 180] to (0.25, 0);
		\end{scope}

		\begin{scope}[shift={(5,4)}]
		\draw[dashed] (.25,0) arc[radius=.25, start angle = 0, end angle = 90];
		\draw[dashed] (0,1.75) arc[radius=.25, start angle = 270, end angle = 360];
		\draw[dashed] (1.75,0) arc[radius=.25, start angle = 180, end angle = 90];
		\draw[dashed] (1.75,2) arc[radius=.25, start angle = 180, end angle = 270];
		\draw[color = red,thick] (.25,0) to (1.75,0);
		\draw[color = red,thick] (0,0.25) to (0,1.75);
		\draw[color = red,thick] (.25,2) to (1.75,2);
		\draw[color = red,thick] (2,0.25) to (2,1.75);

		\draw[color = blue,thick] (0,1.2) to[in=260,out=10] (.5,2);
		\draw[color = blue,thick] (0.5,0) to (.65,2);
		\draw[color = blue,thick] (0.65,0) to (0.8,2);
		
		\draw[color = blue,thick] (2,1.2) to[in=280,out=170] (1.5,2);
		\draw[color = blue,thick] (1.5,0) to (1.35,2);
		\draw[color = blue,thick] (1.35,0) to (1.2,2);
		
		\draw[color = blue,thick] (.8,0) to[in=180,out=80] (1,0.4) to[in=100,out=0] (1.2,0);
		\fill[pattern=north west lines, pattern color=green,opacity=0.8] (0.8,0) to[in=180,out=80] (1,0.4) to[in=100,out=0] (1.2,0) to cycle;
		\end{scope}
		\begin{scope}[shift={(5,2)}]
		\draw[dashed] (.25,0) arc[radius=.25, start angle = 0, end angle = 90];
		\draw[dashed] (0,1.75) arc[radius=.25, start angle = 270, end angle = 360];
		\draw[dashed] (1.75,0) arc[radius=.25, start angle = 180, end angle = 90];
		\draw[dashed] (1.75,2) arc[radius=.25, start angle = 180, end angle = 270];
		\draw[color = red,thick] (.25,0) to (1.75,0);
		\draw[color = red,thick] (0,0.25) to (0,1.75);
		\draw[color = red,thick] (.25,2) to (1.75,2);
		\draw[color = red,thick] (2,0.25) to (2,1.75);
		\draw[color = blue,thick] (0,1.2) to[in=260,out=10] (.5,2);
		\draw[color = blue,thick] (0.5,0) to (.65,2);
		\draw[color = blue,thick] (0.65,0) to (0.8,2);
		
		\draw[color = blue,thick] (2,1.2) to[in=280,out=170] (1.5,2);
		\draw[color = blue,thick] (1.5,0) to (1.35,2);
		\draw[color = blue,thick] (1.35,0) to (1.2,2);
		
		\draw[color = blue,thick] (.8,0) to[in=180,out=80] (1,0.4) to[in=100,out=0] (1.2,0);
		\fill[pattern=north west lines, pattern color=green,opacity=0.8] (0.65,0) to (0.8, 2) to (1.2,2) to (1.35,0) to (0.65,0);
		\end{scope}
		\begin{scope}[shift={(5,0)}]
		\draw[dashed] (.25,0) arc[radius=.25, start angle = 0, end angle = 90];
		\draw[dashed] (0,1.75) arc[radius=.25, start angle = 270, end angle = 360];
		\draw[dashed] (1.75,0) arc[radius=.25, start angle = 180, end angle = 90];
		\draw[dashed] (1.75,2) arc[radius=.25, start angle = 180, end angle = 270];
		\draw[color = red,thick] (.25,0) to (1.75,0);
		\draw[color = red,thick] (0,0.25) to (0,1.75);
		\draw[color = red,thick] (.25,2) to (1.75,2);
		\draw[color = red,thick] (2,0.25) to (2,1.75);
		\draw[color = blue,thick] (0,1.2) to[in=260,out=10] (.5,2);
		\draw[color = blue,thick] (0.5,0) to (.65,2);
		\draw[color = blue,thick] (0.65,0) to (0.8,2);
		
		\draw[color = blue,thick] (2,1.2) to[in=280,out=170] (1.5,2);
		\draw[color = blue,thick] (1.5,0) to (1.35,2);
		\draw[color = blue,thick] (1.35,0) to (1.2,2);
		
		\draw[color = blue,thick] (.8,0) to[in=180,out=80] (1,0.4) to[in=100,out=0] (1.2,0);
		
		\fill[fill=purple,opacity=0.2] (0.5,0) to (0.65, 2) to (1.35,2) to (1.5,0) to (0.5,0);
		\fill[fill=purple,opacity=0.2] (1.5, 0) to (1.75, 0) arc[radius=.25, start angle = 180, end angle = 90] to (2, 1.75) arc[radius=.25, start angle = 270, end angle=180] to (1.35, 2);
		\end{scope}
		\begin{scope}[shift={(7,0)}]
		\draw[dashed] (.25,0) arc[radius=.25, start angle = 0, end angle = 90];
		\draw[dashed] (0,1.75) arc[radius=.25, start angle = 270, end angle = 360];
		\draw[dashed] (1.75,0) arc[radius=.25, start angle = 180, end angle = 90];
		\draw[dashed] (1.75,2) arc[radius=.25, start angle = 180, end angle = 270];
		\draw[color = red,thick] (.25,0) to (1.75,0);
		\draw[color = red,thick] (0,0.25) to (0,1.75);
		\draw[color = red,thick] (.25,2) to (1.75,2);
		\draw[color = red,thick] (2,0.25) to (2,1.75);
		\draw[color = blue,thick] (0,1.2) to[in=260,out=10] (.5,2);
		\draw[color = blue,thick] (0.5,0) to (.65,2);
		\draw[color = blue,thick] (0.65,0) to (0.8,2);
		
		\draw[color = blue,thick] (2,1.2) to[in=280,out=170] (1.5,2);
		\draw[color = blue,thick] (1.5,0) to (1.35,2);
		\draw[color = blue,thick] (1.35,0) to (1.2,2);
		
		\draw[color = blue,thick] (.8,0) to[in=180,out=80] (1,0.4) to[in=100,out=0] (1.2,0);

		\fill[fill=purple, opacity=0.2] (0.25,0) arc[radius=.25, start angle = 0, end angle = 90] to (0, 1.75) arc[radius=.25, start angle = 270, end angle = 360] to (1.75, 2) arc[radius=.25, start angle = 180, end angle = 270] to (2, 0.25) 
		arc[radius=.25, start angle = 90, end angle = 180] to (0.25, 0);
		\end{scope}
	\end{tikzpicture}
	\caption [Immersed disks representing a cancellation in \ref{term:cp3} of Theorem~\ref{thm:cable-p-ainf}.] {\textbf{Immersed disks representing a cancellation in \ref{term:cp3} of Theorem~\ref{thm:cable-p-ainf}}. In solid purple is an immersed disk for $\Gamma_1$ and in hatched green is an immersed disk for $\Gamma_2$, in the notation of \ref{term:cp3}. The composition of operations on left cancels against the composition of operations on the right.}
	\label{fig:cable-4-disks}
\end{figure}
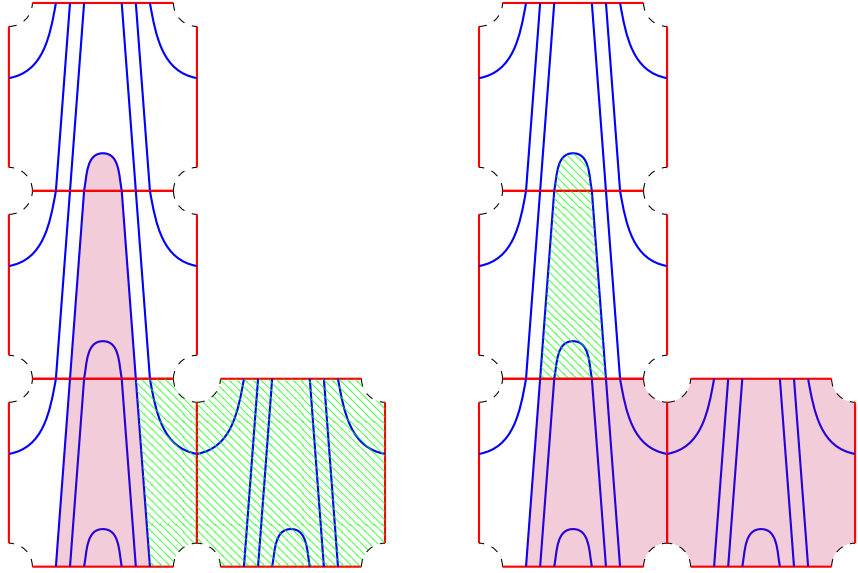

	\end{enumerate}
\end{proof}

\subsection{Other \texorpdfstring{$(1,1)$}{(1,1)}-patterns}
\label{sec:11-patterns}
\begin{figure}
	\centering
		\begin{tikzpicture}[scale=1.25]
			\begin{scope}[shift={(0,0)}]
				\draw[dashed] (.25,0) arc[radius=.25, start angle = 0, end angle = 90];
				\draw[dashed] (0,1.75) arc[radius=.25, start angle = 270, end angle = 360];
				\draw[dashed] (1.75,0) arc[radius=.25, start angle = 180, end angle = 90];
				\draw[dashed] (1.75,2) arc[radius=.25, start angle = 180, end angle = 270];
				\draw[color = red,thick] (.25,0) to (1.75,0);
				\draw[color = red,thick] (0,0.25) to (0,1.75);
				\draw[color = red,thick] (.25,2) to (1.75,2);
				\draw[color = red,thick] (2,0.25) to (2,1.75);
		
				\draw[color = blue,thick] (0,1.375) to[in=270,out=0] (0.25+0.3,2);
				\draw[color = blue,thick] (0,1) to[in=270,out=0] (0.25+2*0.3,2);
				\draw[color = blue,thick] (0,0.625) to[in=90,out=0] (0.25+3*0.3,0);
				
				\draw[color = blue, thick] (0.25+0.3,0) to [in=180,out=90] (0.25+0.45, 0.25) to[out=0,in=90] (0.25+2*0.3,0);

				\draw[color = blue,thick] (0.25+3*0.3,2) to[in=180,out=270] (2,1);
				\draw[color = blue,thick] (0.25+4*0.3,2) to[in=180,out=270] (2,1.375);
				\draw[color = blue,thick] (0.25+4*0.3,0) to[in=180,out=90] (2,0.625);

				\fill[fill=purple,opacity=0.2] (0.25,0) arc[radius=.25, start angle = 0, end angle = 90] to (0, 1) to[in=270,out=0] (0.25 + 2 * 0.3, 2) to (1.75,2) arc[radius=.25, start angle = 180, end angle = 270] to (2, 0.25) 
				arc[radius=.25, start angle = 90, end angle = 180] -- cycle;
				\fill[fill=purple,opacity=0.2] (0,1) to[in=270,out=0] (0.25+2*0.3,2) to (0.25+0.3,2) to[in=0,out=270] (0,1.375) -- cycle;
			\end{scope}
			\begin{scope}[shift={(0,2)}]
				\draw[dashed] (.25,0) arc[radius=.25, start angle = 0, end angle = 90];
				\draw[dashed] (0,1.75) arc[radius=.25, start angle = 270, end angle = 360];
				\draw[dashed] (1.75,0) arc[radius=.25, start angle = 180, end angle = 90];
				\draw[dashed] (1.75,2) arc[radius=.25, start angle = 180, end angle = 270];
				\draw[color = red,thick] (.25,0) to (1.75,0);
				\draw[color = red,thick] (0,0.25) to (0,1.75);
				\draw[color = red,thick] (.25,2) to (1.75,2);
				\draw[color = red,thick] (2,0.25) to (2,1.75);
		
				\draw[color = blue,thick] (0,1.375) to[in=270,out=0] (0.25+0.3,2);
				\draw[color = blue,thick] (0,1) to[in=270,out=0] (0.25+2*0.3,2);
				\draw[color = blue,thick] (0,0.625) to[in=90,out=0] (0.25+3*0.3,0);
				
				\draw[color = blue, thick] (0.25+0.3,0) to [in=180,out=90] (0.25+0.45, 0.25) to[out=0,in=90] (0.25+2*0.3,0);
				\fill[pattern=north west lines, pattern color=green,opacity=0.8] (0.25+0.3,0) to [in=180,out=90] (0.25+0.45, 0.25) to[out=0,in=90] (0.25+2*0.3,0) -- cycle;
				
				\draw[color = blue,thick] (0.25+3*0.3,2) to[in=180,out=270] (2,1);
				\draw[color = blue,thick] (0.25+4*0.3,2) to[in=180,out=270] (2,1.375);
				\draw[color = blue,thick] (0.25+4*0.3,0) to[in=180,out=90] (2,0.625);
			\end{scope}

			\begin{scope}[shift={(4,0)}]
				\draw[dashed] (.25,0) arc[radius=.25, start angle = 0, end angle = 90];
				\draw[dashed] (0,1.75) arc[radius=.25, start angle = 270, end angle = 360];
				\draw[dashed] (1.75,0) arc[radius=.25, start angle = 180, end angle = 90];
				\draw[dashed] (1.75,2) arc[radius=.25, start angle = 180, end angle = 270];
				\draw[color = red,thick] (.25,0) to (1.75,0);
				\draw[color = red,thick] (0,0.25) to (0,1.75);
				\draw[color = red,thick] (.25,2) to (1.75,2);
				\draw[color = red,thick] (2,0.25) to (2,1.75);
		
				\draw[color = blue,thick] (0,1.375) to[in=270,out=0] (0.25+0.3,2);
				\draw[color = blue,thick] (0,1) to[in=270,out=0] (0.25+2*0.3,2);
				\draw[color = blue,thick] (0,0.625) to[in=90,out=0] (0.25+3*0.3,0);
				
				\draw[color = blue, thick] (0.25+0.3,0) to [in=180,out=90] (0.25+0.45, 0.25) to[out=0,in=90] (0.25+2*0.3,0);

				\draw[color = blue,thick] (0.25+3*0.3,2) to[in=180,out=270] (2,1);
				\draw[color = blue,thick] (0.25+4*0.3,2) to[in=180,out=270] (2,1.375);
				\draw[color = blue,thick] (0.25+4*0.3,0) to[in=180,out=90] (2,0.625);

				\fill[fill=purple,opacity=0.2] (0.25+4*0.3,2) to (1.75,2) arc[radius=.25, start angle = 180, end angle = 270] to (2,1.375) to[in=270,out=180] (0.25+4*0.3,2);
			\end{scope}
			\begin{scope}[shift={(4,2)}]
				\draw[dashed] (.25,0) arc[radius=.25, start angle = 0, end angle = 90];
				\draw[dashed] (0,1.75) arc[radius=.25, start angle = 270, end angle = 360];
				\draw[dashed] (1.75,0) arc[radius=.25, start angle = 180, end angle = 90];
				\draw[dashed] (1.75,2) arc[radius=.25, start angle = 180, end angle = 270];
				\draw[color = red,thick] (.25,0) to (1.75,0);
				\draw[color = red,thick] (0,0.25) to (0,1.75);
				\draw[color = red,thick] (.25,2) to (1.75,2);
				\draw[color = red,thick] (2,0.25) to (2,1.75);
		
				\draw[color = blue,thick] (0,1.375) to[in=270,out=0] (0.25+0.3,2);
				\draw[color = blue,thick] (0,1) to[in=270,out=0] (0.25+2*0.3,2);
				\draw[color = blue,thick] (0,0.625) to[in=90,out=0] (0.25+3*0.3,0);
				
				\draw[color = blue, thick] (0.25+0.3,0) to [in=180,out=90] (0.25+0.45, 0.25) to[out=0,in=90] (0.25+2*0.3,0);
				\fill[fill=purple,opacity=0.2] (0.25+4*0.3,0) to (1.75,0) arc[radius=.25, start angle = 180, end angle = 90] to (2,0.625) to[in=90,out=180] (0.25+4*0.3,0);
				
				\draw[color = blue,thick] (0.25+3*0.3,2) to[in=180,out=270] (2,1);
				\draw[color = blue,thick] (0.25+4*0.3,2) to[in=180,out=270] (2,1.375);
				\draw[color = blue,thick] (0.25+4*0.3,0) to[in=180,out=90] (2,0.625);
			\end{scope}
			\begin{scope}[shift={(6,0)}]
				\draw[dashed] (.25,0) arc[radius=.25, start angle = 0, end angle = 90];
				\draw[dashed] (0,1.75) arc[radius=.25, start angle = 270, end angle = 360];
				\draw[dashed] (1.75,0) arc[radius=.25, start angle = 180, end angle = 90];
				\draw[dashed] (1.75,2) arc[radius=.25, start angle = 180, end angle = 270];
				\draw[color = red,thick] (.25,0) to (1.75,0);
				\draw[color = red,thick] (0,0.25) to (0,1.75);
				\draw[color = red,thick] (.25,2) to (1.75,2);
				\draw[color = red,thick] (2,0.25) to (2,1.75);
		
				\draw[color = blue,thick] (0,1.375) to[in=270,out=0] (0.25+0.3,2);
				\draw[color = blue,thick] (0,1) to[in=270,out=0] (0.25+2*0.3,2);
				\draw[color = blue,thick] (0,0.625) to[in=90,out=0] (0.25+3*0.3,0);
				
				\draw[color = blue, thick] (0.25+0.3,0) to [in=180,out=90] (0.25+0.45, 0.25) to[out=0,in=90] (0.25+2*0.3,0);

				\draw[color = blue,thick] (0.25+3*0.3,2) to[in=180,out=270] (2,1);
				\draw[color = blue,thick] (0.25+4*0.3,2) to[in=180,out=270] (2,1.375);
				\draw[color = blue,thick] (0.25+4*0.3,0) to[in=180,out=90] (2,0.625);

				\fill[fill=purple,opacity=0.2] (0,1.75) arc[radius=.25, start angle = 270, end angle = 360] to (0.25+0.3,2) to[in=0,out=270] (0,1.375) to (0,1.75);
			\end{scope}
			\begin{scope}[shift={(6,2)}]
				\draw[dashed] (.25,0) arc[radius=.25, start angle = 0, end angle = 90];
				\draw[dashed] (0,1.75) arc[radius=.25, start angle = 270, end angle = 360];
				\draw[dashed] (1.75,0) arc[radius=.25, start angle = 180, end angle = 90];
				\draw[dashed] (1.75,2) arc[radius=.25, start angle = 180, end angle = 270];
				\draw[color = red,thick] (.25,0) to (1.75,0);
				\draw[color = red,thick] (0,0.25) to (0,1.75);
				\draw[color = red,thick] (.25,2) to (1.75,2);
				\draw[color = red,thick] (2,0.25) to (2,1.75);
		
				\draw[color = blue,thick] (0,1.375) to[in=270,out=0] (0.25+0.3,2);
				\draw[color = blue,thick] (0,1) to[in=270,out=0] (0.25+2*0.3,2);
				\draw[color = blue,thick] (0,0.625) to[in=90,out=0] (0.25+3*0.3,0);
				
				\draw[color = blue, thick] (0.25+0.3,0) to [in=180,out=90] (0.25+0.45, 0.25) to[out=0,in=90] (0.25+2*0.3,0);

				\draw [line width = 1.5pt] (0,0.625) to (0,0.25);

				\fill[fill=purple,opacity=0.2] (.25,0) arc[radius=.25, start angle = 0, end angle = 90] to (0,0.625) to[in=90,out=0] (0.25+3*0.3,0) to (0.25+2*0.3,0) to[in=0,out=90] (0.25+0.45, 0.25) to[in=90,out=180] (0.25+0.3,0) to (0.25,0);

				\fill[fill=purple,opacity=0.2] (0,0.625) to[in=90,out=0] (0.25+3*0.3,0) to (1.75,0) arc[radius=.25, start angle = 180, end angle = 90] to (2,1.75) arc[radius=.25, start angle = 270, end angle = 180] to (0.25, 2) arc[radius=.25, start angle = 360, end angle = 270] -- cycle;

				\draw[color = blue,thick] (0.25+3*0.3,2) to[in=180,out=270] (2,1);
				\draw[color = blue,thick] (0.25+4*0.3,2) to[in=180,out=270] (2,1.375);
				\draw[color = blue,thick] (0.25+4*0.3,0) to[in=180,out=90] (2,0.625);
			\end{scope}


			\begin{scope}[shift={(0,5)}]
				\draw[dashed] (.25,0) arc[radius=.25, start angle = 0, end angle = 90];
				\draw[dashed] (0,1.75) arc[radius=.25, start angle = 270, end angle = 360];
				\draw[dashed] (1.75,0) arc[radius=.25, start angle = 180, end angle = 90];
				\draw[dashed] (1.75,2) arc[radius=.25, start angle = 180, end angle = 270];
				\draw[color = red,thick] (.25,0) to (1.75,0);
				\draw[color = red,thick] (0,0.25) to (0,1.75);
				\draw[color = red,thick] (.25,2) to (1.75,2);
				\draw[color = red,thick] (2,0.25) to (2,1.75);
		
				\draw[color = blue,thick] (0,1.375) to[in=270,out=0] (0.25+0.3,2);
				\draw[color = blue,thick] (0,1) to[in=270,out=0] (0.25+2*0.3,2);
				\draw[color = blue,thick] (0,0.625) to[in=90,out=0] (0.25+3*0.3,0);
				
				\draw[color = blue, thick] (0.25+0.3,0) to [in=180,out=90] (0.25+0.45, 0.25) to[out=0,in=90] (0.25+2*0.3,0);

				\draw[color = blue,thick] (0.25+3*0.3,2) to[in=180,out=270] (2,1);
				\draw[color = blue,thick] (0.25+4*0.3,2) to[in=180,out=270] (2,1.375);
				\draw[color = blue,thick] (0.25+4*0.3,0) to[in=180,out=90] (2,0.625);

				\fill[pattern=north west lines, pattern color=green,opacity=0.8] (0.25,0) arc[radius=.25, start angle = 0, end angle = 90] to (0, 1) to[in=270,out=0] (0.25 + 2 * 0.3, 2) to (1.75,2) arc[radius=.25, start angle = 180, end angle = 270] to (2, 0.25) 
				arc[radius=.25, start angle = 90, end angle = 180] -- cycle;
				\fill[fill=purple,opacity=0.2] (0,1) to[in=270,out=0] (0.25+2*0.3,2) to (0.25+0.3,2) to[in=0,out=270] (0,1.375) -- cycle;
			\end{scope}
			\begin{scope}[shift={(0,7)}]
				\draw[dashed] (.25,0) arc[radius=.25, start angle = 0, end angle = 90];
				\draw[dashed] (0,1.75) arc[radius=.25, start angle = 270, end angle = 360];
				\draw[dashed] (1.75,0) arc[radius=.25, start angle = 180, end angle = 90];
				\draw[dashed] (1.75,2) arc[radius=.25, start angle = 180, end angle = 270];
				\draw[color = red,thick] (.25,0) to (1.75,0);
				\draw[color = red,thick] (0,0.25) to (0,1.75);
				\draw[color = red,thick] (.25,2) to (1.75,2);
				\draw[color = red,thick] (2,0.25) to (2,1.75);
		
				\draw[color = blue,thick] (0,1.375) to[in=270,out=0] (0.25+0.3,2);
				\draw[color = blue,thick] (0,1) to[in=270,out=0] (0.25+2*0.3,2);
				\draw[color = blue,thick] (0,0.625) to[in=90,out=0] (0.25+3*0.3,0);
				
				\draw[color = blue, thick] (0.25+0.3,0) to [in=180,out=90] (0.25+0.45, 0.25) to[out=0,in=90] (0.25+2*0.3,0);
				\fill[fill=purple,opacity=0.2] (0.25+0.3,0) to [in=180,out=90] (0.25+0.45, 0.25) to[out=0,in=90] (0.25+2*0.3,0) -- cycle;
				
				\draw[color = blue,thick] (0.25+3*0.3,2) to[in=180,out=270] (2,1);
				\draw[color = blue,thick] (0.25+4*0.3,2) to[in=180,out=270] (2,1.375);
				\draw[color = blue,thick] (0.25+4*0.3,0) to[in=180,out=90] (2,0.625);
			\end{scope}

			\begin{scope}[shift={(4,5)}]
				\draw[dashed] (.25,0) arc[radius=.25, start angle = 0, end angle = 90];
				\draw[dashed] (0,1.75) arc[radius=.25, start angle = 270, end angle = 360];
				\draw[dashed] (1.75,0) arc[radius=.25, start angle = 180, end angle = 90];
				\draw[dashed] (1.75,2) arc[radius=.25, start angle = 180, end angle = 270];
				\draw[color = red,thick] (.25,0) to (1.75,0);
				\draw[color = red,thick] (0,0.25) to (0,1.75);
				\draw[color = red,thick] (.25,2) to (1.75,2);
				\draw[color = red,thick] (2,0.25) to (2,1.75);
		
				\draw[color = blue,thick] (0,1.375) to[in=270,out=0] (0.25+0.3,2);
				\draw[color = blue,thick] (0,1) to[in=270,out=0] (0.25+2*0.3,2);
				\draw[color = blue,thick] (0,0.625) to[in=90,out=0] (0.25+3*0.3,0);
				
				\draw[color = blue, thick] (0.25+0.3,0) to [in=180,out=90] (0.25+0.45, 0.25) to[out=0,in=90] (0.25+2*0.3,0);

				\draw[color = blue,thick] (0.25+3*0.3,2) to[in=180,out=270] (2,1);
				\draw[color = blue,thick] (0.25+4*0.3,2) to[in=180,out=270] (2,1.375);
				\draw[color = blue,thick] (0.25+4*0.3,0) to[in=180,out=90] (2,0.625);

				\fill[fill=purple,opacity=0.2] (0.25+4*0.3,2) to (1.75,2) arc[radius=.25, start angle = 180, end angle = 270] to (2,1.375) to[in=270,out=180] (0.25+4*0.3,2);
			\end{scope}
			\begin{scope}[shift={(4,7)}]
				\draw[dashed] (.25,0) arc[radius=.25, start angle = 0, end angle = 90];
				\draw[dashed] (0,1.75) arc[radius=.25, start angle = 270, end angle = 360];
				\draw[dashed] (1.75,0) arc[radius=.25, start angle = 180, end angle = 90];
				\draw[dashed] (1.75,2) arc[radius=.25, start angle = 180, end angle = 270];
				\draw[color = red,thick] (.25,0) to (1.75,0);
				\draw[color = red,thick] (0,0.25) to (0,1.75);
				\draw[color = red,thick] (.25,2) to (1.75,2);
				\draw[color = red,thick] (2,0.25) to (2,1.75);
		
				\draw[color = blue,thick] (0,1.375) to[in=270,out=0] (0.25+0.3,2);
				\draw[color = blue,thick] (0,1) to[in=270,out=0] (0.25+2*0.3,2);
				\draw[color = blue,thick] (0,0.625) to[in=90,out=0] (0.25+3*0.3,0);
				
				\draw[color = blue, thick] (0.25+0.3,0) to [in=180,out=90] (0.25+0.45, 0.25) to[out=0,in=90] (0.25+2*0.3,0);
				\fill[fill=purple,opacity=0.2] (0.25+4*0.3,0) to (1.75,0) arc[radius=.25, start angle = 180, end angle = 90] to (2,0.625) to[in=90,out=180] (0.25+4*0.3,0);
				
				\draw[color = blue,thick] (0.25+3*0.3,2) to[in=180,out=270] (2,1);
				\draw[color = blue,thick] (0.25+4*0.3,2) to[in=180,out=270] (2,1.375);
				\draw[color = blue,thick] (0.25+4*0.3,0) to[in=180,out=90] (2,0.625);
			\end{scope}
			\begin{scope}[shift={(6,5)}]
				\draw[dashed] (.25,0) arc[radius=.25, start angle = 0, end angle = 90];
				\draw[dashed] (0,1.75) arc[radius=.25, start angle = 270, end angle = 360];
				\draw[dashed] (1.75,0) arc[radius=.25, start angle = 180, end angle = 90];
				\draw[dashed] (1.75,2) arc[radius=.25, start angle = 180, end angle = 270];
				\draw[color = red,thick] (.25,0) to (1.75,0);
				\draw[color = red,thick] (0,0.25) to (0,1.75);
				\draw[color = red,thick] (.25,2) to (1.75,2);
				\draw[color = red,thick] (2,0.25) to (2,1.75);
		
				\draw[color = blue,thick] (0,1.375) to[in=270,out=0] (0.25+0.3,2);
				\draw[color = blue,thick] (0,1) to[in=270,out=0] (0.25+2*0.3,2);
				\draw[color = blue,thick] (0,0.625) to[in=90,out=0] (0.25+3*0.3,0);
				
				\draw[color = blue, thick] (0.25+0.3,0) to [in=180,out=90] (0.25+0.45, 0.25) to[out=0,in=90] (0.25+2*0.3,0);

				\draw[color = blue,thick] (0.25+3*0.3,2) to[in=180,out=270] (2,1);
				\draw[color = blue,thick] (0.25+4*0.3,2) to[in=180,out=270] (2,1.375);
				\draw[color = blue,thick] (0.25+4*0.3,0) to[in=180,out=90] (2,0.625);

				\fill[fill=purple,opacity=0.2] (0,1.75) arc[radius=.25, start angle = 270, end angle = 360] to (0.25+0.3,2) to[in=0,out=270] (0,1.375) to (0,1.75);
			\end{scope}
			\begin{scope}[shift={(6,7)}]
				\draw[dashed] (.25,0) arc[radius=.25, start angle = 0, end angle = 90];
				\draw[dashed] (0,1.75) arc[radius=.25, start angle = 270, end angle = 360];
				\draw[dashed] (1.75,0) arc[radius=.25, start angle = 180, end angle = 90];
				\draw[dashed] (1.75,2) arc[radius=.25, start angle = 180, end angle = 270];
				\draw[color = red,thick] (.25,0) to (1.75,0);
				\draw[color = red,thick] (0,0.25) to (0,1.75);
				\draw[color = red,thick] (.25,2) to (1.75,2);
				\draw[color = red,thick] (2,0.25) to (2,1.75);
		
				\draw[color = blue,thick] (0,1.375) to[in=270,out=0] (0.25+0.3,2);
				\draw[color = blue,thick] (0,1) to[in=270,out=0] (0.25+2*0.3,2);
				\draw[color = blue,thick] (0,0.625) to[in=90,out=0] (0.25+3*0.3,0);
				
				\draw[color = blue, thick] (0.25+0.3,0) to [in=180,out=90] (0.25+0.45, 0.25) to[out=0,in=90] (0.25+2*0.3,0);

				\fill[fill=purple,opacity=0.2] (.25,0) arc[radius=.25, start angle = 0, end angle = 90] to (0,0.625) to[in=90,out=0] (0.25+3*0.3,0) to (0.25+2*0.3,0) to[in=0,out=90] (0.25+0.45, 0.25) to[in=90,out=180] (0.25+0.3,0) to (0.25,0);

				\fill[pattern=north west lines, pattern color=green,opacity=0.8] (0,0.625) to[in=90,out=0] (0.25+3*0.3,0) to (1.75,0) arc[radius=.25, start angle = 180, end angle = 90] to (2,1.75) arc[radius=.25, start angle = 270, end angle = 180] to (0.25, 2) arc[radius=.25, start angle = 360, end angle = 270] -- cycle;

				\draw[color = blue,thick] (0.25+3*0.3,2) to[in=180,out=270] (2,1);
				\draw[color = blue,thick] (0.25+4*0.3,2) to[in=180,out=270] (2,1.375);
				\draw[color = blue,thick] (0.25+4*0.3,0) to[in=180,out=90] (2,0.625);
			\end{scope}
		\end{tikzpicture}	\caption [Immersed disks representing cancellations for the Whitehead double.] {\textbf{Immersed disks representing cancellations for the Whitehead double.} The two illustrations on the left cancel against each other, and are analogous to the cancellation of Figure~\ref{fig:cable-4-disks}. The illustration on the right provides for a slightly new kind of cancellation; in particular, the disk given by red on the top represents a weight $1$ zero-input operation. The black line on the bottom indicates a cut in the disk.}
		\label{fig:whitehead-double}
\end{figure}
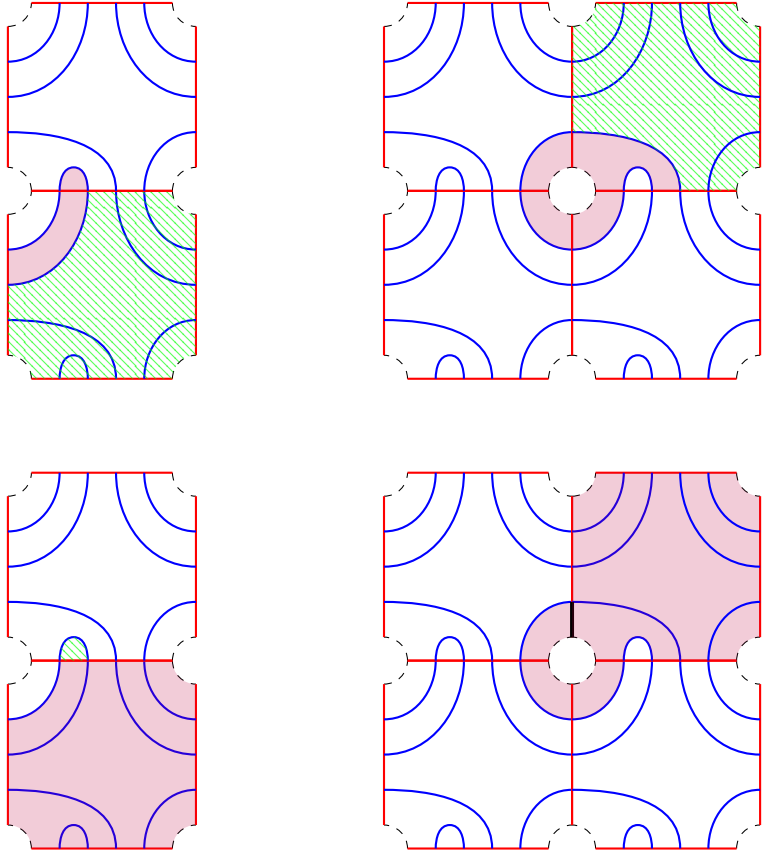

A similar analysis should extend to other satellite operations where the pattern knot admits a genus one Heegaard diagram; that is, other $(1,1)$-pattern knots. The schema of Theorems~\ref{prop:solid-torus-ainf},~\ref{prop:cable-2-ainf},~and~\ref{thm:cable-p-ainf} should extend to proving the $\Ainf$ structure relations for the associated weighted type A modules. In particular, all terms except \ref{term:c3} provide for cancellations that should extend to other bordered Heegaard diagrams. To extend \ref{term:c3} requires a careful analysis of the structure of the module tiling patterns. In Figure~\ref{fig:whitehead-double}, we illustrate two immersed disks (on the left) that provide for a cancellation akin to that of Figure~\ref{fig:cable-4-disks}; and two immersed disks (on the right) that provide for a new type of cancellation.

\section{Algebraic properties of \texorpdfstring{$\CFAm$}{CFA minus}}
\label{ch:algprop}

Now that we have our weighted $\Ainf$-modules constructed, we consider their algebraic properties. In Section~\ref{sec:boundedness}, we prove a useful boundedness property that allows us to take well-defined tensor products. In Section~\ref{sec:uniqueness}, we characterize how the \emph{minus} variant of the modules relate to the \emph{hat} variant.

\subsection{Gradings}
We recall from Section~\ref{sec:gradings} the gradings of algebra elements given by the intermediate grading group $G$, and compute certain additional ones.
\begin{equation} \label{eqn:alg-grading-table}
\begin{gathered}
\begin{aligned}
	\gr(\rho_1) &= \left(-\tfrac{1}{2}; \tfrac{1}{2}, -\tfrac{1}{2}\right)
	\qquad &\gr(\rho_{12}) &= \left(-\tfrac{1}{2}; 1, 0\right)
	\qquad &\gr(\rho_{123}) &= \left(-\tfrac{1}{2}; \tfrac{1}{2}, \tfrac{1}{2}\right) \\
	\gr(\rho_2) &= \left(-\tfrac{1}{2}; \tfrac{1}{2}, \tfrac{1}{2}\right) 
	\qquad &\gr(\rho_{23}) &= \left(-\tfrac{1}{2}; 0, 1\right) 
	\qquad &\gr(\rho_{234}) &= \left(-\tfrac{3}{2}; -\tfrac{1}{2}, \tfrac{1}{2}\right)  \\
	\gr(\rho_3) &= \left(-\tfrac{1}{2}; -\tfrac{1}{2}, \tfrac{1}{2}\right)
	\qquad &\gr(\rho_{34}) &= \left(-\tfrac{3}{2}; -1, 0\right) 
	\qquad &\gr(\rho_{341}) &= \left(-\tfrac{3}{2}; -\tfrac{1}{2}, -\tfrac{1}{2}\right)  \\
	\gr(\rho_4) &= \left(-\tfrac{3}{2}; -\tfrac{1}{2}, -\tfrac{1}{2}\right)
	\qquad &\gr(\rho_{41}) &= \left(-\tfrac{3}{2}; 0, -1\right)
	\qquad &\gr(\rho_{412}) &= \left(-\tfrac{3}{2}; \tfrac{1}{2}, -\tfrac{1}{2}\right)
\end{aligned} \\
\gr(\rho_{1234}) = \gr(\rho_{2341}) = \gr(\rho_{3412}) = \gr(\rho_{4123}) = \gr(U) = (-2; 0, 0).
\end{gathered}
\end{equation}
This renders the computation of the grading of any homogeneous algebra element tractable. We factorize our algebra element as $U^n a' a$, where $a'$ is a Reeb element with length a multiple of four, and $a$ is a Reeb element with length at most three. If $|U^n a'| = 4k$, then $\gr(U^n a') = (-2k; 0, 0)$. The possible gradings for $a$ are recorded in Formula~(\ref{eqn:alg-grading-table}), and in particular we observe that if $a \neq 1$, $\gr(a) = (m; \epsilon_1, \epsilon_2)$, where $m \in \{-1/2, -3/2\}$, $\epsilon_1, \epsilon_2 \in \{0, \pm 1/2, \pm 1\}$ and $|\epsilon_1| + |\epsilon_2| = 1$. The grading of our element is 
	\begin{equation}
		\label{form:grading-fact}
		\gr(U^n a' a) = (-2k + m; \epsilon_1, \epsilon_2).	
	\end{equation}

As we remarked in Section~\ref{sec:gradings}, the gradings on $\CFAm$ work exactly as they do for $\CFAa$. 

\subsubsection{The solid torus}

We record briefly the gradings for $\CFAm(\HD_{st})$ for use in Section~\ref{sec:uniqueness}. The module $\CFAm(\HD_{st})$ has the grading set $\lambda \gr'(\rho_2) \gr'(\rho_1) \backslash G'$. The operation gives $m^0_3(x, \rho_2, \rho_1) = x$ gives us the indeterminacy in the grading set. We can let $\gr'(x) = e$. Using Equation~\ref{eqn:alg-grading-table}, we can compute the grading indeterminacy to be $(-1/2, 1, 0)$ in $G$.

\subsubsection{The \texorpdfstring{$(p,1)$}{(p,1)}-cable}
The gradings for $\CFAa(\CD_p)$ have been computed in \cite{Petkova2013}, with a familiar presentation in \cite{O17}. Recall that $\CFAm(\CD_p)$ is a module over the algebra $\Algm^{U,V}$. The variable $V$ tracks the intersection with the basepoint $z$ and serves the role that $U$ previously served. With a slight abuse of notation, when considering modules over $\Algm^{U,V}$ we let $\gr'(V) = (-1; 1, 1, 1, 1)$ and $\gr'(U) = e$.

To encode the $U$ grading, we consider the enhanced grading group $\widetilde{G}' = G' \times \ZZ$. The grading of the algebra element $U$ is given by $\gre'(U) = (e, 1)$. The weighted type $A$ module $\CFAm(\CD_p)$ has an enhanced grading with this group. To construct this, we let $\gre'(B) = (\gr'(B), n_w(B) - n_z(B))$, where $\gr'(B)$ is the grading in Equation~(\ref{eqn:grading-domain}). This is analogous to the construction in the \emph{hat} flavor of bordered Floer homology, and we refer the reader to \cite[Section~11.4]{LOT1} for details.

\begin{lemma}
\label{lem:cable-p-gr}
The weighted $\Ainf$-module $\CFAm(\CD_p)$ has the grading set $\lambda \gre'(U)^{p} \gre'(\rho_3) \gre'(\rho_2) \backslash \widetilde{G}'$. The gradings are specified by
\begin{align*}
\gre'(x) &= e \\
\gre'(b_i) &= \lambda^{p - i} \gre'(U)^{i} \gre'(\rho_1) (\gre'(\rho_2) \gre'(\rho_1))^{p-i-1} \\
\gre'(c_i) &= \lambda^{p - i - 1} \gre'(\rho_1) (\gre'(\rho_2) \gre'(\rho_1))^{p-i-1},
\end{align*}
for all $1 \leq i \leq p - 1$.
\end{lemma}
\begin{proof}
	The gradings of $b_i$ and $c_i$ follow immediately from Figure~\ref{fig:cable-p-graph} and Equation~(\ref{eqn:module-grading}). In particular, the operations
	\[ m^0_{2 + i}(x, \overbrace{\rho_{12}, \dots, \rho_{12}}^{i}, \rho_1) = c_{p-i-1} \]
	for $0 \leq i \leq p - 2$ give the gradings of $c_i$.
	The operations $m^0_1(b_i) = U^i c_i$ imply that $\lambda^{-1} \gre'(b_i) = \gre'(U)^i \gre'(c_i)$.

	The operation $m^0_3(x, \rho_3, \rho_2) = U^p x$ gives us the indeterminacy in the grading set.
\end{proof}

The value of the gradings in $G$ are straightforward to compute using Formula~(\ref{eqn:alg-grading-table}). The grading indeterminacy is given by $(-1/2; 0, 1)$. The gradings of $b_i$ and $c_i$ for all $1 \leq i \leq p - 1$ is given by
\begin{equation}
\label{eqn:cable-grading}
\gr(b_{p-i}) = \left(\tfrac{1}{2}; -\tfrac{1}{2} + i, -\tfrac{1}{2}\right) \qquad \qquad \gr(c_{p-i})	= \left(-\tfrac{1}{2}; -\tfrac{1}{2} + i, -\tfrac{1}{2}\right) 
\end{equation}

\subsection{Boundedness}
\label{sec:boundedness}

As we noted in Section~\ref{sec:pairing-theorems}, we desire our modules to be \emph{filtered bonsai} for the tensor products in the pairing theorems to be well-defined. This is guaranteed when the bordered Heegaard diagram $\HD$ is admissible. However, our bordered Heegaard diagrams $\CD_p$ are not admissible (but instead provincially admissible).

\begin{proposition} The modules $\CFAc(\CD_p)$ are filtered
	bonsai with respect to the filtration given by total $U$ plus $V$ power.
\end{proposition}
\begin{proof} We wish to show that all quotients $\CFAc(\CD_p) / \Filt^n \CFAc(\CD_p)$ are bonsai; equivalently, any operation of sufficiently large dimension has output in $\Filt^n \CFAc(\CD_p)$.
	It thus suffices to show that for any fixed integer $n$,
	there is a bound on $w + \sum |a_i|$ for which there exists a non-zero operation on $\CFAc(\CD_p)$ of the form
	\[m^w_{1 + n}(\x, a_1, \dots, a_n) = U^k V^l \y,\]
	where the total power $k + l = n$ (c.f. \cite[Lemma~7.5]{LOT:torus-mod}).
	
	Consider how the three moves to generate new operations affect both the total $U$ plus $V$
	power and the quantity to bound, $w + \sum |a_i|$.

	Move (1) yields from operations $m^w_{1+n}(\x, a_1, \dots, a_n) = U^k V^l \y$ and $m^{w'}_{1+m}(\y, a_1', \dots, a_m') = U^{k'}V^{l'}$ the operation $m^{w+w'}_{m+n}(\x, a_1, \dots, a_n a_1', \dots, a_m') = U^{k+k'} V^{l+l'} \y$. The total $U$ plus $V$ power, and the quantity $w + \sum |a_i|$, equals the sum of the quantities from
	the two prior operations. 
	Move (2) yields from the operation $m^w_{1+n}(\x, a_1, \dots, a_n) = U^k V^l \y$ the operation $m^w_{3+n}(\x, a_1, \dots, a_i \rho_{j+1}, \rho_j, \rho_{j-1}, \rho_{j-2} a_{i+1}, \dots, a_n) = U^{k+1} V^{l+1} \y$. This increases the total power $U$ plus $V$ power by two, and
	and the quantity $w + \sum |a_i|$ by four. Move (3) yields from the operation $m^w_{1+n}(\x, a_1, \dots, a_n) = U^k V^l \y$ the operation $m^{w+1}_{n-1}(\x, a_1, \dots, a_{i-2}, a_{i-1}a_{i+1}, a_{i+2}, \dots, a_n) = U^k V^l \y$, where $|a_i| = 4$. This leaves the total $U$ plus $V$ power unaffected, and
	decreases the quantity $w + \sum |a_i|$ by three.

	Suppose that for a fixed $n$, there is no such bound on
	$w + \sum |a_i|$. Choose generators $\x$ and $\y$ such that there exists
	a sequence of operations from $\x$ to $U^k V^l y$ where $w + \sum |a_i|$
	is unbounded and $k + l = n$. We note that move (2) increases the total $U$ plus $V$ power, while
	move (3) decreases the quantity $w + \sum |a_i|$. To then generate such a sequence
	of operations starting with the operations in Figure~\ref{fig:cable-p-graph}, we must use move (1) on operations where one of the two has total $U$ plus $V$ power zero
	and $w + \sum |a_i| > 0$, and we must do so an unbounded number of times.

	The operations with total $U$ plus $V$ power zero for our modules are those given by
	the operations in Figure~\ref{fig:cable-p-graph} with total $U$ plus $V$ power zero, and their closure under move (1).
	These operations, however, form no directed cycles. Thus, they
	cannot be performed an unbounded number of times.
\end{proof}

\subsection{Uniqueness}
\label{sec:uniqueness}

Given a graded weighted $\Ainf$-module $\fModule$ constructed as an extension of the operations on $\CFAa(\HD)$, the operations on $\CFAm(\HD)$ seem to be forced by the $\Ainf$-relations and the gradings. For example, in the case of the $0$-framed solid torus, the operation $m^0_3(x, \rho_4, \rho_3) = Ux$ follows from the $\Ainf$-relation for the inputs $x, \rho_4, \rho_3, \rho_2, \rho_1$, to cancel against the non-zero term $m^0_2(x, \mu^0_4(\rho_4, \rho_3, \rho_2, \rho_1))$. Similarly, the operation $m^0_3(x, \rho_4, \rho_{34}, \rho_4) = U^2 x$ as the result of a move (1) follows from the $\Ainf$-relation for the inputs $x, \rho_4, \rho_3, \rho_4, \rho_3$ to cancel against the non-zero term $m^0_3(m^0_3(x, \rho_4, \rho_3), \rho_4, \rho_3)$. There are analogous ways to \emph{force} the existence of operations corresponding to move (2) and move (3), using the terms they cancelled against in Theorem~\ref{prop:solid-torus-ainf}. That there are no extraneous operations seems to be forced by the gradings.

In this section, we prove a certain uniqueness property for the weighted $\Ainf$-modules we have constructed to capture this notion. We see that any graded weighted $\Ainf$-module constructed as an extension of the operations on $\CFAa(\CD_p)$ has isomorphic tensor product with the dualizing bimodule $\CFDDm(\Id)$ to that of $\CFAm(\CD_p)$. In the other words, the associated type D module is uniquely determined. As with the \emph{hat} flavor of bordered Floer homology \cite{LOT2}, we expect the bimodule $\CFDDm(\Id)$ to be quasi-invertible; if so, this would determine such a weighted extension to be homotopy equivalent to $\CFAm(\CD_p)$.

Recall that the dualizing bimodule $\CFDDm(\Id)$ is the type DD bimodule over $(\Algm^{U=1}, \Algm)$ with generators $\iota_0 \otimes \iota_0$ and $\iota_1 \otimes \iota_1$, and differential $\delta^1$ given by
\begin{equation}
	\label{eqn:dualizing-bimodule}
	\begin{aligned}
	\delta^1(\iota_0 \tensor \iota_0) &= (\rho_1 \tensor \rho_3 + \rho_3 \tensor \rho_1 + \rho_{123} \tensor \rho_{123} + \rho_{341} \tensor \rho_{341}) \tensor (\iota_1 \tensor \iota_1) \\
	\delta^1(\iota_1 \tensor \iota_1) &= (\rho_2 \tensor \rho_2 + \rho_4 \tensor \rho_4 + \rho_{234} \tensor \rho_{412} + \rho_{412} \tensor \rho_{234}) \tensor (\iota_0 \tensor \iota_0) \\
	\end{aligned}
\end{equation}

Further recall that we have fixed our choice of a weighted algebra $\AsDiag$, and a weighted module diagonal primitive $\TrPMDiag$ compatible with $\AsDiag$. \cite[Proposition~7.7]{LOT:torus-mod} proves that our results are independent of the choice of the weighted module diagonal primitive $\TrPMDiag$.

\subsubsection{The solid torus}

As a warm-up, we consider the case of the $0$-framed solid torus. To start with, the following lemma computes $\CFDm(\HD_{st})$.

\begin{lemma} \label{lem:solid-torus-cfd} The tensor product $\CFDm(\HD_{st}) = \CFAm_{U=1}(\HD_{st}) \boxtimes_{\Algm^{U=1}} \CFDDm(\Id)$ is as illustrated in Figure~\ref{fig:solid-torus-cfd}.
\begin{figure}[h!tbp]
\centering
\begin{tikzpicture}
	\node at (0,0) (x) {$\circ$};
	\node at (-1,0) (lphant) {};
	\node at (1,0) (rphant) {};
	\draw[->, bend left=80] (x) to node[inner sep=2ex,left]
	{\lab{\rho_{41}}} (lphant) to (x);
	\draw[->, bend right=80] (x) to node[inner sep=2ex,right]
	{\lab{\rho_{23}}} (rphant) to (x);
\end{tikzpicture}
	\caption[The tensor product \texorpdfstring{$\CFDm(\HD_{st}) = \CFAm_{U=1}(\HD_{st}) \boxtimes_{\Algm^{U=1}} \CFDDm(\Id)$}{CFD minus of the 0-framed solid torus equalling CFA minus of the 0-framed solid torus tensor CFDD minus of identity}.]
	{\textbf{The tensor product $\CFDm(\HD_{st}) = \CFAm_{U=1}(\HD_{st}) \boxtimes_{\Algm^{U=1}} \CFDDm(\Id)$.}}
	\label{fig:solid-torus-cfd}
\end{figure}
\end{lemma}
\begin{proof}
	We abbreviate the generator $x \otimes (\iota_1 \otimes \iota_1)$ of the tensor product $\CFDm(\HD_{st})$ as $x$. Both terms in $\delta^1(x) = \rho_{41} \tensor x + \rho_{23} \tensor x$ follow from the chosen weighted module diagonal primitive $\TrMPrim^{3,0}$. On the sequence of type DD module outputs $(\rho_4 \tensor \rho_4), (\rho_3 \tensor \rho_1)$ and $(\rho_2 \tensor \rho_2), (\rho_1 \tensor \rho_3)$, we apply a $m^0_3$ on the left, and a $\mu_2^0$ on the right. Neither term is cancelled by any higher module diagonal primitive terms, as there is a term
	\[ \mu^1_0 \tensor x = \rho_{2341} \tensor x + \rho_{4123} \tensor x, \]
	in the structure equation for $\CFDm(\HD_{st})$, and they cancel it.
	
	That there are no more terms in the differential follows from an argument using the gradings. The indeterminacy in the grading for $\CFDm(\HD_{st})$ is given by $\lambda\gr(\rho_{23}) = (1/2; 0, 1)$. Consider any term $a'' \tensor x$ in the differential, and factorize $a''$ as $U^n a' a$ such that Formula~(\ref{form:grading-fact}) holds true and $\gr(a'') = (-2k + m; \epsilon_1, \epsilon_2)$. As there is no indeterminacy in the first $\SpinC$-component, we know that $\epsilon_2 = 0$ and so $a = \rho_{23}$ or $a = \rho_{41}$, or $a = 1$. If the length $|U^n a'| > 0$ (so that $k > 0$) or $a = 1$, the relative Maslov components do not line up. So, $k = 0$, and $a''$ is either $\rho_{23}$ or $\rho_{41}$.
	
	We remark that the previous paragraph shows that $|a''| = 2$ for any term $a'' \tensor x$ in the differential on the tensor product. For such a term of length two to arise in Equation~(\ref{eqn:m-dd-tensor}) as the output of a $\mu(T)$, the operation must be a $\mu^0_2$. The only way $\rho_{23}$ or $\rho_{41}$ arises from a $\mu^0_2$ is via the sequence of type DD module outputs $(\rho_4 \tensor \rho_4), (\rho_3 \tensor \rho_1)$ and $(\rho_2 \tensor \rho_2), (\rho_1 \tensor \rho_3)$. This provides an alternative reason for why the terms $\rho_{41} \tensor x$ and $\rho_{23} \tensor x$ are not cancelled by any higher module diagonal primitive terms.
	\end{proof}

With that model computation, we can prove our uniqueness result.

\begin{proposition}
Let $\fModule$ be a graded weighted $\Ainf$-module over $\Algm$ such that
$\widehat{M} = M / UM$, viewed as an $\Ainf$-module over $\widehat{\Alg}$, is the graded $\Ainf$-module $\CFAm(\HD_{st}) / U\CFAm(\HD_{st}) = \CFAa(\HD_{st})$.

Then, $\fModule_{U = 1} \boxtimes_{\Algm^{U=1}}\CFDDm(\Id)$ is isomorphic to $\CFAm_{U=1}(\HD_{st}) \boxtimes_{\Algm^{U=1}} \CFDDm(\Id) = \CFDm(\HD_{st})$ as graded weighted type D modules.
\end{proposition}
\begin{proof}
	Let $D$ be $\CFAm_{U=1}(\HD_{st}) \boxtimes_{\Algm^{U=1}} \CFDDm(\Id)$. We denote its differential by $\delta^1_D$ and its generator by $\widetilde{x}$. Let $\widetilde{m}$ be the module operation on $\fModule$. The gradings on $D$ are identified with the gradings on $\CFDm(\HD_{st})$ by construction. Using the same line of reasoning as in the proof of Lemma~\ref{lem:solid-torus-cfd}, the gradings on $D$ force any term in the differential $\delta^1_D$ to be of the form $\rho_{23} \tensor \wt x$ and $\rho_{41} \tensor \wt x$. Both these terms must appear in the differential to cancel the term $\mu^1_0 \tensor \wt x$, and so $D \cong \CFDm(\HD_{st})$.
	
	We remark that we could show that $\wt{m}^0_3(\wt x, \rho_{4}, \rho_{3}) = U \wt x$ by considering the weight $0$ $5$-input $\Ainf$ relation with inputs $\wt x, \rho_4, \rho_3, \rho_2, \rho_1$. The only term that cancel the non-zero term
	$\wt{m}^0_2(\wt x, \mu^0_4(\rho_4, \rho_3, \rho_2, \rho_1))$
	is the term
	$\wt{m}^0_3(\wt{m}^0_3(\wt x, \rho_4, \rho_3), \rho_2, \rho_1)$; any other algebra operation on the inputs is zero, and the only composition of module operations compatible with the idempotents is either $\wt{m}^0_1(\wt{m}^0_5(\wt x, \rho_4, \rho_3, \rho_2, \rho_1))$ or $\wt{m}^0_5(\wt{m}^0_1(\wt x), \rho_4, \rho_3, \rho_2, \rho_1)$, where neither module operation would be compatible with the gradings.
\end{proof}

\subsubsection{The \texorpdfstring{$(p,1)$}{(p,1)}-cable}

The results follow for the $(p, 1)$-cable with a similar argument, where the grading on the type D module eliminates most possibilities.
	We start by characterizing what terms can appear in the differential for $\CFDm(\CD_p)$ using its gradings. This involves a lot of case work; an impatient reader can skip ahead to Proposition~\ref{prop:cable-p-cfd}, where we summarize the results.
	
	The grading indeterminacy for $\CFDm(\CD_p)$ is given by $\lambda \gr(\rho_{12}) = (1/2; 1, 0)$. Analogous to Formula~(\ref{eqn:cable-grading}), the gradings on $\CFDm(\CD_p)$ are given by $\gr(x) = (0; 0, 0)$ and
 \[ \gr(b_{p-i}) = \left(\tfrac{1}{2}; \tfrac{1}{2}, \tfrac{1}{2} - i\right) \qquad \qquad \gr(c_{p-i})	= \left(-\tfrac{1}{2}; \tfrac{1}{2}, \tfrac{1}{2} - i\right) \]
	
	\textbf{Between the generators $b_i$.} First, consider any term $a'' \tensor b_{p-i}$ in the differential $\delta^1(b_{p-j})$. Factorize $a''$ as $U^n a' a$ such that Formula~(\ref{form:grading-fact}) holds true and $\gr(a'') = (-2k + m; \epsilon_1, \epsilon_2)$. For there to be such a term in the differential,
	\begin{align*} 
	\gr(a'' \tensor b_{p-i}) = &\lambda^{-1} \gr(b_{p-j})  \\
	\implies &\lambda \left(\tfrac{1}{2}; \tfrac{1}{2}, \tfrac{1}{2} - j\right)^{-1} (-2k + m; \epsilon_1, \epsilon_2) \left(\tfrac{1}{2}; \tfrac{1}{2}, \tfrac{1}{2} - i\right) (\tfrac{l}{2}; l, 0) = (0; 0, 0), \numberthis \label{eqn:b-b-grading} 
	\end{align*}
	where the factor $(l/2; l, 0)$ comes from the ambiguity in the grading set, for some $l \in \ZZ$. Comparing the second $\SpinC$-component, observe that $j - i = -\epsilon_2$, so either $\epsilon_2 = 0$ or $\epsilon_2 = \pm 1$. Suppose first that $\epsilon_2 = 0$. Let us compute Equation~(\ref{eqn:b-b-grading}),
	\begin{align*} 
	& &\lambda \left(-\tfrac{1}{2}; -\tfrac{1}{2}, i - \tfrac{1}{2}\right) (-2k + m; \epsilon_1, \epsilon_2) \left(\tfrac{1}{2}; \tfrac{1}{2}, \tfrac{1}{2} - i\right) (\tfrac{l}{2}; l, 0) &= (0; 0, 0) \\
	&\implies &\lambda^{1 - 2k + m + l/2} \left(0; -\tfrac{1}{2}, i - \tfrac{1}{2}\right) (0; \epsilon_1, 0) \left(0; \tfrac{1}{2}, \tfrac{1}{2} - i\right) (0; l, 0) &= (0; 0, 0) \\
	&\implies &\lambda^{1 - 2k + m + il} \left(0; -\tfrac{1}{2}, i - \tfrac{1}{2}\right) (0; \epsilon_1, 0) \left(0; l + \tfrac{1}{2}, \tfrac{1}{2} - i\right)  &= (0; 0, 0) \\
	&\implies &\lambda^{1 - 2k + m + il + \epsilon_1/2 - \epsilon_1 i} \left(0; -\tfrac{1}{2}, i - \tfrac{1}{2}\right) \left(0; \epsilon_1 + l + \tfrac{1}{2}, \tfrac{1}{2} - i\right)  &= (0; 0, 0) \\
	&\implies &\lambda^{1 - 2k + m + \epsilon_1 - 2 \epsilon_1 i + l/2 } (0; \epsilon_1 + l, 0) &= (0; 0, 0) \numberthis \label{eqn:b-b-e2-zero}
	\end{align*}
	Let's consider the case when $\epsilon_1 = 0$. When both $\epsilon_1 = \epsilon_2 = 0$, then $m = 0$. There is no solution to Equation~(\ref{eqn:b-b-e2-zero}), as the Maslov component of the left-hand side is $1 - 2k$. When $\epsilon_1 = 1$ the Maslov component $m$ of $a$ (from our factorization $a'' = U^n a' a$) is $-1/2$. The Maslov component in Equation~(\ref{eqn:b-b-e2-zero}) says $1 - 2k - 2i = 0$. There is again no solution. When $\epsilon_2 = -1$, $m = -3/2$. The Maslov component in Equation~(\ref{eqn:b-b-e2-zero}) says $-1 - 2k +2i  = 0$ for which there is no solution.
	
	Consider the case when $\epsilon_2 = 1$, so $\epsilon_1 = 0$, $m = -1/2$, and $a = \rho_{23}$. Note that $i = j + 1$. Equation~(\ref{eqn:b-b-grading}) evaluates to
	\[ \left(-2k + \tfrac{l}{2}; l, 0\right) = (0; 0, 0). \]
	So $l = 0$, and $k = 0$. This implies that the only permissible term in this case is
	$\rho_{23} \tensor b_{p-(j+1)}$ in the differential $\delta^1(b_{p-j})$. The case $\epsilon_2 = -1$, so $\epsilon_1 = 0$, $m = -3/2$, and $a = \rho_{41}$ is entirely analogous, and we get the constraint that $k = 0$ again. The only permissible term in that case is $\rho_{23} \tensor b_{p-(j-1)}$ in the differential $\delta^1(b_{p-j})$.
	
	\textbf{Between the generators $c_i$.} This case is entirely analogous to the case where we considered permissible terms between the generators $b_i$. The permissible terms are $\rho_{41} \tensor c_{p-(j-1)}$ and $\rho_{23} \tensor c_{p-(j+1)}$ in the differential $\delta^1(c_{p-j})$.
	
	\textbf{From the generators $b_i$ to $c_i$.} Consider any term $a'' \tensor c_{p-i}$ in the differential $\delta^1(b_{p-j})$. Factorize $a''$ as $U^n a' a$ such that Formula~(\ref{form:grading-fact}) holds true and $\gr(a'') = (-2k + m; \epsilon_1, \epsilon_2)$. For there to be such a term in the differential,
	\begin{align*} 
	\gr(a'' \tensor c_{p-i}) = &\lambda^{-1} \gr(b_{p-j})  \\
	\implies &\lambda \left(\tfrac{1}{2}; \tfrac{1}{2}, \tfrac{1}{2} - j\right)^{-1} (-2k + m; \epsilon_1, \epsilon_2) \left(-\tfrac{1}{2}; \tfrac{1}{2}, \tfrac{1}{2} - i\right) (\tfrac{l}{2}; l, 0) = (0; 0, 0), \numberthis \label{eqn:b-c-grading} 
	\end{align*}
	for some $l \in \ZZ$. Comparing the second $\SpinC$-component, observe that $j - i = -\epsilon_2$, so either $\epsilon_2 = 0$ or $\epsilon_2 = \pm 1$. Suppose first that $\epsilon_2 = 0$. Equation~(\ref{eqn:b-c-grading}) evaluates to
	\[ \left(-2k + m + \epsilon_1 - 2 \epsilon_1 i + \tfrac{l}{2}; \epsilon_1 + l, 0\right) = (0; 0, 0). \]
	Comparing this to Equation~(\ref{eqn:b-b-grading}) tells us that when $\epsilon_1 = 0$ so $a = 1$, we get that $k = 0$, and there is a permissible term $1 \tensor c_{p-i}$ in the differential $\delta^1(b_{p-i})$. When $\epsilon_1 = 1$, $a$ is forced to be $\rho_{12}$; when $\epsilon_1 = -1$, $a$ is forced to be $\rho_{34}$. Neither have compatible idempotents. (There are elements with compatible gradings, but we do not spell these out.)
	
	Consider the case when $\epsilon_2 = 1$, so $\epsilon_1 = 0$, $m = -1/2$, and $a = \rho_{23}$. Note that $i = j + 1$. Equation~(\ref{eqn:b-c-grading}) evaluates to
	\[ \left(-2k + \tfrac{l}{2} - 1; l, 0\right) = (0; 0, 0). \]
	This has no solutions. The case $\epsilon_2 = -1$, so $\epsilon_1 = 0$, $m = -3/2$, and $a = \rho_{41}$ is entirely analogous, with no solutions.

	\textbf{From the generators $c_i$ to $b_i$.} Consider any term $a'' \tensor b_{p-i}$ in the differential $\delta^1(c_{p-j})$. Factorize $a''$ as $U^n a' a$ such that Formula~(\ref{form:grading-fact}) holds true and $\gr(a'') = (-2k + m; \epsilon_1, \epsilon_2)$. For there to be such a term in the differential,
	\begin{align*} 
	\gr(a'' \tensor b_{p-i}) = &\lambda^{-1} \gr(c_{p-j})  \\
	\implies &\lambda \left(-\tfrac{1}{2}; \tfrac{1}{2}, \tfrac{1}{2} - j\right)^{-1} (-2k + m; \epsilon_1, \epsilon_2) \left(\tfrac{1}{2}; \tfrac{1}{2}, \tfrac{1}{2} - i\right) (\tfrac{l}{2}; l, 0) = (0; 0, 0), \numberthis \label{eqn:c-b-grading} 
	\end{align*}
	for some $l \in \ZZ$. Comparing the second $\SpinC$-component, observe that $j - i = -\epsilon_2$, so either $\epsilon_2 = 0$ or $\epsilon_2 = \pm 1$. Suppose first that $\epsilon_2 = 0$. Equation~(\ref{eqn:c-b-grading}) evaluates to
	\[ \left(-2k + m + \epsilon_1 - 2 \epsilon_1 i + \tfrac{l}{2} + 2; \epsilon_1 + l, 0\right) = (0; 0, 0). \]
	When $\epsilon_1 = 0$ so $a = 1$, we get that $2k - 2 = 0$. There is a permissible term $U^n a' \tensor b_{p-i}$ in the differential $\delta^1(c_{p-i})$, where either $n = 1$ and $a' = 1$ or $n = 0$ and $|a'| = 4$. When $\epsilon_1 = 1$, $a$ is forced to be $\rho_{12}$; when $\epsilon_1 = -1$, $a$ is forced to be $\rho_{34}$. Neither have compatible idempotents.
	
	Consider the case when $\epsilon_2 = 1$, so $\epsilon_1 = 0$, $m = -1/2$, and $a = \rho_{23}$. Note that $i = j + 1$. Equation~(\ref{eqn:b-c-grading}) evaluates to
	\[ \left(-2k + \tfrac{l}{2} + 1; l, 0\right) = (0; 0, 0). \]
	This has no solutions. The case $\epsilon_2 = -1$, so $\epsilon_1 = 0$, $m = -3/2$, and $a = \rho_{41}$ is entirely analogous, with no solutions.
	
	\textbf{From the generator $x$ to $b_i$.} Consider any term $a'' \tensor b_{p-i}$ in the differential $\delta^1 (x)$.  Factorize $a''$ as $U^n a' a$ such that Formula~(\ref{form:grading-fact}) holds true and $\gr(a'') = (-2k + m; \epsilon_1, \epsilon_2)$. For there to be such a term in the differential,
	\begin{align*} 
	\gr(a'' \tensor b_{p-i}) = &\lambda^{-1} \gr(x)  \\
	\implies &\lambda (-2k + m; \epsilon_1, \epsilon_2) \left(\tfrac{1}{2}; \tfrac{1}{2}, \tfrac{1}{2} - i\right) (\tfrac{l}{2}; l, 0) = (0; 0, 0), \numberthis \label{eqn:x-b-grading} 
	\end{align*}
	for some $l \in \ZZ$. Comparing the second $\SpinC$-component, observe that $i - 1/2 = \epsilon_2$ and as $i \geq 1$, $\epsilon_2$ is forced to be $1/2$ and $i$ is forced to be one. The possibilities for $a$ are $\rho_2, \rho_3, \rho_{123}$, or $\rho_{234}$; of these, only $\rho_3$ and $\rho_{123}$ have compatible idempotents. Let us first consider $a = \rho_3$, in which case Equation~(\ref{eqn:x-b-grading}) simplifies to
	$\left(-2k + \tfrac{l}{2} + 1; l, 0\right) = (0; 0, 0).$
	This has no solutions. In the case $a = \rho_{123}$, Equation~(\ref{eqn:x-b-grading}) simplifies to
	$\left(-2k + \tfrac{l}{2}; l, 0\right) = (0; 0, 0).$
	So, the only permissible term in this case is the term $\rho_{123} \tensor b_{p-1}$ in the differential $\delta^1(x)$.

	\textbf{From the generator $b_i$ to $x$.} Consider any term $a'' \tensor x$ in the differential $\delta^1 (b_{p-j})$.  Factorize $a''$ as $U^n a' a$ such that Formula~(\ref{form:grading-fact}) holds true and $\gr(a'') = (-2k + m; \epsilon_1, \epsilon_2)$. For there to be such a term in the differential,
	\begin{align*} 
	\gr(a'' \tensor x) = &\lambda^{-1} \gr(b_{p-j})  \\
	\implies &\lambda \left(\tfrac{1}{2}; \tfrac{1}{2}, \tfrac{1}{2} - j\right)^{-1} (-2k + m; \epsilon_1, \epsilon_2) (\tfrac{l}{2}; l, 0) = (0; 0, 0), \numberthis \label{eqn:b-x-grading} 
	\end{align*}
	for some $l \in \ZZ$. Comparing the second $\SpinC$-component, observe that $j - 1/2 = -\epsilon_2$ and as $j \geq 1$, $\epsilon_2$ is forced to be $-1/2$ and $j$ is forced to be one. The possibilities for $a$ are $\rho_1, \rho_4, \rho_{341}$, or $\rho_{412}$; of these, only $\rho_4$ and $\rho_{412}$ have compatible idempotents. In the case $a = \rho_4$,  Equation~(\ref{eqn:b-x-grading}) simplifies to
	$\left(-2k + \tfrac{l-1}{2}; l-1, 0\right) = (0; 0, 0)$. The only permissible term in this case is the term $\rho_{4} \tensor x$ in the differential $\delta^1(b_{p-1})$.
	 In the case $a = \rho_{412}$, Equation~(\ref{eqn:b-x-grading}) simplifies to
	$\left(-2k + \tfrac{l}{2} - 1; l, 0\right) = (0; 0, 0).$
	This has no solutions.

	\textbf{From the generator $x$ to $c_i$; or from the generator $c_i$ to $x$.} The analysis in these cases are entirely analogous to the previous two cases. The only permissible terms in this case are the terms $\rho_{3} \tensor c_{p-1}$ in the differential $\delta^1(x)$; and $\rho_{412} \tensor x$ in the differential $\delta^1(c_{p-1})$.

\begin{proposition} \label{prop:cable-p-cfd}
	The tensor product $\CFDm(\CD_{p}) = \CFAm_{U,V=1}(\CD_{p}) \boxtimes_{\Algm^{U=1}} \CFDDm(\Id)$ is as illustrated in Figure~\ref{fig:cable-p-cfd}.
	\begin{figure}[h!]
	\centering
	\begin{tikzpicture}[scale=2]
		\node at (0,0) (x) {$x$};
		\node at (2,1) (cpm1) {$c_{p-1}$};
		\node at (3,1) (cpm2) {$c_{p-2}$};
		\node at (4,1) (dotsup1) {$\dots$};
		\node at (5,1) (ci) {$c_i$};
		\node at (6,1) (dotsup2) {$\dots$};
		\node at (7,1) (c1) {$c_1$};
		
		\node at (2,-1) (bpm1) {$b_{p-1}$};
		\node at (3,-1) (bpm2) {$b_{p-2}$};
		\node at (4,-1) (dotsdown1) {$\dots$};
		\node at (5,-1) (bi) {$b_i$};
		\node at (6,-1) (dotsdown2) {$\dots$};
		\node at (7,-1) (b1) {$b_1$};
		
		\draw[->, bend right=15] (x) to node[below]{\lab{\rho_3}} (cpm1);
		\draw[->, bend right=15] (cpm1) to node[above,sloped]{\lab{\rho_{412}}} (x);
		\draw[->, bend left=80] (x) to node[left]{\lab{\rho_{12}}} (-.5,0) to (x);
		\draw[->] (bpm1) to node[right]{} (cpm1);
		\draw[->, bend right=15] (x) to node[below, sloped]{\lab{\rho_{123}}} (bpm1);
		\draw[->, bend right=15] (bpm1) to node[above]{\lab{\rho_4}} (x);

		\draw[->, bend left=15](cpm1) to node[above]{\lab{\rho_{23}}} (cpm2);
		\draw[->, bend left=15](cpm2) to node[below]{\lab{\rho_{41}}} (cpm1);
		\draw[->, bend right=15](bpm1) to node[below]{\lab{\rho_{23}}} (bpm2);
		\draw[->, bend right=15](bpm2) to node[above]{\lab{\rho_{41}}} (bpm1);

		\draw[->] (bpm2) to node[right]{} (cpm2);

		\draw[->, bend left=15](cpm2) to node[above]{\lab{\rho_{23}}} (dotsup1);
		\draw[->, bend left=15](dotsup1) to node[below]{\lab{\rho_{41}}} (cpm2);
		\draw[->, bend right=15](bpm2) to node[below]{\lab{\rho_{23}}} (dotsdown1);
		\draw[->, bend right=15](dotsdown1) to node[above]{\lab{\rho_{41}}} (bpm2);

		\draw[->] (bi) to node[right]{} (ci);

		\draw[->, bend left=15](dotsup1) to node[above]{\lab{\rho_{23}}} (ci);
		\draw[->, bend left=15](ci) to node[below]{\lab{\rho_{41}}} (dotsup1);
		\draw[->, bend right=15](dotsdown1) to node[below]{\lab{\rho_{23}}} (bi);
		\draw[->, bend right=15](bi) to node[above]{\lab{\rho_{41}}} (dotsdown1);

		\draw[->, bend left=15](ci) to node[above]{\lab{\rho_{23}}} (dotsup2);
		\draw[->, bend left=15](dotsup2) to node[below]{\lab{\rho_{41}}} (ci);
		\draw[->, bend right=15](bi) to node[below]{\lab{\rho_{23}}} (dotsdown2);
		\draw[->, bend right=15](dotsdown2) to node[above]{\lab{\rho_{41}}} (bi);

		\draw[->, bend left=15](dotsup2) to node[above]{\lab{\rho_{23}}} (c1);
		\draw[->, bend left=15](c1) to node[below]{\lab{\rho_{41}}} (dotsup2);
		\draw[->, bend right=15](dotsdown2) to node[below]{\lab{\rho_{23}}} (b1);
		\draw[->, bend right=15](b1) to node[above]{\lab{\rho_{41}}} (dotsdown2);

		\draw[->, bend right=15] (b1) to node[left]{} (c1);
		\draw[->, bend right=15] (c1) to node[above,sloped]{\lab{\rho_{2341}}} (b1);
	\end{tikzpicture}
		\caption [The tensor product \texorpdfstring{$\CFDm(\CD_p) = \CFAm_{U,V=1}(\CD_p) \boxtimes_{\Algm^{U=1}} \CFDDm(\Id)$}{CFD minus of CDp equalling CFA minus of CDp tensor CFDD minus of identity}.] {\textbf{The tensor product $\CFDm(\CD_p) = \CFAm_{U,V=1}(\CD_p) \boxtimes_{\Algm^{U=1}} \CFDDm(\Id)$.}}
	\label{fig:cable-p-cfd}
	\end{figure}
\end{proposition}
\begin{proof}
	All terms in the differential with the corresponding Reeb element of length less than four in Figure~\ref{fig:cable-p-cfd} follow from our choice of $\TrMPrim^{3,0}$ and $\TrMPrim^{4,0}$, and the correspondingly labelled operations in Figure~\ref{fig:cable-p-graph}. In particular, we note that terms with length three Reeb elements follow from the first pair of trees of $\TrMPrim^{4,0}$, and every other summand in the weighted module diagonal primitive is zero as $\mu^0_3 = 0$.
	
	The sole way a length two Reeb element arises in the differential of the tensor product is with a $\mu^0_2$, and so such terms do not get cancelled by any higher module diagonal primitive terms. The relevant length three Reeb elements arise from a binary tree of $\mu^0_2$ (as they indeed do in our differential), or from the type DD module outputs $\rho_{123} \tensor \rho_{123}$ or $\rho_{234} \tensor \rho_{412}$. The latter do not lead to cancellations as the gradings preclude there being any operations $m_2(x, \rho_{123})$ or $m_2(c_{p-1}, \rho_{234})$.
	
	The term $\delta^1(c_1) = \rho_{2341} \tensor b_1$ follows from a compatible choice of $\TrMPrim^{5,0}$. Instead, we argue it follows from the structure relations for $\CFDm(\CD_p)$. Indeed, there are terms $\rho_{2341} \tensor b_1$ and $\rho_{2341} \tensor c_1$ in the structure relations that arise from $\mu^1_0 \tensor b_1$ and $\mu^1_0 \tensor c_1$, and the only term permissible with the gradings that can cancel them is $\delta^1(c_1)$ having a term $\rho_{2341} \tensor b_1$.
	
	We recap our analysis of the possible terms in the differential for $\CFDm(\CD_p)$. Except for the terms illustrated in Figure~\ref{fig:cable-p-cfd}, the only possible terms are in the differential $\delta^1(c_{p-i})$ and of the form $a \tensor b_{p-i}$ for $a \in \{U, \rho_{2341}, \rho_{4123}\}$. All three would introduce corresponding terms $a \tensor c_{p-i}$ and $a \tensor b_{p-i}$ in the structure relations that do not cancel. Thus, they do not occur in the differential.
\end{proof}

\begin{theorem}
Let $\fModule$ be a graded weighted $\Ainf$-module over $\Algm^{U,V}$ such that
$\widehat{M} = M / VM$, viewed as an $\Ainf$-module over $\widehat{\Alg}[U]$, is the graded $\Ainf$-module $\CFAm(\CD_{p}) / V\CFAm(\CD_{p}) = \CFAa(\CD_{p})$.

Then, $\fModule_{U,V = 1} \boxtimes_{\Algm^{U=1}}\CFDDm(\Id)$ is isomorphic to $\CFAm_{U,V=1}(\CD_{p}) \boxtimes_{\Algm^{U=1}} \CFDDm(\Id) = \CFDm(\CD_p)$ as graded weighted type D modules.
\end{theorem}
\begin{proof}
	Let $\wt m$ be the module operation on $\fModule$, and let $\wt{x}, \wt{b}_i, \wt{c}_i$ be its generators, for $1 \leq i \leq p - 1$. We make use of the relative $U$-grading on $\fModule$. From Lemma~\ref{lem:cable-p-gr}, we know the ambiguity in the $U$-grading (for both $\CFAm(\CD_p)$ and $\CFAa(\CD_p)$) is $p$. The module operations respect the $U$-grading.
	
	Consider the weight $0$ $5$-input $\Ainf$ relation for $\fModule$ with inputs $\wt{b}_{p-i}, \rho_2, \rho_1, \rho_4, \rho_3$ where $i < p - 1$. In this $\Ainf$ relation, there is a non-zero term
	$\wt{m}^0_2(\wt{b}_{p-i}, \mu^0_4(\rho_2, \rho_1, \rho_4, \rho_3)) = UV \wt{b}_{p-i}$ that must be cancelled by another term.
	Any other algebra operation on the inputs is zero, so the cancelling term must be a composition of module operations.
	
	First consider the composition of module operations $\wt{m}^0_1(\wt{m}^0_5(\wt{b}_{p-i}, \rho_2, \rho_1, \rho_4, \rho_3))$ or $\wt{m}^0_5(\wt{m}^0_1(\wt{b}_{p-i}), \rho_2, \rho_1, \rho_4, \rho_3)$. A simple analysis of the gradings in $G$ shows that the output of both $\wt{m}^0_5(\wt{b}_{p-i}, \rho_2, \rho_1, \rho_4, \rho_3)$ and $\wt{m}^0_1(\wt{b}_{p-i})$ must be a scalar multiple of $\wt{c}_{p-i}$. Such a composition of module operations then picks up a power of $U^{p-i}$ to respect the $U$-grading, where $p-i > 1$. This cannot cancel against $UV \wt{b}_{p-i}$. The other possibilities for the composition of module operations are $\wt{m}^0_4(\wt{m}^0_2(\wt{b}_{p-i}, \rho_2), \rho_1, \rho_4, \rho_3)$ and $\wt{m}^0_2(\wt{m}^0_4(\wt{b}_{p-i}, \rho_2, \rho_1, \rho_4), \rho_3)$. The idempotents force the output of the inner module operations to be a scalar multiple of $\wt{x}$, but the gradings in $G$ do not allow for such operations.

	The only non-zero term that remains is the composition $\wt{m}^0_3(m^0_3(\wt{b}_{p-i}, \rho_2, \rho_1), \rho_4, \rho_3)$, where we know $m^0_3(\wt{b}_{p-i}, \rho_2, \rho_1) = U \wt{b}_{p-(i+1)}$. To cancel against $UV\wt{b}_{p-i}$, the $\Ainf$ relation necessitates that
	$\wt{m}^0_3(\wt{b}_{p-(i+1)}, \rho_4, \rho_3)$ has a term $V \wt{b}_p$.
	An entirely analogous argument shows that $\wt{m}^0_3(\wt{c}_{p-(i+1)}, \rho_4, \rho_3)$ has a term $UV \wt{c}_{p-i}$.
	
	We are now ready to force things.
	
	Let $D$ be the tensor product $\fModule_{U,V = 1} \boxtimes_{\Algm^{U=1}}\CFDDm(\Id)$, and let $\delta^1_D$ be the differential on it. All the terms in Figure~\ref{fig:cable-p-cfd} that do not have a factor of $\rho_4$ follow from the operations on $\CFAa(\CD_p)$, and our choice of $\TrMPrim^{3,0}$ and $\TrMPrim^{4,0}$. The differential terms labelled with $\rho_{41}$ follow from our analysis above. None of these get cancelled by any other terms in the differential for precisely the same length and grading reasons as in Proposition~\ref{prop:cable-p-cfd}.
	
	The remaining terms in the differential are forced by the type D structure relation for $D$.
	
	The terms $\rho_{1234} \tensor \wt{x}$ and $\rho_{3412} \tensor \wt{x}$ arising from $\mu^1_0 \tensor \wt{x}$ in the structure relation necessitate $\rho_{412} \tensor \wt{x}$ and $\rho_4 \tensor \wt{x}$ in $\delta^1_D(\wt{c}_{p-1})$ and $\delta^1_D(\wt{b}_{p-1})$ respectively.
	
	The only grading permissible terms left to analyze are terms in $\delta^1_D(\wt{c}_{p-i})$. But these, we saw in the proof of Proposition~\ref{prop:cable-p-cfd}, are forced by the structure relation.
\end{proof}

\section{\texorpdfstring{$\CFDm$}{CFD minus} and tensor products}
\label{ch:cfd}

In this section, we take our modules out for a test drive and compute a sample tensor product. We remark that computation of any tensor product with a weighted type D module is a finite task that can be performed by a computer. In particular, the differential on the tensor product preserves the $U$- and $V$-grading. The enumeration of all operations on our $\Ainf$-modules that output at most a certain power of $U$ and $V$ is a tractable task; c.f. Section~\ref{sec:boundedness}.

\begin{figure}
	\centering
	\begin{tikzpicture}
		\node at (0,0) (a1) {$\bullet$};
		\node at (2,0) (a2) {$\bullet$};
		\node at (0,-2) (a3) {$\bullet$};
		\node at (2,-2) (a4) {$\bullet$};
		\draw[->] (a2) to (a1);
		\draw[->] (a2) to (a4);
		\draw[->] (a1) to (a3);
		\draw[->] (a4) to (a3);
	\end{tikzpicture}
	\qquad
	\begin{tikzpicture}
		\node at (4,0) (x1) {$x_1$};
		\node at (0,0) (x2) {$x_2$};
		\node at (0,-4) (x3) {$x_3$};
		\node at (4,-4) (x4) {$x_4$};

		\node at (2,0) (y1) {$y_1$};
		\node at (0,-2) (y2) {$y_2$};
		\node at (2,-4) (y3) {$y_3$};
		\node at (4,-2) (y4) {$y_4$};
	
		\draw[->, bend right=20] (x1) to node[above]{\lab{\rho_3}} (y1);
		\draw[->, bend right=20] (y1) to node[below]{\lab{\rho_{412}}} (x1);
		
		\draw[->, bend right=20] (y1) to node[above]{\lab{\rho_2}} (x2);
		\draw[->, bend right=20] (x2) to node[below]{\lab{\rho_{341}}} (y1);
		
		\draw[->, bend right=20] (x2) to node[left]{\lab{\rho_1}} (y2);
		\draw[->, bend right=20] (y2) to node[right]{\lab{\rho_{234}}} (x2);
		
		\draw[->, bend right=20] (x3) to node[right]{\lab{\rho_{123}}} (y2);
		\draw[->, bend right=20] (y2) to node[left]{\lab{\rho_{4}}} (x3);

		\draw[->, bend left=20] (x4) to node[below]{\lab{\rho_3}} (y3);
		\draw[->, bend left=20] (y3) to node[above]{\lab{\rho_{412}}} (x4);

		\draw[->, bend left=20] (y3) to node[below]{\lab{\rho_2}} (x3);
		\draw[->, bend left=20] (x3) to node[above]{\lab{\rho_{341}}} (y3);

		\draw[->, bend left=20] (x1) to node[right]{\lab{\rho_1}} (y4);
		\draw[->, bend left=20] (y4) to node[left]{\lab{\rho_{234}}} (x1);

		\draw[->, bend left=20] (x4) to node[left]{\lab{\rho_{123}}} (y4);
		\draw[->, bend left=20] (y4) to node[right]{\lab{\rho_{4}}} (x4);
	\end{tikzpicture}
	\qquad
	\begin{tikzpicture}[xscale=2, yscale=2]
		\draw[thin, gray] (-1,-1) grid (1,1);
		\node at (0.12, 0) (x1) {$x_1$};
		\node at (0, 1) (x2) {$x_2$};
		\node at (-0.12, -0) (x3) {$x_3$};
		\node at (0, -1) (x4) {$x_4$};
		
		\node at (-0.5,.5) (y1) {$y_1$};
		\node at (.5,.5) (y2) {$y_2$};
		\node at (-0.5,-0.5) (y3) {$y_3$};
		\node at (0.5,-0.5) (y4) {$y_4$};
	
		\draw[->] (x1) to node[left]{\lab{\rho_3}} (y1);

		\draw[->] (y1) to node[left]{\lab{\rho_2}} (x2);
		
		\draw[->] (x2) to node[right]{\lab{\rho_1}} (y2);

		\draw[->] (y2) to node[right]{\lab{\rho_{4}}} (x3);

		\draw[->] (y3) to node[left]{\lab{\rho_{412}}} (x4);

		\draw[->] (x3) to node[left]{\lab{\rho_{341}}} (y3);

		\draw[->] (y4) to node[right]{\lab{\rho_{234}}} (x1);
		
		\draw[->] (x4) to node[right]{\lab{\rho_{123}}} (y4);
	\end{tikzpicture}
	\caption[A $1\times 1$ rectangle and the associated type $D$
	module.]{\textbf{A $1\times 1$ rectangle and the associated weighted type $D$
			module.} On the right is illustrated the relative $\SpinC$-gradings.}
	\label{fig:11-rectangle-type-D}
\end{figure}
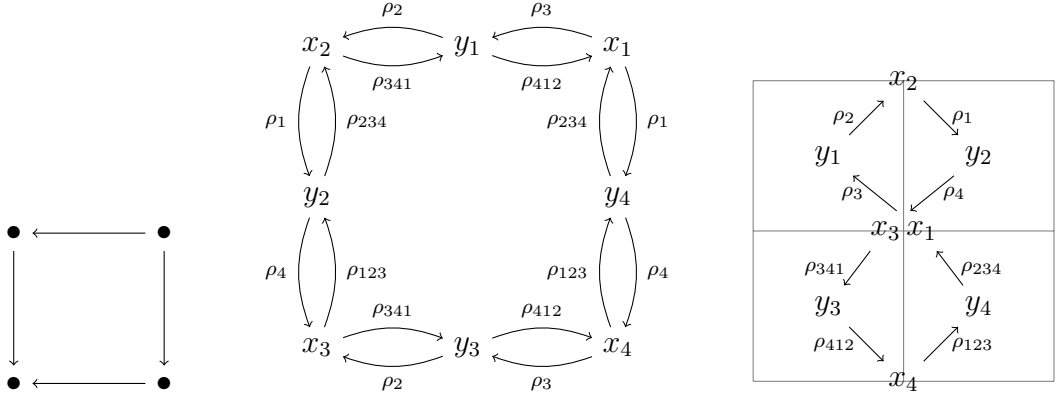

\begin{lemma}
	The type D structure of a $1 \times 1$ rectangle is given as in
	Figure~\ref{fig:11-rectangle-type-D}.
\end{lemma}
\begin{proof}
	Proposition~\ref{prop:type-d-structure}, along with the formula to compute $\CFDa$ from $\CFKm$ \cite[Theorem~11.27]{LOT1}, provides us all terms in the differential in Figure~\ref{fig:11-rectangle-type-D} that do not have a factor of $\rho_4$. The terms that do feature a factor of $\rho_4$ are forced by the weighted type D structure relation, to cancel terms in $\mu_0^1 \tensor x_i$ and $\mu_0^1 \tensor y_i$.
	
	That there are no other terms follows from a simple analysis of the gradings, c.f. Section~\ref{sec:uniqueness}. We illustrate the relative $\SpinC$-gradings on the right in Figure~\ref{fig:11-rectangle-type-D}; in particular, the $\SpinC$-gradings necessitate that there are no terms in the differential between generators that are not adjacent. The Maslov gradings show that adjacent generators do not have further terms between them.
\end{proof}

In \cite[Section~9]{LOT:torus-mod}, Lipshitz, Ozsv\'ath, and Thurston illustrate further examples of constructing $\CFDm$ from $\CFKm$.

The weighted type D module of Figure~\ref{fig:11-rectangle-type-D} is a summand of the type D module of the figure-eight knot. The knot Floer complex of $(p,-1)$-cables of the figure-eight knot have been computed by Hom, Kang,  Park, and Stoffregen using the immersed curve invariants of Hanselman \cite{Hom2022}. In particular, they compute the $UV = 0$ complex and show it lifts uniquely. We compute the tensor product of $\CFAm(\CD_p)$ with the weighted type D module of Figure~\ref{fig:11-rectangle-type-D}, and remark that it matches precisely the computation in \cite{Hom2022}.

\begin{proposition} The tensor product of $\CFAc(\CD_p)$ with
	the weighted type D module of a $1 \times 1$ rectangle (as in
	Figure~\ref{fig:11-rectangle-type-D})
	is given by the complex in Figure~\ref{fig:11-rectangle-cable-cfk}.
\begin{figure}
	\centering
	\subfloat[]{\begin{tikzpicture}[scale=1.25]
		\node at (0,0)(b1)  {$b_1y_1$};
		\node at (2,0)(c1)  {$c_1y_1$};
		\node at (0,-2)(b3)  {$b_1y_3$};
		\node at (2,-2)(c3)  {$c_1y_3$};
		\draw[->](b1) to node[above]{$U$} (c1);
		\draw[->](b3) to node[above]{$U$}(c3);
		\draw[->](b1) to node[right]{$V^p$} (b3);
		\draw[->](c1) to node[right]{$V^p$} (c3);
	\end{tikzpicture} \label{subfl:cable-11-start}} \qquad
	\subfloat[] {\begin{tikzpicture}[scale=1.25]
		\node at (0,0) (xx4) {$xx_4$};
		\node at (0,2) (by4) {$b_{p-1} y_4$};
		\node at (2,0) (xx3) {$xx_3$};
		\node at (2,2) (cy4) {$c_{p-1} y_4$};
		\node at (4,2) (xx1) {$xx_1$};
		\node at (4,0) (by2) {$b_{p-1} y_2$};
		\node at (6,0) (cy2) {$c_{p-1} y_2$};
		\node at (6,2) (xx2) {$xx_2$};
		
		\draw[->] (xx1) to node[above] {$U^p$} (xx2);
		\draw[->] (xx1) to node[right] {$U$} (by2);
		\draw[->] (xx1) to node {} (cy4);
		\draw[->] (xx2) to node {} (cy2);
		\draw[->] (by2) to node[above] {$U^{p-1}$} (cy2);
		\draw[->] (by2) to node[above] {$V$} (xx3);
		\draw[->] (xx4) to node[above] {$U^p$} (xx3);
		\draw[->] (cy4) to node[right] {$UV$} (xx3);
		\draw[->] (by4) to node[above] {$U^{p-1}$} (cy4);
		\draw[->] (by4) to node[right] {$V$} (xx4);
	\end{tikzpicture} \label{subfl:cable-11-end}} \\
	\bigskip
	\subfloat[] {\begin{tikzpicture}[scale=1.25]
		\node at (0,0) (b2y3) {$b_{k+1}y_3$};
		\node at (0,2) (b1y4) {$b_k y_4$};
		\node at (2,0) (c2y3) {$c_{k+1}y_3$};
		\node at (2,2) (c1y4) {$c_k y_4$};
		\node at (4,2) (b2y1) {$b_{k+1}y_1$};
		\node at (4,0) (b1y2) {$b_k y_2$};
		\node at (6,0) (c1y2) {$c_k y_2$};
		\node at (6,2) (c2y1) {$c_{k+1} y_1$};
		
		\draw[->] (b1y4) to node[right] {$V$} (b2y3);
		\draw[->] (b1y4) to node[above] {$U^k$} (c1y4);
		\draw[->] (b2y3) to node[above] {$U^{k+1}$} (c2y3);
		\draw[->] (c1y4) to node[right] {$UV$} (c2y3);
		\draw[->] (b2y1) to node[above] {$V^{p-k-1}$} (c1y4);
		\draw[->] (b2y1) to node[right] {$U$} (b1y2);
		\draw[->] (b1y2) to node[above] {$V^{p-k}$} (c2y3);
		\draw[->] (b2y1) to node[above] {$U^{k+1}$} (c2y1);
		\draw[->] (b1y2) to node[above] {$U^k$} (c1y2);
		\draw[->] (c2y1) to node {} (c1y2);
	\end{tikzpicture} \label{subfl:cable-11-mid}}
	\caption[Summands of the tensor product of $\CFAc(\CD_p)$ with the type D structure
	of a $1 \times 1$ rectangle as in Figure~\ref{fig:11-rectangle-type-D}.]
	{\textbf{Summands of the tensor product of $\CFAc(\CD_p)$ with the type D structure
			of a $1 \times 1$ rectangle as in Figure~\ref{fig:11-rectangle-type-D}.}
			The tensor product for $p = 2$ is
			given by the two summands in the top, where $b_1 = b$ and $c_1 = c$. For $p > 2$,
			there is one additional summand of the bottom form for each
			$1 \leq k \leq p - 2$.}
	\label{fig:11-rectangle-cable-cfk}
\end{figure}
\end{proposition}
\begin{proof} The chord sequences that contribute to the tensor product are as
	follows, labelled along
	the corresponding arrows from Figure~\ref{fig:11-rectangle-cable-cfk}. In a box, we note the weight of corresponding operations in cases with positive weight.

	\begin{center}
	\begin{tikzpicture}[scale=1.5]
		\node at (0,0)(b1)  {$b_1y_1$};
		\node at (2,0)(c1)  {$c_1y_1$};
		\node at (0,-2)(b3)  {$b_1y_3$};
		\node at (2,-2)(c3)  {$c_1y_3$};
		\draw[->](b1) to node[above]{} (c1);
		\draw[->](b3) to node[above]{}(c3);
		\draw[->](b1) to node[sloped,anchor=center,above]{\lab{\rho_{412} \otimes \rho_1 \otimes \rho_4 \otimes \rho_3}} (b3);
		\draw[->](b1) to node[sloped,anchor=center,below]{\tikz \node[rectangle,draw] {\lab{p-2}};} (b3);
		\draw[->](c1) to node[sloped,anchor=center,above,rotate=180]{\lab{\rho_2 \otimes \rho_1 \otimes \rho_4 \otimes \rho_{341}}} (c3);
		\draw[->](c1) to node[sloped,anchor=center,below,rotate=180]{\tikz \node[rectangle,draw] {\lab{p-2}};} (c3);
	\end{tikzpicture} \qquad
	\begin{tikzpicture}[scale=1.5]
		\node at (0,0) (xx4) {$xx_4$};
		\node at (0,2) (by4) {$b_{p-1} y_4$};
		\node at (2,0) (xx3) {$xx_3$};
		\node at (2,2) (cy4) {$c_{p-1} y_4$};
		\node at (4,2) (xx1) {$xx_1$};
		\node at (4,0) (by2) {$b_{p-1} y_2$};
		\node at (6,0) (cy2) {$c_{p-1} y_2$};
		\node at (6,2) (xx2) {$xx_2$};
		
		\draw[->] (xx1) to node[sloped,anchor=center,above] {\lab{\rho_3 \otimes \rho_2}} (xx2);
		\draw[->] (xx1) to node[sloped,anchor=center,above] {\lab{\rho_3 \otimes \rho_2 \otimes \rho_1}} (by2);
		\draw[->] (xx1) to node[sloped,anchor=center,above] {\lab{\rho_1}} (cy4);
		\draw[->] (xx2) to node[sloped,anchor=center,rotate=180,above] {\lab{\rho_1}} (cy2);
		\draw[->] (by2) to node[above] {} (cy2);
		\draw[->] (by2) to node[sloped,anchor=center,above] {\lab{\rho_4}} (xx3);
		\draw[->] (xx4) to node[sloped,anchor=center,above] {\lab{\rho_3 \otimes \rho_2}} (xx3);
		\draw[->] (cy4) to node[sloped,anchor=center,rotate=180,above] {\lab{\rho_4 \otimes \rho_3 \otimes \rho_2}} (xx3);
		\draw[->] (by4) to node {} (cy4);
		\draw[->] (by4) to node[sloped,anchor=center,above] {\lab{\rho_4}} (xx4);
	\end{tikzpicture} \\
	\bigskip
	\begin{tikzpicture}[scale=1.5]
		\node at (0,0) (b2y3) {$b_{k+1}y_3$};
		\node at (0,2) (b1y4) {$b_k y_4$};
		\node at (2,0) (c2y3) {$c_{k+1}y_3$};
		\node at (2,2) (c1y4) {$c_k y_4$};
		\node at (4,2) (b2y1) {$b_{k+1}y_1$};
		\node at (4,0) (b1y2) {$b_k y_2$};
		\node at (6,0) (c1y2) {$c_k y_2$};
		\node at (6,2) (c2y1) {$c_{k+1} y_1$};
		
		\draw[->] (b1y4) to node[sloped,anchor=center,above] {\lab{\rho_4 \otimes \rho_3}} (b2y3);
		\draw[->] (b1y4) to node {} (c1y4);
		\draw[->] (b2y3) to node {} (c2y3);
		\draw[->] (c1y4) to node[sloped,anchor=center,rotate=180,above] {\lab{\rho_4 \otimes \rho_3}} (c2y3);
		\draw[->] (b2y1) to node[sloped,anchor=center,above] {\lab{\rho_{412} \otimes \rho_1}} (c1y4);
		\draw[->] (b2y1) to node[sloped,anchor=center,below] {\tikz \node[rectangle,draw] {\lab{p-k-2}};} (c1y4);
		\draw[->] (b2y1) to node[sloped,anchor=center,rotate=180,above] {\lab{\rho_2 \otimes \rho_1}} (b1y2);
		\draw[->] (b1y2) to node[sloped,anchor=center,above] {\lab{\rho_4 \otimes \rho_{341}}} (c2y3);
		\draw[->] (b1y2) to node[sloped,anchor=center,below] {\tikz \node[rectangle,draw] {\lab{p-k-2}};} (c2y3);
		\draw[->] (b2y1) to node {} (c2y1);
		\draw[->] (b1y2) to node {} (c1y2);
		\draw[->] (c2y1) to node[sloped,anchor=center,rotate=180,above] {\lab{\rho_2 \otimes \rho_1}} (c1y2);
	\end{tikzpicture}
	\end{center}
	
	We will see that no other chord sequence from the type D structure contributes
	to the tensor product. First, any chord sequence
	$a_1 \otimes \dots \otimes a_i \otimes \dots \otimes a_k$ that arises from a non-zero
	operation on the type A module has the property that $a_i a_{i+1} = 0$ for any $i$.
	This limits the chord sequences of interest from the type D structure as those
	obtained from a clockwise or counter-clockwise walk on the graph in
	Figure~\ref{fig:11-rectangle-type-D}.
	
	Next, consider any such chord sequence of length $\geq 5$ that has as a subsequence
	$\rho_3 \otimes \rho_2 \otimes \rho_1 \otimes \rho_4$ or
	$\rho_1 \otimes \rho_4 \otimes \rho_3 \otimes \rho_2$. This cannot arise
	from a tiling pattern as it would necessitate a vertex that meets more than
	four non-blue faces. This limits us to only finitely many chord sequences
	from the type D structure, and we can check that none except the ones
	listed above correspond to a disk.

	\begin{center}
	\begin{tikzpicture}
		\foreach \i/\j in {0/0,0/2,0/4,2/0,2/4,4/0,4/2,4/4} {
			\begin{scope}[shift={(\i,\j)}]
			\draw[dashed] (.25,0) arc[radius=.25, start angle = 0, end angle = 90];
			\draw[dashed] (0,1.75) arc[radius=.25, start angle = 270, end angle = 360];
			\draw[dashed] (1.75,0) arc[radius=.25, start angle = 180, end angle = 90];
			\draw[dashed] (1.75,2) arc[radius=.25, start angle = 180, end angle = 270];
			\draw[color = red,thick] (.25,0) to (1.75,0);
			\draw[color = red,thick] (0,0.25) to (0,1.75);
			\draw[color = red,thick] (.25,2) to (1.75,2);
			\draw[color = red,thick] (2,0.25) to (2,1.75);
			\draw[color = blue,thick] (0,1) to[in=270,out=0] (.667,2);
			\draw[color = blue,thick] (2,1) to[in=270,out=180] (1.333,2);
			\draw[color = blue,thick] (.667,0) to[in=180,out=90] (1,0.75) to[in=90,out=0] (1.333,0);
			\end{scope}
		};
		
		\draw [line width = 4pt, draw=purple, opacity=0.4]
		(6,3.25) -- (6,3.75) arc (270:180:0.25) -- (4.25,4)
		arc (360:270:0.25) -- (4,2.25) arc (90:-180:0.25)
		-- (2.25, 2) arc (0:-270:0.25) -- (2,3.75)
		arc (270:0:0.25) -- (3.75,4) arc (180:-90:0.25)
		-- (4,2.25) arc (90:0:0.25);
		\draw [line width = 4pt, draw=gray, opacity=0.6] 
		(4.25,4) arc (0:-90:0.25) -- (4,2.25) arc (90:0:0.25)
		-- (5.75, 2) arc (180:90:0.25)
		-- (6,2.75);
	\end{tikzpicture}
	\end{center}

	As a demonstrative example, consider the chord sequence 
	$\rho_1 \otimes \rho_4 \otimes \rho_{341} \otimes \rho_{412} \otimes
	\rho_{123} \otimes \rho_{342} \otimes \rho_3 \otimes \rho_2$. Highlighted
	above is a lift of the red curve with this chord sequence;
	starting with the overlap (in gray), the curve lifts to a different copy
	of the square. The blue curve does not bound a disk with the highlighted
	curve. A consideration of the gradings also shows that no such chord sequence corresponds to an operation.
\end{proof}

The computation of the corresponding tensor product for the \emph{hat} flavor has been performed by Petkova \cite{Petkova2013}; in particular, the gradings are computed in Section~5 and remain precisely the same.

\renewbibmacro{in:}{}
\printbibliography

\end{document}